\documentclass[11pt,twoside,reqno]{amsart}
\usepackage{amsmath, amssymb, amsthm, graphicx,marvosym}
\usepackage{mathrsfs}
\usepackage{lipsum}
\usepackage{bbm}
\usepackage{enumitem}
\usepackage{mathtools}
\usepackage{xcolor}
\usepackage{graphicx}
\usepackage{tensor}
\usepackage{hyperref}
\usepackage{enumitem}
\usepackage{cancel}
\usepackage[nobysame]{amsrefs}[11pts]
\usepackage{epigraph}

\usepackage{color}
\usepackage{geometry}
\newtheorem{theorem}{Theorem}[section]
\newtheorem{definition}{Definition}[section]
\newtheorem{lemma}{Lemma}[section]
\newtheorem{remark}{Remark}[section]

\numberwithin{equation}{section}
\numberwithin{figure}{section}

\newtheoremstyle{case}{}{}{}{}{}{:}{ }{}
\theoremstyle{case}

\makeatletter \@addtoreset{equation}{section} \makeatother

\hypersetup{
	colorlinks,
	citecolor=black,
	filecolor=black,
	linkcolor=blue,
	urlcolor=blue
}

\numberwithin{equation}{section}
\numberwithin{figure}{section}
\allowdisplaybreaks[4]
\renewcommand{\a}{\alpha}

\renewcommand{\O}{\Omega}

\newcommand{\supp}{{\mathrm{supp}}}
\newcommand{\loc}{{\mathrm{loc}}}

\renewcommand{\div}{{\mathrm{div}}}

\newcommand{\dx}{\,\mathrm dx}
\newcommand{\dt}{\,\mathrm dt}
\newcommand{\dy}{\,\mathrm dy}
\newcommand{\ds}{\,\mathrm ds}
\newcommand{\dz}{\,\mathrm dz}
\newcommand{\dr}{\,\mathrm dr}

\DeclareMathOperator*{\esssup}{ess\,sup}

\newcommand{\si}{\sigma}
\newcommand{\ga}{\gamma}

\newcommand\reallywidehat[1]{%
  \savestack{\tmpbox}{\stretchto{%
    \scaleto{%
      \scalerel*[\widthof{\ensuremath{#1}}]
      {\kern-.6pt\bigwedge\kern-.6pt}%
      {\rule[-\textheight/2]{1ex}{\textheight}}%
    }{\textheight}%
  }{0.5ex}}%
  \stackon[1pt]{#1}{\tmpbox}%
}

\begin{document}

\title[Non-Formation of Vacuum States]{Global Weak Solutions and Non-Formation of Vacuum States for   Compressible Navier--Stokes--Korteweg  Equations}

\date{\today} 

\author{Yachun Li}
\address[Yachun Li]{School of Mathematical Sciences, CMA--Shanghai, and MOE--LSC, Shanghai Jiao Tong 
University, \\
Shanghai 200240, P. R. China;}
\email{\tt ycli@sjtu.edu.cn}

\author{Fangzhe Qin}
\address[Fangzhe Qin]{School of Mathematical Sciences, Shanghai Jiao Tong University, Shanghai 200240, P.R.China;}
\email{\tt qfzqfz19991219@sjtu.edu.cn}

\author{Shengguo Zhu }
\address[Shengguo Zhu]{School of Mathematical Sciences, CMA--Shanghai, and MOE--LSC, Shanghai Jiao Tong 
University, \\
Shanghai 200240, P. R. China.} \email{\tt   zhushengguo@sjtu.edu.cn}

\begin{abstract}

A long-standing open problem in the mathematical theory of   multidimensional  fluid dynamics  in viscous compressible flow is whether vacuum can form in finite time from initial data with density  bounded away from zero.
The underlying assumption in the derivation of systems of Navier--Stokes type from physical principles is that the fluid is non-dilute and can be described as a continuum. A negative answer to this problem  would support the self-consistency
of the continuum assumption for the Navier--Stokes-type  models. This issue is particularly important when the viscosity coefficients are constant, since, from a physical standpoint, it is unclear how to meaningfully define the fluid velocity in vacuum regions, which in turn makes it difficult to interpret the associated viscous stresses.
 In this paper, for the  barotropic  compressible fluid with capillary effects, when  the coefficients of  viscosity   and capillarity  are all constants, 
we  establish the global existence of weak solutions of the Cauchy problem for 
the corresponding  Navier--Stokes--Korteweg equations  with initial data that are large,  radially symmetric, and away from the vacuum in two and three spatial dimensions. The failure of the well-known  BD entropy estimate poses a substantial challenge to establishing this existence theory.
Moreover,  we rule out finite-time vacuum formation for the two-dimensional radially symmetric weak solutions constructed here. This extends, for the first time, the well-known Hoff--Smoller non-vacuum-formation theory (Commun. Math. Phys. 216, 255–276, 2001) from one-dimensional compressible Navier--Stokes equations with constant viscosity coefficients  to two-dimensional radially symmetric flows with capillarity.
\end{abstract}




\maketitle
\setcounter{tocdepth}{1}
\tableofcontents

\section{Introduction}\label{1}

In  the mathematical theory of  the  multidimensional (M-D)  fluid dynamics  in viscous compressible flow, 
the question  of vacuum formation is of fundamental importance from the perspectives of both  modeling and well-posedness theory with large data, which can be roughly stated as follows: 
\begin{equation}\label{assume00}
\begin{split}
&\textit{Do there exist   solutions of the systems of Navier--Stokes type that exhibit vacuum in}\\
&\textit{finite time when the initial density is  bounded away from zero?}
\end{split}
\end{equation}
As far as we know, when the viscosity coefficients are constants, there  is   no definitive answer until now in M-D  space, even when the initial data have  some form of symmetry; see Bresch--Jabin \cite{BJ}, 
 Chen--Huang--Zhu \cite{CHZ},   Feireisl \cite{fu3}, Feireisl--Novotn\'{y}--Petzeltov\'{a} \cite{fu1}, Hoff \cite{H2}, Hoff--Jenssen \cite{HJ},  Lions \cite{L2}, Jiang--Zhang \cite{JZ}, Plotnikov--Weigant \cite{Pw} and the references therein. It is worth noting  that,  in the view of physics, there is no reasonable way to define the fluid velocity and the  associated viscous stresses in the vacuum region. Actually,  if we try to model the  viscous  compressible  flows with  vacuum states by the  Navier-Stokes type equations   with constant viscosity coefficients, 
 several  counterintuitive behaviors of solutions have been observed, which include  the non-conservation of momentum \cite{zz2}, the failure of continuous dependence on the initial data \cite{hoffserre}, among others.

In this paper, for the  barotropic  compressible fluid with internal capillarity, when  the coefficients of  viscosity   and capillarity  are all constants, 
we establish the global existence of weak solutions  for general   radially symmetric initial data with positive and bounded density for the corresponding   compressible Navier--Stokes--Korteweg equations (\textbf{CNSK}) 
in $\mathbb{R}^n$ ($n=2$ or $3$):
\begin{equation}\label{1.1}
\begin{cases}
\rho_t+\div (\rho{U})=0 , \\
(\rho{{U}})_t+\div (\rho{{U}}\otimes{{U}})+\nabla{P}=\div \mathbb T+\div\mathbb{K}.  \\
\end{cases} 
\end{equation}
 with the initial  data:
\begin{equation}\label{eqs:CauchyInit}
(\rho,U)(0,\boldsymbol{x})=(\rho_0,U_0)(\boldsymbol{x})=(\rho_0(|\boldsymbol{x}|),u_0(|\boldsymbol{x}|)\frac{\boldsymbol{x}}{|\boldsymbol{x}|}) 
\qquad \text{for $\,\boldsymbol{x} \in \mathbb{R}^n$},
\end{equation}
where  $t\geq 0$ and   $\boldsymbol{x}=(x_1,\cdots\!,x_n)^\top\in \mathbb{R}^n$   denote the time  and  Eulerian spatial coordinates, respectively,  $\rho\ge 0$  is  the density of the fluid and the initial density 
$\rho_0$
here is strictly positive and bounded,   
$U=(U_1,\cdots\!,U_n)^\top\in \mathbb{R}^n$ denotes  the  velocity of fluid,  and $P$ is  the pressure. For the polytropic gas,  the constitutive relation is given by $P=A\rho^{\gamma}$, where  $\gamma\geq1$ is a parameter,  and $A>0$ is a  given constant.
The viscous  stress tensor $\mathbb{T}$ and the Korteweg  stress tensor $\mathbb{K}$ are given by
\begin{equation}\label{1.4}
\begin{split}
           \mathbb{T}&=2\mu D(U)+\lambda \div U\mathbb{I}_n,\\
          \mathbb{K}&=\frac{\kappa}{2}\big(\Delta(\rho^2)-|\nabla\rho|^2\big) \mathbb{I}_n-\kappa\nabla\rho\otimes\nabla\rho,
\end{split}
\end{equation}
where $D(U)=\dfrac{\nabla U+(\nabla U)^T}{2}$ is the deformation tensor and  $\mathbb I_n$ is the $n\times n$ identity matrix.  Moreover, $\mu$ denotes the shear viscosity  coefficient, $\lambda+\dfrac{2}{n}\mu$ denotes the bulk viscosity coefficient, and $\kappa>0$ is the capillarity coefficient, which are all constants and  satisfy 
\begin{equation}\label{physical viscosity assumption}
    \mu>0,\quad \lambda+\dfrac{2}{n}\mu\ge 0,\quad \kappa>0.
\end{equation}
It is worth emphasizing that, by  using  a 
momentum–impulse logarithmic contradiction
argument to deal with the coordinates singularity at the origin, we also prove that   for the two-dimensional (2-D) weak solutions with strictly positive initial density constructed here, there is no formation of vacuum states in finite time. Such a feature has not been observed in the study on  weak solutions of the compressible Navier--Stokes equations (\textbf{CNS}); see \cite{BJ,CHZ,DZ,fu3,fu1,H2,HJ,HS,L2,JZ,Pw,XY} and the references therein.

\textbf{CNSK} is   used to  model the motion of viscous compressible   fluids with  capillarity. The model was  first introduced by Van der Waals \cite{Van-der-Waals} and Korteweg \cite{KO}, and then derived rigorously by Dunn--Serrin \cite{DS} via the second gradient theory. This model can also describe the phase transition phenomena (see \cite{Anderson,SE}).
Note that when $\kappa=0$, system \eqref{1.1} is reduced to the classical   \textbf{CNS} for barotropic flows.

There is a vast literature on the global well-posedness of solutions to \textbf{CNS} 
when the viscosity coefficients are constants and the initial density is strictly positive. The one-dimensional (1-D) problem  has been studied extensively; 
see \cite{CHT, HS,Kanel, KN,KS, Denis,ZA1} and the references therein. 
In the M-D case, 
Matsumura--Nishida \cite{MN} first established the global existence of three-dimensional (3-D) classical solutions for initial data close to a non-vacuum equilibrium in  $H^s(\mathbb{R}^n)$. 
We also refer the reader to  Danchin \cite{danchin} for  
global strong solutions with small initial data in some Besov spaces of $\mathbb{R}^n$ 
for $n\geq 2$.

Despite these significant achievements mentioned above, for \textbf{CNS} with constant viscosity coefficients, 
the progress on the global existence of M-D large solutions  with initial data bounded away from vacuum remains very limited  due to the possible implosion or  formation of vacuum states  inside the fluids. Under the assumption of  radial symmetry, 
Hoff \cite{H2} established the global existence of weak solutions with strictly positive initial density  to  the Cauchy problem of \textbf{CNS}
for  isothermal flow in $\mathbb{R}^n$ ($n\geq 2$), and the corresponding study for M-D non-isentropic flows can be found in  Hoff--Jenssen \cite{HJ} and  Chen--Huang--Zhu \cite{CHZ}.
It is worth pointing out that,  due to the lack of the lower bound estimates of the density,  the analysis in \cite{CHZ,H2,HJ} allows the possibility that vacuum  emerges at the origin.
Hoff--Smoller \cite{HS} proved that the 1-D weak solutions of \textbf{CNS} do not exhibit vacuum states, provided that no vacuum states are present initially. Later,  a sufficient 
condition on the regularity of the velocity excluding the formation of vacuum states in finite time  is
given by Xin-Yuan \cite{XY}. Moreover, in  Devault--Gremaud--Jenssen \cite{Jenssen2007}, a numerical result indicates that there are weak solutions with strictly positive initial density  to \textbf{CNS} in two and three spatial dimensions, which develop vacuum  when the initial data are discontinuous and sufficiently large. 
In summary, for M-D \textbf{CNS}  with constant
viscosity coefficients and the initial density  bounded away  from zero{\rm :} 
\begin{itemize}
\item   under radial symmetry, it remains unclear  whether the   Hoff-type weak solutions constructed in \cite{CHZ,H2,HJ}  can develop vacuum  in finite time;
\smallskip
\item Without the symmetry assumption,   we know very little about the dynamics of the vacuum, and then the global  existence of weak solutions  with general  data. 
\end{itemize}
On the other hand,   significant progress on the global existence of M-D weak solutions with general data  has been made when vacuum is allowed in the initial data; see 
Bresch--Jabin \cite{BJ},  Feireisl--Novotn\'{y}--Petzeltov\'{a} \cite{fu1},  Lions \cite{L2}, Jiang--Zhang \cite{JZ}, Plotnikov--Weigant \cite{Pw} and the references therein.

\par When  the  internal capillarity is considered,   the motion of  viscous  compressible fluids can be  governed by  \textbf{CNSK}.
For the  flow  with density-dependent  viscosity  and capillarity  coefficients, substantial progress have been made 
in the  well-posedness theory in 1-D space; 
see \cite{ABS,BH,CH,GL} and the references therein.
In M-D space, by introducing the well-known BD entropy estimates,  Bresch--Desjardins--Lin \cite{BDL} established the global existence of weak solutions with finite energy  to the barotropic \textbf{CNSK}  on the torus or strip domains in two and three spatial  dimensions.  Recently, Antonelli--Spirito \cite{AS1} established the global  existence of  3-D weak solutions with finite energy on the torus via  some subtle truncation arguments, and the  solutions constructed here   hold in  the whole domain they considered.
More,  when $\inf \rho_0>0$,  Haspot \cite{Has3} first established  the global well-posedness  of  strong solutions  with general 
data  to the quantum \textbf{CNS} in  $\mathbb{R}^n(n\geq2)$. We also refer the reader  to 
\cite{AM2,JG,BGV,VV,Has2} for  other important progress.

Moreover,  when the coefficients of the viscosity and capillary are all constants,   \textbf{CNSK} also  has very wide and important applications in  both  physical modeling and  
numerical simulations;  see \cite{MM,GN} and the references therein. When the initial data is away from vacuum,  some key achievements in the well-posedness theory also have been obtained for this case.   Hattori--Li \cite{HL1} established the local well-posedness of strong solutions to the  Cauchy problem, and  then  proved  the corresponding  global existence with  small  data in \cite{HL2}. Later, in $\mathbb{R}^n$ ($n\geq 2$),
Danchin--Desjardins \cite{DD} proved the global existence of the unique strong solution with small initial data in some Besov spaces to the isothermal  flow.
Chen--Zhao \cite{CZ} established the  nonlinear stability of stationary solutions in $H^s$ to the Cauchy problem for the 3-D non-isentropic flow. We also refer the reader to Hou--Peng--Zhu \cite{HPZ1, HPZ2} for M-D global well-posedness of strong solutions  with small data and  vacuum. 

\par However, in stark contrast,  for the M-D flow with constant viscosity and capillary coefficients, the BD entropy estimates fail, which,  along with  the strong nonlinearities in  both the pressure term $\nabla P$ and  the Korteweg stress  term $\text{div}\mathbb{K}$, makes the study on the global well-posedness theory with large data challenging. 
 To the best of our knowledge, the global existence problem for weak or strong solutions to \textbf{CNSK} with large M-D initial data remains largely open, even under symmetry assumptions. This is due, in part, to the limited understanding of the associated vacuum formation problem.
 
In this paper,  
we  establish the global existence of weak solutions to  the Cauchy problem of  
system \eqref{1.1}  with initial data that are large,  radially symmetric, and away from the vacuum in $\mathbb{R}^n$ ($n=2,3$). This partially solves the open problem  proposed by Bresch--Gisclon--Lacroix-Violet on page 981 of \cite{BGV}. Moreover, we will prove that, in two spatial dimensions,  the solutions  we constructed here do not allow the formation of vacuum states in finite time. This extends, for the first time, the well-known Hoff--Smoller non-vacuum-formation theory in \cite{HS} from 1-D \textbf{CNS}   to 2-D radially symmetric flows with  internal capillarity. Specifically, we consider the radially  symmetric solutions of the form
\begin{equation}\label{solutionform}
    (\rho,U)(t,\boldsymbol{x})=\left(\rho(t,|\boldsymbol{x}|),u(t,|\boldsymbol{x}|)\dfrac{\boldsymbol{x}}{|\boldsymbol{x}|}\right)
\end{equation}
of system \eqref{1.1} with \eqref{1.4}--\eqref{physical viscosity assumption} in $\mathbb{R}^n$ ($n=2,3$)  with Cauchy  data \eqref{eqs:CauchyInit}.
Therefore, by  the radial coordinate variable $r=|\boldsymbol{x}|$, problem  \eqref{1.1}--\eqref{physical viscosity assumption} can be reformulated as
\begin{equation}\label{1.6}
    \begin{cases}
        \rho_t+(\rho u)_r+\dfrac{m\rho u}{r}=0,\\
   \rho(u_t+uu_r )+(A\rho^\gamma)_r=(2\mu+\lambda)u_{rr}+\kappa\rho\rho_{rrr}+(2\mu+\lambda)\left(\dfrac{mu}{r}\right)_r+\kappa\rho\left(\dfrac{m\rho_r}{r}\right)_r,
    \end{cases}
\end{equation}
in $I_T=[0,T]\times[0,\infty)$ with  $m=n-1$, and 
\begin{align}
    &(\rho, u)(0,r)=(\rho_0,u_0)(r)\qquad \ \ \text{for} \quad r\in [0,\infty),\label{radilinitial}\\
    &u(t,0)=0,\quad\rho_r(t,0)=0\qquad \text{for} \quad t\in [0,T].\label{bd condition}
\end{align}
The boundary conditions \eqref{bd condition} are derived from the continuity of $x\rightarrow (U,\nabla  \rho)(t,x)$ at $|\boldsymbol{x}|=0$.

Due to the compressibility of fluids, 
our analysis encounters two major obstacles:
\begin{itemize}
\item the possible cavitation inside the fluids in finite time.
\item the possible implosion inside the fluids in finite time.
\end{itemize}
Overcoming these difficulties is particularly challenging because of several inherent issues:
\begin{itemize}
\item the coordinate singularity at the origin;
\smallskip
\item the degeneracy in the time evolution, arising  from the possible formation of vacuum states;
\smallskip
\item the failure of the BD entropy estimate;
\smallskip
\item the strong nonlinearities in  both the pressure term $\nabla P$ and  the Korteweg stress  term $\text{div}\mathbb{K}$.
\end{itemize}
Thus, new ideas  \eqref{1.1} will be  provided in our anlysis.
Next, before stating  the   main results, we first introduce  the definition of  weak solutions to   problem \eqref{1.1}-\eqref{physical viscosity assumption} as follows.
\begin{definition}\label{defni}
A vector function   $(\rho,U)(t,\boldsymbol{x})$ is called a  weak solution of the Cauchy problem \eqref{1.1}-\eqref{physical viscosity assumption} if it satisfies the following conditions: \\
\smallskip
\rm{(i)} There exist  a H$\ddot{{\rm o}}$lder continuous curve $\underline r(t): [0,T]\rightarrow [0,\infty)$ and a positive constant $C_0=C_0(U_0,\rho_0,A, \\
\gamma, \mu,\lambda, \kappa,n)$ such that 
    \[
    \lim_{t\rightarrow0}\underline r(t)=0,\qquad0\le\underline r(t)\le C_0\text{ for all }t\in[0,T].
    \]
We define $\mathcal{F}:=\{(t,x)\in[0,T]\times \mathbb{R}^n:|\boldsymbol{x}|>\underline r(t)\}$ as the fluid region and $\mathring{\mathcal{F}}=\mathcal{F}\cap(0,T)$ as an open set. \\
\smallskip
\rm{(ii)} $\rho-1\in L^\infty([0,T];H^1(\mathbb{R}^n))$, $\sqrt{\rho}U\in L^\infty([0,T];L^2(\mathbb{R}^n))$,  $\nabla U\in L^2([0,T];L^2(\mathbb{R}^n))$; $\rho\in L^\infty_{loc}(\mathcal{F})$, $\rho=0 $ for a.e.  in $\{(t,r)\in(0,T]\times\mathbb{R}^n\backslash \mathcal{F}\}$; and  $(U,\nabla\rho)$ is locally H$\ddot{{\rm o}}$lder continuous in $\mathcal{F}\backslash\{t=0\}$. \\
\smallskip
\rm{(iii)}
For each $\Phi\in  C^1([0,T];C^1_c(\mathbb{R}^n))$ and for all $t\in[0,T]$,
    \begin{equation*}
\int_{\mathbb{R}^n}\rho(t,\boldsymbol{x})\Phi(t,\boldsymbol{x})\mathrm{d}\boldsymbol{x}=\int_{\mathbb{R}^n}\rho_0(\boldsymbol{x})\Phi(0,\boldsymbol{x})\mathrm{d}\boldsymbol{x}+\int^t_0\int_{\mathbb{R}^n}(\rho\Phi_t+\rho U\cdot\nabla\Phi)\mathrm{d}\boldsymbol{x}\ds.
    \end{equation*}
\smallskip
\rm{(iv)}
For each $\Phi\in C^1([0,T];C^2_c(\mathbb{R}^n);\mathbb{R}^n)$ and for all $t\in[0,T]$,
    \begin{align*}
\int_{\mathbb{R}^n}\rho(t,\boldsymbol{x})U(t,\boldsymbol{x})\cdot\Phi(t,\boldsymbol{x})\mathrm{d}\boldsymbol{x}&=\int_{\mathbb{R}^n}\rho_0(\boldsymbol{x})U_0(\boldsymbol{x})\cdot\Phi(0,\boldsymbol{x})\mathrm{d}\boldsymbol{x}+\int^t_0\int_{\mathbb{R}^n}(\rho U\cdot\Phi_t+\rho U\cdot(U\cdot\nabla\Phi))\mathrm{d}\boldsymbol{x}\ds\\
        &-\int^t_0\int_{\mathbb{R}^n}(-A\rho^\gamma \div\Phi+  \mu\nabla U:\nabla\Phi+(\mu+\lambda)\div U\div\Phi)\mathrm{d}\boldsymbol{x}\ds\\
        &+\kappa\int^t_0\int_{\mathbb{R}^n}\rho\nabla\rho\cdot\nabla\div\Phi\dx\ds+\dfrac{\kappa}{2}\int^t_0\int_{\mathbb{R}^{n}}|\nabla\rho|^2\div\Phi\mathrm{d}\boldsymbol{x}\ds\\
        &+\kappa\int^t_0\int_{\mathbb{R}^n}\nabla\rho\otimes\nabla\rho:\nabla\Phi\mathrm{d}\boldsymbol{x}\ds.
    \end{align*}
 
\end{definition}

Our first theorem establishes the global existence of radially symmetric weak solutions to the Cauchy problem \eqref{1.1}--\eqref{physical viscosity assumption} in two-dimensions and three-dimensions. 
\begin{theorem}\label{2.4}
    Let $n=2$ or $3$, $\gamma>1$,  and let the radially symmetric initial data  satisfy
\begin{align}\label{2.5}
&
\widetilde{C}^{-1}\le \rho_0(r)\le \widetilde{C},
\quad \text{for  } a.e.~r\in[0,\infty),                                      \notag\\
&\int_0^\infty
\left\{
\rho_0\left(\dfrac12 |u_0|^2\right)
+G(\rho_0)
+\left|(\rho_0-1,(\rho_0)_r)\right|^2
\right\} r^m\,\dr
\le \widetilde{C},
\end{align}
for some given constant $\widetilde{C}>1$, where 
\[
G(\zeta):=\dfrac{A}{\gamma-1}\zeta^\gamma-\dfrac{A\gamma}{\gamma-1}\zeta+A,
\qquad
\psi(\zeta):=\dfrac{A}{\gamma-1}\left(\zeta^{1-\gamma}+(\gamma-1)\zeta-\gamma\right).
\]
Then, for each $T>0$, there exists a global radially symmetric weak solution
$(\rho,U)(t,x)$ of  problem \eqref{1.1}-\eqref{physical viscosity assumption} in the sense of \textnormal{Definition \ref{defni}} such that the following properties hold:
\begin{itemize}
    \item There exists a positive constant $C(T)=C(T,\tilde{C},A,\gamma, \mu, \lambda,\kappa,n)$ such that
    \begin{align*}
        \esssup_{t\in[0,T]}\int_{\mathbb{R}^n}(\rho|U|^2+G(\rho)+|\nabla\rho|^2)\mathrm{d}\boldsymbol{x}+\int^T_0\int_{\mathbb{R}^n}|\nabla U|^2\mathrm{d}\boldsymbol{x}\dt\le C(T).
    \end{align*}
    \smallskip
    \item There exists a continuous function $(t,y)\mapsto\widetilde{r}(t,y):[0,T]\times[0,\infty)\rightarrow[0,\infty)$ such that $y\mapsto\widetilde{r}(t,y)$ is strictly increasing. The function $t\mapsto\underline{r}(t)=\lim_{y\rightarrow0}\widetilde{r}(t,y)$ is H$\ddot{o}$lder continuous on $[0,T]$. Moreover, there exists a positive constant $C_0=C_0(\widetilde{C},A,\gamma,\mu,\lambda,\kappa,n)$ such that
    \begin{align*}
        \begin{cases}
        \displaystyle
            \int_{\{\underline{r}(t)<|\boldsymbol{x}|<\widetilde{r}(t,y)\}}\rho(t,|\boldsymbol{x}|)\mathrm{d}\boldsymbol{x}=y|S^{n-1}|\quad&\text{ for }a.e. (t,y)\in[0,T]\times (0,\infty),\\
            \displaystyle
            \widetilde{r}(t,y)=\widetilde{r}_0(y)+\int^t_0 u(s,\widetilde{r}(s,y))\ds\quad&\text{ for }a.e. (t,y)\in[0,T]\times (0,\infty),\\
            \displaystyle
            ny\psi^{-1}_-\left(\dfrac{C_0}{y}\right)\le\widetilde{r}^n(t,y)\le C_0(1+y)\quad&\text{ for all }(t,y)\in[0,T]\times (0,\infty),\\
        \end{cases}
    \end{align*}
    where $u(t,|\boldsymbol{x}|):=U(t,\boldsymbol{x})\cdot\dfrac{\boldsymbol{x}}{|\boldsymbol{x}|}$, $\psi^{-1}_-$ is the left inverse branch of $\psi$, $|S^{n-1}|$ is the surface measure of unit sphere, and function $y\mapsto\widetilde{r}_0(y)$ is implicitly defined by
    \begin{align*}
        y=\int^{\widetilde{r}_0(y)}_0\rho_0(r)r^m\dr\quad\text{ for each }y\in[0,\infty).
    \end{align*}
    In addition, the speed of the vacuum interface  satisfies 
     \begin{align*}
         \dfrac{\mathrm{d}\underline{r}(t)} {\mathrm{d}t}=u(t,\underline{r}(t)) \text{ for} \ \text{a.e. } \ t\in \{t\in(0,T):\underline{r}(t)>0\}.
     \end{align*}
    \item For each  $\varepsilon\in(0,1)$, there exists a  constant $C(\varepsilon,T):=C(\varepsilon, T,\tilde{C}, A,\gamma,\mu,\lambda,n,\kappa)>0$ such that
    \begin{align*}
        \begin{cases}
            |\widetilde{r}(t,y_1)-\widetilde{r}(t,y_2)|\le C(\varepsilon,T)|y_1-y_2|\quad&\text{ for all }(t,y_1,y_2)\in[0,T]\times[\varepsilon,\infty)^2,\\
            |\widetilde{r}(t_1,y)-\widetilde{r}(t_2,y)|\le C(\varepsilon,T)|t^{\frac34}_1-t^{\frac34}_2|\quad&\text{ for all }(t_1,t_2,y)\in[0,T]^2\times[\varepsilon,\infty).
        \end{cases}
    \end{align*}
    Moreover, let $\si(t)=\min\{1,t\}$ and $B^c(r):=\{\boldsymbol{x}\in \mathbb{R}^n:|\boldsymbol{x}|\ge r\}$. Then
    \begin{align*}
        C^{-1}(\varepsilon,T)\le\rho(\boldsymbol{x},t)\le C(\varepsilon,T),\quad |(U(\boldsymbol{x},t),\nabla\rho(\boldsymbol{x},t))|\le \sigma^{-\frac14}(t)C(\varepsilon,T),
    \end{align*}
    for a.e.  $(\boldsymbol{x},t)$ such that $t\in[0,T]$ and $x\in B^c(\widetilde{r}(t,\varepsilon))$, and 
    \begin{align*}
        &\sup_{t\in[0,T]}\int_{B^c(\widetilde{r}(t,\varepsilon))}\left|\left(\rho-1,\sigma^{\frac12}(|U|^2,\nabla U,\nabla^{2}\rho,|\nabla\rho|^2\right)\right|^2\mathrm{d}\boldsymbol{x}\\
        &+\int^T_0\int_{B^c(\widetilde{r}(t,\varepsilon))}\left|\left(\nabla U,\nabla^2\rho,\rho_t,|\nabla\rho|^2,\sigma^{\frac12}(U_t,\nabla^2 U,\nabla^3\rho,\nabla\rho_t)\right)\right|^2\mathrm{d}\boldsymbol{x}\dt\le C(\varepsilon,T).
    \end{align*}
    In addition, for all $t\in(0,T]$ and $\boldsymbol{x},\boldsymbol{y}\in B^c(\widetilde{r}(t,\varepsilon))$, 
    \begin{align*}
        \sigma^{\frac12}(t)|U(t,\boldsymbol{x})-U(t,\boldsymbol{y})|+\sigma^{\frac12}(t)|\nabla\rho(t,\boldsymbol{x})-\nabla\rho(t,\boldsymbol{y})|\le C(\varepsilon,T)|\boldsymbol{x}-\boldsymbol{y}|^{\frac12},
    \end{align*}
    and for all $0<t_1<t_2\le T$ and $|\boldsymbol{x}|\ge\sup_{t\in[t_1,t_2]}\widetilde{r}(t,\varepsilon)$,
    \begin{align*}
        \sigma^{\frac12}(t_1)|U(t_1,\boldsymbol{x})-U(t_2,\boldsymbol{x})|+\sigma^{\frac12}(t_1)|\nabla\rho(t_1,\boldsymbol{x})-\nabla\rho(t_2,\boldsymbol{x})|\le C(\varepsilon,T)|t_1-t_2|^{\frac14}.
    \end{align*}
\end{itemize}
\end{theorem}

Our  second theorem  on the non-formation of vacuum states in  $\mathbb{R}^2$ can be stated as follows.
\begin{theorem}\label{th-nonvacuum}
 Let $n=2$, $\gamma>1$,  and the radially symmetric initial data $(\rho_0,u_0)$  satisfy \eqref{2.5}. Then for each $T>0$, let 
$(\rho,U)(t,\boldsymbol{x})$ be the  radially symmetric weak solution in $[0,T]\times \mathbb{R}^2$  to the Cauchy problem \eqref{1.1}-\eqref{physical viscosity assumption} obtained in {\rm Theorem \ref{2.4}}. Then $\underline{r}(t)=0$ for each $t\in[0,T]$, and $\rho$ satisfies: for $0<t \leq T$,
\begin{align*}
        &\min\{C^{-1}(T,p)|\boldsymbol{x}|^{\frac{2p}{1-p}},C^{-1}(T)\}\le\rho(t,\boldsymbol{x})\ \text{ for}\ |\boldsymbol{x}|>0,\ p\in\left(\frac{2}{3},1\right)\\
        &\rho(t,\boldsymbol{x})\le C(T)\text{ for  }\ |\boldsymbol{x}|\ge1,\quad\rho(t,\boldsymbol{x})\le C(T)\left(1+\log\dfrac{1}{|\boldsymbol{x}|}\right)\text{ for} \ 0<|\boldsymbol{x}|<1,
    \end{align*}
for some  constant $C(T)=C(T,\tilde{C},A,\gamma, \mu, \lambda,\kappa,n)>0$,
which means that the vacuum will not form in any non-empty open set in $\mathbb{R}^2$ in finite time. 
\end{theorem}

\begin{remark}
    To the best of our knowledge, this is the first result, stated in  {\rm Theorem \ref{2.4}}, concerning the global existence of weak solutions  to   \textbf{CNSK} with  large  data   for the flow with constant viscosity coefficients and internal capillarity in M-D space, even  under radial symmetry. Our main contribution here is to develop a  careful compactness argument    to deal with the strong nonlinearities in  both the pressure term $\nabla P$ and  the Korteweg stress  term $\text{div}\mathbb{K}$. \end{remark}

\begin{remark}
For the viscous compressible  flow with constant viscosity coefficients, 
the non-vacuum formation conclusion stated in {\rm Theorem \ref{th-nonvacuum}}  extends, for the first time,  the well-known Hoff--Smoller \cite{HS} non-vacuum-formation theory established in 2001  from 1-D \textbf{CNS} flow   to 2-D radially symmetric flows with internal capillarity.  The key point in our analysis is to introduce    a 
momentum–impulse logarithmic contradiction
argument to deal with the coordinate singularity at the origin.\end{remark}

\begin{remark}
   For the question of vacuum formation in finite time, in stark contrast with the 2-D flow,  the corresponding study  for the 3-D flow seems formidable, since the  momentum-impulse logarithmic contradiction argument we introduced here fails in three spatial dimensions. Specifically, if we assume that vacuum forms in some open set  in finite time,  in radial Eulerian coordinates,  the impulse imparted by viscous stresses exhibits $\log\dfrac{1}{\underline{r}(s)}$ growth near the vacuum region, while  
 the available  estimates in 3-D space  can only show us that  the available estimates yield only a polynomial bound in $\underline {r}^{-\a}$ for the momentum with $\a>0$, which does not lead to the desired contradiction.\end{remark}
\begin{remark}
With suitable modifications, when the initial density is bounded away from zero and the coefficients of viscosity and capillarity are all constants,  the methodology developed in this paper can be applied to proving the global existence 
of radially symmetric weak solutions   
to the  Cauchy problem of the barotropic {\rm\textbf{CNSK}} 
with the isothermal pressure law $P=A\rho$
in $\mathbb{R}^n$  for {\rm $n=2$, $3$},  and the related non-formation of vacuum states in $\mathbb{R}^2$, which will be addressed 
in our forthcoming paper \cite{LQZ}.
\end{remark}

\par The rest of this  paper is organized as follows: In \S \ref{2},
the main strategy and the overall scheme of the proof are described. The Cauchy problem considered can be regarded as the limit of a series of initial-boundary value problems that are formulated in finite annular regions. For such approximation problems, in \S \ref{3}, we derive some \textit{a priori} estimates that are independent of both the inner and outer radii of the annuli considered in the radially symmetric Lagrangian coordinates. Based on these uniform  estimates,  in \S \ref{4}, we establish  global weak solutions of exterior problems in Eulerian coordinates, and in \S \ref{5}, global weak solutions of the Cauchy problem in the domain containing the origin in Eulerian coordinates is obtained. In \S \ref{6}, the non-formation of vacuum states in $\mathbb{R}^2$ is proved via some  momentum-impulse logarithmic contradiction argument, and,  furthermore, the lower bound  of the density $\rho$  near the origin  is characterized. Finally,  some basic lemmas  are collected in Appendix A.

\section{Reformulation and Main Strategy}\label{2}
 Throughout this paper, we use the following notation for five types of positive generic constants
that are all larger than $1$ and may be different at different occurrence (without loss of generality)
to indicate their specific 
dependence on various parameters:
 \begin{align*}
     &C= C(A,\gamma,\mu,\lambda,\kappa,n),\quad C_0=C_0(C,\widetilde{C}),\quad C(T)=C(T,C_0),\\
     &C(\varepsilon,T)=C(\varepsilon,T,C_0),\quad C(a,T)=C(a,T,C_0)\text{ for each }a\in(0,1)\text{ and }\varepsilon\in(0,1).  
 \end{align*}
 Moreover, we  also denote  $\beta=2\mu+\lambda$.

 \subsection{Reformulation}\label{Reformulation}
 The outline of the proof for  Theorem \ref{2.4} on the global existence of weak solutions to  the Cauchy problem  \eqref{1.1}--\eqref{physical viscosity assumption} can be stated as follows: 
\begin{itemize}

\item [\rm{(i)}] first, we  regard the Cauchy problem  \eqref{1.1}--\eqref{physical viscosity assumption}  as the limit of  a series of exterior problems \eqref{approximate problem a}-\eqref{bd condition-ex} formulated in $[a,\infty)$ in  radially symmetric  Eulerian coordinates with $a>0$;

\item [\rm{(ii)}] then  we regard exterior problems \eqref{approximate problem a}--\eqref{bd condition-ex}  as the limit of a series of free boundary value problems (\textbf{FBVP}) \eqref{FBVP}--\eqref{bd condition-annuli} posed on finite annuli in radially symmetric     Eulerian coordinates;

\item [\rm{(iii)}] by introducing Lagrangian coordinates, we can reduce the study of \textbf{FBVP} \eqref{FBVP}--\eqref{bd condition-annuli} to that of the corresponding initial-boundary value problem  (\textbf{IBVP}) \eqref{approximate problem}--\eqref{initial boundary condition} on  fixed finite balls;

\item [\rm{(iv)}] for such 
approximate \textbf{IBVP} \eqref{approximate problem}--\eqref{initial boundary condition}, we  derive careful {\it a-priori} estimates  that are independent 
of both the inner radius  and outer radius  of the annuli  
considered in the radially  symmetric Lagrangian coordinates;

\item [\rm{(v)}]  we convert back into Eulerian coordinates and let the outer radii of \textbf{FBVP}  \eqref{FBVP}--\eqref{bd condition-annuli}   tend to infinity to obtain the global existence of weak solutions to exterior problems \eqref{approximate problem a}-\eqref{bd condition-ex};
\item [\rm{(vi)}]
 we let the inner radii of exterior problems \eqref{approximate problem a}-\eqref{bd condition-ex} tend to zero to obtain the global weak solution of the Cauchy problem \eqref{1.1}--\eqref{physical viscosity assumption}. 
\end{itemize}
Next we show three different reformulations required in our analysis.
 \subsubsection{ Exterior problems in  radially symmetric  Eulerian coordinates}
Let $(\rho_0,U_0)(\boldsymbol{x})$ be   initial data specified  in Theorem \ref{2.4}. First, we mollify $(\rho_0,U_0)$ by the standard mollifier $j_a$, denoted by $$(\overline{\rho}^a_0,\overline{U}^a_0):=(j_a*\rho_0,j_a* U_0).$$ Thus $(\rho^a_0,U^a_0)\in C^\infty(\mathbb{R}^n)$. Next,  let $\chi\in C^\infty:[0,\infty)\rightarrow[0,1]$ be such that $\chi(\varsigma)=0$ if $\varsigma\le1$, $\chi(\varsigma)=1$ if $\varsigma\ge 2$, and  $|\chi'(\varsigma)|\le2$ in  $[0,\infty)$. 
For each $a\in(0,1)$, denote $\chi^a(r)=\chi\left(\dfrac{r}{a}\right)$ by $r=|\boldsymbol{x}|$.
At last, denote 
$$(\rho^a_0, U^a_0)(\boldsymbol{x}):=(\rho^a_0(r),  u^a_0(r)\frac{\boldsymbol{x}}{|\boldsymbol{x}|})=\big(\chi^a(|x|)(\overline{\rho}^a_0-1)+1,\chi^a(|x|)\overline{U}_0\big).$$
According to  Lemma \ref{mollify initial data still spher} in Appendix A, $(\rho_0,U_0)$ is still spherical symmetry, and then 
\begin{align}\label{trun inital data}
    \rho^a_0(r)=\chi^a(r)(\overline{\rho}^a_0(r)-1)+1,\quad u^a_0(r)=\chi^a(r)\overline{u}_0(r)\text{ for all }r\ge0.
\end{align}
It follows from the initial conditions on $(\rho_0,U_0)$ in Theorem \ref{2.4}  that  there exists a constant $\widetilde{C}>0$ independent of $a\in(0,1)$ such that
\begin{align}\label{3.3}
    &(\overline{\rho}^a_0,\overline{u}^a_0)(r)\in C^\infty([0,\infty)),\quad\widetilde{C}^{-1}\le\overline{\rho}^a_0\le\widetilde{C}\text{   for  }a.e.~r\ge0,\nonumber\\
    &(\overline{\rho}^a_0,\overline{u}^a_0)(r)\rightarrow(\rho_0,u_0)\quad\text{ as }~a \rightarrow0\text{ for  }a.e. \varsigma\in(0,\infty).
\end{align}
One has $(\rho^a_0,u^a_0)(r)=(1,0)$ for $r\le a$ and $\left(\rho^a_0,u^a_0\right)(r)=(\overline{\rho}^a_0,\overline{u}^a_0)$ for $r\ge 2a$.
We note  that the data $ (\rho^a_0,u^a_0)   $  preserve the uniform bound \eqref{2.5} (see Lemmas \ref{lemma3.2} and  \ref{lemma 3.4}), which is  independent of $a\in(0,1)$.

Now we introduce the  exterior problems in  radially symmetric  Eulerian coordinates in our proof.
 \begin{equation}\label{approximate problem a}
    \begin{cases}
        \rho_t+(\rho u)_r+\dfrac{m\rho u}{r}=0,\\
   \rho(u_t+uu_r )+(A\rho^\gamma)_r=(2\mu+\lambda)u_{rr}+\kappa\rho\rho_{rrr}+(2\mu+\lambda)\left(\dfrac{mu}{r}\right)_r+\kappa\rho\left(\dfrac{m\rho_r}{r}\right)_r,\\
    \end{cases}
\end{equation}
in $I_T=[0,T]\times[a,\infty)$ with  $m=n-1$ for some positive constant $a$, and 
\begin{align}
    &(\rho, u)(0,r)=(\rho^a_0,u^a_0)(r)\qquad \ \ \text{for} \quad r\in [a,\infty),\label{radilinitial-ex}\\
    &u(t,a)=0,\quad\rho_r(t,a)=0\qquad \text{for} \quad t\in [0,T].\label{bd condition-ex}
\end{align}

\subsubsection{{\rm \textbf{IBVP}} posed in  finite balls in radially symmetric  Lagrangian coordinates}
First, let $\chi\in C^\infty:\mathbb{R}\rightarrow[0,1]$  be such that $\chi(\varsigma)=1$ if $\varsigma\le0$ and $\chi(\varsigma)=0$ if $\varsigma\ge1$.

Second, we define $h(r)=\int^r_a\rho^a_0(r)r^m\dr$. Since $0<\rho^a_0(r)<\infty$, by the inverse function theorem, there exists a smooth, bijective, strictly increasing map $\widetilde{r}^a_0(h):[0,\infty)\rightarrow[a,\infty)$, which satisfies 
$$
h(\widetilde{r}^a_0(y))=\int^{ \widetilde{r}^a_0(y)  }_a\rho^a_0(r)r^m\dr=y.
$$
Next  let  $(\rho^a_0,u^a_0)$ be the data given in \eqref{trun inital data}, and   $\varphi^{a,b}_0(r):=\chi\left(\dfrac{2r-\widetilde{r}^a_0(b)}{\widetilde{r}^a_0(b)}\right)$ for some constants  $a\in(0,1)$ and $b\in  \mathbb{N}$, which satisfies 
\begin{align*}
    \varphi^{a,b}_0=1\text{ for } r\in\left[a,\dfrac{\widetilde{r}^a_0(b)}{2}\right],\qquad\varphi^{a,b}_0=0\text{ for } r\in[\widetilde{r}^a_0(b),\infty).
\end{align*}
Then we can truncate the initial data  for each $r\in[a,\infty)$ as follows:
\begin{align}\label{ trun ca ini data}
    (\rho^{a,b}_0,u^{a,b}_0)(r)=\left(\varphi^{a,b}_0(\rho^a_0-1)+1,\varphi^{a,b}_0u^a_0\right).
\end{align}

Now, by denoting 
\begin{align}\label{approximate particle path}
        P=\dfrac{A}{(\widetilde{v}^{a,b})^\gamma},\quad\widetilde{r}^{a,b}(t,x)=\left(a^n+n\int^x_0\widetilde{v}^{a,b}(t,s)\ds\right)^{1/n},
\end{align}
we are ready to introduce the approximation problem posed in finite balls in radially symmetric Lagrangian coordinates:
 \begin{align}\label{approximate problem}
\left\{
    \begin{array}{ll}
 \partial_t\widetilde{v}^{a,b}&=\partial_x\bigl((\widetilde{r}^{^{a,b}})^m\widetilde{u}^{a,b}\bigl),\\
          \partial_t\widetilde{u}^{a,b}&=\left(\widetilde{r}^{a,b}\right)^m\partial_x\left(\beta\dfrac{\partial_x\bigl((\widetilde{r}^{^{a,b}})^m\widetilde{u}^{a,b}\bigl)}{\widetilde{v}^{a,b}}-\kappa\partial_x\left(\dfrac{(\widetilde{r}^{a,b})^{2m}\partial_x\widetilde{v}^{a,b}}{(\widetilde{v}^{a,b})^5}\right)-\dfrac{5\kappa}{2}\dfrac{(\widetilde{r}^{a,b})^{2m}|\partial_x\widetilde{v}^{a,b}|^2}{(\widetilde{v}^{a,b})^6}\right)
         \\
         &-(\widetilde{r}^{a,b})^m\partial_xP(\widetilde{v}^{a,b})-m\kappa\dfrac{(\widetilde{r}^{a,b})^{2m-1}|\partial_x\widetilde{v}^{a,b}|^2}{(\widetilde{v}^{a,b})^5},
    \end{array}
    \right.
\end{align}
in $[0,T]\times [0,b]$ with  $m=n-1$ for some constants  $a\in(0,1)$ and $b\in  \mathbb{N}$, and 
\begin{align}\label{initial boundary condition}
    \begin{cases}
        \widetilde{u}^{a,b}=\partial_x\widetilde{v}^{a,b}=0 \text{ at }x=0,b,\\
        \left(\widetilde{v}^{a,b},\widetilde{u}^{a,b}\right)(0,x)=\left(\widetilde{v}^{a,b}_0,\widetilde{u}^{a,b}_0\right)(x):=\left((\rho^{a,b}_0)^{-1},u^{a,b}_0)\right)(\widetilde{r}^a_0(x)).
    \end{cases}
\end{align}


 \subsubsection{{\rm \textbf{FBVP}} posed on finite annuli in  radially symmetric  Eulerian coordinates}
First,  let $(\widetilde v^{a,b},\widetilde u^{a,b})$
be a solution to  problem \eqref{approximate problem}--\eqref{initial boundary condition}, and $\widetilde v^{a,b}=(\rho^{a,b})^{-1}$ is the so-called specific volume for the corresponding density function $\rho^{a,b}$. 
Then   we can  introduce the following Lagrangian coordinate transformation:
\begin{equation}\label{mass cordi}
    x(t,r)=\int^r_a\rho^{a,b}(t,s)s^m\ds, \quad t(t,r)=t,
\end{equation}
and $(t,x)$ is the so-called Lagrangian coordinates.
When  $\rho^{a,b}(t,r)$ is bounded away  from zero, the  coordinate transformation \eqref{mass cordi} is invertible. We define 
$
(\widetilde{r}^{a,b}(t,x),t) 
$ be the corresponding inverse transformation of \eqref{mass cordi}, which  naturally satisfies, for any $y\geq 0$, 
\begin{equation}\label{mass cordi for r}
    x(t,\widetilde{r}^{a,b}(t,y))=\int^{ \widetilde{r}^{a,b}(t,y) }_a\rho^{a,b}(t,s)s^m\ds=y,\quad t(t,\widetilde{r}^{a}(t,y))=t,
\end{equation}
and $\widetilde{r}^{a,b}(t,x)$ can be given in $\eqref{approximate particle path}_2$.
\par By virtue of \eqref{approximate particle path}, we have $\widetilde{r}^{a,b}(0,x)=\left( a^n+n\int^x_0\widetilde{v}^{a,b}_0(y)\dy\right)^{\frac1n}$ for $x\in[0,b]$. We denote $\widetilde{r}^{a,b}_0(x)=\widetilde{r}^{a,b}(0,x)$. Since $\partial_x\widetilde{r}^{a,b}_0=\dfrac{\widetilde{v}^{a,b}_0}{(\widetilde{r}^{a,b}_0)^m}>0$, by virtue of the inverse function theorem, we have
\begin{align*}
    (\widetilde{r}^{a,b}_0)^{-1}(r):[a,\widetilde{r}^{a,b}_0(b)]\rightarrow[0,b].
\end{align*}
Thus, we have $\dfrac{r^n-a^n}{n}=\displaystyle\int^{(\widetilde{r}^{a,b}_0)^{-1}(r)}_0\widetilde{v}^{a,b}_0(y)\dy$ for each $r\in[a,\widetilde{r}^{a,b}_0(b)]$. Now we define 
\begin{align*}
    {\widetilde{\rho}^{a,b}_0}(r):=\rho^{a,b}_0\left(\widetilde{r}^a_0\left((\widetilde{r}^{a,b}_0)^{-1}(r)\right)\right),\quad{\widetilde{u}^{a,b}_0}(r):=u^{a,b}_0\left(\widetilde{r}^a_0\left((\widetilde{r}^{a,b}_0)^{-1}(r)\right)\right).
\end{align*}
Next,  based on the transformations \eqref{mass cordi}-\eqref{mass cordi for r}, we can reformulate the approximate \textbf{IBVP} \eqref{approximate problem}--\eqref{initial boundary condition}
in Lagrangian coordinates as a \textbf{FBVP}
posed on a finite annulus in Eulerian coordinates: 
\begin{equation}\label{FBVP}
    \begin{cases}
        \rho^{a,b}_t+(\rho^{a,b} u^{a,b})_r+\dfrac{m\rho^{a,b} u^{a,b} }{r}=0,\\
   \rho^{a,b} (u^{a,b}_t+u^{a,b} u^{a,b}_r )+(A(\rho^{a,b})^\gamma)_r=(2\mu+\lambda)u^{a,b}_{rr}+\kappa\rho^{a,b} \rho^{a,b}_{rrr}\\
   \qquad \qquad \qquad \qquad \qquad \qquad \qquad \quad  +(2\mu+\lambda)\left(\dfrac{mu^{a,b} }{r}\right)_r+\kappa\rho^{a,b}\left(\dfrac{m\rho^{a,b}_r}{r}\right)_r,\\
    \end{cases}
\end{equation}
in $[0,T]\times [a,\widetilde{r}^{a,b}(t,b)]$ with  $m=n-1$ for some constants  $a\in(0,1)$ and $b\in  \mathbb{N}$, and 
\begin{align}
    &(\rho^{a,b}, u^{a,b})(0,r)=({\widetilde{\rho}^{a,b}_0}(r),{\widetilde{u}^{a,b}_0}(r))\qquad \ \ \text{for} \quad r\in [a,\widetilde{r}^{a,b}(0,b)],\label{radilinitial-annuli}\\
    &u^{a,b}(t,a)= u^{a,b}(t,\widetilde{r}^{a,b}(t,b))    =\rho^{a,b}_r(t,a)=\rho^{a,b}_r(t,\widetilde{r}^{a,b}(t,b)   )=0\qquad \text{for} \quad t\in [0,T].\label{bd condition-annuli}
\end{align}

\subsection{Main Strategies }
We now  briefly outline the difficulties encountered in the proof for the global existence of  weak solutions and the non-formation of vacuum states, and also show the corresponding strategies that can be used  to overcome them.

\subsubsection{The uniform upper and lower bound estimates  of $\rho$}\label{The lower bound of rho in main strategy} 
Compared with  \textbf{CNS}, since  \textbf{CNSK}  contains second-order density derivatives in the Korteweg stress tensor,  we can not apply the same argument for deriving the lower and upper bounds of $\rho$  used in  \cite{CHZ,H2,HJ}  directly to our case.  To deal with this issue, we introduce the spatial averaging operator, the domain segmentation  argument, and a maximum value principle argument. First, for a sufficiently large outer radii $b$ in  Lagrangian coordinates, by virtue of fundamental  energy estimates, we can control the upper bound of  $v=\dfrac{1}{\rho}$ in  $[x_*,b]$ for some $x_*>0$; see Lemma \ref{upper bound of v}. The main difficulty is obtain such bound in  $[0,x_*]$. Second, making use of   the logarithmic structure of $v$, we introduce  an auxiliary function: 
\begin{align*}
    Y(t,x)=S_1(v(t,x))+M(t,x)+\delta x,
\end{align*}
where  $\delta>0$ is a suitably large constant, and 
the quantity \(M\) is a linear combination of the spatial average of the modified velocity \(u/r^m\) and that of \(\log v\).
The quantity $S_1(v) $ is given by 
\begin{align*}
     S_1(v)=\int_{\mathscr V}^{v}S'(s)\,ds,
    \qquad
   S'(v)=\frac{\beta}{2v}
   \left(
        1+\sqrt{1-\frac{4\kappa}{\beta^2v^3}}
   \right),
\end{align*}
for $v>\mathscr V$ with some carefully chosen constant $\mathscr V$.
 This construction allows us to cancel the high order derivatives of $\rho$ arising from the Korteweg stress in the time  evolution equation for $Y$; see Lemma \ref{upper bound of v}. 
 
 Roughly speaking, the upper bound of $Y$ controls $\log v$ in the  region $\{v>\mathscr{V}\}\cap[0,x_*]$, which implies that  the upper bound  of $Y$ gives the desired estimates on  $v$ in $[0,x_*]$.  
We then apply the  maximum value principle argument to the time evolution equation for  $Y$, and then  obtain the upper bound of $v$. Specifically, we choose a upper barrier with respect to the corresponding  differential inequality, initial data, and boundary estimates. If $Y$ cross this upper barrier, the inequality of $Y_t$ would lead to a contradiction. Consequently, $Y$ can not cross the barrier, yielding an upper bound for $v$ in $[0,x_*]$, which gives  the upper bound of $v$ for $x\in[0,b]$, and hence the lower bound for $\rho$; see Lemma \ref{upper bound of v}.
    \smallskip

    \subsubsection{Compactness theory}\label{Passing limit to the nonlinear Korteweg stress in sec 2}

A major difficulty in constructing weak solutions is to identify the weak limit of the quadratic capillarity term $\nabla\rho\otimes\nabla\rho$. Under radial symmetry, this tensor becomes  $|\rho_r|^2$ for $r>0$. The available higher-order estimates are not uniform as the inner radius $a\downarrow0$. Consequently, the weak limit of $|\rho_r^a|^2$ may differ from $|\rho_r|^2$ by a defect measure $\mathfrak d$.

To analyze this limit, we consider separately the fluid region ${(t,r):r>\underline r(t)}$, the possibly  vacuum ball $\mathcal V:={(t,r):r<\underline r(t)}$, and the interface ${(t,r):r=\underline r(t)}$. We start from the weak formulation of the momentum equation for the exterior problem \eqref{approximate problem a}–\eqref{bd condition-ex}, obtained in Lemma~\ref{momentum equation} by letting the outer radius $b\to\infty$ with $a>0$ fixed. For every admissible test function $\Psi$ satisfying $\Psi(t,\cdot)\in C_c^\infty([a,\infty))$ and $\Psi(a,t)=0$, the following identity holds:
    \begin{align*}
        \int^\infty_a(\rho^au^a)(t,r)\Psi(t,r) r^m\dr&=\int^\infty_a(\rho^a_0u^a_0)(r)\Psi (0,r)r^m\dr+\int^t_0\int^\infty_a\rho^au^a(\Psi_t+u^a\Psi_r)r^m\dr\ds\\
        &+\int^t_0\int^\infty_a\left(A(\rho^a)^\gamma-\beta\left(u^a_r+\dfrac{mu^a}{r}\right)\right)\left(\Psi_r+\dfrac{m\Psi}{r}\right)r^m\dr\ds\\
        &+\kappa\int^t_0\int^\infty_a\rho^a\rho^a_r\left(\Psi_r+\dfrac{m\Psi}{r}\right)_rr^m\dr\ds\\
        &+\dfrac{\kappa}{2}\int^t_0\int^\infty_a|\rho^a_r|^2(3\Psi_r+\dfrac{m\Psi}{r}) r^m\dr\ds.
    \end{align*}
There is no available estimate to deal with the last term. By virtue of Lemma \ref{last entropy estimate}, $\sup_{t\in[0,T]}\int^\infty_0|\rho^a_r|^2r^m\dr\le C(T)$. So with the help of the Banach-Alaoglu theorem and the Riesz representation theorem, there exists a non-negative finite Radon measure $\nu$ such that 
    \begin{align*}
        \int_Q\phi\mathrm{d}\nu=\lim_{n\rightarrow\infty}\iint_Q|\rho^{a_j}_r|^2\phi r^m\dr\ds\text{ for all }Q\Subset(0,T)\times(0,\infty), \phi\in C^\infty_0(Q).
    \end{align*}
    Thus, $\sqrt{\phi}\rho^{a_j}_r\rightharpoonup\sqrt{\phi}\rho_r$, where $\rho_r$ is specified in Lemma \ref{last convergence}. Then we define the defect measure by $\mathfrak{d}(E)=\nu(E)-\iint_E|\rho_r|^2\dr\ds$, where $E\subset Q$ is a Borel set. The defect measure $\mathfrak{d}$ is actually the weak* limit of the difference of $|\rho^a_r-\rho_r|^2$, \textit{i.e.},
    \begin{align*}
        |\rho^a_r-\rho_r|^2\dr\ds\overset{*}{\rightharpoonup}\mathfrak{d}\text{ in }\mathcal{M}_{\loc}((0,T)\times(0,\infty)),
    \end{align*}
    where $\mathcal{M}$ is the dual space of $C^\infty_0$. The strong convergence established in the fluid region implies that the defect measure is equal to zero in the fluid region $\mathscr{F}$. It therefore remains to prove that it equals to zero in the vacuum region and on the vacuum interface. However, the high-order regularity (independent of inner radii $a$) in these regions prevents a direct application of the compactness arguments. 
    To prove $\mathfrak{d}=0$ on the vacuum interface, we introduce the disintegration theorem for measures to obtain 
    \begin{align*}
        \beta\left(u_r+\dfrac{mu}{r}\right)_r=\dfrac{\kappa}{2}\left(3\mathfrak{d}_r+\dfrac{2m\mathfrak{d}}{r}\right)\text{ in }\mathcal{D}'(\mathcal{V}).
    \end{align*}
By the mollification and truncation arguments, it can be shown for each $\eta\in C^\infty_c(t_1,t_2)$,
    \begin{align*}
    \int_{\{(t,r):t\in(t_1,t_2),r=\underline{r}(t)\}}\eta(t)r^m\mathrm{d}\mathfrak{d}=0,
\end{align*}
Next, we have to prove $\mathfrak{d}=0$ in the vacuum region. To overcome this difficulty, we combine mollification and truncation arguments, arguments of total variation and bounded variation, the Hardy’s inequality; see the proof
of Lemma \ref{momentum equation} for more details.
    \smallskip
    \subsubsection{Non-formation of vacuum states}\label{Non-formation of vacuum states in sec2}
We rule out the formation of vacuum states by the momentum-impulse logarithmic contradiction.  Since the capillary defect  vanished in the vacuum region, then 
\begin{align}
    \left(u_r+\frac{u}{r}\right)_r=0
    \qquad\text{in }\mathcal{D}'(\mathcal{V}).
\end{align}
One has  $u(t,r)=c(t)r$ for $0<r<\underline{r}(t)$. Then, according to  the fact that 
$\underline{r}'(t)=u(t,\underline{r}(t))$, it can be proved that the unweighted viscous impulse is actually equal to $2\beta\log\dfrac{\underline{r}(t)}{\underline{r}(s)}$, where $[s,t]\Subset \{\tau:\underline{r}(\tau)>0\}$. The corresponding truncated
unweighted momentum $\mathscr{M}_{\chi}$ satisfies
\begin{align}
    |\mathscr{M}_{\chi}|
    \le
    C\left(
        1+\log\frac{\underline{r}(t)}{\underline{r}(s)}
    \right)^{3/4};
\end{align}
see Lemma \ref{Ruling out the possible vacuum} for the detailed construction of $\mathscr{M}_\chi$. Furthermore, the remaining terms in the formula of $(\mathscr{M}_\chi)_t$
are bounded by $C(T)$. Since $\{\tau:\underline{r}(\tau)>0\}$ is an open set so that there exists a union of at most countable open interval such that 
\begin{align*}
\{\tau:\underline{r}(\tau)>0\}=\bigcup^\infty_{i=1}(t^i_1,t^i_2).
\end{align*}
Now we fix $t\in(t^i_1,t^i_2)$ and let $s_j$ tend to $t^i_1$, hence $\underline{r}(s_j)\to0$. The momentum-impulse contradiction means the unweighted viscous impulse tends to infinity faster than
the unweighted momentum. This contradiction shows that the vacuum region $\mathcal{V}$ is empty.

\subsection{Mollification and truncation of the initial data in the Eulerian coordinates}\label{3.1}
\begin{lemma}\label{lemma3.2}
    The mollified and truncated initial data $(\rho^a_0,u^a_0)(r)$ satisfy
    \item[(i).]$(\rho^a_0,u^a_0)\in C^\infty([0,\infty)),\quad u^a_0(a)=(\rho^a_0)_r(a)=0$.
    \item[(ii).]   There exists  $C_0>0$ such that 
    \begin{align}\label{3.4}
        &C^{-1}_0\le \rho^a_0(r)\le C_0\text{ for }\forall (a,r)\in (0,1)\times[0,\infty),\nonumber\\
        &\sup_{a\in(0,1)}\int^\infty_a\left(\rho^a_0|u^a_0|^2+G(\rho^a_0)+|((\rho^a_0)_r,\rho^a_0-1)|^2\right)r^m\dr\le C_0,
    \end{align}
    where $C_0$ is independent of $a\in(0,1)$.
    \item[(iii).] As $a\rightarrow0$, we have
    \begin{align}\label{3.5}
        &\Vert \left(\rho^a_0-1\right)-\left(\rho_0-1\right)\Vert_{L^2([0,\infty);r^m\dr)}\rightarrow0,\quad(\rho^a_0,u^a_0)(r)\rightarrow(\rho_0,u_0)(r) \text{ for  } a.e.  r\in [0,\infty).
    \end{align}
\end{lemma}
\begin{proof}
We divide the proof into three  steps.

\textbf{1. Smoothness, upper and lower bounds, and almost every convergence.}
    It is straightforward to check $(\rho^a_0,u^a_0)$ is smooth by the properties of mollifier. By the construction of the mollified and truncated initial data $\rho^a_0$, the $\eqref{3.4}_1$ holds as long as we choose $C_0$ larger than $\widetilde{C}$. Furthermore, since $j_a*\rho_0\rightarrow\rho_0,\chi^a\rightarrow1$ a.e., $\rho^a_0\rightarrow\rho_0,u^a_0\rightarrow u_0$.
    
\textbf{2. The bound of the mollified and truncated initial data.} For $\eqref{3.4}_1$, by Minkowski's inequality, we have
    \begin{align*}
        \Vert(\rho^a_0-1)\Vert^2_{L^2([0,\infty);r^m\dr)}&\le\int_{\mathbb{R}^n}|j_a*(\rho_0-1)|^2\dx\\
        &\le\int_{\mathbb{R}^n}|\rho_0-1|^2\dx=|S^{n-1}|\int^\infty_0|\rho_0-1|^2 r^m\dr\le 2^{m+2}\widetilde{C}.
    \end{align*}
    Next for the term $(\rho^a_0)_r$, by the construction of $\rho^a_0$ we have 
    \begin{align*}
        (\rho^a_0)_r=\chi^a(r)(\overline{\rho}^a_0)_r+(\chi^a(r))'(\bar \rho^a_0-1),
    \end{align*}
    hence
    \begin{align*}
        \Vert(\rho^a_0)_r\Vert_{L^2([0,\infty);r^m\dr)}&=\Vert\chi^a(r)(\overline{\rho}^a_0)_r+(\chi^a(r))'(\overline{\rho}^a_0-1)\Vert_{L^2([0,\infty);r^m\dr)}\\
        &\leq\Vert\chi^a(r)(\overline{\rho}^a_0)_r\Vert_{L^2([0,\infty);r^m\dr)}+\Vert(\chi^a(r))'(\overline{\rho}^a_0-1)\Vert_{L^2([0,\infty);r^m\dr)}\\
        &=I_1+I_2,
    \end{align*}
    For $I_1$,
    \begin{align*}
        (I_1)^2=\Vert\chi^a(r)(\overline{\rho}^a_0)_r\Vert^2_{L^2([0,\infty);r^m\dr)}&\le \int^\infty_0|(\overline{\rho}^a_0)_r|^2r^m\dr=\int_{\mathbb{R}^n}|j_a*\nabla{}\rho_0|^2\dx\\
        &\le|S^{n-1}|\int^\infty_0|(\rho_0)_r|^2r^m\dr\le 2^{m+2}\widetilde{C}.
    \end{align*}
    For $I_2$, it follows from the same argument that
    \begin{align*}
        (I_2)^2&=\Vert(\chi^a(r))'(\overline{\rho}^a_0-1)\Vert^2_{L^2([0,\infty);r^m\dr)}\le \dfrac{4}{a^2}\int^{2a}_a|\overline{\rho}^a_0-1|^2r^m\dr\le 2^{m+2}\widetilde{C},
    \end{align*}
    Combining the above, we have $\Vert(\rho^a_0)_r\Vert^2_{L^2([0,\infty);r^m\dr)}\le 2^{m+3}\widetilde{C}$. Next, by the Jensen inequality and the facts that  $G(\varsigma)$  are convex and  $G(1)=0$, one has  $G(\rho^a_0)\le \chi^aG(\overline{\rho}^a_0)\le j_a*G(\rho_0)$. Thus
    \begin{align*}
        \int^\infty_aG(\rho^a_0)r^m\dr\le2^{m+3}\widetilde{C}.
    \end{align*}
    Next since $|\cdot|^2$ is a convex function, we can derive the following estimate
    \begin{align*}
        \int^\infty_a\left(\rho^a_0|u^a_0|^2\right)r^m\dr&\le \widetilde{C}\int^\infty_0\chi^a|\overline{u}^a_0|^2r^m\dr\le \widetilde{C}\int_{\mathbb{R}^n}|j_a*U_0|^2\mathrm{d}\boldsymbol{x}\\
        &\le \widetilde{C}\int_{\mathbb{R}^n}|U_0|^2\mathrm{d}\boldsymbol{x}\le \widetilde{C}^2|S^{n-1}|\int^\infty_0\rho_0|u_0|^2r^m\dr\le \widetilde{C}^32^{m+2}.
    \end{align*}
    Setting $C_0=2^{m+2}\widetilde{C}^3>\widetilde{C}$, we complete the proof of $\eqref{3.4}_2$.
    
    \textbf{3. The convergence of $L^2$ norm.} For the term $\rho-1$, since
    \begin{align*}
        \Vert(\rho^a_0-1)-(\rho_0-1)\Vert_{L^2([0,\infty);r^mdr)}\le C_0\Vert\rho^a_0-\rho_0\Vert_{L^2([0,\infty);r^mdr)},
    \end{align*}
    we just need to prove $\Vert\rho^a_0-\rho_0\Vert_{L^2([0,\infty);r^mdr)}\rightarrow0$ as $a\rightarrow0$. By the triangle inequality,
    \begin{align*}
        \Vert\rho^a_0-\rho_0\Vert_{L^2([0,\infty);r^mdr)}\le\Vert\chi^a(\bar\rho^a_0-\rho_0)\Vert_{L^2([0,\infty);r^mdr)}+\Vert(1-\chi^a)(\rho_0-1)\Vert_{L^2([0,\infty);r^mdr)}.
    \end{align*}
Since $\chi^a\rightarrow1$ as $a\rightarrow0$ almost everywhere and $\Vert(\rho_0-1)\Vert_{L^2([0,\infty);r^mdr)}\le C_0$, we have the following convergence via the dominated convergence theorem
\begin{align*}
    \Vert(1-\chi^a)(\rho_0-1)\Vert_{L^2([0,\infty);r^mdr)}\rightarrow0 \text{ as } a\rightarrow0.
\end{align*}
Next, Since $\Vert\chi^a(\bar\rho^a_0-\rho_0)\Vert_{L^2([0,\infty);r^mdr)}\le C\Vert j_a*\rho_0-\rho_0\Vert_{L^2(\mathbb{R}^n)}\rightarrow0$, we obtain the desired conclusion.

\end{proof}
\par It is obvious that  $r\rightarrow\displaystyle\int^r_a\rho^a_0s^m\ds$ is strictly increasing. Therefore, by the inverse function theorem, there exists a strictly increasing, bijective, smooth function $\widetilde{r}^a_0(x):[0,\infty)\rightarrow[a,\infty)$ such that 
\begin{align}\label{initial x}
    x=\displaystyle\int^{\widetilde{r}^a_0(x)}_a\rho^a_0r^m\dr \text{ for } \forall x\in[0,\infty).
\end{align}
From now on, we use $\partial_x$ to denote the spatial derivatives of $x$ in Lagrangian coordinates.
\begin{lemma}\label{lemma 3.3}
For all $x\in[0,\infty)$,
\begin{align}\label{initial par path}
    \partial_x \widetilde r^a_0(x)
=\dfrac{1}{
\bigl(\widetilde r^a_0(x)\bigr)^m
\rho^a_0\bigl(\widetilde r^a_0(x)\bigr)},
\qquad\widetilde r^a_0(x)=\left(a^n+n\int_0^x\dfrac{1}{\rho^a_0\bigl(\widetilde r^a_0(y)\bigr)}
\dy
\right)^{\frac1n}.
\end{align}
\begin{proof}
    By \textnormal{Lemma \ref{lemma3.2}}, we have $\rho^a_0(r)\in C^\infty$ and
$C_0^{-1}\le \rho^a_0(r)\le C_0$, thus
\begin{align*}
\lim_{r\to \widetilde r^a_0(x)}
\dfrac{1}{r-\widetilde r^a_0(x)}
\int_{\widetilde r^a_0(x)}^{r}
\rho^a_0(\zeta)\varsigma^m\,\mathrm{d}\varsigma=
\bigl(\widetilde r^a_0(x)\bigr)^m
\rho^a_0\bigl(\widetilde r^a_0(x)\bigr).
\end{align*}
Since
$x\rightarrow \widetilde r^a_0(x)$ is smooth, for $\forall x\in[0,\infty)$ we have
\begin{align*}
\lim_{y\to x}\dfrac{1}{\widetilde r^a_0(y)-\widetilde r^a_0(x)}\int_{\widetilde r^a_0(x)}^{\widetilde r^a_0(y)}\rho^a_0(\zeta)\zeta^m\,\mathrm{d}\zeta=
\bigl(\widetilde r^a_0(x)\bigr)^m
\rho^a_0\bigl(\widetilde r^a_0(x)\bigr).
\end{align*}
Since $x=\displaystyle\int^{\widetilde{r}^a_0(x)}_a\rho^a_0r^m\dr$, we have
\begin{align*}
\dfrac{\widetilde r^a_0(y)-\widetilde r^a_0(x)}{y-x}&\left(\dfrac{1}{\widetilde r^a_0(y)-\widetilde r^a_0(x)}\int_{\widetilde r^a_0(x)}^{\widetilde r^a_0(y)}
\rho^a_0(r)r^m\,\dr
\right)=\dfrac{1}{y-x}
\int_{\widetilde r^a_0(x)}^{\widetilde r^a_0(y)}\rho^a_0(r)r^m\,\dr=1 .
\end{align*}
Thus combining the above results,
\begin{align*}
    \lim_{y\rightarrow x}\dfrac{\widetilde r^a_0(y)-\widetilde r^a_0(x)}{y-x}&=\lim_{y\rightarrow x}\left(\dfrac{1}{\widetilde r^a_0(y)-\widetilde r^a_0(x)}\int_{\widetilde r^a_0(x)}^{\widetilde r^a_0(y)}\rho^a_0(r)r^m\,\dr\right)^{-1}\\
&=\dfrac{1}{\bigl(\widetilde r^a_0(x)\bigr)^m\rho^a_0\bigl(\widetilde r^a_0(x)\bigr)}.
\end{align*}
where we use $0<\rho^a_0(\widetilde{r}^a_0(y))(\widetilde{r}^a_0(y))^m<\infty$. Thus $\eqref{initial par path}_1$ is proved. Lemma $\eqref{initial par path}_2$ can be proved by multiplying $(\widetilde{r}^a_0)^m$ to $\partial_x \widetilde r^a_0(x)
=\dfrac{1}{
\bigl(\widetilde r^a_0(x)\bigr)^m
\rho^a_0\bigl(\widetilde r^a_0(x)\bigr)}$ and integrating in $y\in[0,x]$.

\end{proof}
\end{lemma}

\begin{remark}\label{remark 3.1}
    \begin{align*}
        &(\rho^{a,b}_0,u^{a,b}_0)(r)=
            (1,0)\text{ for all }r\in[\widetilde{r}^a_0(b),\infty),\quad(\rho^{a,b}_0,u^{a,b}_0)=(\rho^a_0,u^a_0)(r)\text{ for }r\in\left[a,\dfrac{\widetilde{r}^a_0(b)}{2}\right),
    \end{align*}
    and $u^{a,b}_0(r)=(\rho^{a,b}_0)_r(r)=0$ for $r=a\text{ and  }\widetilde{r}^a_0(b)$. 
    \par Moreover, one has  $(\rho^{a,b}_0,u^{a,b}_0)(r)\in C^\infty([a,\infty))$. 
\end{remark}
The next Lemma shows that the truncated initial data satisfy the same bound constructed before.
\begin{lemma}\label{lemma 3.4}
    Let  $(\rho^{a,b}_0,u^{a,b}_0)(r)$ be the initial data constructed above. Then there exists a positive constant $C_0$ such that
    \begin{align*}
        &C^{-1}_0\le \rho^{a,b}_0\le C_0\text{ for all } (r,b)\in[a,\infty)\times (\mathscr{B},\infty),\\
        &\sup_{a\in(0,1)}\sup_{b>\mathscr{B}}\int^\infty_a \left\{\left(\rho^{a,b}_0|u^{a,b}_0|^2+G(\rho^{a,b}_0)\right)+\left|(\rho^{a,b}_0)_r\right|^2+|\rho^{a,b}_0-1|^2\right\}r^m\dr\le C_0.
    \end{align*}
    \begin{proof}
       The results follow directly from Lemma \ref{lemma3.2} and the construction.
    \end{proof}
\end{lemma}
By Lemma \ref{lemma 3.3} and Remark \ref{remark 3.1}, it follows that
\begin{align*}
\left(\widetilde{v}^{a,b}_0,\widetilde{u}^{a,b}_0\right)(x)\in C^\infty([0,\infty))
,\quad\widetilde{u}^{a,b}_0=\partial_x(\widetilde{v}^{a,b}_0)=0 \text{ at }x=0,b.
\end{align*}

Furthermore, we have
\begin{lemma}\label{energy estimates in L}
    There exists a positive constant $C_0$ such that
    \begin{align}
        \sup_{a\in(0,1)}\sup_{b>\mathscr{B}}&\int^b_0\left(\dfrac{1}{2}|\widetilde{u}^{a,b}_0|^2+\psi(\widetilde{v}^{a,b}_0)+\left|\widetilde{v}^{a,b}_0-1\right|^2\right)\dx\le C_0.
    \end{align}
    \begin{proof}
        First,  it implies from Lemma  \ref{lemma 3.3} that,
        for each $x\in[0,\infty)$,
        \begin{align*}
            0<C^{-1}_0(\widetilde{r}^a_0(x))^{-m}\le|\partial_x\widetilde{r}^a_0(x)|\le(\widetilde{r}^a_0(x))^{-m}(\rho^a_0(\widetilde{r}^a_0))^{-1}\le C_0a^{-m}<\infty.
        \end{align*}
        Then for any locally integrable function $f(r):[a,\infty)\rightarrow \mathbb{R}$, we have
        \begin{align*}
            \int^b_0 f(\widetilde{r}^a_0(x))\dx=\int^{\widetilde{r}^a_0(b)}_af(r)\rho^a_0(r)r^m\dr.
        \end{align*}
        We provide the details only for the first term; the remaining terms, except for the one involving the convex function, can be treated by the same argument.
        \begin{align*}
            \int^b_0|\widetilde{u}^{a,b}_0|^2\dx=\int^{\widetilde{r}^a_0(b)}_a\dfrac{\rho^a_0}{\rho^{a,b}_0}\rho^{a,b}_0|\widetilde{u}^{a,b}_0|^2r^m\dr\le C_0.
        \end{align*}
        For the convex function term, since $\psi\left(v\right)\dfrac{1}{v}=G\left(\dfrac{1}{v}\right)$, 
        \begin{align*}
            \int^b_0\psi(\widetilde{v}^{a,b}_0)\dx=\int^{\widetilde{r}_0^{a}(b)}_a\dfrac{\rho^a_0}{\rho^{a,b}_0}G(\rho^{a,b}_0)r^m\dr\le C_0.
        \end{align*}

    \end{proof}
\end{lemma}
\begin{remark}\label{ particle path remark}
    The equation of approximate problems \eqref{approximate problem} is derived from 
    \begin{align}
        x(t,r)=\int^r_a \rho(t,r)r^m\dr,\quad t(t,r)=t
    \end{align}
    The corresponding solutions in Lagrangian coordinates are defined as
    \begin{align*}
        \left(v(t,x),u(t,x)\right)=\left((\rho(t,r(t,x)))^{-1},u(t,r(t,x))\right).
    \end{align*}
    hence the equation  could be transformed from \eqref{1.6} to \eqref{approximate problem} with the following relation:
    \begin{align}
\dfrac{\partial(x,t)}{\partial(r,t)}
&=\begin{pmatrix}
 \rho r^m & -\rho u r^m \\
0 & 1
\end{pmatrix},
\qquad\dfrac{\partial(r,t)}{\partial(x,t)}=\begin{pmatrix}
v r^{-m} & u \\
0 & 1
\end{pmatrix},
\end{align}
and for particle path $r(t,x)$, we have 
\begin{align}\label{constructure of particle path}
&r(t,x)
=r(0,x)+\int_0^t u(s,x)\ds, \quad r(0,x)
:=\left(a^n+n\int_0^x \rho_0^{-1}(r_0(y))\,\dy
\right)^{1/n},\nonumber\\
&r^n(t,x)=a^n+n\int_0^x v(t,y)\dy.
\end{align}
Once $0<\rho<\infty$ has been established, the coordiante map  $x\rightarrow r(t,x)$ is bijective and has a positive Jacobian determinant. We could transform back to the Eulerian coordinates, this is the motivation.
\par Specifically, for the initial data of \textbf{IBVP} \eqref{approximate problem}, it holds that
\begin{align*}
    \widetilde{r}^{a,b}(0,x)=\left(a^n+n\int^x_0\dfrac{1}{\rho^{a,b}_0(\widetilde{r}^a_0(y))}dy\right)^{\frac1n}.
\end{align*}
With the help of the upper and lower bound of the truncated $\rho^{a,b}$, we have
\begin{align*}
    a^n+\dfrac{nx}{C_0}\le(\widetilde{r}^{a,b}(0,x))^n\le a^n+nC_0x.
\end{align*}
Thus, for the mollified initial particle path and the truncated particle path, we have
\begin{align*}
    &\dfrac{(\widetilde{r}^{a,b}(0,x))^n}{(\widetilde{r}^a_0(x))^n}\le\dfrac{a^n+nC_0x}{a^n+\dfrac{nx}{C_0}}\le\dfrac{C^2_0a^n+nC_0x}{a^n+\dfrac{nx}{C_0}}=C^2_0,\quad\dfrac{(\widetilde{r}^{a,b}(0,x))^n}{(\widetilde{r}^a_0(x))^n}\ge\dfrac{a^n+\dfrac{nx}{C_0}}{a^n+nC_0x}\ge C^{-2}_0.
\end{align*}
This relation is important to obtain the truncated first order derivative.
\end{remark}
Based on the analysis in Remark \ref{ particle path remark}, we have
\begin{lemma}
There exists a positive constant $C_0$ such that
    \begin{align}
        \int^b_0|(\widetilde{r}^{a,b}_0)^m\partial_x\widetilde{v}^{a,b}_0|^2\dx\le C_0.
    \end{align}
    \begin{proof}
\begin{align*}
    \int^b_0|(\widetilde{r}^{a,b}_0)^m\partial_x\widetilde{v}^{a,b}_0|^2\dx&\le C_0^{\frac{4m}{n}}\int^{\widetilde{r}^{a}_0(b)}_a\dfrac{(\rho^{a,b}_r)^2}{(\rho^{a,b})^4\rho^a_0}r^m\dr\le C_0,
\end{align*}
        where we have used $C^{-1}_0\le\rho^{a,b}_0\le C_0$ mentioned in Lemma \ref{lemma 3.4}.
        \end{proof}
\end{lemma}

The local well-posedness of strong solutions to \textbf{IBVP} \eqref{approximate problem} can be proved by a linearization argument. The argument is standard; see \cite{HL1, PB}. Therefore, we only state the results and omit details here.
\begin{theorem}
For  the smooth initial data   $(\widetilde{v}^{a,b}_0,\widetilde{u}^{a,b}_0)(x)$ constructed above for each $x\in[0,b]$,  there exist a $T_*>0$ and a  unique strong solution $(\widetilde{v}^{a,b},\widetilde{u}^{a,b})(t,x)$ in  $[0,T_*]\times [0,b]$ to  \eqref{approximate problem}--\eqref{initial boundary condition} such that
\begin{align}
    &\widetilde{v}^{a,b}\in C([0,T_*];H^{3}(0,b))\cap L^2(0,T_*;H^4(0,b)),\quad C^{-1}(a,T_*)\le \widetilde{v}^{a,b}\le C(a,T_*),\nonumber\\
     &\widetilde{u}^{a,b}\in C([0,T_*];H^{2}(0,b))\cap L^2(0,T_*;H^3(0,b)),\label{abregularity}\\
    &\partial_t\widetilde{v}^{a,b}\in C([0,T_*];H^{1}(0,b)),\quad \partial_t\widetilde{u}^{a,b}\in C([0,T_*];L^2(0,b)).\nonumber
\end{align}
\end{theorem}

\section{Uniform A Priori Estimates }\label{3}
In this section, we omit the superscripts of approximation parameters$(a,b)\in(0,1)\times \mathbb{N}$ for simplicity. Also, we denote the solution pair $(v,u)\equiv(\widetilde{v}^{a,b},\widetilde{u}^{a,b})$.
\subsection{Main results}
Now  we construct global-in-time strong solutions to the approximate problems \eqref{approximate problem}--\eqref{initial boundary condition}. First, we specify the lower bound of $b$ as follows:
\begin{align*}
    b>4+2\max\left\{\dfrac{64(C_P+\kappa C_\mathscr{A}+C_{rest})^2}{A^2},4+4\max\left\{\dfrac{C_0}{\psi(\frac12)},\dfrac{2R^n_*}{n}\right\}\right\}:=\mathscr{B},
\end{align*}
where $R^n_*,C_P,C_\mathscr{A},C_{rest}$ will be specified later in our proof (see Lemma \ref{lower bound of v}).

\begin{theorem}\label{global existence of strong solutions}
Let the initial data and boundary conditions be as specified in \eqref{initial boundary condition}. Then for any $T>0$, there  exists a  unique strong solution $(v^{a,b},u^{a,b})$ in $[0,T]\times [0,b]$ to  problem \eqref{approximate problem} for each $b>\mathscr{B}$, which satisfies \eqref{abregularity} with $T_*$ replaced  by $T$. Moreover the following estimates hold for each $T>0$:
    \begin{enumerate}
        \item There exists a positive constant $C_0$ such that
    \begin{align}
    \sup_{t\in[0,T]}&\int^b_0\left(|u|^2+\psi(v)+\dfrac{\kappa r^{2m}|\partial_x v|^2}{2v^5}\right)\dx+\beta\int^T_0\int^b_0\dfrac{|\partial_x(r^mu)|^2}{v}\dx\ds\le C_0.
    \end{align}
    \item For each $a\in(0,1)$, there exists a positive constant $C(a,T)$ such that $C^{-1}(a,T)\le v(t,x)\le C(a,T)$ for all $(t,x)\in[0,T]\times[0,b]$.
    \item There exists a positive constant $C(a,T)$ such that
    \begin{align*}
    &\sup_{t\in[0,T]}\int^b_0\left(|v-1|^2+\dfrac{|r^m\partial_xu|^2}{v}+\left|\partial_x\left(\dfrac{r^m\partial_xv}{v^{\frac52}}\right)\right|^2\dfrac{r^{2m}}{v}+|u|^4+\left|\dfrac{r^m\partial_xv}{v^{\frac52}}\right|^4\right)\dx\\
    &+\int^T_0\int^b_0\left(\dfrac{r^{4m}|\partial_{xx}v|^2}{v^7}+\dfrac{r^{4m}|\partial_{x}v|^4}{v^9}+\dfrac{r^{2m-2}|\partial_{x}v|^2}{v^5}\right)\dx\dt\\
    &+\int^T_0\int^b_0 \left(|\partial_tv|^2+|r^m\partial_{xt}v|^2+|\partial_t u|^2+\dfrac{r^{4m}}{v^2}\left(|\partial_{xx}u|^2+\left|\partial_{xx}\left(\dfrac{r^m\partial_xv}{v^{\frac52}}\right)\right|^2\right)\right)\dx\dt\le C(a,T).
    \end{align*}
    \item There exists a positive constant $C(a,T)$ such that
    \begin{align*}
        &\sup_{t\in[0,T]}\int^b_0\left(|\partial_tu|^2+\left|\partial_t\left(\dfrac{r^m\partial_xv}{v^{\frac52}}\right)\right|^2\right)\dx+\int^T_0\int^b_0\left(|\partial_{xt}u|^2+\left|\partial_{xt}\left(\dfrac{r^m\partial_xv}{v^{\frac52}}\right)\right|^2\right)\dx\dt\\
        &+\sup_{t\in[0,T]}\left(\Vert v-1\Vert^2_{H^3_x(0,b)}+\Vert u\Vert^2_{H^2(0,b)}\right)+\int^T_0\left(\Vert v-1\Vert^2_{H^4(0,b)}+\Vert u\Vert^2_{H^3(0,b)}\right)\dt\le C(a,T)
    \end{align*}
    \end{enumerate}
\end{theorem}
\begin{remark}
    It is worth noting that the estimates established  above allow us to establish the global existence of strong solutions to \eqref{approximate problem}--\eqref{initial boundary condition}. The proof of \textnormal{Theorem} \ref{global existence of strong solutions} is similar to the proof of \textnormal{Theorem} \ref{global existence of strong solutions independent of a}, except for the last estimate. The last estimate and \eqref{approximate problem} are crucial to ensure the spatial-time continuity of $Y_{xx}$ and $Y_t$, see \textnormal{Lemma} \ref{lower bound of v}. Since the last estimate can be proved by a standard energy method, we omit its proof here for simplicity. In \textnormal{Theorem} \ref{global existence of strong solutions}, high-order estimates  are allowed to depend on $a$ and the corresponding high-order norms of initial data \eqref{initial boundary condition}  are finite. This global existence of approximate solutions allows us to derive uniform estimates independent of $a$. The time-weighted  estimates presented below are used to establish some uniform estimates independent of $a$. 
\end{remark}
     In addition, we obtain some uniform estimates independent of inner radius $a$,
    \begin{theorem}\label{global existence of strong solutions independent of a}
The global-in-time strong solution $(v^{a,b},u^{a,b})$ satisfies the following estimates for each $T>0$:
    \begin{enumerate}
        \item There exists a positive constant $C_0$ such that
    \begin{align}
    \sup_{t\in[0,T]}&\int^b_0\left(|u|^2+\psi(v)+\dfrac{\kappa r^{2m}|\partial_x v|^2}{2v^5}\right)\dx+\beta\int^T_0\int^b_0\dfrac{|\partial_x(r^mu)|^2}{v}\dx\ds\le C_0.
    \end{align}
    \item For each $\varepsilon>0$, there exists a positive constant $C(\varepsilon,T)$ such that $C^{-1}(\varepsilon,T)\le v(t,x)\le C(\varepsilon,T)$ for all $(t,x)\in[0,T]\times[\varepsilon,b]$.
    \item There exists a positive constant $C(\varepsilon,T)$ such that
    \begin{align*}
    &\sup_{t\in[0,T]}\int^b_\varepsilon\left(|v-1|^2+\sigma\dfrac{|r^m\partial_xu|^2}{v}+\sigma\left|\partial_x\left(\dfrac{r^m\partial_xv}{v^{\frac52}}\right)\right|^2\dfrac{r^{2m}}{v}+\sigma|u|^4+\sigma\left|\dfrac{r^m\partial_xv}{v^{\frac52}}\right|^4\right)\dx\\
    &+\int^T_0\int^b_\varepsilon\left(\dfrac{r^{4m}|\partial_{xx}v|^2}{v^7}+\dfrac{r^{4m}|\partial_{x}v|^4}{v^9}+\dfrac{r^{2m-2}|\partial_{x}v|^2}{v^5}\right)\dx\dt\\
    &+\int^T_0\int^b_\varepsilon \left(|\partial_tv|^2+\sigma|r^m\partial_{xt}v|^2+\sigma|\partial_t u|^2+\sigma\dfrac{r^{4m}}{v^2}\left(|\partial_{xx}u|^2+\left|\partial_{xx}\left(\dfrac{r^m\partial_xv}{v^{\frac52}}\right)\right|^2\right)\right)\dx\dt\le C(\varepsilon,T).
    \end{align*}
    \item $\sigma^{\frac14}(t)(|u(t,x)|+|r^m\partial_xv(t,x)|)\le C(\varepsilon,T)$ for $x\ge\varepsilon$.
    \end{enumerate}
\end{theorem}
The proof of \textnormal{Theorem} \ref{global existence of strong solutions independent of a} is divided into several Lemmas in the rest of the section. 
\subsection{Basic Energy Estimate}\label{sec4.1}
Here we give  some basic energy   estimates.
\begin{lemma}\label{entropy estimate}
There exists a positive constant $C_0$ such that
    \begin{align}
        \sup_{t\in[0,T]}&\int^b_0\left(\dfrac{1}{2}|u|^2+\psi(v)+\dfrac{\kappa r^{2m}|\partial_x v|^2}{2v^5}\right)\dx+\beta\int^T_0\int^b_0\dfrac{|\partial_x(r^mu)|^2}{v}\dx\ds\le C_0.
    \end{align}
    \begin{proof}
        The proof is finished by multiplying  the equation $\eqref{approximate problem}_2$ by $u$ and integrating by parts. 
    \end{proof}
\end{lemma}
Now we give the upper and lower bounds of the particle path.
\begin{lemma}\label{particle path for ab}
    Let $\psi^{-1}_-(\cdot)$ be the left branch inverse of $\psi$. Then 
    \begin{align}
        a^n+nx\psi^{-1}_-\left(\dfrac{C_0}{x}\right)\le\left[\widetilde{r}^{a,b}(t,x)\right]^n\le C_0(1+x)\text{ for }(t,x)\in[0,T]\times[0,b].
    \end{align}
    \begin{proof}
        We divide the proof into 4 steps.\\
        \textbf{1. } Since $\psi(1)=0$ and $\psi$ is strictly convex in $(0,\infty)$, there exists a strictly decreasing left inverse function $\psi^{-1}_-:[0,\infty)\rightarrow(0,1]$  and a strictly increasing right inverse function $\psi^{-1}_+:[0,\infty)\rightarrow[1,\infty)$. By the Jensen inequality, we conclude that
        \begin{align}\label{left branch}
            \psi\left(\dfrac{[\widetilde{r}^{a,b}(t,x)]^n-a^n}{nx}\right)=\psi\left(\dfrac{1}{x}\int^x_0\widetilde{v}^{a,b}(t,y)\dy\right)\le\dfrac{1}{x}\int^x_0\psi\left(\widetilde{v}^{a,b}(t,y)\right)\dy\le\dfrac{C_0}{x}.
        \end{align}
      Then, we split the proof into two cases: $0<\dfrac{[\widetilde{r}^{a,b}(t,x)]^n-a^n}{nx}\le1$ and $\dfrac{[\widetilde{r}^{a,b}(t,x)]^n-a^n}{nx}\ge1$.\\
        \textbf{2.} For the first case $0<\dfrac{[\widetilde{r}^{a,b}(t,x)]^n-a^n}{nx}\le1$,  we have 
        \begin{align*}
            \dfrac{[\widetilde{r}^{a,b}(t,x)]^n-a^n}{nx}=\psi^{-1}_-\left(\psi\left(\dfrac{[\widetilde{r}^{a,b}(t,x)]^n-a^n}{nx}\right)\right).
        \end{align*}
        According to the strictly decreasing of $\psi^{-1}_-$ and \eqref{left branch}, it follows that
        \begin{align*}
            \dfrac{[\widetilde{r}^{a,b}(t,x)]^n-a^n}{nx}=\psi^{-1}_-\left(\psi\left(\dfrac{[\widetilde{r}^{a,b}(t,x)]^n-a^n}{nx}\right)\right)\ge\psi^{-1}_-\left(\dfrac{C_0}{x}\right),
        \end{align*}
        on the other hand, since $\dfrac{[\widetilde{r}^{a,b}(t,x)]^n-a^n}{nx}\le1$, the upper bound is also proved.\\ 
        \textbf{3.} For the second case: $\dfrac{[\widetilde{r}^{a,b}(t,x)]^n-a^n}{nx}\ge1$. Since $\psi^{-1}_-(z)\le1$ for $\forall z\in[0,\infty)$, it follows that
        \begin{align*}
            [\widetilde{r}^{a,b}(t,x)]^n\ge a^n+nx\ge a^n+nx\psi^{-1}_-\left(\dfrac{C_0}{x}\right).
        \end{align*}
        For the upper bound in this case, it holds that
        \begin{align*}
            \int^x_0\psi(v)\dx=\dfrac{A}{\gamma-1}\int^x_0\left(v^{1-\gamma}+(\gamma-1)v-\gamma\right)\dx=\dfrac{A}{n}[r^n(t,x)-a^n]+\dfrac{A}{\gamma-1}\int^x_0v^{1-\gamma}\dx-\dfrac{A\gamma x}{\gamma-1}.
        \end{align*}
        Combining the above we have
        \begin{align*}
            r^n(t,x)\le C_0(1+x)\text{ for each }x\in[0,b].
        \end{align*}
        \textbf{4. } For $x=0$, based on the fact $\psi^{-1}_-\le 1$, we extend $x\psi^{-1}_-\left(\dfrac{C_0}{x}\right)$ to $x=0$ by $\lim_{x\rightarrow0}x\psi\left(\dfrac{C_0}{x}\right)=0$, thus the upper and lower bound of $\widetilde{r}(t,0)$ can be proved straightly.
        \par The proof is completed now.
    \end{proof}
\end{lemma}
Then we can obtain estimates at some point by the mean-value theorem.
\begin{lemma}\label{point estimate}
    For $\forall i=1, \cdots,b$, and for any fixed $t\in[0,T]$, there exists $B_i(t)$ such that
    \begin{align*}
        C^{-1}_0\le v(t,B_i(t))\le C_0.
    \end{align*}
    \begin{proof}
        Since $\psi$ is convex, it follows  that
    \begin{align*}
        \psi\left(\int^i_{i-1}v(t,x)\dx\right)\le\int^i_{i-1}\psi\left(v(t,x)\right)\dx\le C_0,
    \end{align*}
    where we have used the Jensen inequality and entropy estimates mentioned in Lemma \ref{entropy estimate}. Furthermore, we obtain that
    \begin{align*}
        0<\psi^{-1}_-\left(C_0\right)\le\int^i_{i-1}v(t,x)\dx\le\psi^{-1}_+\left(C_0\right)<\infty.
    \end{align*}
    Now by the mean value theorem, there exists a point $B_i(t)\in[i-1,i]$, for each fixed $t\in[0,T]$, such that
    \begin{align*}
        0<\psi^{-1}_-\left(C_0\right)\le v(t,B_i(t))\le\psi^{-1}_+\left(C_0\right)<\infty.
    \end{align*}
    Thus this proof is finished.
    \end{proof}
\end{lemma}

\subsection{Upper and Lower Bounds of density}\label{point wise lower bound of v}
Now we give the proof for the lower bound for $v$.
\begin{lemma}\label{lower bound of v}
For any  $(t,\varepsilon)\in(0,1)\times[0,T]$, there exist  positive constants $C_0,C(\varepsilon,T),C(T)$ such that
\begin{align}\label{vupperestimates}
       & C^{-1}(T)\le v(t,x) \quad \quad \text{ for } x\ge 1,\nonumber\\
       & C^{-1}(\varepsilon,T)\le v(t,x)\quad \text{ for } x\ge\varepsilon,  \\
    &C^{-1}_0a^{\frac{2m}{3}}\le v(t,x)\quad  \ \ \text{ for } x\in[0,b].\nonumber
    \end{align}
    \begin{proof}
        The proof is  based on the Sobolev embedding and we only show the proof for the last lower bound. For  $\forall(t,x)\in[0,T]\times[0,b]$, there exists some $B_i(t)$ satisfying $|x-B_i(t)|\le1$ such that
        \begin{align*}
           \left| \left(\dfrac{1}{v(t,x)}\right)^{\frac{3}{2}}-\left(\dfrac{1}{v(t,B_i(t))}\right)^{\frac{3}{2}}\right|&\le\int^{\max\{x, B_i(t)\}}_{\min\{x, B_i(t)\}}\left|\partial_y\left(\dfrac{1}{v(t,y)}\right)^{\frac32}\right|dy\\
            &\le C\left(\int^b_0\left|\dfrac{r^m\partial_y v}{v^{\frac{5}{2}}}\right|^2dy\right)^{\frac12}|x-B_i(t)|^{\frac12}\sup_{y\in[0,b]}\dfrac{1}{r^m(t,y)}\\
            &\le C_0\sup_{y\in[0,b]}\dfrac{1}{r^m(t,y)}.
        \end{align*}
        Based on the lower bound of particle path $r(t,x)$, it holds that
        \begin{align*}
            r^n(t,x)\ge a^n+nx\psi^{-1}_-\left(\dfrac{C_0}{x}\right)\ge a^n,
        \end{align*}
        Thus one has \eqref{vupperestimates}.
    \end{proof}
    \begin{remark}
        Note that the bound: $C^{-1}_0a^{\frac{2m}{3}}\le v(t,x)$ for  $x\in[0,b]$ is crucial later to extend the $\rho^a$ defined in $[a,\infty)$ as a constant to the undefined domain $(0,a)$.
    \end{remark}
Establishing an upper bound for $v$ requires a more delicate argument. Since higher-order derivatives appear in integrating by parts, Kazhikhov’s method is no longer directly applicable. To overcome this difficulty, we introduce a modified effective viscous flux and construct auxiliary functions to cancel the higher-order derivative terms. We then combine integral averaging with a maximum principle argument to obtain the desired upper bound.
\begin{lemma}\label{upper bound of v}
    For $\forall\varepsilon\in(0,1)$, there exists  positive constants $C(T),C(\varepsilon,T), C(a,T)$ such that the following estimates hold for all $t\in[0,T]$\\
\begin{align*}
        \begin{cases}
        v(t,x) \le C(T)&\text{ for } x\ge 1,\\
         v(t,x)\le C(\varepsilon,T)&\text{ for } x\ge\varepsilon, \\
     v(t,x)\le C(a,T)&\text{ for } x\in[0,b].
    \end{cases}
\end{align*}
    \begin{proof}
        We divide the proof into 10 steps.\\
        \textbf{1. The upper bound away from the origin.} In this step, $C_0$ is a fixed positive constant. Since $\psi(v)\backsim v$ as $v\rightarrow\infty$, there exists $\mathscr{V}_0$ such that $\psi(v)\ge\dfrac{A}{2}v$ for each $v\ge\mathscr{V}_0$.
    Now we define two positive constants $R^n_*>1+\max\left\{2^n,\dfrac{2nC_0}{A(1-\frac{1}{2^n})},\left(2^{m+1}C_0\mathscr{V}_0\sqrt{\dfrac{2}{A}}\right)^{\frac nm}\right\}$, $x_*=1+\max\left\{\dfrac{C_0}{\psi(\frac{1}{2})},\dfrac{2R^n_*}{n}\right\}$, and
    \begin{align*}
        K(v)=\int^v_{\mathscr{V}_0}\dfrac{\sqrt{\psi(s)}}{s^{\frac{5}{2}}}ds.
    \end{align*} 
    For  $\forall x\ge x_*$, there must exist a point $x_0>0$ such that $r(t,x_0)=\dfrac{r(t,x)}{2}$. If $v\ge\mathscr{V}_0$ everywhere on the region $[x_0,x]$, then
    \begin{align*}
        \dfrac{A}{2}\int^x_{x_0}v\dy\le\int^b_0\psi(v)\dy\le C_0,
    \end{align*}
    Thus    
    \begin{align*}
        \dfrac{A}{2n}r^n(t,x)(1-\dfrac{1}{2^n})=\dfrac{A}{2n}[r^n(t,x)-r^n(t,x_0)]\le C_0,
    \end{align*}

Thus the contradiction appears and there must exist a point $y\in[x_0,x]$ such that
\begin{align*}
    r(t,x_0)\le r(t,y)\le r(t,x),\quad v(t,y)\le \mathscr{V}_0.
\end{align*}
Now by the Newton--Leibniz formula
\begin{align*}
    K(v(t,x))-K(v(t,y))&=\int^x_y\partial_zK(v(t,z))\dz=\int^x_yK'(v(t,z))\partial_zv \dz\\
    &=\int^x_y\dfrac{\sqrt{\psi(v)}}{v^{\frac{5}{2}}}\partial_zv \dz
    \le C_0\sup_{z\in[y,x]}\dfrac{1}{r^m(t,z)}.
\end{align*}
Since $K(v(t,y))\le0$, we have $K(v(t,x))\le C_0 \sup_{z\in[y,x]}\dfrac{1}{r^m(t,z)}$. Next since $r^n(t,x)>R^n_*$ and
\begin{align*}
    K(2\mathscr{V}_0)=\int^{2\mathscr{V}_0}_{\mathscr{V}_0}\sqrt{\dfrac{\psi(s)}{s}} \dfrac{1}{s^2}\ds\ge\sqrt{\dfrac{A}{2}}\int^{2\mathscr{V}_0}_{\mathscr{V}_0}\dfrac{1}{s^2}\ds=\sqrt{\dfrac{A}{2}}\dfrac{1}{2\mathscr{V}_0}.
\end{align*}
We have 
\begin{align*}
    \sup_{z\in[y,x]}\dfrac{C_0}{r^m(t,z)}\le \dfrac{C_02^m}{R^m_*}<\sqrt{\dfrac{A}{2}}\dfrac{1}{2\mathscr{V}_0}.
\end{align*}
Thus $K(2\mathscr{V}_0)>K(v(t,x))$. By virtue of the monotonicity of $K$, we have $v(t,x)<2\mathscr{V}_0$  for each $x\ge x_*$. Then we just need to prove the upper bound of $v$ in $[0,x_*]$.\\
\textbf{2.Reformulation of the momentum equation.} In this step, we introduce the averaging operator and reformulate the momentum equation. Now we define the effective viscous flux $\mathscr{E}$
\begin{align*}
    \mathscr{E}=\beta\dfrac{\partial_tv}{v}-(P(v)-P(1))-\kappa\partial_x(\dfrac{r^{2m}\partial_xv}{v^5})-\dfrac{5\kappa}{2}\dfrac{r^{2m}|\partial_xv|^2}{v^6}.
\end{align*}
Thus, it holds that
\begin{align}\label{effective flux equation}
    \partial_x\mathscr{E}=\partial_t\left(\dfrac{u}{r^m}\right)+\underline{mr^{-n}u^2+m\kappa r^{m-1}\dfrac{|\partial_xv|^2}{v^5}}_{:=\clubsuit},
\end{align}
where $\clubsuit\ge0$. 
\par For  $\forall s\in(0,1)$, let $L=L_0+x$ and $y_s=x+Ls$. Now we define the cut-off function $\omega$ by 
\begin{align*}
    0\le\omega\in C^\infty_c\left(\dfrac{1}{4},\dfrac{3}{4}\right)\text{ with }\int^1_0\omega(s)ds=1,
\end{align*}
where $L_0=2+\max\left\{4\left(1+\max\left\{\dfrac{C_0}{\psi(\frac{1}{2})},\dfrac{2R^n_*}{n}\right\}\right),\left(\dfrac{8(C_P+\kappa C_{\mathscr{A}}+C_{rest})}{A}\right)^2\right\}$ and $C_P,C_{\mathscr{A}},C_{rest}$ will be specified later. Let $\displaystyle W(s)=\int^1_s\omega(\sigma)\mathrm{d}\sigma$. We define the averaging operator by
\begin{align*}
    \mathscr{A}f(t,x)=\int^1_0\omega(s)f(y_s,t)ds.
\end{align*}
Thus for each $x\in[0,x_*]$, we integrate \eqref{effective flux equation} from $x$ to $y_s$, then multiply by $\omega(s)$ and integrate from $0$ to $1$ to obtain
\begin{align*}
    \mathscr{A}\mathscr{E}-\mathscr{E}&=L\int^1_0 W(s)\partial_x \mathscr{E}(t,y_\sigma)\mathrm{d}\sigma\\
    &=L\int^1_0W(s)\partial_t\left(\dfrac{u}{r^m}\right)(t,y_s)\ds+L\int^1_0W(s)\clubsuit(t,y_s)\ds\\
    &=\partial_t\left(L\int^1_0W(s)\dfrac{u}{r^m}\ds\right)+L\int^1_0W(s)\clubsuit(t,y_s)\ds
\end{align*}
Denote 
$$M=L\displaystyle\int^1_0W(s)\left(\dfrac{u}{r^m}\right)(t,y_s)ds-\beta\mathscr{A}\log v.$$
Then we have the reformulated effective viscous flux equation
\begin{align}\label{reformulated effective flux equation}
    \partial_t\left(\beta\log v+M\right)&=-L\int^1_0W(s)\clubsuit(t,y_s)\ds-\mathscr{A}(P(v)-P(1))+(P(v)-P(1))+\kappa\partial_x\left(\dfrac{r^{2m}\partial_xv}{v^5}\right)\nonumber\\
    &+\dfrac{5\kappa}{2}\left(\dfrac{r^{2m}|\partial_xv|^2}{v^6}\right)-\kappa\mathscr{A}\left(\partial_x\left(\dfrac{r^{2m}\partial_xv}{v^5}\right)\right)-\dfrac{5\kappa}{2}\mathscr{A}\left(\dfrac{r^{2m}|\partial_xv|^2}{v^6}\right).
\end{align}
\textbf{3. Estimates for $M$.} We first estimate  $\displaystyle L\int^1_0W(s)\left(\dfrac{u}{r^m}\right)(t,y_s)\ds$
\begin{align*}
    \left|L\int^1_0W(s)\left(\dfrac{u}{r^m}\right)(t,y_s)\ds\right|=\left|\int^{x+L}_xW\left(\dfrac{y-x}{L}\right)\dfrac{u}{r^m}(t,y)\dy\right|\le a^{-m}\sqrt{L}C^{\frac12}_0.
\end{align*}
Based on the fact that $C^{-1}_0 a^{\frac{2m}{3}}\le v$,
it follows that $\psi(v)=\bar\psi(e^z)=\dfrac{A}{\gamma-1}\left(e^{z(1-\gamma)}+(\gamma-1)e^z-\gamma\right)$, where $z=\log v$. Thus $\bar\psi(0)=\bar\psi'(0)=0$. Furthermore, it holds that $\bar\psi''(z)=A(e^z+(\gamma-1)e^{z(1-\gamma)})$. Since $\bar\psi''(z)\ge1$ as $z\ge0$ and $\bar\psi''(z)\ge(\gamma-1)$ as $z\le0$, thus for each $z\in\mathbb{R}$, we have $\bar\psi(z)\ge \dfrac{A}{2}\min\{1,\gamma-1\}z^2$. Returning  to $v$, we have
\begin{align*}
    \psi(v)\ge\dfrac{A}{2}\min\{1,\gamma-1\}|\log v|^2.
\end{align*}
Thus for the second term of $M$, we have
\begin{align*}
    |\mathscr{A}(\log v)|&=\left|\int^1_0\omega(s)\log v(t,x+Ls)\ds\right|
    =\dfrac{1}{L}\left|\int^{x+L}_x\omega\left(\dfrac{y-x}{L}\right)\log v(t,y)\dy\right|\\
    &\le \Vert \omega \Vert_{L^2} L^{-\frac{1}{2}}\left(\int^b_0|\log v(t,y)|^2\dy\right)^{\frac12}
    \le C\dfrac{\left(\int^b_0\psi(v)\dx\right)^{\frac12}}{\min\{1,\gamma-1\}\sqrt{L}}
    \le C^{\frac12}_0L^{-\frac{1}{2}}.
\end{align*}
Combining the above estimates we obtain that for each $t\in[0,T]$,
\begin{align*}
    \Vert M\Vert_{L^\infty([0,x_*])}\le C_0^{\frac12}L^{\frac{1}{2}}a^{-m}+ C^{\frac12}_0L_0^{-\frac12}.
\end{align*}

\textbf{4. Estimates for other average integration terms.} In this step we estimate other average integration terms on the right hand side of \eqref{reformulated effective flux equation} and the constant $C_0>0$ will increase with each occurrence. First, since $L_0>1$, $y_s>\dfrac{1}{4}$. Thus we have $v\ge C^{-1}\left(\dfrac{1}{4},T\right) $, and the lower bound is independent of $a$. By the scale analysis and L'Hôpital's rule, we have the following limit near the singularity point $v=1$.
\begin{align*}
    \lim_{v\rightarrow1}\dfrac{|P(v)-P(1)|^2}{\psi(v)}=2A\gamma.
\end{align*}
Next since $|P(v)-P(1)|^2\rightarrow A^2$ as $v\rightarrow\infty$, we have
\begin{align*}
    |P(v)-P(1)|^2\le C_0\psi(v) \text{ for each } v\ge C^{-1}_0.
\end{align*}
Thus,
\begin{align*}
    |\mathscr{A}(P(v)-P(1))|^2&\le\int^1_0\omega(s)|\left(P(v)-P(1)\right)(t,y_s)|^2\ds\le C_0\int^1_0\omega(s)\psi(v)(t,y_s)\ds\\
    &\le \dfrac{C_0}{L_0}\int^{x+L}_x\psi(v)(t,y)\dy\le\dfrac{C^2_P}{L_0},
\end{align*}
where $C_P$ is a fixed positive constant. With the help of Taylor expansion, we have the following relation for each $h$ and $s\in\supp\omega$, 
\begin{align*}
0\le\omega(s+h)&=\omega(s)+\omega'(s)h+\dfrac{\omega''(s+\theta h)}{2}h^2
\le\omega(s)+\omega'(s)h+\dfrac{C_\omega}{2}h^2,
\end{align*}
where we choose $C_\omega=\max\{4\Vert\omega''\Vert_{L^\infty},4\Vert\omega'\Vert_{L^\infty}\}$ sufficiently large. If $h=-\dfrac{\omega'(s)}{C_\omega}$, then
\begin{align*}
    0\le\omega(s)-\dfrac{|\omega'(s)|^2}{2C_\omega}.
\end{align*}
As a result, $\dfrac{|\omega'(s)|^2}{\omega}\le 2 C_\omega$ for $s\in\supp\omega$. Based on this fact, for the third term in \eqref{reformulated effective flux equation}, it holds from integrating by parts that
\begin{align*}
    \left|\mathscr{A}\left(\partial_x\left(\dfrac{r^{2m}\partial_xv}{v^5}\right)\right)\right|&=\left|\int^1_0\omega(s)\partial_y\left(\dfrac{r^{2m}\partial_xv}{v^5}\right)(t,y_s)\ds\right|
    =\left|\int^1_{0}\omega(s)\partial_s\left(\dfrac{r^{2m}\partial_xv}{v^5}\right)(t,y_s)\dfrac{\partial s}{\partial y}\ds\right|\\
    &=\left|-\dfrac{1}{L}\int^1_{0}\omega'(s)\dfrac{r^{2m}\partial_xv}{v^5}(t,y_s)\ds\right|\\
    &\le\dfrac{5}{2}\int^1_0 \omega(s)\dfrac{r^{2m}|\partial_xv|^2}{v^6}\ds+\dfrac{C_0}{L^2}\int_{\supp~\omega}\dfrac{|\omega'(s)|^2}{\omega(s)} \dfrac{r^{2m}}{v^4}\ds.
\end{align*}
For the second term in the right hand side of the last inequality,
\begin{align*}
    \dfrac{C_0}{L^2}\int_{\supp~\omega}\dfrac{|\omega'(s)|^2}{\omega(s)} \dfrac{r^{2m}}{v^4}\ds=\dfrac{C_\omega C_0}{L^2}\int^{\frac{3}{4}}_{\frac{1}{4}}\dfrac{r^{2m}}{v^4}(t,y_s)\ds\le C_\mathscr{A}L^{-\frac2n},
\end{align*}
where we used the upper bound of Lemma \ref{particle path for ab}, $x_*$ is the constant independent of $a$ and $b$, and  $C_\mathscr{A}$ is a fixed positive constant depending on $C_0$.\\

\textbf{5. Estimates of derivatives of $M$.} In this step we give the upper bound of $M_x$ and $M_{xx}$. Since $\partial_xL=1,\partial_xy_s=1+s$, 
\begin{align*}
    \partial_x\left(L\int^1_0W(s)\dfrac{u}{r^m}(t,y_s)\ds\right)=\left(-\dfrac{u}{r^m}\right)(t,x)+\int^1_0\omega(s)(1+s)\dfrac{u}{r^m}(t,y_s)\ds.
\end{align*}
Next, 
\begin{align*}
    \partial_x[\mathscr{A}(\log v)]=-\dfrac{1}{L}\int^1_0\left(\omega'(s)(1+s)+\omega(s)\right)\log v(t,y_s)\ds
\end{align*}
Combining the above two equalities we have
\begin{align*}
    M_x=\left(-\dfrac{u}{r^m}\right)(x)+\underline{\int^1_0\omega(s)(1+s)\dfrac{u}{r^m}(t,y_s)\ds+\dfrac{\beta}{L}\int^1_0\left(\omega'(s)(1+s)+\omega(s)\right)\log v(t,y_s)\ds}_{:=M_1},
\end{align*}
By the same argument mentioned above, 
\begin{align*}
    &\left|\int^1_0\omega(s)(1+s)\dfrac{u}{r^m}(t,y_s)\ds\right|\le C_ML_0^{-\frac{1}{2}-\frac{m}{n}},\quad\left|\dfrac{\beta}{L}\int^1_0\left(\omega'(s)(1+s)+\omega(s)\right)\log v(t,y_s)\ds\right|\le C_ML_0^{-\frac{3}{2}}.
\end{align*}
where $C_M$ is a fixed positive constant depending on $C_0$. Then the supremum bound is obtained: $\Vert M_1\Vert_{L^\infty([0,x_*])}\le C_ML_0^{-\frac12-\frac mn}$. Now we differentiate $M_x$,
\begin{align*}
    M_{xx}=-\partial_x(\dfrac{u}{r^m}(x))+\underline{\partial_x\left(\int^1_0\omega(s)(1+s)\dfrac{u}{r^m}(t,y_s)\ds+\dfrac{\beta}{L}\int^1_0\left(\omega'(s)(1+s)+\omega(s)\right)\log v(t,y_s)\ds\right)}_{:=M_2}
\end{align*}
By the same argument mentioned above
\begin{align*}
    &\left|\partial_x\left(\int^1_0\omega(s)(1+s)\dfrac{u}{r^m}(t,y_s)\ds+\dfrac{\beta}{L}\int^1_0\left(\omega'(s)(1+s)+\omega(s)\right)\log v(t,y_s)\ds\right)\right|\\
    &\le\left|-\dfrac{1}{L}\int^1_0\left(\omega'(s)(1+s)^2+2\omega(s)(1+s)\right)\dfrac{u}{r^m}(t,y_s)\ds\right|\\
    &+\dfrac{\beta}{L^2}\left|\int^1_0\left\{(1+s)(4\omega'+\omega''(1+s))+2\omega(s)\right\}\log v(t,y_s)\ds\right|\\
    &\le C_M(L_0^{-\frac32-\frac mn}+L_0^{-\frac52}).
\end{align*}
Thus for $M_2$, it holds that $ \Vert M_2\Vert_{L^\infty([0,x_*])}\le C_ML^{-\frac{3}{2}-\frac mn}_0$.\\
\textbf{6. The definition of auxiliary functions.}
Now, we define an auxiliary function called algebraic cancellation function
\begin{align*}
    S_1(v)=\left\{\begin{array}{ll}
         \displaystyle\int^v_{\mathscr{V}}S'(s)\ds \quad&\text{ for } v>\mathscr{V},  \\
          S'(\mathscr{V})(v-\mathscr{V})+\dfrac{1}{2}S''(\mathscr{V})(v-\mathscr{V})^2\quad&\text{ for }v\le \mathscr{V},
    \end{array}\right.
\end{align*}
where the function $S'(v)$ is defined as follows
\begin{align*}
    S'(v)=\dfrac{\beta}{2v}\left(1+\sqrt{1-\dfrac{4\kappa}{\beta^2v^3}}\right)\text{ for }v\ge\max\left\{4^\gamma,\left(\dfrac{16\kappa}{\beta^2}\right)^{\frac13},2\mathscr{V}_0\right\}:=\mathscr{V}
\end{align*}
Since $v\ge C^{-1}_0 a^{\frac{2m}{3}}$, it holds that $|S_1|\le C(a,T)$ for each $v\in[C^{-1}_0a^{-\frac{2m}{3}},\mathscr{V}]$. It is worth noting that
\begin{align*}
    (S')^2-\dfrac{\beta}{v}S'+\dfrac{\kappa}{v^5}=0.
\end{align*}
This formula is crucial to cancel the high-order derivatives.\\
\textbf{7. The Maximum principle.}
Now we define $Y(t,x)=S_1(v(t,x))+M(t,x)+\delta x$ for $x\in[0,x_*]$, where $\delta=2C_ML_0^{-\frac{1}{2}-\frac{m}{n}}$ is a fixed constant. Since $Y_x(t,0)=M_x(t,0)+\delta$, $Y(t,0)$ cannot be the maximum value in $[0,x_*]$. For each fixed $t\in[0,T]$, there exists a $x_0(t)$ such that $Y(x_0,t)=\max_{x\in[0,x_*]}Y(t,x)$, and we consider two cases according to the value of $v(x_0(t),t)$. 
\begin{itemize}
    \item[1. ] $v(x_0(t),t)\le\mathscr{V}$. It is straightforward to check that $\Vert Y\Vert_{L^\infty([0,x_*])}\le C(a,T)$.\\
    \item[2. ] $v(x_0(t))>\mathscr{V}$. Since $Y_x(x_0(t),t)=0,Y_{xx}(x_0(t),t)\le0$, then
    \begin{align*}
        \dfrac{u}{r^m}=S'(v)\partial_xv+M_1+\delta,\quad
        \partial_x\left(\dfrac{u}{r^m}\right)\ge \partial_x\left(S'\partial_xv\right)+M_2\quad \text{ at }x_0(t).
    \end{align*}
    Combining the mass equation in \eqref{approximate problem},
    \begin{align*}
        \partial_tv&=r^{2m}\partial_x\left(\dfrac{u}{r^m}\right)+2mr^{m-1}v\dfrac{u}{r^m}\\
        &\ge2mr^{m-1}v\dfrac{u}{r^m}+r^m\partial_x\left(S'(v)\partial_xv\right)+r^mM_2\\
        &=\partial_x\left(r^{2m}S'(v)\partial_xv\right)+2mr^{m-1}v(M_1+\delta)+r^{2m}M_2 \quad\text{ at }x_0(t).
    \end{align*}
    Thus 
    \begin{align*}
        Y_t&=\left(S'(v)-\dfrac{\beta}{v}\right)\partial_tv+\partial_t\left(\beta\log v+M\right)\\
        &=\left(S'(v)-\dfrac{\beta}{v}\right)\partial_tv-\mathscr{A}(P(v)-P(1))-\kappa\mathscr{A}\left(\partial_x\left(\dfrac{r^{2m}\partial_xv}{v^5}\right)\right)-\dfrac{5\kappa}{2}\mathscr{A}\left(\dfrac{r^{2m}|\partial_xv|^2}{v^6}\right)\nonumber\\
    &+(P(v)-P(1))+\kappa\partial_x\left(\dfrac{r^{2m}\partial_xv}{v^5}\right)+\dfrac{5\kappa}{2}\left(\dfrac{r^{2m}|\partial_xv|^2}{v^6}\right)-L\int^1_0W(s)\clubsuit(t,y_s)\ds\\
    &\le \left(S'(v)-\dfrac{\beta}{v}\right)\partial_x\left(r^{2m}S'(v)\partial_xv\right)+\left(S'(v)-\dfrac{\beta}{v}\right)\left(2mr^{m-1}v(M_1+\delta)+r^{2m}M_2\right)\\
    &-\mathscr{A}(P(v)-P(1))-\kappa\mathscr{A}\left(\partial_x\left(\dfrac{r^{2m}\partial_xv}{v^5}\right)\right)-\dfrac{5\kappa}{2}\mathscr{A}\left(\dfrac{r^{2m}|\partial_xv|^2}{v^6}\right)\\
    &+(P(v)-P(1))+\kappa\partial_x\left(\dfrac{r^{2m}\partial_xv}{v^5}\right)+\dfrac{5\kappa}{2}\left(\dfrac{r^{2m}|\partial_xv|^2}{v^6}\right)\quad\text{ at }x_0(t),
    \end{align*}
    where we used $\left(S'(v)-\dfrac{\beta}{v}\right)\le0$. Since 
    \begin{align*}
        \left|S'(v)-\dfrac{\beta}{v}\right|=\left|\dfrac{\beta}{2v}\left(\sqrt{1-\dfrac{4\kappa}{\beta^2v^3}}-1\right)\right|\le\dfrac{C_0}{v}.
    \end{align*}
    Thus
    \begin{align*}
        \left|\left(S'(v)-\dfrac{\beta}{v}\right)\left(2mr^{m-1}v(M_1+\delta)+r^{2m}M_2\right)\right|\le C_{rest}L_0^{-\frac{1}{2}-\frac mn},
    \end{align*}
    where $C_{rest}$ is a fixed positive constant depending on $C_0$. Combining  step 4 and step 5 and the fact derived from $v(x_0(t),t)\ge\mathscr{V}$ that $P(v)-P(1)\le-\dfrac{P(1)}{2}$, we obtain that
    \begin{align*}
        Y_t\le-\dfrac{P(1)}{4} \text{ at }x_0(t).
    \end{align*}
\end{itemize} 
\textbf{8. Contradictions.} Based on the analysis mentioned before, $0<x_0(t)\le x_*$, and we choose a positive constant $\mathscr C$ such that  
\begin{align*}
    \mathscr{C}>\max\left\{\max_{t\in[0,T]}Y(x_*,t),\max_{x\in[0,x_*]}Y(0,x),\max_{C^{-1}_0a^{\frac{2m}{3}}\le v\le\mathscr{V}}Y(t,x)\right\}.
\end{align*}
 Assuming that there exists $t_1$ such that $\max_{x\in[0,x_*]}Y(t_1,x)>\mathscr{C}$. Let $x_1$ be a corresponding spatial maximizer,  $i.e.~Y(t_1,x_1)>\mathscr{C}$. Thus, there must exist a $t_2$ such that $\max_{x\in[0,x_*]}Y(t_2,x)=\mathscr{C}$ and $\max_{x\in[0,x_*]}Y(t,x)<\max_{x\in[0,x_*]}Y(t_2,x)$ for each $t<t_2$ and we let $x_2\in(0,x_*)$ be the corresponding spatial maximizer. It is straightforward to check that 
\begin{align*}
    v(t_2,x_2)>\mathscr{V},\quad Y_x(t_2,x_2)=0,\quad Y_{xx}(t_2,x_2)\le0,\quad Y_t(t_2,x_2)\ge0.
\end{align*}
However, this  contradicts $Y_t(t_2,x_2)<0$. Thus $\max_{x\in[0,x_*]}Y(t,x)\le \mathscr{C}\le C(a,T)$ for each $t\in[0,T]$. Thus $\max_{x\in[0,x_*]}S_1(v(t,x))\le C(a,T)$.\\
\textbf{9. The upper bound of $v$ in $[0,x_*]$} For high specific volume region $\{v\ge\mathscr{V}\}$, since $\dfrac{\beta}{2v}\le S'(v)\le\dfrac{\beta}{v}$, we integrate this relation from $\mathscr{V}$ to $v$ to obtain that
\begin{align*}
    \beta\log v(t,x)\le C(a,T) \text{ for each } x\in[0,x_*].
\end{align*}
\textbf{10. The upper bound of $v$ in $[\varepsilon,b]$ }
For  $\forall\varepsilon>0$, we choose $O\ge\varepsilon$ and consider the interval $I_O=(\dfrac{O}{2},3O),J_O=[O,2O]$. Define $\theta^O(x)=\dfrac{4\left(x-\dfrac{O}{2}\right)(3O-x)}{\dfrac{25}{4}O^2}$ for each $x\in I_O$.
Since $\theta^O$ has positive lower bound in $J_O$, now we define $L_O=L_0+O$. The corresponding averaging operator is defined by
\begin{align*}
    \mathscr{A}_Of(t,x)=\int^1_0\omega(s)f(t,y^O_s)\ds,
\end{align*}
where $y^O_s:=x+L_Os$. Then, by the same argument, 
\begin{align*}
    &\partial_t\left(L_O\int^1_0W(s)\dfrac{u}{r^m}(t,y^O_s)\ds-\beta\mathscr{A}_O(\log v)+\beta\log v\right)\\
    &=-L_O\int^1_0W(s)\clubsuit(t,y^O_s)\ds-\mathscr{A}_O(P(v)-P(1))-\kappa\mathscr{A}_O\left(\partial_x\left(\dfrac{r^{2m}\partial_xv}{v^5}\right)\right)-\dfrac{5\kappa}{2}\mathscr{A}_O\left(\dfrac{r^{2m}|\partial_xv|^2}{v^6}\right)\nonumber\\
    &+(P(v)-P(1))+\kappa\partial_x\left(\dfrac{r^{2m}\partial_xv}{v^5}\right)+\dfrac{5\kappa}{2}\left(\dfrac{r^{2m}|\partial_xv|^2}{v^6}\right).
\end{align*}
Now we define $M^O=L_O\displaystyle\int^1_0W(s)\dfrac{u}{r^m}(t,y^O_s)\ds-\beta\mathscr{A}_O(\log v)$. In addition, $M^O_1,M^O_2$ are the corresponding residual terms of the first order and the second order derivative to $M^O$. By a similar argument, it holds that 
\begin{align*}
    &\left\Vert M^O\right\Vert_{L^\infty([\frac{O}{2},3O])}\le C(\varepsilon,T),\quad\left\Vert r^{m-1}M^O_1\right\Vert_{L^\infty([\frac{O}{2},3O])}\le C(\varepsilon,T),\quad\left\Vert r^{2m}M^O_2\right\Vert_{L^\infty([\frac{O}{2},3O])}\le C(\varepsilon,T).
\end{align*}
The modified auxiliary function $Y^O$ is defined by
\begin{align*}
    Y^O(t,x)=S_1(v(t,x))+M^O+\beta\log\theta^O
\end{align*}
Since $x\rightarrow\dfrac{O}{2}\text{ or }3O$, $Y^O\rightarrow-\infty$, the maximum value point $x_0(t)$ is attained at an interior point of $I_O$ for each fixed $t_0\in[0,t]$. By a similar argument, at $x_0(t)$, we have
\begin{align*}
       &\dfrac{u}{r^m}=S'(v)\partial_xv+M^O_1+\beta\partial_x(\log \theta^O),\quad\partial_x\left(\dfrac{u}{r^m}\right)\ge\partial_x\left(S'(v)\partial_xv\right)+M^O_2+\beta\partial_{xx}(\log \theta^O),\\
       &\partial_tY^O\le \left(S'(v)-\dfrac{\beta}{v}\right)\left(2mr^{m-1}v(M^O_1+\beta\partial_x\log\theta^O)+r^{2m}(M^O_2+\beta\partial_{xx}\log\theta^O)\right)\\
    &-\mathscr{A}_O(P(v)-P(1))-\kappa\mathscr{A}_O\left(\partial_x\left(\dfrac{r^{2m}\partial_xv}{v^5}\right)\right)-\dfrac{5\kappa}{2}\mathscr{A}_O\left(\dfrac{r^{2m}|\partial_xv|^2}{v^6}\right)\\
    &+(P(v)-P(1))+\kappa\partial_x\left(\dfrac{r^{2m}\partial_xv}{v^5}\right)+\dfrac{5\kappa}{2}\left(\dfrac{r^{2m}|\partial_xv|^2}{v^6}\right)+\left(S'(v)-\dfrac{\beta}{v}\right)\partial_x\left(r^{2m}S'(v)\partial_xv\right).
\end{align*}
Now, for the function $\max_{x\in[\frac{O}{2},3O]}Y(t,x)$, we assume that this function takes the maximum $\mathscr{C}_\varepsilon$ at $t=t_1$, and the corresponding spatial point is $x_1$ where the $\mathscr{C}_\varepsilon>\max_{[C^{-1}(\varepsilon,T),\mathscr{V}]}Y(t,x) $ is some constant dependent on $\varepsilon$ and will be fixed later. It holds that
\begin{align*}
    \mathscr{C}_\varepsilon&=S(v)+M^O+\beta\log \theta^O \le\beta\log v+\beta\log \theta^O+C(\varepsilon,T)\quad\text{ at }(t_1,x_1).
\end{align*}
Thus $v\theta^O\ge \exp\left(\dfrac{\mathscr{C}_\varepsilon-C(\varepsilon,T)}{\beta}\right)$. It also holds that $ \dfrac{1}{v}\le C(\varepsilon,T)\exp\left(-\dfrac{\mathscr{C}_\varepsilon}{\beta}\right)\theta^O$. As a result, we have
\begin{align*}
    \left(S'(v)-\dfrac{\beta}{v}\right)\left(2mr^{m-1}v\beta\partial_x\log\theta^O+r^{2m}\beta\partial_{xx}\log\theta^O\right)\le C(\varepsilon,T)\exp\left(\dfrac{-2\mathscr{C}_\varepsilon}{\beta}\right).
\end{align*}
We choose $\mathscr{C}_\varepsilon$ sufficiently large such that $C(\varepsilon,T)\exp\left(\dfrac{-2\mathscr{C}_\varepsilon}{\beta}\right)\le\dfrac{P(1)}{8}$. By a similar argument, we can obtain that
\begin{align*}
    v(t,x)\le C(\varepsilon,T) \text{ for } (t,x)\in J_O\times[0,T].
\end{align*}
Now let $I_k=[2^k\varepsilon,2^{k+1}\varepsilon]$ for each $k=0,\cdots,n,\cdots$, there exists a $k_0\in \mathbb{N}_+$ such that $[\varepsilon,x_*]\subset\bigcup^{k_0}_{k=0}I_{k}$. Thus we have $v\le C(\varepsilon,T)$ for $x\in[\varepsilon,b]$.
\par The proof of Lemma \ref{upper bound of v} is finished now.
    \end{proof}
    \end{lemma}
\subsection{Estimates Independent of the Inner Radii}\label{Estimates independent of inner radii}
Now we start to establish some estimates independent of the inner radii which are crucial for the global existence of the Cauchy problem \eqref{1.6}--\eqref{bd condition}.
Motivated by Hoff \cite{H2} and Chen--Huang--Zhu \cite{CHZ}, we define the following spatial cut-off function $g_\varepsilon(x)$ by
\begin{align}
g_\varepsilon(x):=
    \begin{cases}
        0&x\le\dfrac{\varepsilon}{2},\\
        \dfrac{8}{\varepsilon^2}(x-\dfrac{\varepsilon}{2})^2&\dfrac{\varepsilon}{2}\le x\le\dfrac{3\varepsilon}{4},\\
        1-\dfrac{8}{\varepsilon^2}(x-\varepsilon)^2&\dfrac{3\varepsilon}{4}\le x\le\varepsilon,\\
        1&x\ge\varepsilon,
    \end{cases}
\end{align}
with the following properties
\begin{align}\label{cut off function for exterior}
    &\supp(g'_\varepsilon)=\left[\dfrac{\varepsilon}{2},\varepsilon\right],\quad0\le g_\varepsilon\le1,\quad g'_\varepsilon\ge0\nonumber,\quad|g'_\varepsilon(x)|^2\le\dfrac{32}{\varepsilon^2}g_\varepsilon(x) \text{ for all }x\in[0,\infty).
\end{align}
In addition, it follows from Lemma \ref{particle path for ab} that, for any $0<\varepsilon<1$,
\begin{align}
    \sup_{x\in[\varepsilon,k]}r^{-1}(t,x) \le C(\varepsilon,T),\sup_{x\in[\frac{\varepsilon}{2},\varepsilon]}r(t,x)\le C(\varepsilon,T).
\end{align}
\begin{lemma}\label{exterior estimate of v-1}
For any fixed $\varepsilon>0$, there exists a positive constant $C(\varepsilon,T)$ such that
    \begin{align}
        \sup_{t\in[0,T]} \int^b_0g_\varepsilon|v-1|^2\dx\le C (\varepsilon,T).
    \end{align}
    \begin{proof}
    It comes from mass equation of \eqref{approximate problem} that
\begin{align*}
    \dfrac{1}{2}\dfrac{\mathrm{d}}{\dt}\int^b_0g_\varepsilon|v-1|^2\dx=\int^b_0g_\varepsilon\partial(r^mu)(v-1)\dx\le C(\varepsilon,T)\left(1+\int^b_0g_\varepsilon|v-1|^2\dx\right).
\end{align*}
    \end{proof}
\end{lemma}
\begin{lemma}\label{u^2}
    For any fixed $\varepsilon\in(0,1)$, there exists a positive constant $C(\varepsilon,T)$ such that
    \begin{align}
        \int^T_0\int^b_\varepsilon|r^m\partial_xu|^2\dx\ds+\int^T_0\sup_{x\in[\varepsilon,b]}u^2(s,x)\ds\le C(\varepsilon,T).
    \end{align}
    \begin{proof}
    Since $\partial_x(r^mu)=\dfrac{mv}{r}u+r^m\partial_xu$, we have
    \begin{align*}
        \int^T_0\int^b_\varepsilon|r^m\partial_xu|^2\dx\ds\le\int^T_0\int^b_\varepsilon|\partial_x(r^mu)|^2\dx\ds+\int^T_0\int^b_\varepsilon\left|\dfrac{mv}{r}u\right|^2\dx\ds\le C(\varepsilon,T).
    \end{align*}
        With the help of Sobolev embedding theorem  $W^{1,1}\hookrightarrow C^0$, we have
        \begin{align}\label{l1 for u}
            \sup_{y\in[\varepsilon,b]}u^2(t,y)\le C\int^b_\varepsilon u^2\dx+C\int^b_\varepsilon|u\partial_xu|\dx\le C_0+\int^b_\varepsilon\left|\dfrac{u}{r^m}\partial_x(r^mu)+mr^{-n}vu^2\right|\dx.
        \end{align}
        Thus integrating for $t$, the proof is finished. 
    \end{proof}
\end{lemma}

Now we obtain the second order derivatives of $v$,
\begin{lemma}\label{second order derivatives of v}
    For each fixed $\varepsilon>0$, there exists a positive constant $C(\varepsilon,T)$ such that
    \begin{align}\label{sc or v}
         &\int_0^T\int_0^bg_\varepsilon
 \left[
 \dfrac{r^{4m}|\partial_{xx} v|^2}{v^7}
 +\dfrac{r^{4m}|\partial_x v|^4}{v^9}
 +\frac{r^{2m-2}|\partial_x v|^2}{v^5}
 \right]\dx\dt
 \le C(\varepsilon,T),\nonumber\\
&\int_0^T\left\|g_\varepsilon r^m\partial_x v\right\|_{L^\infty[0,b]}^4\dt
 \le C(\varepsilon,T).
    \end{align}
    \begin{proof}
        Multiplying $g_\varepsilon\dfrac{r^m\partial_xv}{v^2}$ and integrating by parts, we have
        \begin{align*}
&\kappa\int_0^bg_\varepsilon\left[
 \dfrac{r^{4m}|\partial_{xx}v|^2}{v^7}
 +4m^2\frac{r^{2m-2}|\partial_{x} v|^2}{v^5}
 +5\frac{r^{4m}|\partial_{x} v|^4}{v^9}\right]\dx\\
 &=-A\gamma\int_0^bg_\varepsilon
        \dfrac{r^{2m}|\partial_x v|^2}{v^{\gamma+3}}\dx+\int_0^bg_\varepsilon\partial_t u\dfrac{r^m\partial_x v}{v^2}\dx+\beta\int_0^b\partial_x\left(g_\varepsilon\dfrac{r^{2m}\partial_x v}{v^2}\right)\dfrac{\partial_x(r^mu)}v\dx\\
&+10m\kappa\int_0^bg_\varepsilon\dfrac{r^{3m-1}(\partial_{x} v)^3}{v^7}\dx-4m\kappa\int_0^bg_\varepsilon\dfrac{r^{3m-1}\partial_{x} v\partial_{xx} v}{v^6}\dx
 +\dfrac{9\kappa}{2}\int_0^bg_\varepsilon\frac{r^{4m}|\partial_{x} v|^2\partial_{xx} v}{v^8}\dx\\
 &-\kappa\int_0^bg_\varepsilon'\left[
 \dfrac{r^{4m}\partial_{x} v\partial_{xx} v}{v^7}
 +2m\dfrac{r^{3m-1}|\partial_{x} v|^2}{v^6}
 -\dfrac52\dfrac{r^{4m}(\partial_{x} v)^3}{v^8}
 \right]\dx:=\sum^{7}_{i=1}I_i.
        \end{align*}
For $I_1+I_4+I_5+I_6$, with the help of Cauchy inequality, it holds that
\begin{align*}
    I_1+I_4+I_5+I_6\le C(\varepsilon,T)+\dfrac{\kappa}{2}\int^b_0g_\varepsilon\dfrac{r^{4m}|\partial_{xx}v|^2}{v^7}\dx+\dfrac{9\kappa}{4}\int^b_0g_\varepsilon\dfrac{r^{4m}|\partial_{x}v|^4}{v^9}\dx.
\end{align*}
For $I_7$, since $\sup_{x\in[\varepsilon,k]}r^{-1}(t,x) \le C(\varepsilon,T),\sup_{x\in[\frac{\varepsilon}{2},\varepsilon]}r(t,x)\le C(\varepsilon,T)$, it holds that
\begin{align*}
    I_7\le\dfrac{\kappa}{8}\int^b_0g_\varepsilon\dfrac{r^{4m}|\partial_{xx}v|^2}{v^7}\dx+\dfrac{\kappa}{8}\int^b_0g_\varepsilon\dfrac{r^{4m}|\partial_{x}v|^4}{v^9}\dx+C(\varepsilon,T).
\end{align*}
Then for $I_4$, by virtue of  Lemma \ref{entropy estimate}, we have
\begin{align*}
    I_4&=\beta\int_0^bg'_\varepsilon\dfrac{r^{2m}\partial_x v}{v^2}\dfrac{\partial_x(r^mu)}v\dx+\beta\int_0^bg_\varepsilon\dfrac{r^{2m}\partial_{xx} v}{v^2}\dfrac{\partial_x(r^mu)}v\dx+2m\beta\int_0^bg_\varepsilon\dfrac{r^{m-1}\partial_x v}{v}\dfrac{\partial_x(r^mu)}v\dx\\
    &-2\beta\int_0^bg_\varepsilon\dfrac{r^{2m}|\partial_x v|^2}{v^3}\dfrac{\partial_x(r^mu)}v\dx\\
    &\le C(\varepsilon,T)\left(1+\int^b_0g_\varepsilon\dfrac{|\partial_x(r^mu)|^2}{v}\dx\right)+\dfrac{\kappa}{8}\int^b_0g_\varepsilon\dfrac{r^{4m}|\partial_{xx}v|^2}{v^7}\dx+\dfrac{\kappa}{8}\int^b_0g_\varepsilon\dfrac{r^{4m}|\partial_{x}v|^4}{v^9}\dx.
\end{align*}
Next for $I_2$,
\begin{align*}
    I_2&=\partial_t\left(\int_0^bg_\varepsilon u\dfrac{r^m\partial_x v}{v^2}\dx\right)-\int^b_0g_\varepsilon u^2\dfrac{mr^{m-1}\partial_xv}{v^2}\dx+2\int^b_0g_\varepsilon u\dfrac{r^m\partial_xv\partial_x(r^mu)}{v^3}\dx\\
    &+\int^b_0g'_\varepsilon u\dfrac{r^m\partial_xv}{v^2}\partial_x(r^mu)\dx+\int^b_0g_\varepsilon \dfrac{|\partial_x(r^mu)|^2}{v^2}\dx-2\int^b_0g_\varepsilon u\dfrac{r^m\partial_xv\partial_x(r^mu)}{v^3}\dx\\
    &:=\sum^{6}_{j=1} I^j_2.
\end{align*}
It is straight to check that
\begin{align*}
    \sum^{6}_{j=2} I^j_2\le C(\varepsilon,T)\left(1+\Vert g_\varepsilon u\Vert^2_{L^\infty}+\int^b_0g_\varepsilon\dfrac{|\partial_x(r^mu)|^2}{v}\dx\right).
\end{align*}
Combining the above results, we have
\begin{align*}
    &\dfrac{\kappa}{4}\int_0^bg_\varepsilon\left[
 \dfrac{r^{4m}|\partial_{xx}v|^2}{v^7}
 +4m^2\frac{r^{2m-2}|\partial_{x} v|^2}{v^5}
 +5\frac{r^{4m}|\partial_{x} v|^4}{v^9}\right]\dx\\
 &\le C(\varepsilon,T)\left(1+\Vert g_\varepsilon u\Vert^2_{L^\infty}\right)+\partial_t\left(\int_0^bg_\varepsilon u\dfrac{r^m\partial_x v}{v^2}\dx\right)+C(\varepsilon,T)\int^b_0g_\varepsilon\dfrac{|\partial_x(r^mu)|^2}{v}\dx.
\end{align*}
Then we integrate the above inequality in $[0,T]$ to obtain the desired result$\eqref{sc or v}_1$. Next, for each $x\in[0,b]$
\begin{align*}
    |g_\varepsilon r^m\partial_xv|^2&=2\int^x_0\partial_x\left(g_\varepsilon r^m\partial_xv\right)g_\varepsilon r^m\partial_xv\dy\\
    &\le C(\varepsilon,T)\Vert g_\varepsilon r^m\partial_xv\Vert_{L^2([0,b])}\Vert\partial_x\left(g_\varepsilon r^m\partial_xv\right)\Vert_{L^2([0,b])}
\end{align*}
Then $\eqref{sc or v}_2$ is proved.
\par The proof of Lemma \ref{second order derivatives of v} is finished.
    \end{proof}
\end{lemma}
Next, we give the second-order estimate of velocity and the third-order estimate of the specific volume, and we let $\sigma(t)=\min\{1,t\}$.
\begin{lemma}\label{lemma third estimates of v}
    For each fixed $\varepsilon>0$, there exists a positive constant $C(\varepsilon,T)$ such that
    \begin{align}\label{third estimates of v}
        &\sup_{t\in[0,T]}\sigma(t)\int^b_0g_\varepsilon\left(\dfrac{|r^m\partial_xu|^2}{v}+\left|\partial_x\left(\dfrac{r^m\partial_xv}{v^{\frac52}}\right)\right|^2\dfrac{r^{2m}}{v}\right)\dx\nonumber\\
        &+\int^T_0\sigma(\tau)\int^b_0g_\varepsilon\dfrac{r^{4m}}{v^2}\left(|\partial_{xx}u|^2+\left|\partial_{xx}\left(\dfrac{r^m\partial_xv}{v^{\frac52}}\right)\right|^2\right)\dx\mathrm{d}\tau\le C(\varepsilon,T).
    \end{align}
    \begin{proof}
    In this lemma, we temporarily let $(\sup_{t\in[0,T]}\sup_{x\in[\frac{\varepsilon}{2},b]}v)^{-1}:=\bar{\mathcal{C}}>0$ and $\inf_{t\in[0,T]}\inf_{x\in[\frac{\varepsilon}{2},b]}v:= \underline{\mathcal{C}}$. It is obvious that
    \begin{align}\label{equation of density gradient}
        \partial_t\left(\dfrac{r^m\partial_xv}{v^{\frac52}}\right)=\dfrac{r^m\partial_{xt}v}{v^{\frac52}}+\dfrac{mur^{{m-1}}\partial_xv}{v^{\frac52}}-\dfrac{5}{2}\dfrac{r^m\partial_xv\partial_tv}{v^{\frac72}}.
    \end{align}
        We multiply  \eqref{equation of density gradient} by  $\partial_x\left(g_\varepsilon\dfrac{r^{2m}}{v}\left(-\kappa\partial_x\left(\dfrac{r^m\partial_x v}{v^{\frac{5}{2}}}\right)+\a\partial_xu\right) \right)$ and multiply the momentum equation by $\partial_x\left(g_\varepsilon\dfrac{r^{2m}}{v}\left(\a\partial_x\left(\dfrac{r^m\partial_x v}{v^{\frac{5}{2}}}\right)-\partial_xu\right) \right)$. Let $\a\le\min\left\{\dfrac{\sqrt{\kappa}}{2},\dfrac{\kappa\bar{\mathcal{C}}^{\frac{3}{2}}}{2\beta},\dfrac{\beta\underline{\mathcal{C}}^\frac32}{4}\right\}$ be a fixed positive constant for each $\varepsilon>0$. Integrating by parts and adding the results together, we obtain
    \begin{align*}
        &\dfrac{1}{2}\dfrac{\mathrm{d}}{\dt}\left\{\int^b_0g_\varepsilon\dfrac{r^{2m}}{v}\left(|\partial_xu|^2+\kappa\left|\partial_x\left(\dfrac{r^m\partial_xv}{v^{\frac{5}{2}}}\right)\right|^2\right)\dx-2\a\int^b_0g_\varepsilon\dfrac{r^{2m}}{v}\partial_xu\partial_x\left(\dfrac{r^m\partial_xv}{v^{\frac{5}{2}}}\right)\dx\right\}\\
        &+\beta\int^b_0g_\varepsilon\dfrac{r^{4m}}{v^2}|\partial_{xx}u|^2\dx+\a\kappa\int^b_0g_\varepsilon\dfrac{r^{4m}}{v^{\frac{7}{2}}}\left|\partial_{xx}\left(\dfrac{r^m\partial_xv}{v^{\frac{5}{2}}}\right)\right|^2\dx-\a\int^b_0g_\varepsilon\dfrac{r^{4m}}{v^{\frac72}}|\partial_{xx}u|^2\dx\\
        &-\a\beta\int^b_0g_\varepsilon\partial_{xx}\left(\dfrac{r^m\partial_xv}{v^{\frac{5}{2}}}\right)\partial_{xx}u\dfrac{r^{4m}}{v^2}\dx\\
        &=\dfrac{1}{2}\int^b_0g_\varepsilon\partial_t\left(\dfrac{r^{2m}}{v}\right)\left(|\partial_xu|^2-2\a\partial_xu\partial_x\left(\dfrac{r^m\partial_xv}{v^{\frac{5}{2}}}\right)+\kappa\left|\partial_x\left(\dfrac{r^m\partial_xv}{v^{\frac{5}{2}}}\right)\right|^2\right)\dx+\sum^{14}_{i=1}I_i,
    \end{align*}
    where $I_i,~i=1,\cdots,14$ are the remaining terms in integration and will be specified later.
Since $\a\le\dfrac{\sqrt{\kappa}}{2}$, for the first term on the left hand side, we have the following inequality for each $t\in[0,T]$
\begin{align}\label{equiv first order}
    \dfrac{1}{2}|\partial_xu|^2+\dfrac{\kappa}{2}\left|\partial_x\left(\dfrac{r^m\partial_xv}{v^{\frac{5}{2}}}\right)\right|^2&\le|\partial_xu|^2+\kappa\left|\partial_x\left(\dfrac{r^m\partial_xv}{v^{\frac{5}{2}}}\right)\right|^2-2\a\partial_xu\partial_x\left(\dfrac{r^m\partial_xv}{v^{\frac{5}{2}}}\right)\\
    &\le2|\partial_xu|^2+2\kappa\left|\partial_x\left(\dfrac{r^m\partial_xv}{v^{\frac{5}{2}}}\right)\right|^2.
\end{align}
For the remaining terms on the left hand side, we let $\widetilde{\theta}=\min\left\{\dfrac{\beta}{4},\dfrac{\alpha\kappa}{2}\mathcal{C}^{\frac{3}{2}}\right\}$. Thus 
\begin{align*}
    &\dfrac{1}{2}\dfrac{\mathrm{d}}{\dt}\left\{\int^b_0g_\varepsilon\dfrac{r^{2m}}{v}\left(|\partial_xu|^2+\kappa\left|\partial_x\left(\dfrac{r^m\partial_xv}{v^{\frac{5}{2}}}\right)\right|^2-2\a\partial_xu\partial_x\left(\dfrac{r^m\partial_xv}{v^{\frac52}}\right)\right)\dx\right\}\\
    &+\widetilde{\theta}\int^b_0g_\varepsilon\dfrac{r^{4m}}{v^2}\left(|\partial_{xx}u|^2+\left|\partial_{xx}\left(\dfrac{r^m\partial_xv}{v^{\frac{5}{2}}}\right)\right|^2\right)\dx\\
    &\le\dfrac{1}{2}\int^b_0g_\varepsilon\partial_t(\dfrac{r^{2m}}{v})\left(|\partial_xu|^2-2\a\partial_xu\partial_x\left(\dfrac{r^m\partial_xv}{v^{\frac{5}{2}}}\right)+\kappa\left|\partial_x\left(\dfrac{r^m\partial_xv}{v^{\frac{5}{2}}}\right)\right|^2\dx\right)+\sum^{14}_{i=1}I_i.
\end{align*}
For the first term on the right hand side, we obtain
\begin{align*}
    &\dfrac{1}{2}\int^b_0g_\varepsilon\partial_t\left(\dfrac{r^{2m}}{v}\right)\left(|\partial_xu|^2-2\a\partial_xu\partial_x\left(\dfrac{r^m\partial_xv}{v^{\frac{5}{2}}}\right)+\kappa\left|\partial_x\left(\dfrac{r^m\partial_xv}{v^{\frac{5}{2}}}\right)\right|^2\right)\dx\\
    &=\dfrac{1}{2}\int^b_0g_\varepsilon\left(\dfrac{mr^{2m-1}u}{v}-\dfrac{r^{3m}\partial_xu}{v^2}\right)\left(|\partial_xu|^2-2\a\partial_xu\partial_x\left(\dfrac{r^m\partial_xv}{v^{\frac{5}{2}}}\right)+\kappa\left|\partial_x\left(\dfrac{r^m\partial_xv}{v^{\frac{5}{2}}}\right)\right|^2\right)\dx\\
    &\le C(\varepsilon,T)\int^b_0g_\varepsilon\dfrac{r^{2m}}{v}\left(|\partial_xu|^2+\kappa\left|\partial_x\left(\dfrac{r^m\partial_xv}{v^{\frac{5}{2}}}\right)\right|^2\right)\dx\\
    &+\dfrac{\widetilde{\theta}}{20}\int^b_0g_\varepsilon\dfrac{r^{4m}}{v^2}\left(|\partial_{xx}u|^2+\left|\partial_{xx}\left(\dfrac{r^m\partial_xv}{v^{\frac{5}{2}}}\right)\right|^2\right)\dx
\end{align*}
Next, for $I_1+I_2+I_7$, it follows that
\begin{align*}
        &I_1+I_2+I_7=-\int^b_0 r^m\partial_xP\partial_x\left[-g_\varepsilon\dfrac{r^{2m}}{v}\partial_xu+\a g_\varepsilon\dfrac{r^{2m}}{v}\partial_x\left(\dfrac{r^m\partial_x v}{v^{\frac{5}{2}}}\right)\right]\dx\\
        &+\int^b_0g_\varepsilon'\dfrac{r^{2m}}{v}\left(\a\partial_x\left(\dfrac{r^m\partial_x v}{v^{\frac{5}{2}}}\right)-\partial_xu\right)\beta\left[\dfrac{r^{2m}}{v}\partial_{xx}u+2mr^{m-1}\partial_xu-\dfrac{r^{2m}\partial_xv\partial_xu}{v^2}-muvr^{-2}\right]\dx\\
        &-\int^b_0g_\varepsilon'\dfrac{r^{2m}}{v}\left(\a\partial_x\left(\dfrac{r^m\partial_x v}{v^{\frac{5}{2}}}\right)-\partial_xu\right)\kappa\bigg[\dfrac{r^{2m}}{v^{\frac52}}\partial_{xx}\left(\dfrac{r^m\partial_xv}{v^{\frac{5}{2}}}\right)+2r^m\partial_x\left(\dfrac{r^m}{v^{\frac52}}\right)\partial_{x}\left(\dfrac{r^m\partial_xv}{v^{\frac{5}{2}}}\right)\\
        &+\partial_{xx}\left(\dfrac{r^m}{v^{\frac52}}\right)\dfrac{r^{2m}\partial_xv}{v^{\frac52}}\bigg]\dx\\
        &\le \dfrac{\widetilde{\theta}}{20}\int^b_0g_\varepsilon\dfrac{r^{4m}}{v^2}\left(|\partial_{xx}u|^2+\left|\partial_{xx}\left(\dfrac{r^m\partial_xv}{v^{\frac{5}{2}}}\right)\right|^2\right)\dx+C(\varepsilon,T)\left(\int^\varepsilon_{\frac{\varepsilon}{2}}|\partial_xu|^2+\left|\partial_x\left(\dfrac{r^m\partial_xv}{v^{\frac52}}\right)\right|^2\dx\right)\\
        &+C(\varepsilon,T)\left(1+\Vert g_\varepsilon u\Vert^2_{L^\infty([0,b])}+\Vert g_\varepsilon r^m\partial_xv\Vert^2_{L^\infty([0,b])}\right)\int^b_0g_\varepsilon\dfrac{r^{2m}}{v}\left(|\partial_xu|^2+\kappa\left|\partial_x\left(\dfrac{r^m\partial_xv}{v^{\frac{5}{2}}}\right)\right|^2\right)\dx\\
        &+C(\varepsilon,T)\left(1+\int^\varepsilon_{\frac{\varepsilon}{2}}\left[|r^m\partial_xv|^2+|r^{m}\partial_xv|^4+|r^{2m}\partial_{xx}v|^2\right]\dx\right).
    \end{align*}
    For $I_4+I_9+I_{10}$, it holds that
    \begin{align*}
        &I_4+I_9+I_{10}=-\int^b_0g_\varepsilon\dfrac{r^{2m}}{v}\partial_{xx}u\left[\beta\left(2mr^{m-1}\partial_xu-\dfrac{r^{2m}\partial_xv\partial_xu}{v^2}-\dfrac{muv}{r^{2}}\right)\right]\dx\\
        &+\kappa\int^b_0g_\varepsilon\dfrac{r^{2m}}{v}\partial_{xx}ur^m\left[2\partial_x\left(\dfrac{r^m}{v^{\frac52}}\right)\partial_x\left(\dfrac{r^m\partial_xv}{v^{\frac52}}\right)+\dfrac{r^m\partial_xv}{v^{\frac52}}\partial_{xx}\left(\dfrac{r^m}{v^{\frac52}}\right)\right]\dx\\
        &-\a\kappa\int^b_0g_\varepsilon\dfrac{r^{2m}}{v}\partial_{xx}\left(\dfrac{r^m\partial_xv}{v^{\frac{5}{2}}}\right)r^m\left[2\partial_x\left(\dfrac{r^m}{v^{\frac52}}\right)\partial_x\left(\dfrac{r^m\partial_xv}{v^{\frac52}}\right)+\dfrac{r^m\partial_xv}{v^{\frac52}}\partial_{xx}\left(\dfrac{r^m}{v^{\frac52}}\right)\right]\dx\\
        &\le \dfrac{\widetilde{\theta}}{20}\int^b_0g_\varepsilon\dfrac{r^{4m}}{v^2}\left(|\partial_{xx}u|^2+\left|\partial_{xx}\left(\dfrac{r^m\partial_xv}{v^{\frac{5}{2}}}\right)\right|^2\right)\dx\\
        &+C(\varepsilon,T)\left(1+\Vert g_\varepsilon r^m\partial_xv\Vert^4_{L^\infty([0,b])}+\int^b_0g_\varepsilon\left(|r^m\partial_xu|^2+|u|^2\right)\dx\right)\\
        &+C(\varepsilon,T)\left(1+\Vert g_\varepsilon u\Vert^2_{L^\infty([0,b])}+\Vert g_\varepsilon r^m\partial_xv\Vert^2_{L^\infty([0,b])}\right)\int^b_0g_\varepsilon\dfrac{r^{2m}}{v}\left(|\partial_xu|^2+\kappa\left|\partial_x\left(\dfrac{r^m\partial_xv}{v^{\frac{5}{2}}}\right)\right|^2\right)\dx.
    \end{align*}
    Next, we estimate $I_3+I_8$ as follows
    \begin{align*}
        I_3+I_8&=\int^b_0g_\varepsilon\partial_x\left(\dfrac{r^{2m}}{v}\right)\left(\a\partial_x\left(\dfrac{r^m\partial_x v}{v^{\frac{5}{2}}}\right)-\partial_xu\right)\beta\left[\dfrac{r^{2m}\partial_{xx}u}{v}+2mr^{m-1}\partial_xu\right]\dx\\
        &+\int^b_0g_\varepsilon\partial_x\left(\dfrac{r^{2m}}{v}\right)\left(\a\partial_x\left(\dfrac{r^m\partial_x v}{v^{\frac{5}{2}}}\right)-\partial_xu\right)\beta\left[-\dfrac{r^{2m}\partial_xv\partial_xu}{v^2}-\dfrac{muv}{r^{2}}\right]\dx\\
        &-\kappa\int^b_0g_\varepsilon\partial_x\left(\dfrac{r^{2m}}{v}\right)\left(\a\partial_x\left(\dfrac{r^m\partial_x v}{v^{\frac{5}{2}}}\right)-\partial_xu\right)\dfrac{r^{2m}}{v^{\frac52}}\partial_{xx}\left(\dfrac{r^m\partial_xv}{v^{\frac{5}{2}}}\right)\dx\\
        &-\kappa\int^b_0g_\varepsilon\partial_x\left(\dfrac{r^{2m}}{v}\right)\left(\a\partial_x\left(\dfrac{r^m\partial_x v}{v^{\frac{5}{2}}}\right)-\partial_xu\right)2r^m\partial_x\left(\dfrac{r^m}{v^{\frac52}}\right)\partial_{x}\left(\dfrac{r^m\partial_xv}{v^{\frac{5}{2}}}\right)\dx\\
        &-\kappa\int^b_0g_\varepsilon\partial_x\left(\dfrac{r^{2m}}{v}\right)\left(\a\partial_x\left(\dfrac{r^m\partial_x v}{v^{\frac{5}{2}}}\right)-\partial_xu\right)\partial_{xx}\left(\dfrac{r^m}{v^{\frac52}}\right)\dfrac{r^{2m}\partial_xv}{v^{\frac52}}\dx\\
        &\le\dfrac{\widetilde{\theta}}{20}\int^b_0g_\varepsilon\dfrac{r^{4m}}{v^2}\left(|\partial_{xx}u|^2+\left|\partial_{xx}\left(\dfrac{r^m\partial_xv}{v^{\frac{5}{2}}}\right)\right|^2\right)\dx\\
        &+C(\varepsilon,T)\left(1+\Vert g_\varepsilon u\Vert^2_{L^\infty([0,b])}+\Vert g_\varepsilon r^m\partial_xv\Vert^2_{L^\infty([0,b])}\right)\int^b_0g_\varepsilon\dfrac{r^{2m}}{v}\left(|\partial_xu|^2+\kappa\left|\partial_x\left(\dfrac{r^m\partial_xv}{v^{\frac{5}{2}}}\right)\right|^2\right)\dx\\
        &+C(\varepsilon,T)\left(1+\Vert g_\varepsilon u\Vert^2_{L^\infty([0,b])}+\Vert g_\varepsilon r^m\partial_xv\Vert^4_{L^\infty([0,b])}\right).
    \end{align*}
    Based on the above estimates, for $I_5+I_6$, we have
\begin{align*}
    &I_5+I_6=\a\beta\int^b_0g_\varepsilon\dfrac{r^{2m}}{v}\partial_{xx}\left(\dfrac{r^m\partial_xv}{v^{\frac{5}{2}}}\right)\left[2mr^{m-1}\partial_xu-\dfrac{r^{2m}\partial_xv\partial_xu}{v^2}-\dfrac{muv}{r^{2}}\right]\dx\\
   &+ \int^b_0g'_\varepsilon\dfrac{r^{2m}}{v}\left(-\kappa\partial_x\left(\dfrac{r^m\partial_x v}{v^{\frac{5}{2}}}\right)+\a\partial_xu\right)\left[\dfrac{r^{2m}\partial_{xx}u}{v^\frac52}+\dfrac{2mr^{m-1}\partial_xu}{v^\frac32}-\dfrac{mu}{v^{\frac12}r^{2}}+\dfrac{mr^{m-1}u\partial_xv}{v^{\frac52}}\right]\dx\\
   &+\int^b_0g_\varepsilon\partial_x\left(\dfrac{r^{2m}}{v}\right)\left(-\kappa\partial_x\left(\dfrac{r^m\partial_x v}{v^{\frac{5}{2}}}\right)+\a\partial_xu\right)\left[\dfrac{r^{2m}\partial_{xx}u}{v^\frac52}+\dfrac{2mr^{m-1}\partial_xu}{v^\frac32}-\dfrac{mu}{v^{\frac12}r^{2}}+\dfrac{mr^{m-1}u\partial_xv}{v^{\frac52}}\right]\dx\\
   &+\int^b_0g_\varepsilon\dfrac{r^{2m}}{v}\left(\a\partial_{xx}u-\kappa\partial_{xx}\left(\dfrac{r^m\partial_xv}{v^{\frac{5}{2}}}\right)\right)\left[\dfrac{2mr^{m-1}\partial_xu}{v^\frac32}-\dfrac{mu}{v^{\frac12}r^{2}}+\dfrac{mr^{m-1}u\partial_xv}{v^{\frac52}}\right]\dx\\
   &\le \dfrac{\widetilde{\theta}}{20}\int^b_0g_\varepsilon\dfrac{r^{4m}}{v^2}\left(|\partial_{xx}u|^2+\left|\partial_{xx}\left(\dfrac{r^m\partial_xv}{v^{\frac{5}{2}}}\right)\right|^2\right)\dx+C(\varepsilon,T)\left(\int^\varepsilon_{\frac{\varepsilon}{2}}|\partial_xu|^2+\left|\partial_x\left(\dfrac{r^m\partial_xv}{v^{\frac52}}\right)\right|^2\dx\right)\\
        &+C(\varepsilon,T)\left(1+\Vert g_\varepsilon r^m\partial_xv\Vert^4_{L^\infty([0,b])}+\int^b_0g_\varepsilon|r^m\partial_xu|^2\dx\right)\\
        &+C(\varepsilon,T)\left(1+\Vert g_\varepsilon u\Vert^2_{L^\infty([0,b])}+\Vert g_\varepsilon r^m\partial_xv\Vert^2_{L^\infty([0,b])}\right)\int^b_0g_\varepsilon\dfrac{r^{2m}}{v}\left(|\partial_xu|^2+\kappa\left|\partial_x\left(\dfrac{r^m\partial_xv}{v^{\frac{5}{2}}}\right)\right|^2\right)\dx.
\end{align*}
   For $I_{11}+I_{12}$, it holds that
   \begin{align*}
&I_{11}+I_{12}=\int^b_0g'_\varepsilon\dfrac{r^{2m}}v{}\left(\a\partial_x\left(\dfrac{r^m\partial_x v}{v^{\frac{5}{2}}}\right)-\partial_xu\right)\left[-\dfrac{5\kappa}{2}r^m\left(2\partial_x\left(\dfrac{\partial_xv r^m}{v^{\frac52}}\right)\dfrac{\partial_xvr^m}{v^{\frac72}}-\dfrac{r^{2m}(\partial_xv)^3}{v^7}\right)\right]\dx\\
       &+\int^b_0g_\varepsilon\partial_x\left(\dfrac{r^{2m}}{v}\right)\left(\a\partial_x\left(\dfrac{r^m\partial_x v}{v^{\frac{5}{2}}}\right)-\partial_xu\right)\left[-\dfrac{5\kappa}{2}r^m\left(2\partial_x\left(\dfrac{\partial_xv r^m}{v^{\frac52}}\right)\dfrac{\partial_xvr^m}{v^{\frac72}}-\dfrac{r^{2m}(\partial_xv)^3}{v^7}\right)\right]\dx\\
       &+\int^b_0g_\varepsilon\dfrac{r^{2m}}{v}\left(\a\partial_{xx}\left(\dfrac{r^m\partial_x v}{v^{\frac{5}{2}}}\right)-\partial_{xx}u\right)\left[-\dfrac{5\kappa}{2}r^m\left(2\partial_x\left(\dfrac{\partial_xv r^m}{v^{\frac52}}\right)\dfrac{\partial_xvr^m}{v^{\frac72}}-\dfrac{r^{2m}(\partial_xv)^3}{v^7}\right)\right]\dx\\
       &+\int^b_0g'_\varepsilon\dfrac{r^{2m}}v{}\left(\a\partial_x\left(\dfrac{r^m\partial_x v}{v^{\frac{5}{2}}}\right)-\partial_xu\right)\left(-m\kappa r^{2m-1}\dfrac{|\partial_x v|^2}{v^5}\right)\dx\\
       &+\int^b_0g_\varepsilon\dfrac{r^{2m}}{v}\left(\a\partial_{xx}\left(\dfrac{r^m\partial_x v}{v^{\frac{5}{2}}}\right)-\partial_{xx}u\right)\left(-m\kappa r^{2m-1}\dfrac{|\partial_x v|^2}{v^5}\right)\dx\\
       &+\int^b_0g_\varepsilon\partial_x\left(\dfrac{r^{2m}}{v}\right)\left(\a\partial_x\left(\dfrac{r^m\partial_x v}{v^{\frac{5}{2}}}\right)-\partial_xu\right)\left(-m\kappa r^{2m-1}\dfrac{|\partial_x v|^2}{v^5}\right)\dx\\
       &\le \dfrac{\widetilde{\theta}}{10}\int^b_0g_\varepsilon\dfrac{r^{4m}}{v^2}\left(|\partial_{xx}u|^2+\left|\partial_{xx}\left(\dfrac{r^m\partial_xv}{v^{\frac{5}{2}}}\right)\right|^2\right)\dx+C(\varepsilon,T)\left(1+\Vert g_\varepsilon r^m\partial_xv\Vert^4_{L^\infty([0,b])}\right).
   \end{align*}
   For the last two terms $I_{13}+I_{14}$, it follows that
   \begin{align*}
       &I_{13}+I_{14}=\int^b_0g'_\varepsilon\dfrac{r^{2m}}{v}\left(-\kappa\partial_x\left(\dfrac{r^m\partial_x v}{v^{\frac{5}{2}}}\right)+\a\partial_xu\right)\left(\dfrac{mur^{m-1}\partial_xv}{v^{\frac52}}-\dfrac{5}{2}\dfrac{r^m\partial_xv}{v^{\frac72}}(r^m\partial_xu+\dfrac{mvu}{r})\right)\dx\\
       &+\int^b_0g_\varepsilon\partial_x\left(\dfrac{r^{2m}}{v}\right)\left(-\kappa\partial_x\left(\dfrac{r^m\partial_x v}{v^{\frac{5}{2}}}\right)+\a\partial_xu\right)\left(\dfrac{mur^{m-1}\partial_xv}{v^{\frac52}}-\dfrac{5}{2}\dfrac{r^m\partial_xv}{v^{\frac72}}(r^m\partial_xu+\dfrac{mvu}{r})\right)\dx\\
       &+\int^b_0g_\varepsilon\dfrac{r^{2m}}{v}\left(-\kappa\partial_{xx}\left(\dfrac{r^m\partial_x v}{v^{\frac{5}{2}}}\right)+\a\partial_{xx}u\right)\left(\dfrac{mur^{m-1}\partial_xv}{v^{\frac52}}-\dfrac{5}{2}\dfrac{r^m\partial_xv}{v^{\frac72}}(r^m\partial_xu+\dfrac{mvu}{r})\right)\dx\\
       &\le\dfrac{\widetilde{\theta}}{10}\int^b_0g_\varepsilon\dfrac{r^{4m}}{v^2}\left(|\partial_{xx}u|^2+\left|\partial_{xx}\left(\dfrac{r^m\partial_xv}{v^{\frac{5}{2}}}\right)\right|^2\right)\dx+C(\varepsilon,T)\left(\int^\varepsilon_{\frac{\varepsilon}{2}}|\partial_xu|^2+\left|\partial_x\left(\dfrac{r^m\partial_xv}{v^{\frac52}}\right)\right|^2\dx\right)\\
       &+C(\varepsilon,T)\left(1+\Vert g_\varepsilon u\Vert^2_{L^\infty([0,b])}+\Vert g_\varepsilon r^m\partial_xv\Vert^4_{L^\infty([0,b])}+\int^b_0g_\varepsilon\left(|r^m\partial_xu|^2+|u|^2\right)\dx\right)\\
        &+C(\varepsilon,T)\left(1+\Vert g_\varepsilon u\Vert^2_{L^\infty([0,b])}+\Vert g_\varepsilon r^m\partial_xv\Vert^2_{L^\infty([0,b])}\right)\int^b_0g_\varepsilon\dfrac{r^{2m}}{v}\left(|\partial_xu|^2+\kappa\left|\partial_x\left(\dfrac{r^m\partial_xv}{v^{\frac{5}{2}}}\right)\right|^2\right)\dx.
   \end{align*}
Combining the above estimates, we obtain
\begin{align*}
    &\dfrac{1}{2}\dfrac{\mathrm{d}}{\dt}\left\{\int^b_0g_\varepsilon\dfrac{r^{2m}}{v}\left(|\partial_xu|^2+\kappa\left|\partial_x\left(\dfrac{r^m\partial_xv}{v^{\frac{5}{2}}}\right)\right|^2-2\a\partial_xu\partial_x\left(\dfrac{r^m\partial_xv}{v^{\frac52}}\right)\right)\dx\right\}\\
    &+\dfrac{\widetilde{\theta}}{2}\int^b_0g_\varepsilon\dfrac{r^{4m}}{v^2}\left(|\partial_{xx}u|^2+\left|\partial_{xx}\left(\dfrac{r^m\partial_xv}{v^{\frac{5}{2}}}\right)\right|^2\right)\dx\\
    &\le C(\varepsilon,T)\left(1+\Vert g_\varepsilon u\Vert^2_{L^\infty([0,b])}+\Vert g_\varepsilon r^m\partial_xv\Vert^2_{L^\infty([0,b])}\right)\int^b_0g_\varepsilon\dfrac{r^{2m}}{v}\left(|\partial_xu|^2+\kappa\left|\partial_x\left(\dfrac{r^m\partial_xv}{v^{\frac{5}{2}}}\right)\right|^2\right)\dx\\
    &+C(\varepsilon,T)\left(1+\Vert g_\varepsilon u\Vert^2_{L^\infty([0,b])}+\Vert g_\varepsilon r^m\partial_xv\Vert^4_{L^\infty([0,b])}+\int^b_0g_\varepsilon\left(|r^m\partial_xu|^2+|u|^2\right)\dx\right)\\
    &+C(\varepsilon,T)\left(1+\int^\varepsilon_{\frac{\varepsilon}{2}}\left[|r^m\partial_xv|^2+|r^{m}\partial_xv|^4+|r^{2m}\partial_{xx}v|^2\right]\dx\right)\\
    &+C(\varepsilon,T)\left(\int^\varepsilon_{\frac{\varepsilon}{2}}|\partial_xu|^2+\left|\partial_x\left(\dfrac{r^m\partial_xv}{v^{\frac52}}\right)\right|^2\dx\right).
\end{align*}
Multiplying the above inequality by $\sigma$ and applying Gronwall's inequality, we obtain the desired estimates in \eqref{third estimates of v} by virtue of \eqref{equiv first order}.
\par The proof of Lemma \ref{lemma third estimates of v} is finished.
        \end{proof}
\end{lemma}
Next the $L^2_tL^2_x$ estimate of $\partial_tu$ is obtained, this estimate is crucial to obtain the strong convergence of velocity.
\begin{lemma}\label{time derivative of velocity}
    For each fixed $\varepsilon>0$, there exists a positive constant $C(\varepsilon,T)$ such that
    \begin{align*}
        \int^T_0\int^b_0g_\varepsilon\sigma(s)|\partial_tu|^2\dx\ds+\sup_{t\in[0,T]}\left(\sigma(t)\int^b_0g_\varepsilon\dfrac{|\partial_x(r^mu)^2|}{v}\dx\right)\le C(\varepsilon,T).
    \end{align*}
    \begin{proof}
        Multiplying $g_\varepsilon\partial_tu$ to the momentum equation \eqref{approximate problem}, it follows that
        \begin{align*}
            \int^b_0g_\varepsilon|\partial_tu|^2\dx&=\int^b_0g_\varepsilon r^m\partial_x\left(\dfrac{\beta\partial_x(r^mu)}{v}\right)\partial_tu\dx-\kappa\int^b_0g_\varepsilon r^m\partial_tu\partial_{xx}\left(\dfrac{r^{2m}\partial_xv}{v^5}\right)\dx\\
            &-\dfrac{5\kappa}{2}\int^b_0g_\varepsilon\partial_tu\partial_x\left(\dfrac{r^{2m}|\partial_xv|^2}{v^6}\right)r^m\dx-m\kappa\int^b_0g_\varepsilon\partial_tu\dfrac{r^{2m-1}|\partial_xv|^2}{v^5}\dx\\
            &-\int^b_0g_\varepsilon\partial_xP\cdot r^m\partial_tu\dx\\
            &:=\sum^{5}_{i=1}I_i.
        \end{align*}
        For $I_1$, we have
        \begin{align*}
            I_1&=-\dfrac{\beta}{2}\left(\dfrac{\mathrm{d}}{\dt}\int^b_0g_\varepsilon\dfrac{|\partial_x(r^mu)|^2}{v}\dx\right)-\dfrac{\beta}{2}\int^b_0g_\varepsilon\dfrac{\left(\partial_x(r^mu)\right)^3}{v^2}\dx\\
            &+\beta\int^b_0g_\varepsilon\dfrac{\partial_x(r^mu)}{v}\partial_x\left(mr^{m-1}u^2\right)\dx-\beta\int^b_0r^m\beta\partial_tug'_\varepsilon\dfrac{\partial_x(r^mu)}{v}\dx\\
            &\le -\dfrac{\beta}{2}\left(\dfrac{\mathrm{d}}{\dt}\int^b_0g_\varepsilon\dfrac{|\partial_x(r^mu)|^2}{v}\dx\right)+C(\varepsilon,T)\left(1+\int^b_0g_\varepsilon|r^m\partial_xu|^2\dx+\int^b_0g_\varepsilon\dfrac{|\partial_x(r^mu)|^2}{v}\dx\right)\\
            &+C(\varepsilon,T)\left(\Vert g_\varepsilon u\Vert^2_{L^\infty([0,b])}+\Vert g_\varepsilon r^m\partial_xu
            \Vert^2_{L^\infty([0,b])}+\Vert g_\varepsilon r^m\partial_xu\Vert^{\frac52}_{L^2([0,b])}\Vert g_\varepsilon r^{2m}\partial_{xx}u\Vert^{\frac12}_{L^2([0,b])}\right)\\
            &+\dfrac{1}{10}\int^b_0g_\varepsilon|\partial_tu|^2\dx.
        \end{align*}
For $I_2$, it holds that
\begin{align*}
    I_2&=-\kappa\int^b_0g_\varepsilon\partial_tu\left(\dfrac{r^{2m}}{v^{\frac52}}\partial_{xx}\left(\dfrac{r^m\partial_xv}{v^{\frac{5}{2}}}\right)+2r^m\partial_x\left(\dfrac{r^m}{v^{\frac52}}\right)\partial_{x}\left(\dfrac{r^m\partial_xv}{v^{\frac{5}{2}}}\right)+\partial_{xx}\left(\dfrac{r^m}{v^{\frac52}}\right)\dfrac{r^{2m}\partial_xv}{v^{\frac52}}\right)\dx\\
    &\le\dfrac{1}{10}\int^b_0g_\varepsilon|\partial_tu|^2\dx+C(\varepsilon,T)\left(1+\Vert g_\varepsilon r^m\partial_xv\Vert^4_{L^\infty([0,b])}+\int^b_0g_\varepsilon\dfrac{r^{4m}}{v^2}\left|\partial_{xx}\left(\dfrac{r^m\partial_xv}{v^{\frac{5}{2}}}\right)\right|^2\dx\right)\\
    &+C(\varepsilon,T)\left(1+\Vert g_\varepsilon r^m\partial_xv\Vert^2_{L^\infty([0,b])}\right)\int^b_0g_\varepsilon\dfrac{r^{2m}}{v}\left|\partial_{x}\left(\dfrac{r^m\partial_xv}{v^{\frac{5}{2}}}\right)\right|^2\dx.
\end{align*}
The estimates of $I_3+I_3+I_5$ are similar to $I_1$ and $I_2$, here we only state the results and leave the proofs to the reader,
\begin{align*}
    I_3+I_3+I_5&\le\dfrac{3}{10}\int^b_0g_\varepsilon|\partial_tu|^2\dx+ C(\varepsilon,T)\left(1+\Vert g_\varepsilon r^m\partial_xv\Vert^4_{L^\infty([0,b])}\right)\\
&+C(\varepsilon,T)\left(1+\Vert g_\varepsilon r^m\partial_xv\Vert^2_{L^\infty([0,b])}\right)\int^b_0g_\varepsilon\dfrac{r^{2m}}{v}\left|\partial_{x}\left(\dfrac{r^m\partial_xv}{v^{\frac{5}{2}}}\right)\right|^2\dx.
\end{align*}
Combining the above estimates, it holds that
\begin{align*}
    &\dfrac{1}{2}\int^b_0g_\varepsilon|\partial_tu|^2\dx+\dfrac{1}{2}\dfrac{\mathrm{d}}{\dt}\int^b_0g_\varepsilon\dfrac{|\partial_x(r^mu)|^2}{v}\dx\\
    &\le C(\varepsilon,T)\left(1+\Vert g_\varepsilon r^m\partial_xv\Vert^2_{L^\infty([0,b])}\right)\int^b_0g_\varepsilon\dfrac{r^{2m}}{v}\left|\partial_{x}\left(\dfrac{r^m\partial_xv}{v^{\frac{5}{2}}}\right)\right|^2\dx+C(\varepsilon,T)\Vert g_\varepsilon r^m\partial_xv\Vert^4_{L^\infty([0,b])}\\
    &+C(\varepsilon,T)\left(1+\int^b_0g_\varepsilon|r^m\partial_xu|^2\dx+\int^b_0g_\varepsilon\dfrac{|\partial_x(r^mu)|^2}{v}\dx\right)\\
    &+C(\varepsilon,T)\left(1+\Vert g_\varepsilon r^m\partial_xv\Vert^2_{L^\infty([0,b])}\right)\int^b_0g_\varepsilon\dfrac{r^{2m}}{v}\left|\partial_{x}\left(\dfrac{r^m\partial_xv}{v^{\frac{5}{2}}}\right)\right|^2\dx\\
            &+C(\varepsilon,T)\left(\Vert g_\varepsilon u\Vert^2_{L^\infty([0,b])}+\Vert g_\varepsilon r^m\partial_xu
            \Vert^2_{L^\infty([0,b])}+\Vert g_\varepsilon r^m\partial_xu\Vert^{\frac52}_{L^2([0,b])}\Vert g_\varepsilon r^{2m}\partial_{xx}u\Vert^{\frac12}_{L^2([0,b])}\right).
\end{align*}
Now multiplying $\sigma(s)$ to the above inequality and integrating in $[0,T]$, the desired results can be obtained by estimates established before, and we only prove the term 
\begin{align*}
    \int^T_0\sigma(t)\Vert g_\varepsilon r^m\partial_xu\Vert^{\frac52}_{L^2([0,b])}\Vert g_\varepsilon r^{2m}\partial_{xx}u\Vert^{\frac12}_{L^2([0,b])}\dt.
\end{align*}
With the help of the Cauchy's inequality, we have
\begin{align*}
    &\int^T_0\sigma(t)\Vert g_\varepsilon r^m\partial_xu\Vert^{\frac52}_{L^2([0,b])}\Vert g_\varepsilon r^{2m}\partial_{xx}u\Vert^{\frac12}_{L^2([0,b])}\dt\\
&\le\int^T_0\sigma(t)\int^b_0g_\varepsilon r^{4m}|\partial_{xx}u|^2\dx+\int^T_0\sigma(t)\left(\int^b_0g_\varepsilon r^{2m}|\partial_xu|^2\dx\right)^{\frac53}\dt\\
    &\le\int^T_0\sigma(t)\int^b_0g_\varepsilon r^{4m}|\partial_{xx}u|^2\dx+\sup_{t\in[0,T]}\left(\sigma(t)\int^b_0g_\varepsilon r^{2m}|\partial_xu|^2\dx\right)^{\frac23}\int^T_0\int^b_0g_\varepsilon r^{2m}|\partial_xu|^2\dt\\
    &\le C(\varepsilon,T).
\end{align*}
\par The proof of Lemma \ref{time derivative of velocity} is finished now.
    \end{proof}
\end{lemma}
In addition, with the help of interpolation inequality and the Sobolev embedding, we have
\begin{lemma}
    For each fixed $\varepsilon>0$, there exists a positive constant $C(\varepsilon,T)$ such that
    \begin{align*}
\sup_{t\in[0,T]}\sigma(t)\left(\int^b_0g_\varepsilon|u|^4\dx+\int^b_0g_\varepsilon\left|\dfrac{r^m\partial_xv}{v^{\frac52}}\right|^4\dx\right)\le C(\varepsilon,T),\quad \sigma^{\frac14}(t)\left(|g_\varepsilon u(t,x)|+|g_\varepsilon r^m\partial_xv|\right)\le C(\varepsilon,T).
    \end{align*}
\end{lemma}
The following lemma is crucial to obtain the strong convergence of $\rho_r$
\begin{lemma}\label{time and space derivatives of v}
     For each fixed $\varepsilon>0$, there exists a positive constant $C(\varepsilon,T)$ such that
     \begin{align}
\int^T_0\int^b_0g_\varepsilon|\partial_tv|^2\dx\ds\le C(\varepsilon,T),\quad\int^T_0\int^b_0\sigma(s)g_\varepsilon|r^m\partial_{xt}v|^2\dx\ds\le C(\varepsilon,T).
     \end{align}
     \begin{proof}
         The first result is derived from mass equation of \eqref{approximate problem}. Differentiating the mass equation with respect to $x$, we obtain
         \begin{align*}
             \partial_{xt} v=\partial_{xx}(r^mu)=r^m\partial_{xx}u+\dfrac{2mv\partial_xu}{r}+\dfrac{m\partial_xvu}{r}-\dfrac{mv^2u}{r^{n+1}}.
         \end{align*}
         Therefore,
         \begin{align*}
             \int^T_0\int^b_0\sigma(s)g_\varepsilon|r^m\partial_{xt}v|^2\dx\ds&\le C(\varepsilon,T)\left[\int^T_0\int^b_0\sigma(s)g_\varepsilon\left(|r^{2m}\partial_{xx}u|^2+|r^m\partial_xu|^2\right)\dx\ds\right]\\
             &+C(\varepsilon,T)\int^T_0\Vert g_\varepsilon u\Vert^2_{L^\infty{(0,b)}}\int^b_0g_\varepsilon\left(|r^m\partial_xv|^2+|u|^2\right)\dx\ds\\
             &\le C(\varepsilon,T).
         \end{align*}
     \end{proof}
\end{lemma}

\section{Global Weak Solutions of  the Exterior Problem}\label{4}
\upshape
In this section, the global existence of weak solutions to the exterior problem is established via the uniform estimates obtained in \S \ref{3}. The estimates derived in Lagrangian coordinates are uniform with respect to the inner and outer radii, while the corresponding estimates in Eulerian coordinates are uniform with respect to both the inner and outer radii. These bounds allow us to pass to the limit as the outer radius tends to infinity, and thereby obtain a global solution.
Let  $(\widetilde{v}^{a,b},\widetilde{u}^{a,b},\widetilde{r}^{a,b})$ be solutions of approximate problems \eqref{approximate problem} for each fixed $(a,b)\in(0,1)\times (\mathscr{B},\infty)$. 

Now we define the coordinate transformations  $(t,r)=\mathcal{O}^{a,b}(s,x)$ for each $(s,x)\in[0,T]\times[0,b]$,
\begin{align*}
    r=\mathcal{O}^{a,b}_1(s,x):=\widetilde{r}^{a,b}(s,x),t=\mathcal{O}^{a,b}_2(s,x):=s.
\end{align*}
The corresponding Jacobian 
$J=\det \dfrac{\partial(\mathcal{O}^{a,b}_1,\mathcal{O}^{a,b}_2)}{\partial(x,s)}=\dfrac{\widetilde{v}^{a,b}}{(\widetilde{r}^{a,b})^{m}}$
satisfies the following bound derived for each $(t,x)\in[0,T]\times[0,b]$  from Lemma \ref{particle path for ab}, Lemma \ref{upper bound of v}, and Lemma \ref{lower bound of v}.
\begin{align}
    C^{-1}(a,T)\left(a^n+C(a,T)x\right)^{-\frac{m}{n}}\le J\le a^{-m}C (a,T).
\end{align}
Then the map $\mathcal{O}^{a,b}$ is a diffeomorphism so that by the inverse function theorem the inverse map $(\mathcal{O}^{a,b})^{-1}$ exists for each $(a,k)\in(0,b)\times (\mathscr{B},\infty)$. Thus we can define $(x^{a,b},s^{a,b})=(\mathcal{O}^{a,b})^{-1}(t,r)$ such that
\begin{align*}
    \widetilde{r}^{a,b}(x^{a,b}\left(t,r),t\right)=r,\quad\widetilde{s}^{a,b}(t,r)=t.
\end{align*}
The corresponding Jacobian is
\begin{align*}
    \dfrac{\partial(r,t)}{\partial(x,s)}=\dfrac{\widetilde{v}^{a,b}}{(\widetilde{r}^{a,b})^m}.
\end{align*}
For any smooth function $\widetilde{\phi}(t,x^{a,b}(t,r)):[0,T]\times[0,b]\rightarrow \mathbb{R}$, we define $\phi(t,r)=\widetilde{\phi}(t,x^{a,b}(t,r))$, hence the corresponding derivatives in Eulerian coordinates are
\begin{align}\label{e derivative and l derivative}
    &\phi_r(t,r)=\dfrac{\partial_x\widetilde{\phi}(x^{a,b}(t,r),t)(\widetilde{r}(x^{a,b}(t,r),t))^m}{\widetilde{v}^{a,b}(x^{a,b}(t,r),t)}=r^m\dfrac{\partial_x\widetilde{\phi}}{\widetilde{v}^{a,b}}(x^{a,b}(t,r),t).\\
    &\phi_t(t,r)=\partial_t\widetilde{\phi}(x^{a,b}(t,r),t)-r^m\left(\dfrac{\widetilde{u}^{a,b}\partial_x\widetilde{\phi}}{\widetilde{v}^{a,b}}(x^{a,b}(t,r),t)\right).
\end{align}
Then we convert the function in Lagrangian coordinates back into Eulerian coordinates by
\begin{align}\label{intermediate step}
    \left(\bar v^{a,b},\bar u^{a,b}\right)(t,r)=\left(\widetilde{v}^{a,b},\widetilde{u}^{a,b}\right)(x^{a,b}(t,r),t).
\end{align}
Next, we define the approximate solutions which extend the \eqref{intermediate step} to the whole exterior domain $[0,T]\times[a,\infty)$.
\begin{align}\label{approximate solution}
    \left(\rho^{a,b},u^{a,b}\right)=\left(\varphi^{a,b}\dfrac{1}{\bar v^{a,b}}+1-\varphi^{a,b},\varphi^{a,b}\bar u^{a,b}\right)
\end{align}
Here  with the smooth cut-off function $\chi\in C^\infty:\mathbb{R}\rightarrow[0,1]$ which satisfies: $\chi(\zeta)=1$ if $\zeta\le0$ and $\chi(\zeta)=0$ if $\zeta\ge1$. In addition, $|\chi'|\le2,|\chi'',\chi'''|\le C$.
\par The cutoff function$\varphi^{a,b}$ is defined by
\begin{align}
    \varphi ^{a,b}(t,r)=\chi\left(\dfrac{2r-\widetilde{r}^{a,b}(t,b)}{\widetilde{r}^{a,b}(t,b)}\right).
\end{align}
By the definition, we have
\begin{align}
    \varphi^{a,b}(t,r)=1 \text{ if } r\in[0,\dfrac{\widetilde{r}^{a,b}(t,b)}{2}],\quad\varphi^{a,b}(t,r)=0 \text{ if r} \in[\widetilde{r}^{a,b}(t,b),\infty).
\end{align}
With Lemma \ref{upper bound of v}, Lemma \ref{lower bound of v}, and the construction of approximate particle path $\widetilde{r}^{a,b}$, we have
\begin{align*}
    \widetilde{r}^{a,b}(t,b)&=\left(a^n+n\int^b_0\widetilde{v}^{a,b}(t,y)\dy\right)^{\frac{1}{n}}\ge \left(n\int^b_{\frac{1}{2}}\widetilde{v}^{a,b}(t,y)\dy\right)^{\frac{1}{n}}\\
    &\ge \left(\dfrac{n}{C(T)}\int^b_{\frac{1}{2}}\dy\right)^{\frac{1}{n}}\ge\dfrac{1}{C(T)}.\\
    \widetilde{r}^{a,b}(t,b)&=\left(a^n+n\int^b_0\widetilde{v}^{a,b}(t,y)\dy\right)^{\frac{1}{n}}\ge \dfrac{b^{\frac{1}{n}}}{C(a,T)}.
\end{align*}
Since 
\begin{align*}
    (\varphi^{a,b})_r(t,r)=\dfrac{2}{\widetilde{r}^{a,b}(t,b)}\chi'\left(\dfrac{2r-\widetilde{r}^{a,b}(t,b)}{\widetilde{r}^{a,b}(t,b)}\right).
\end{align*}
Combining the above estimates, we have
\begin{align}\label{cut off varphi}
    |(\varphi^{a,b})_r(t,r)|\le\min\{C(T),C(a,T)b^{-\frac{1}{n}}\},\quad\varphi^{a,b}_t(t,r)=0 \text{ for all }(t,r)\in[0,T]\times[a,\infty).
\end{align}
\subsection{Weak Solutions to the Exterior Domain}
\label{Weak solutions to the exterior domain}
\par Now we define the weak solutions $(\rho^a,u^a)$ to the exterior approximate problems on $I_{a,T}:=[0,T]\times[a,\infty)$ for a fixed $a\in(0,1)$, subject to the following boundary conditions
\begin{align*}
    u^a(t,a)=\rho^a_r(t,a)=0\text{ for all }t\in[0,T].
\end{align*}
The initial data for the  approximation problems $(\rho^a_0,u^a_0)$ have been constructed in \S \ref{3}.

\begin{definition}\label{5.1}
    Let 
    \begin{align*}
        &\mathscr{H}^a:=\bigg\{\Psi\in C^\infty([0,T]\times[a,\infty)):\exists~ N>a~s.t.~\Psi(t,r)=0\text{ for all } r\ge N  \bigg\}\\
        &\mathscr{H}_0^a:=\bigg\{\Psi\in\mathscr{H}^a:\Psi(t,a)=0~for~all ~t\in[0,T] \bigg\}.
    \end{align*}
    Then we say $(\rho^a,u^a)$ is a weak solution of exterior problems  if the following weak forms hold:
  \begin{enumerate}
      \item  For $\forall\Psi\in\mathscr{H}^a$,
    \begin{equation*}
        \int^\infty_a\rho^a(t,r)\Psi (t,r)r^m\dr= \int^\infty_a\rho_0^a(r)\Psi (0,r)r^m\dr+\int^t_0\int^\infty_a\big(\rho^a\Psi_t +\rho^au^a\Psi_r\big)r^m\dr\ds.
    \end{equation*}
    \item For $\forall \Psi\in \mathscr{H}_0^a$,
    \begin{align*}
        \int^\infty_a(\rho^au^a)(t,r)\Psi(t,r) r^m\dr&=\int^\infty_a(\rho^a_0u^a_0)(r)\Psi (0,r)r^m\dr+\int^t_0\int^\infty_a\rho^au^a(\Psi_t+u^a\Psi_r)r^m\dr\ds\\
        &+\int^t_0\int^\infty_a\left(A(\rho^a)^\gamma-\beta\left(u^a_r+\dfrac{mu^a}{r}\right)\right)\left(\Psi_r+\dfrac{m\Psi}{r}\right)r^m\dr\ds\\
        &-\kappa\int^t_0\int^\infty_a\rho^a(\rho^{a}_{rr}+\dfrac{m\rho^a_r}{r})\left(\Psi_r+\dfrac{m\Psi}{r}\right)r^m\dr\ds\\
        &-\kappa\int^t_0\int^\infty_a\rho^a_r(\rho^{a}_{rr}+\dfrac{m\rho^a_r}{r})\Psi r^m\dr\ds,
    \end{align*}
    where $\beta=2\mu+\lambda$.
  \end{enumerate}
\end{definition}
\begin{theorem}\label{exterior existence}
    For fixed $a\in(0,1)$, let $(\rho^a_0,u^a_0)$ be the initial data constructed in \S \ref{3}. Then, for each $T>0$, the exterior problem admits a global-in-time weak solution $(\rho^a,u^a)$ in the sense of \textnormal{Definition \ref{5.1}} with the following properties: 
    \begin{enumerate}
        \item The entropy estimates
    \begin{align*}
        \esssup_{t\in(0,T)}\int^{\infty}_a\left\{\dfrac{\rho^a}{2}|u^a|^2+G(\rho^a)+\dfrac{\kappa}{2}|\rho^a_r|^2\right\}r^m\dr+\beta\int^T_0\int^\infty_a\left|u^a_r+\dfrac{mu^a}{r}\right|^2r^m\dr\ds\le C(T).
    \end{align*}
    \item  For  $a.e. ~(t,r)\in[0,T]\times[a,\infty)$, we have
    \begin{equation*}
        C^{-1}(a,T)\le\rho^a\le C (a,T),\quad|u^a|+|\rho^a_r|\le C(a,T)\sigma^{-1/4}(t).
    \end{equation*}
    \item  There exists a continuous function $(t,x)\rightarrow\widetilde{r}^a(t,x)$ such that $x\rightarrow \widetilde{r}^a(t,x)$ is strictly increasing for each $t\in[0,T]$  and satisfies 
    \begin{align*}
        \begin{cases}
            x=\displaystyle\int^{\widetilde{r}^a(t,x)}_a\rho^ar^m\dr\qquad&\text{for }a.e.  (t,x)\in [0,T]\times[0,\infty),\\
            \widetilde{r}^a(t,x)=\widetilde{r}^a_0(t,x)+\displaystyle\int^t_0u^a(s,\widetilde{r}^a(s,x))\ds\qquad&\text{for }a.e.  (t,x)\in [0,T]\times[0,\infty),\\
            nx\psi^{-1}_{-}\left(\dfrac{C_0}{x}\right)\le(\widetilde{r}^a(t,x))^n\le C_0(1+x)\qquad&\text{for }a.e.  (t,x)\in [0,T]\times[0,\infty),
        \end{cases}
    \end{align*}
    where the function $x\rightarrow \widetilde{r}^a_0(x):[0,\infty)\rightarrow[a,\infty)$ is defined  by
    \begin{align*}
        x=\int^{\widetilde{r}^a_0(x)}_a\rho^a_0r^m\dr \text{  for }\forall x\in [0,\infty).
    \end{align*}
    \item For  each $\varepsilon\in(0,1]$, the particle path satisfies the $H\ddot{o}lder$ continuity for both $x$ and $t$:
    \begin{align*}
        \begin{cases}
            |\widetilde{r}^a(t,x_1)-\widetilde{r}^a(t,x_2)|\le C(\varepsilon,T)|x_1-x_2|\quad&\text{ for }(t,x_1,x_2)\in[0,T]\times[\varepsilon,\infty)^2,\\
            |\widetilde{r}^a(t_1,x)-\widetilde{r}^a(t_2,x)|\le C(\varepsilon,T)|t_1^{\frac{3}{4}}-t_2^{\frac{3}{4}}|\quad&\text{ for }(t_1,t_2,x)\in[0,T]^2\times[\varepsilon,\infty).
        \end{cases}
    \end{align*}
    \item For the solution $(\rho^a,u^a)$ in the  region $[\widetilde{r}^a(t,\varepsilon),\infty)\times[0,T]$, we have the following bounds:
    \begin{align*}
        C^{-1}(\varepsilon,T)\le \rho^a\le C(\varepsilon,T),\quad\sigma^{1/4}(t)|\rho^a_r|+|u^a|\le C(\varepsilon,T).
    \end{align*}
    Moreover, the following $H\ddot{o}lder$ estimates hold: 
    \begin{align*}
        \begin{cases}
            \sigma^{1/2}(t)|u^a(t,r_1)-u^a(t,r_2)|+|\rho^a(t,r_1)-\rho^a(t,r_2)|+\sigma^{1/2}(t)|\rho^a_r(t,r_1)-\rho^a_r(t,r_2)|\le C(\varepsilon,T)|r_1-r_2|^{\frac{1}{2}}\\
            \text{ for all } ~(t,r_1),(t,r_2)\in[0,T]\times[\widetilde{r}^a(t,\varepsilon),\infty),\\\\
            \sigma^{1/2}(t_1)|u^a(t_1,r)-u^a(t_2,r)|+|\rho^a(t_1,r)-\rho^a(t_2,r)|+\sigma^{1/2}(t_1)|\rho^a_r(t_1,r)-\rho^a_r(t_2,r)|\le C(\varepsilon,T)|t_1-t_2|^{\frac{1}{4}}\\\text{ for all} ~0<t_1<t_2\le T,\underset{t_1\le t\le t_2}{\sup}\widetilde{r}^a(t,\varepsilon)\le r.
        \end{cases}
        \end{align*}
        \item The following estimates hold
        \begin{align*}
            &\sup_{t\in[0,T]}\int^\infty_{a}\left\vert \left(\rho^a-1,\sigma^{\frac12}|u^a|^2,\sigma^{\frac12}|\rho^a_r|^2,\sigma^{\frac12}\rho^a_{rr},\sqrt{\sigma}u^a_r\right)\right\vert^2r^m\dr\\
  &+\int^T_0\int^\infty_{a}\left|\left(\rho^a_t,\sqrt{\sigma}\rho^a_{rt},\rho^a_{rr},u^a_r,|\rho^a_r|^2,\sigma^{\frac12}u^a_t,\sigma^{\frac12}\left(u^a_{rr},\rho^a_{rrr}\right)\right)\right|^2r^m\dr\ds\le C(a,T),\\
  &\sup_{t\in[0,T]}\int^\infty_{\widetilde{r}^a(t,\varepsilon)}\left\vert \left(\rho^a-1,\sigma^{\frac12}|u^a|^2,\sigma^{\frac12}|\rho^a_r|^2,\sigma^{\frac12}\rho^a_{rr},\sqrt{\sigma}u^a_r\right)\right\vert^2r^m\dr\\
  &+\int^T_0\int^\infty_{\widetilde{r}^a(t,\varepsilon)}\left|\left(\rho^a_t,\sqrt{\sigma}\rho^a_{rt},\rho^a_{rr},u^a_r,|\rho^a_r|^2,\sigma^{\frac12}u^a_t,\sigma^{\frac12}\left(u^a_{rr},\rho^a_{rrr}\right)\right)\right|^2r^m\dr\ds\le C(\varepsilon,T).
        \end{align*}
\end{enumerate}
\end{theorem}
\subsection{Estimates in Eulerian coordinates}\label{Estimates in Eulerian coordinate}
In this subsection, we obtain some estimates in Eulerian coordinates via the coordinate transformations.
\begin{lemma}\label{Entropy estimate of ab}
    There exists a constant $ N(a):=N(T,C_0,a)$  such that the following estimates hold for all $b\ge N(a)$:
    \begin{align}
        \sup_{t\in[0,T]}&\int^\infty_a\left\{\dfrac{\rho^{a,b}}{2}|u^{a,b}|^2+G(\rho^{a,b})+\kappa|\rho^{a,b}_r|^2\right\}r^m\dr+\beta\int^T_0\int^\infty_a\left|u^{a,b}_r+\dfrac{mu^{a,b}}{r}\right|^2r^m\dr\ds\le C(T).
    \end{align}
    \begin{proof}
        \begin{align*}
            &\int^\infty_a\left\{\dfrac{\rho^{a,b}}{2}|u^{a,b}|^2+G(\rho^{a,b})+\dfrac{\kappa}{2}|\rho^{a,b}_r|^2\right\}r^m\dr\\
            &\le\int^{\widetilde{r}^{a,b}(t,b)}_a\left(\dfrac{1}{2}|\bar u^{a,b}|^2+\bar v^{a,b} G\left((\bar v^{a,b})^{-1}\right)\right)(\bar v^{a,b})^{-1}r^m\dr+\int^{\widetilde{r}^{a,b}(t,b)}_a(1-\varphi^{a,b})\dfrac{1}{2}|\bar u^{a,b}|^2r^m\dr\\
            &+\kappa\Vert\varphi_r\Vert_{L^\infty}^2\int^{\widetilde{r}^{a,b}(t,b)}_a|\bar\rho^{a,b}-1|^2r^m\dr+\kappa\int^b_0\dfrac{|r^m\partial_xv|^2}{v^5}\dx:=\sum^4_{i=1}I_i.
        \end{align*}
        For $I_1$, it follows that
        \begin{align*}
            \int^{\widetilde{r}^{a,b}(t,b)}_a\left(\dfrac{1}{2}|\bar u^{a,b}|^2+\bar v^{a,b} G\left((\bar v^{a,b})^{-1}\right)\right)(\bar v^{a,b})^{-1}r^m\dr=\int^b_0\left(\dfrac{1}{2}|\widetilde{u}^{a,b}|^2+\psi(\widetilde{v}^{a,b})\right)\dx\le C_0,
        \end{align*}
        where we have used $yG\left(\dfrac{1}{y}\right)=\psi(y)$ and Lemma \ref{entropy estimate}.
        Now by Lemma \ref{particle path for ab}, Lemma \ref{upper bound of v}, and Lemma \ref{lower bound of v} for each fixed $a\in(0,1)$, we have
        \begin{align*}
            \left(a^n+\frac{n}{C(a,T)}x\right)^{\frac{1}{n}}\le\widetilde{r}^{a,b}(t,x)\le\left(a^n+nC(a,T)x\right)^{\frac{1}{n}}.
        \end{align*}
        Thus, if $b\ge\max\{2^nC^2(a,T)+\dfrac{2^n-1}{n}a^nC(a,T),\mathscr{B}\}$ then $\widetilde{r}^{a,b}(t,b)\ge2\widetilde{r}^{a,b}(t,1)$. Combining the construction of $\varphi^{a,b}$, for $I_2$, we have
        \begin{align*}
            I_2&=\int^{\widetilde{r}^{a,b}(t,b)}_a(1-\varphi^{a,b})\dfrac{1}{2}|\bar u^{a,b}|^2r^m\dr=\int^b_0\dfrac{\left(1-\varphi^{a,b}(\widetilde{r}^{a,b}(t,x),t)\right)}{2}|\widetilde{u}^{a,b}|^2\widetilde{v}^{a,b}\dx\\
           &\le\sup_{x\in[1,b]}\widetilde{v}^{a,b}\int^b_1\dfrac{1}{2}|\widetilde{u}^{a,b}|^2\dx\le C(T).
        \end{align*}
        The proof of $I_3+I_4$ is similar to the above and we leave this to the reader. Since $u^{a,b}=\varphi^{a,b}\bar u^{a,b}$, if we choose $b$ sufficiently large such that $b\ge\max\bigg\{\dfrac{2^n-1}{n}a^n C(a,T)+2^nC^2(a,T),C^{\frac{n}{2}}(a,T)C^{\frac{n}{2}}(T),\mathscr{B}\bigg\}$, it holds that
        \begin{align*}
            &\int^T_0\int^\infty_a\left|u^{a,b}_r+\dfrac{mu^{a,b}}{r}\right|^2r^m\dr\ds\\
            &\le2\Vert\varphi\Vert^2_{L^\infty}\int^T_0\int^{\widetilde{r}^{a,b}(t,b)}_{\widetilde{r}^{a,b}(t,1)}|\bar u^{a,b}|^2r^m\dr\ds+2\int^T_0\int^\infty_a\left|\bar u^{a,b}_r+\dfrac{m\bar u^{a,b}}{r}\right|^2r^m\dr\ds\\
            &\le \int^T_0\int^b_1|u|^2v\dx\ds+\int^T_0\int^b_0\dfrac{|\partial_x(r^mu)|^2}{v}\dx\ds\le C(T).
        \end{align*}
        \par The proof of Lemma \ref{Entropy estimate of ab} is completed.
    \end{proof}
\end{lemma}
\begin{lemma}\label{estimate in euler coordinate }
    For each fixed $a\in(0,1)$, there exists a positive constant $C(a,T)$ such that
    \begin{align*}
        &\sup_{b>\mathscr{B}}\sup_{t\in[0,T]}\int^\infty_{a}\left\vert \left(\rho^{a,b}-1,\sigma^{\frac12}|u^{a,b}|^2,\sigma^{\frac12}|\rho^{a,b}_r|^2,\sigma^{\frac12}\rho^{a,b}_{rr},\sqrt{\sigma}u^{a,b}_r\right)\right\vert^2r^m\dr\\
  &+\sup_{b>\mathscr{B}}\int^T_0\int^\infty_{a}\left|\left(\rho^{a,b}_t,\sqrt{\sigma}\rho^{a,b}_{rt},\rho^{a,b}_{rr},u^{a,b}_r,|\rho^{a,b}_r|^2,\sigma^{\frac12}u^{a,b}_t,\sigma^{\frac12}\left(u^{a,b}_{rr},\rho^{a,b}_{rrr}\right)\right)\right|^2r^m\dr\ds\le C(a,T).
    \end{align*}
    and we have the following point-wise bounds for $(t,r)\in[0,T]\times[a,\infty)$
    \begin{align*}
        C^{-1}(a,T)\le\rho^{a,b}\le C(a,T),\quad
        \sigma^{\frac{1}{4}}(t)\left(|u^{a,b}(t,r)|+|\rho_r^{a,b}(t,r)|\right)\le C(a,T).
    \end{align*}
    Moreover, for each fixed $\varepsilon\in(0,1)$, there exists a positive constant $C(\varepsilon,T)$ such that
 \begin{align*}
        &\sup_{b\ge\mathscr{B}}\sup_{t\in[0,T]}\int^\infty_{\widetilde{r}^{a,b}(t,\varepsilon)}\left\vert \left(\rho^{a,b}-1,\sigma^{\frac12}|u^{a,b}|^2,\sigma^{\frac12}|\rho^{a,b}_r|^2,\sigma^{\frac12}\rho^{a,b}_{rr},\sqrt{\sigma}u^{a,b}_r\right)\right\vert^2r^m\dr\\
  &+\sup_{b>\mathscr{B}}\int^T_0\int^\infty_{\widetilde{r}^{a,b}(t,\varepsilon)}\left|\left(\rho^{a,b}_t,\sqrt{\sigma}\rho^{a,b}_{rt},\rho^{a,b}_{rr},u^{a,b}_r,|\rho^{a,b}_r|^2,\sigma^{\frac12}u^{a,b}_t,\sigma^{\frac12}\left(u^{a,b}_{rr},\rho^{a,b}_{rrr}\right)\right)\right|^2r^m\dr\ds\le C(\varepsilon,T).
    \end{align*}
    and we have the following point-wise bound for $(t,r)\in[0,T]\times[\widetilde{r}^{a,b}(t,\varepsilon),\infty)$
    \begin{align*}
        C^{-1}(\varepsilon,T)\le\rho^{a,b}\le C(\varepsilon,T),\quad
        \sigma^{\frac{1}{4}}(t)\left(|u^{a,b}(t,r)|+|\rho_r^{a,b}(t,r)|\right)\le C(\varepsilon,T).
    \end{align*}
    \begin{proof}
        Now we suppress the superscripts and just prove the case for $[0,T]\times[a,\infty)$, since the proof is similar for $[0,T]\times[\widetilde{r}(t,\varepsilon),\infty)$.
        \begin{align*}
            \int^\infty_{a}|\rho-1|^2r^m\dr\le\int^{\widetilde{r}(k,t)}_a|\dfrac{1}{\bar v}-1|^2r^m\dr=\int^b_0\dfrac{|\widetilde{v}-1|^2}{\widetilde{v}}\dx\le C(a,T),
        \end{align*}
        where we used Lemma \ref{upper bound of v}, Lemma \ref{lower bound of v}, and Lemma \ref{exterior estimate of v-1}. With the same argument
        \begin{align*}
            \int^\infty_a \sigma(t)|u_r|^2r^m\dr&\le \int^{\widetilde{r}(k,t)}_{a}\sigma(t)|\bar u_r|^2r^m\dr+C\int^{\widetilde{r}(k,t)}_{a}\sigma(t)|\bar u|^2r^m\dr\\
            &\le C(a,T).
        \end{align*}
        Since $\rho_r=\varphi_r\left(\dfrac{1}{\bar v}-1\right)-\varphi\dfrac{\bar{v}_r}{\bar v^2}$,
        \begin{align*}
            \int^\infty_a \sigma(t)|\rho_r|^4r^m\dr&\le C(a,T) \left(\int^b_0\sigma(t)|\partial_xvr^m|^4\dx+C\int^b_{0}\sigma(t)|v-1|^4\dx\right)\\
            &\le C(a,T).
        \end{align*}
        Since $\rho_{rr}=\chi''\dfrac{4}{\widetilde{r}(t,b)^2}\left(\dfrac{1-\widetilde{v}}{\widetilde{v}}\right)-\chi'\dfrac{4}{\widetilde{r}(t,b)}\dfrac{\partial_x\widetilde{v} \widetilde{r}^m}{\widetilde{v}^3}+\chi\left(-\partial_x\left(\dfrac{\widetilde{r}^m\partial_x\widetilde{v}}{\widetilde{v}^{\frac52}}\right)\widetilde{r}^m\widetilde{v}^{-\frac32}+\dfrac{1}{2}\dfrac{|\partial_x\widetilde{v} \widetilde{r}^m|^2}{v^5}\right)$, we have
        \begin{align*}
            \sigma(t)\int^\infty_a|\rho_{rr}|^2r^m\dr&\le C(a,T)\left[\int^b_0|v-1|^2\dx+\sigma(t)\int^b_0\left(|r^m\partial_xv|^4\dx+\dfrac{r^{2m}}{v}\left|\partial_x\left(\dfrac{r^m\partial_xv}{v^{\frac52}}\right)\right|^2\right)\dx\right]\\
            &\le C(a,T).
        \end{align*}
        With the same argument, the remaining estimates follow by similar arguments.
  \end{proof}
\end{lemma}
\subsection{The Strong Convergence of Solutions}\label{The strong convergence of solutions}
From now on, we omit the parameter $a$ in this section to emphasize our goal of taking limit for the outer radii. We denote 
\begin{align}\label{simplified}
    (\rho^{b},u^{b}):=(\rho^{a,b},u^{a,b}),\quad(\bar v^{b},\bar u^{b}):=(\bar v^{a,b},\bar u^{a,b}),\quad (\widetilde{v}^{b},\widetilde{u}^{b},\widetilde{r}^b):=(\widetilde{v}^{a,b},\widetilde{u}^{a,b},\widetilde{r}^{a,b}).
\end{align}
\begin{lemma}\label{particle path limit of k}
    There exists a  H$\ddot{o}$lder continuous function $\widetilde{r}$ such that
    \begin{enumerate}
        \item  There exists a subsequence which is still denoted by $\widetilde{r}^b$ such that for any compact subset $\Omega\Subset[0,T]\times[0,\infty)$, 
    \begin{align}
        \lim_{b\rightarrow\infty}\sup_{(t,x)\in\Omega}|\widetilde{r}(t,x)-\bar r^b(t,x)|=0,
    \end{align}
    where
    \begin{align}
        \bar r^b(t,x)=\begin{cases}
            \widetilde{r}^b(t,x)&\qquad(t,x)\in[0,T]\times[0,b],\\
            x-b+ \widetilde{r}^b(t,b)&\qquad(t,x)\in[0,T]\times[b,\infty).
        \end{cases}
    \end{align}
    \item For  $\forall\varepsilon>0$, we have the following H$\ddot{o}$lder continuity:
    \begin{align}
        |\widetilde{r}(t,x_1)-\widetilde{r}(t,x_2)|\le C(\varepsilon,T)|x_1-x_2| \qquad\text{for }&\forall(t,x_1,x_2)\in[0,T]\times[\varepsilon,\infty)^2,\nonumber\\
        |\widetilde{r}(t_1,x)-\widetilde{r}(t_2,x)|\le C(\varepsilon,T)|t_1^{\frac{3}{4}}-t_2^{\frac{3}{4}}|\qquad\text{for }&\forall(t_1,t_2,x)\in[0,T]^2\times[\varepsilon,\infty).
    \end{align}
    \item The limit particle path satisfies the following bounds for each $(t,x)\in[0,T]\times[0,\infty)$
    \begin{align}
        nx\psi^{-1}_-(\dfrac{C_0}{x})\le\widetilde{r}^n(t,x)\le C_0(1+x).
    \end{align}
    \item $x\rightarrow\widetilde{r}(t,x)$ is strictly increasing for each $t\in[0,T]$.
    \end{enumerate}
    \begin{proof}
        The compactness result stated in Lemma \ref{2.1} is used in our proof. Now, for each integer $K\ge \mathscr{B}$, we define the corresponding $\mathbb{D}_K:=[0,T]\times[0,K]\subset[0,T]\times[0,\infty)$. Choosing $(t,x_1),(t,x_2)\in\mathbb{D}_K$, if $b\ge K$, we have $\bar r^b(t,x_i)=\widetilde{r}^b(t,x_i)$ for $i=1,2$. Thus, by the mean value theorem 
        \begin{align*}
            |\bar r^{b}(t,x_1)-\bar r^{b}(t,x_2)|\le|x_1-x_2|\sup_{x\in[0,b]}|\partial_x\widetilde{r}^b(t,x)|\le C(a,T)|x_1-x_2|,
        \end{align*}
        where we used  \eqref{constructure of particle path}, Lemma \ref{particle path for ab}, and Lemma \ref{upper bound of v}.
        For the case $b<K$, if $x_1,x_2\in[0,b]$ , we can obtain the same result using the above procedure. If $x_1\in[0,b],x_2\in[b,K]$, by the construction of $\bar{r}^b$, we have
        \begin{align*}
            |\bar r^{b}(t,x_1)-\bar r^{b}(t,x_2)|&\le|\bar r^{b}(t,x_1)-\bar r^{b}(t,b)|+|x_2-b|\\
            &\le C(a,T)(b-x_1)+x_2-b\le C(a,T)|x_2-x_1|,
        \end{align*}
        and if $x_1,x_2\in[b,K]$
        \begin{align*}
            |\bar r^{b}(t,x_1)-\bar r^{b}(t,x_2)|=|x_1-x_2|\le C(a,T)|x_1-x_2|.
        \end{align*}
        Thus, combining the above for each $K$ we have the H$\ddot{o}$lder continuity of $x$
        \begin{align*}
            |\bar r^{b}(t,x_1)-\bar r^{b}(t,x_2)|\le C(a,T)|x_1-x_2|.
        \end{align*}
        A direct calculation gives for each $b> \mathscr{B}$
        \begin{align*}
            |\bar r^{b}(t_1,x)-\bar r^{b}(t_2,x)|=\left|\int^{t_2}_{t_1}\bar u^b\ds\right|\le C(a,T)\int^{t_2}_{t_1}s^{-\frac{1}{4}}\ds\le C(a,T)t^{-\frac{1}{4}}_*|t_1-t_2|,   
        \end{align*}
        where $t_*=\min\{1,t_1,t_2\}$,$t_1,t_2\in(0,T]$. If $0=t_1<t_2\le T$, we have
        \begin{align*}
            |\bar r^{b}(0,x)-\bar r^{b}(t_2,x)|\le C(a,T)t^{\frac{3}{4}}_2\le C(a,T)t^{-\frac{1}{4}}_2|t_2-0|.
        \end{align*}
        Thus, in summary, for $\forall t_1,t_2\in[0,T],x\in[0,b]$, there exists a $t_{**}$ such that
        \begin{align*}
            &|\bar r^{b}(t_1,x)-\bar r^{b}(t_2,x)|\le C(a,T)t^{-\frac{1}{4}}_{**}|t_2-t_1|, \quad |\bar r^{b}(t_1,x)-\bar r^{b}(t_2,x)|\le C(a,T)|t^{\frac{3}{4}}_2-t^{\frac{3}{4}}_1|.
        \end{align*}
        By a similar analysis, there are some results for the region $[0,T]\times[\varepsilon,b]$.
        \begin{align*}
            &|\bar r^{b}(t,x_1)-\bar r^{b}(t,x_2)|\le C(\varepsilon,T)|x_1-x_2|\qquad\forall(t,x_1),(t,x_2)\in[0,T]\times[\varepsilon,K],\\
            &|\bar r^{b}(t_1,x)-\bar r^{b}(t_2,x)|\le C(\varepsilon,T)|t^{\frac{3}{4}}_2-t^{\frac{3}{4}}_1| \qquad\forall(t_1,x),(t_2,x)\in[0,T]\times[\varepsilon,K].
        \end{align*}
        Furthermore, for the point-wise bound for $\forall b>\mathscr{B}$, we have
        \begin{align*}
            \sup_{(t,x)\in \mathbb{D}_K}\bar r^b(t,x)\le \left(a^n+C(a,T)K\right)^{\frac{1}{n}}+K,
        \end{align*}
        where Lemma \ref{upper bound of v} is used.
        Now we apply Lemma \ref{2.1} for $\bar r^b$ to obtain a continuous function $\widetilde{r}(t,x):[0,T]\times[0,\infty)\rightarrow[a,\infty)$ such that for all compact subsets $\Omega\Subset[0,T]\times[0,\infty)$
        \begin{align*}
        \lim_{b\rightarrow\infty}\sup_{(t,x)\in\Omega}|\widetilde{r}(t,x)-\bar r^b(t,x)|=0.
    \end{align*}
    Based on this convergence, the point-wise bounds mentioned in Lemma $\ref{particle path limit of k}_{2,3}$ are proved directly. We next prove strict monotonicity in $x$, first, we assume $0<x_1<x_2<\infty$, and for all $b\ge\max\{2^nx_2C^2(a,T)+\dfrac{2^n-1}{n}C(a,T),\mathscr{B}\}$, we have
    \begin{align*}
        2\widetilde{r}^b(t,x_1)\le2\widetilde{r}^b(t,x_2)\le\widetilde{r}^b(t,b).
    \end{align*}
    Thus, 
    \begin{align*}
        x_2-x_1&=\int^{x_2}_{x_1}\dx=\int^{\widetilde{r}^k(t,x_2)}_{\widetilde{r}^b(t,x_1)}\dfrac{1}{\bar v^b}r^m\dr=\int^{\widetilde{r}^b(t,x_2)}_{\widetilde{r}^b(t,x_1)}\rho^br^m\dr\\
        &\le\dfrac{C(x_1,T)}{n}\big[(\widetilde{r}^b)^n(t,x_2)-(\widetilde{r}^b)^n(t,x_1)\big],
    \end{align*}
    where we used the bound $\rho^b\le C(x_1,T)$ for $r\in[\widetilde{r}^b(t,x_1),\infty)$ and this bound have been established in Lemma \ref{estimate in euler coordinate }.
    If $0=x_1<x_2<\infty$, the proof is similar to the above. Passing the limit to $b$ and the strong convergence of particle path, the proof of Lemma \ref{particle path limit of k} is finished.
    \end{proof}
\end{lemma}
Next, we establish the strong convergence of $(\rho^b,u^b,\rho^b_r)$.
\begin{lemma}\label{strong convergence of solution}
    There exist both a continuous function $\rho:[0,T]\times[a,\infty)\rightarrow [0,\infty)$, $(u,\rho_r):(0,T]\times[a,\infty)\rightarrow \mathbb{R}\times\mathbb{R}$ and a subsequence which is still denoted by $(\rho^b,u^b,\rho^b_r)$ such that 
    \begin{enumerate}
        \item  For any compact subset $\Omega_1\Subset[0,T]\times[a,\infty)$ and $\O_2\Subset(0,T]\times[a,\infty)$
    \begin{align}
        &\lim_{b\rightarrow\infty}\sup_{(t,r)\in \Omega_1}|\rho^b-\rho|(t,r)=0,\\
        &\lim_{b\rightarrow\infty}\sup_{(t,r)\in \Omega_2}|(\rho^b_r-\rho_r,u^b-u)|(t,r)=0.
    \end{align}
    \item For  $\forall\varepsilon>0$ if $t\in[0,T]$ and $\widetilde{r}(t,\varepsilon)\le r_1,r_2$ then
    \begin{align*}
        &\sigma^{\frac{1}{2}}(t)|u(t,r_1)-u(t,r_2)|+\sigma^{\frac12}(t)|\rho_r(t,r_1)-\rho_r(t,r_2)|+|\rho(t,r_1)-\rho(t,r_2)|\le C(\varepsilon,T)|r_1-r_2|^{\frac{1}{2}},
    \end{align*}
    and for all $0<t_1\le t_2\le T$ and $\sup_{t\in[t_1,t_2]}\widetilde{r}(t,\varepsilon)\le r$ we  have
    \begin{align*}
        \sigma^{\frac{1}{2}}(t_1)|u(t_1,r)-u(t_2,r)|+\sigma^{\frac12}(t_1)|\rho_r(t_1,r)-\rho_r(t_2,r)|+|\rho(t_1,r)-\rho(t_2,r)|\le C(\varepsilon,T)|t_1-t_2|^{\frac{1}{4}}.
    \end{align*}
    \item For  $\forall a\in(0,1)$, there exists a positive constant $C(a,T)$ such that
    \begin{align*}
       &|u(t,r)|+|\rho_r(t,r)|\le C(a,T)\sigma^{-\frac{1}{4}}(t), \text{ for }(t,r)\in(0,T]\times[a,\infty),\\
       &C^{-1}(a,T)\le \rho(t,r)\le C(a,T) \text{ for }(t,r)\in[0,T]\times[a,\infty).
    \end{align*}
    For  $\forall\varepsilon\in(0,1)$, there exists a positive constant $C(\varepsilon,T)$ such that
    \begin{align*}
        &|u(t,r)|+|\rho_r(t,r)|\le C(\varepsilon,T)\sigma^{-\frac{1}{4}}(t) \text{ for }(t,r)\in(0,T]\times[\widetilde{r}(t,\varepsilon),\infty),\\
       &C^{-1}(\varepsilon,T)\le \rho(t,r)\le C(\varepsilon,T) \text{ for }(t,r)\in[0,T]\times[\widetilde{r}(t,\varepsilon),\infty).
    \end{align*}
     \end{enumerate}
    \begin{proof}
        We only prove the convergence of $\rho$, since $\rho_r,u$ can be proved by a similar argument of $\rho$. 
        First, for each $(t,r_1,r_2)\in[0,T]\times[a,\infty)^2$
        \begin{align}\label{rho kr minus rhok r}
            |\rho^b(t,r_1)-\rho^b(t,r_2)|&=\left|\int^{r_2}_{r_1}\rho^b_r\dr\right|\le\left(\int^{r_2}_{r_1}|\rho^b_r|^2r^m\dr\right)^{\frac{1}{2}}\left(\int^{r_2}_{r_1}r^{-m}\dr\right)^{\frac{1}{2}}\nonumber\\
            &\le r^{-\frac{m}{2}}_1(r_2-r_1)^{\frac{1}{2}}C(T).
        \end{align}
        Now for $r\in[a,\infty),0\le t_1< t_2\le T$, we let $h=\sqrt{t_2-t_1}$,
        \begin{align*}
            \left|\dfrac{1}{h}\int^{r+h}_r(\rho^b(t_1,\zeta)-\rho^b(t_2,\zeta)\mathrm{d}\zeta\right|=|\rho^b(t_1,\xi)-\rho^b(t_2,\xi)|,
        \end{align*}
        where we used the mean-value theorem. By the triangle inequality, we have
        \begin{align*}
            |\rho^b(t_1,r)-\rho^b(t_2,r)|-|\rho^b(t_1,\xi)-\rho^b(t_2,\xi)|&\le\left|\int^\xi_r(\rho^b_r(t_1,\zeta)-\rho^b_r(t_2,\zeta)\mathrm{d}\zeta\right|\\
            &\le|\xi-r|^{\frac{1}{2}}r^{-\frac{m}{2}}\sum^2_{i=1}\left(\int^\xi_r|\rho^b_r(t_i,\zeta)|^2\zeta^m\mathrm{d}\zeta\right)^\frac{1}{2}\\
            &\le C(T)a^{-\frac{m}{2}}\sqrt{h}\le C(a,T)\sqrt{h}.
        \end{align*}
        Thus, it follows that
        \begin{align}\label{rho kt minus rhok t}
            |\rho^b(t_1,r)-\rho^b(t_2,r)|&\le C(a,T)\sqrt{h}+\dfrac{1}{h}\int^{r+h}_r\int^{t_2}_{t_1}|\rho^b_t|\dz\mathrm{d}\zeta\le C(a,T)\sqrt{h}+\dfrac{1}{h}\int^{t_2}_{t_1}\int^{r+h}_r|\rho^b_t|d\zeta \dz\nonumber\\
            &\le C(a,T)\sqrt{h}+\dfrac{1}{h}\left(\int^{t_2}_{t_1}\int^{r+h}_rd\zeta \dz\right)^{\frac{1}{2}}\left(\int^{t_2}_{t_1}\int^\infty_a|\rho^b_t|^2r^m\dr\dz\right)^{\frac{1}{2}}a^{-\frac{m}{2}}\nonumber\\
            &\le C(a,T)\sqrt{h}\le C(a,T)|t_2-t_1|^{\frac{1}{4}}.
        \end{align}
        Combining the point-wise bounds $C^{-1}(a,T)\le\rho^b\le C(a,T)$ and Lemma \ref{2.1}, we have the strong convergence result: there exists a continuous function $\rho\ge0$ such that for each compact subset $\Omega_1\Subset[0,T]\times[a,\infty)$
        \begin{align*}
            \lim_{b\rightarrow\infty}\sup_{(t,r)\in\Omega_1}|\rho^b-\rho|=0,
        \end{align*}
        furthermore, it satisfies the following estimate
        \begin{align*}
            &|\rho(t,r_1)-\rho(t,r_2)|\le C(a,T)|r_1-r_2|^{\frac12},\quad|\rho(t_1,r)-\rho(t_2,r)|\le C(a,T)|t_1-t_2|^\frac{1}{4},\\
           & C^{-1}(a,T)\le\rho\le C(a,T).
        \end{align*}
        Since $x\rightarrow\widetilde{r}(t,x)$ is strictly increasing, $d:=\inf_{t\in[0,T]}\{\widetilde{r}(t,\varepsilon)-\widetilde{r}(t,\dfrac{\varepsilon}{2})\}>0$, where we have used the continuity $(t,x)\rightarrow\widetilde{r}(t,x)$ and $[0,T]$ is a compact subset. Hence, there exists $N_\varepsilon$ such that if $b\ge N_\varepsilon$ we have
        \begin{align*}
            \sup_{t\in[0,T]}|\widetilde{r}^b(t,\dfrac{\varepsilon}{2})-\widetilde{r}(t,\dfrac{\varepsilon}{2})|\le\dfrac{d}{2}.
        \end{align*}
        This implies that if $b\ge N_\varepsilon,t\in[0,T]$, then $\widetilde{r}^b(t,\dfrac{\varepsilon}{2})<\widetilde{r}(t,\varepsilon)$. Now we choose $r_1,r_2\ge\widetilde{r}(t,\varepsilon)$, $t\in[0,T]$, thus $r_1,r_2>\widetilde{r}^b(t,\dfrac{\varepsilon}{2})$. By the same calculation of \eqref{rho kr minus rhok r}, we have
        \begin{align*}
            |\rho^b(t,r_2)-\rho^b(t,r_1)|\le C(\varepsilon,T)|r_2-r_1|^{\frac{1}{2}}.
        \end{align*}
        Also, by the same calculation of \eqref{rho kt minus rhok t}, we have the following estimate for $0\le t_1< t_2\le T$ and $r\ge\sup_{t\in[t_1,t_2]}\widetilde{r}(t,\varepsilon)$
        \begin{align*}
            |\rho^b(t_1,r)-\rho^b(t_2,r)|\le C(\varepsilon,T)|t_2-t_1|^{\frac{1}{4}},
        \end{align*}
        where we have used $r>\widetilde{r}^b(t,\dfrac{\varepsilon}{2})\ge C^{-1}(\varepsilon)$. Moreover, we have the point-wise bound for each $t\in[0,T]$ and $r\ge\widetilde{r}(t,\varepsilon)$
        \begin{align*}
            C^{-1}(\varepsilon,T)\le\rho^b\le C(\varepsilon,T).
        \end{align*}
        We pass $b\rightarrow\infty$ to obtain  parts (2) and (3) of Lemma $\ref{strong convergence of solution}$.
        \par The proof of Lemma \ref{strong convergence of solution} is finished.
    \end{proof}
\end{lemma}
Next, we prove  the conservation of mass.
\begin{lemma}\label{x}
    \begin{align}
        &x=\int^{\widetilde{r}(t,x)}_a\rho(t,r)r^m\dr\quad\text{ for each }t\in [0,T] \text{ and each }x>0.\nonumber\\
        &\widetilde{r}(t,x)=\widetilde{r}^a_0(x)+\int^t_0u(s,\widetilde{r}(s,x))\ds\quad\text{ for all }(t,x)\in[0,T]\times[a,\infty).
    \end{align}
    \begin{proof}
        \textbf{1.} For any fixed $(t,x)\in[0,T]\times[0,\infty)$, we choose $b$ sufficiently large such that 
        \begin{align*}
            b\ge \max\left\{\frac{2^n-1}{n }C(a,T)a^n+2^n C(a,T)^2x,\mathscr{B}\right\}
        \end{align*}
        then $\widetilde{r}^b(t,b)\ge2\widetilde{r}^b(t,x)$. Then we integrate $\rho^br^m$ in $[a,\widetilde{r}^b(t,x)]$ to obtain
        \begin{align*}
            \int^{\widetilde{r}^b(t,x)}_a\rho^b(t,r)r^m\dr=\int^{\widetilde{r}^b(t,x)}_a\dfrac{1}{\bar v^b}r^m\dr=x,
        \end{align*}
        where we used $\rho^b(t,r)=\dfrac{1}{\bar v^b}$ for $b\ge\frac{2^n-1}{n }C(a,T)a^n+2^n C(a,T)^2x$.
        Thus,
        \begin{align*}
            x-\int^{\widetilde{r}(t,x)}_a\rho(t,r)r^m\dr&=\int^{\widetilde{r}^b(t,x)}_a\rho^b(t,r)r^m\dr-\int^{\widetilde{r}(t,x)}_a\rho(t,r)r^m\dr\\
            &=\int^{\widetilde{r}^b(t,x)}_{\widetilde{r}(t,x)}\rho^b(t,r)r^m\dr+\int^{\widetilde{r}(t,x)}_a\big[\rho^b(t,r)-\rho(t,r)\big]r^m\dr\\
            &=\sum^2_{i=1}I_i
        \end{align*}
For $I_1$, with the help of the point-wise bound of $\rho^b$ mentioned in Lemma \ref{estimate in euler coordinate } and the strong convergence of particle path mentioned in Lemma \ref{particle path limit of k}, we have
\begin{align*}
    |I_1|=\left|\int^{\widetilde{r}(t,x)}_{\widetilde{r}^b(t,x)}\rho^b(t,r)r^m\dr\right|&\le C(a,T)\frac{1}{n}\sup_{t\in[0,T]}\big|\widetilde{r}^n(t,x)-(\widetilde{r}^b(t,x))^n\big|\rightarrow0.
\end{align*}
as $b\rightarrow\infty$. For $I_2$, it follows from Lemma \ref{strong convergence of solution} and the dominated convergence theorem that
\begin{align*}
   I_2\rightarrow0 \text{ as } b\rightarrow0. 
\end{align*}
 The first part of Lemma \ref{x} is proved. \\
 \textbf{2. } Notice  that for each $(t,x)\in[0,T]\times[0,b]$,
\begin{align*}
    \widetilde{r}^b(t,x)=\widetilde{r}^b(0,x)+\int^t_0\widetilde{u}^b(s,x)\ds,
\end{align*}
where $\widetilde{r}^b(0,x)$ is defined by \eqref{approximate particle path}
\begin{equation*}
    \widetilde{r}^b(0,x)=\bigg(a^n+n\int^x_0\frac{1}{\rho^{a,b}_0(\widetilde{r}^a_0(y))}\dy\bigg)^{\frac{1}{n}}.
\end{equation*}
For any fixed $(t,x)\in[0,T]\times[0,\infty)$, if $b\ge\max\left\{ \dfrac{2^n-1}{n }C(a,T)a^n+2^n C(a,T)^2x,\mathscr{B}\right\}$, then 
\begin{align*}
    u^b(t,\widetilde{r}^b(t,x))=\widetilde{u}^b(t,x).
\end{align*}
In addition,
\begin{align}\label{rk}
    \widetilde{r}^b(t,x)=\widetilde{r}^b(0,x)+\int^t_0\widetilde{u}^b(s,\widetilde{r}^b(s,x))\ds
\end{align}
If we choose $b$ to be sufficiently large as before, we can obtain $2\widetilde{r}^a_0(x)\le\widetilde{r}^a_0(b)$.
Also, we have
\begin{align*}
    \widetilde{v}^{a,b}_0(x)=\dfrac{1}{\rho^{a,b}_0(\widetilde{r}^a_0(x))}=\dfrac{1}{\rho^a_0(\widetilde{r}^a_0(x))\varphi^{a,b}_0(\widetilde{r}^a_0(x))+1-\varphi^{a,b}_0(\widetilde{r}^a_0(x))}=\dfrac{1}{\rho^a_0(\widetilde{r}^a_0(x))}.
\end{align*}
Thus, we choose $b$ as before, and it follows that
\begin{align}\label{r=0}
\widetilde{r}^b(0,x)=\bigg(a^n+n\int^x_0\frac{1}{\rho^{a}_0(\widetilde{r}^a_0(y))}\dy\bigg)^{\frac{1}{n}}=\widetilde{r}^a_0(x).
\end{align}
Furthermore, by virtue of Lemma \ref{strong convergence of solution}, it holds for each $(s,x)\in[0,T]\times[0,\infty)$
\begin{align*}
    |u^b(s,\widetilde{r}^b(s,x))-u^b(s,\widetilde{r}(s,x))|\le\sigma^{-\frac{1}{4}}(s)C(a,T)|\widetilde{r}^b(s,x)-\widetilde{r}(s,x)|^{\frac12}.
\end{align*}
Thus, we have
\begin{align*}
    |u^b(s,\widetilde{r}^b(s,x))-u(s,\widetilde{r}(s,x))|&\le|u^b(s,\widetilde{r}^b(s,x))-u^b(s,\widetilde{r}(s,x))|+|u^b(s,\widetilde{r}(s,x))-u(s,\widetilde{r}(s,x))|\\
    &\le\sigma^{-\frac{1}{4}}(s)C(a,T)|\widetilde{r}^b(s,x)-\widetilde{r}(s,x)|+|u^b(s,\widetilde{r}(s,x))-u(s,\widetilde{r}(s,x))|.
\end{align*}
For $\forall t>0$, it follows that
\begin{align*}
    \left|\int^t_0\big[u^b(s,\widetilde{r}^b(s,x))-u(s,\widetilde{r}(s,x))\big]\right|\ds&\le C(a,T)\sup_{s\in[0,T]}|\widetilde{r}^b(s,x)-\widetilde{r}(s,x)|\int^t_0\sigma^{-\frac{1}{4}}(s)\ds\\
    &+\left|\int^t_0u^b(s,\widetilde{r}(s,x))-u(s,\widetilde{r}(s,x))\ds\right|.
\end{align*}
The first term tends to zero by virtue of Lemma \ref{particle path limit of k}. For the second term, since $u^b(t,r)\rightarrow u(t,r)$ as $b\rightarrow\infty$ for  $a.e. ~t\in[0,T]$ and $u^b,u\in L^1(0,T)$, we can apply the dominated convergence theorem to obtain 
\begin{align*}
    \left|\int^t_0u^b(s,\widetilde{r}(s,x))-u(s,\widetilde{r}(s,x))\ds\right|\rightarrow0\text{ as  } b \rightarrow\infty.
\end{align*}
Thus, by passing the limit to \eqref{rk} for each $t>0$, we have
\begin{align*}
    \widetilde{r}(t,x)=\widetilde{r}^a_0(x)+\int^t_0u(s,\widetilde{r}(s,x))\ds.
\end{align*}
The proof for $t=0$ is straightforward from $\left|\displaystyle\int^t_0(u^b,u)\ds\right|\le C(a,T)t^{\frac{3}{4}}$ and \eqref{r=0}.
\par The proof of Lemma \ref{x} is finished.
    \end{proof}
\end{lemma}

Now based on the convergence results constructed before, we have 
\begin{lemma}\label{weak convergence}
There exists a weak convergence subsequence (still denoted by $b$) $(\rho^k,u^k)$ such that
    \begin{align*}
        &\bigg(\rho^b-1,\sqrt{\sigma}u^b_r,\sigma^{\frac14}u^b,\sqrt{\sigma}\rho^b_{rr},\sigma^{\frac14}\rho^b_r\bigg)\overset{*}{\rightharpoonup}(\rho-1,\sqrt{\sigma}u_r,\sigma^{\frac14}u,\sqrt{\sigma}\rho_{rr},\sigma^{\frac14}\rho_r)\text{ in }L^\infty(0,T;L^2([a,\infty);r^m\dr)),\\
        &\bigg(u^b_r,\sqrt{\sigma}\rho^b_{rt},\rho^b_{rr},\bigg)\rightharpoonup\bigg(u_r,\sqrt{\sigma}\rho_{rt},\rho_{rr},\bigg)\text{ in } L^2(0,T;L^2((a,\infty);r^m\dr)),
        \\&\bigg(\sqrt{\sigma}u^b_t,|\rho^b_r|^2,\rho^b_t,\sqrt{\sigma}u^b_{rr},\sqrt{\sigma}\rho^b_{rrr}\bigg)\rightharpoonup\bigg(\sqrt{\sigma}u_t,|\rho_r|^2,\rho_t,\sqrt{\sigma}u_{rr},\sqrt{\sigma}\rho_{rrr}\bigg)\text{ in }L^2(0,T;L^2([a,\infty);r^m\dr)).
    \end{align*}
    \begin{proof}
        We only prove $\rho^k-1$  since other terms can be proved with the same argument.  From Lemma \ref{estimate in euler coordinate },
        \begin{align*}
            \sup_{b\ge \mathscr{B}}\sup_{t\in[0,T]}\int^\infty_a|\rho^b-1|^2r^m\dr\le C(a,T),
        \end{align*}
        there exists a subsequence (still denoted by $b$) and a function $f\in L^\infty(0,T;L^2([a,\infty);\dr))$ such that
        \begin{align*}
            r^{\frac{m}{2}}(\rho^b-1)\rightharpoonup f \text{ in }L^\infty(0,T;L^2([a,\infty)).
        \end{align*}
        Hence, for $\forall\Phi\in C^\infty_c([0,T]\times[a,\infty))$,  we have
        \begin{align*}
            \int^T_0\int^\infty_a(\rho^b-1)r^{\frac{m}{2}}\Phi \dr\ds\rightarrow  \int^T_0\int^\infty_af\Phi \dr\ds \text{ as }k\rightarrow\infty.
        \end{align*}
        By Lemma \ref{strong convergence of solution} and the dominated convergence theorem, we have
        \begin{align*} 
            \int^T_0\int^\infty_a(\rho^b-1)r^{\frac{m}{2}}\Phi \dr\ds\rightarrow \int^T_0\int^\infty_a(\rho-1)r^{\frac{m}{2}}\Phi \dr\ds.
        \end{align*}
       Thus, $f=(\rho-1)r^{\frac{m}{2}}$ for  $a.e. ~(t,r)\in[0,T]\times[a,\infty)$.  
       \par The proof is finished.
    \end{proof}
\end{lemma}
Now we have the following estimate from \textnormal{Lemma} \ref{weak convergence}.
\begin{lemma}\label{outer estimate}
For each fixed $a\in(0,1)$, there exists a positive constant $C(a,T)$ such that
    \begin{align}\label{estimate a}
        &\sup_{t\in[0,T]}\int^\infty_a\left(|\rho-1|^2+\sigma |u|^4+\sigma|u_r|^2+\sigma|\rho_{rr}|^2+\sigma|\rho_r|^4\right)r^m\dr\nonumber\\
        &+\int^T_0\int^\infty_a\left(|u_r|^2+|\rho_{rr}|^2+|\rho_r|^4+\sigma|(u_t,u_{rr},\rho_{rrr},\rho_{tr})|^2+|\rho_t|^2\right)r^m\dr\ds\nonumber\\
        &\le C(a,T).
    \end{align}
    For each fixed $\varepsilon\in(0,1)$, there exists a positive constant $C(\varepsilon,T)$ such that
        \begin{align}\label{estimate epsilon}
        &\sup_{t\in[0,T]}\int^\infty_{\widetilde{r}(t,\varepsilon)}\left(|\rho-1|^2+\sigma |u|^4+\sigma|u_r|^2+\sigma|\rho_{rr}|^2+\sigma|\rho_r|^4\right)r^m\dr\nonumber\\
        &+\int^T_0\int^\infty_{\widetilde{r}(t,\varepsilon)}\left(|u_r|^2+|\rho_{rr}|^2+|\rho_r|^4+\sigma|(u_t,u_{rr},\rho_{rrr},\rho_{tr})|^2+|\rho_t|^2\right)r^m\dr\ds\nonumber\\
        &\le C(\varepsilon,T).
    \end{align}
    \begin{proof}
        By virtue of Lemma \ref{weak convergence} and the weakly lower semi-continuity of Sobolev norm,  \eqref{estimate a} hold, we only need to prove for  \eqref{estimate epsilon}. We only prove for the term $\sigma|\rho_r|^4$ and other terms can be proved by the same argument. For  $\forall \Phi\in L^1(0,T;L^{\frac43}([a,\infty);r^m\dr)$,
        \begin{align*}
            &\left|\int^T_0\int^\infty_{\widetilde{r}^b(t,\varepsilon)}\rho^b_r\Phi r^m\dr\ds-\int^T_0\int^\infty_{\widetilde{r}(t,\varepsilon)}\rho_r\Phi r^m\dr\ds\right|\\
            &\le\left|\int^T_0\int^\infty_{\widetilde{r}(t,\varepsilon)}(\rho^b_r-\rho_r)\Phi r^m\dr\ds\right|+\int^T_0\int^{\widetilde{r}(t,\varepsilon)}_{\widetilde{r}^b(t,\varepsilon)}|\Phi\rho^b_r|r^m\dr\ds:=\sum^2_{i=1}I_i
        \end{align*}
        By virtue of the dominated convergence theorem and Lemma \ref{strong convergence of solution},
        \begin{align*}
            I_1\rightarrow0 \text{ as }b\rightarrow\infty,\quad I_2\rightarrow0 \text{ as }b\rightarrow\infty.
        \end{align*}
        Hence, it holds that $\si^{\frac14}\chi^b_\varepsilon\rho^b_r\rightharpoonup\si^\frac14\chi_\varepsilon\rho_r$ in $L^\infty(0,T;L^4([a,\infty);r^m\dr)$, where $\chi_\varepsilon$ and $\chi^b_\varepsilon$ is defined by
        \begin{align*}
            &\chi_\varepsilon(t,r)=\begin{cases}
                1&\qquad (t,r)\in[0,T]\times[\widetilde{r}(t,\varepsilon),\infty),\\
                0&\qquad\text{ otherwise.}
            \end{cases}
            &\chi^b_\varepsilon(t,r)=\begin{cases}
                1&\quad (t,r)\in[0,T]\times[\widetilde{r}^b(t,\varepsilon),\infty),\\
                0&\quad\text{ otherwise.}
            \end{cases}
        \end{align*}
        Thus with the help of the weakly lower semi-continuity of Sobolev norm and the interpolation 
        \begin{align*}
\Vert\sigma^{\frac14}\chi_\varepsilon\rho_r\Vert_{L^\infty((0,T);L^4([a,\infty);r^m\dr)}\le C(\varepsilon,T).
        \end{align*}
        \par The proof of Lemma \ref{outer estimate} is finished.
    \end{proof}
\end{lemma}
\subsection{Global Weak Solutions to the Exterior Problems }\label{Global weak solution of exterior problem}
In  this section, we establish the global existence of a weak solution to the exterior problem for each fixed $a\in(0,1)$. We still use the simplified notation \eqref{simplified} to emphasize the limit of the outer radii $b$. We start with the entropy estimates.
\begin{lemma}\label{limit entropy}
     For any $a\in(0,1)$, there exists $C(T)>0$ such that 
    \begin{align}\label{limit entropy estimate of exterior problem}
        &\esssup_{t\in[0,T]}\int^\infty_a\left(\dfrac{\rho}{2}|u|^2+G(\rho)+\kappa|\rho_r|^2\right)r^m\dr+\beta\int^T_0\int^\infty_a\left|u_r+\dfrac{mu}{r}\right|^2r^m\dr\ds\le C(T).
    \end{align}
    \begin{proof}
        First, let the test function $ \Phi_L\in C^\infty_c([a,\infty))$  satisfy $\Phi_L(r)=1$ for $a\le r\le L-1$, and $\Phi_L(r)=0$ for $ r\ge L$. By the construction of $\Phi_L$ and Lemma \ref{estimate in euler coordinate },  we have 
    \begin{align*}
        &\left |\int^\infty_a\rho^b(u^b)^2\Phi_Lr^m\dr-\int^\infty_a\rho u^2\Phi_Lr^m\dr\right|\\
        &\le \left|\int^\infty_a(\rho-\rho^b)(u^b)^2\Phi_Lr^m\dr\right|+\left|\int^\infty_a\rho\big(u^2-(u^b)^2\big)r^m\Phi_L\dr\right|:=\mathcal{I}_1+\mathcal{I}_2.
    \end{align*}
    we can derive $\mathcal{I}_1,\mathcal{I}_2\rightarrow0$ as $b\rightarrow\infty $ via the dominated convergence theorem and Lemma \ref{strong convergence of solution}. Thus for each $L>0$ we have 
    \begin{align*}
        \int^\infty_a\rho u^2\Phi_Lr^m\dr\le C(T),
    \end{align*}
    where we used the weakly lower-semi continuity of Sobolev norm and the dominated convergence theorem. Since $\rho u^2\ge0$, it follows from the monotone convergence theorem that 
    \begin{align*}
        \int^\infty_a\rho u^2r^m\dr\le C(T).
    \end{align*}
    Since $|\rho^b_r|+|\rho_r|\le C(a,T)$ and $\rho^b_r(t,r)\rightarrow\rho_r(t,r)$, with the same argument, we have
    \begin{align*}
        \int^\infty_a|\rho_r|^2r^m\dr\le C(T).
    \end{align*}
    Since $\rho^b\rightarrow\rho$ and $C^{-1}(a,T)\le\rho+\rho^b\le C(a,T)$, it follows from the dominated convergence theorem and monotone convergence theorem that 
    \begin{align*}
        \int^\infty_aG(\rho)r^m\dr\le C(T).
    \end{align*}
    Next, let $\Phi\in L^2(0,T;L^2([a,\infty);r^m\dr))$,
    \begin{align*}
        &\int^T_0\int^\infty_a\left\{\left(u^b_r+\dfrac{mu^b}{r}\right)-\left(u_r+\dfrac{mu}{r}\right)\right\}\Phi r^m\dr\ds\\
        &=\int^T_0\int^\infty_a(u^b_r-u_r)\Phi r^m\dr\ds+\int^T_0\int^\infty_a\dfrac{m}{r}(u^b-u)r^m\Phi\dr\ds\\
        &:=I_3+I_4.
    \end{align*}
    Since $u^b_r\rightharpoonup u_r$ in $L^2(0,T;L^2([a,\infty);r^m\dr))$, $I_3\rightarrow0$ as $b\rightarrow\infty$. Since $u^b\rightarrow u$, $I_4\rightarrow0$.
    \par The proof of Lemma \ref{limit entropy} is finished.
    \end{proof}
    \end{lemma}
Now we prove the pair $(\rho,u)$ satisfies the mass equation  and the momentum equation in the sense of distribution.
\begin{lemma}\label{weak solution of exterior problem}
    $(\rho,u)$ is a weak solution  stated in \textnormal{Definition} \ref{5.1}.
    \begin{proof}
        The proof is divided into 2 steps. \\
        \textbf{1. The weak form of the mass equation.} For any $\Psi \in \mathscr{H}^a$ there exists a fixed $L>0,L\in\mathbb{N}$ such that 
        \begin{align*}
            \Psi(t,r)=0 \text{ if }r>L \text{ and } t\in[0,T].
        \end{align*}
        Thus if we choose $b$ sufficiently large such that $b\ge\max\left\{2^nL^n\max\{C_0,C(a,T)\},\mathscr{B}\right\}$ it holds that
        \begin{align*}
            \widetilde{r}^b(t,b)>2L,\quad\widetilde{r}^a_0(b)>2L.
        \end{align*}
         If $(t,r)\in[0,T]\times[a,L]$ and $b$ sufficiently large, we have
        \begin{align*}
            &\rho^b(t,r)=\dfrac{1}{\bar v^b(t,r)},\quad u^b(t,r)=\bar u^b(t,r),\\
&\rho^{a,b}_0=\varphi^{a,b}_0\rho^a_0+1-\varphi^{a,b}_0=\rho^a_0,\quad u^{a,b}_0=\varphi^{a,b}_0u^a_0=u^a_0.
        \end{align*}
        Furthermore, for $b\ge\max\left\{2^nL^n\max\{C_0,C(a,T)\},\mathscr{B}\right\}$ and $x\in[0,\infty)$ if $\widetilde{r}^a_0(x)\in\supp(\Psi(\cdot,0))$,  it follows that
        \begin{align*}
            \widetilde{r}^a_0(x)=\bigg(a^n+n\int^x_0\dfrac{1}{\rho^a_0(\widetilde{r}^a_0(y))}\dy\bigg)^{\frac{1}{n}}=\bigg(a^n+n\int^x_0\widetilde{v}^{a,b}_0\dy\bigg)^{\frac{1}{n}}=\widetilde{r}^k(0,x).
        \end{align*}
        Thus if $b$ has lower bound stated before, for continuity equation, it follows that 
        \begin{align*}
            \int^\infty_a\rho^b(t,r)\Psi(t,r)r^m\dr-\int^\infty_a\rho^a_0\Psi(0,r)r^m\dr&=\int^{\widetilde{r}^b(t,b)}_a\dfrac{1}{\bar v^b}(t,r)\Psi(t,r) r^m\dr-\int^{\widetilde{r}^a_0(b)}_a\rho^a_0\Psi(0,r)r^m\dr\\
            &=\int^b_0\Psi(\widetilde{r}^b(t,x),t)\dx-\int^b_0\Psi(\widetilde{r}^b(0,x),0)\dx\\
            &=\int^t_0\int^b_0(u^b\Psi_r+\Psi_t)(\widetilde{r}^b(s,x),s)\dx\ds\\
            &=\int^t_0\int^{\widetilde{r}^b(t,b)}_a\dfrac{1}{\bar v^b}(\bar u^b\Psi_r+\Psi_t)r^m\dr\ds\\
            &=\int^t_0\int^\infty_a(\rho^bu^b\Psi_r+\rho^b\Psi_t)r^m\dr\ds.
        \end{align*}
        Now we apply Lemma \ref{strong convergence of solution} and the dominated convergence theorem to the first term of the left hand side also to the first and the second term of the right hand side. For details, since $\rho^b\rightarrow\rho $ and $C^{-1}(a,T)\le\rho^b,\rho\le C(a,T)$, with the help of dominated convergence theorem, we have 
        \begin{align*}
            \int^\infty_a\rho^b(t,r)\Psi(t,r)r^m\dr\rightarrow\int^\infty_a\rho(t,r)\Psi(t,r)r^m\dr
        \end{align*}
        we can prove the convergence of the right hand side by the same argument. Thus we obtain the Definition $\ref{5.1}_1$
        \begin{align*}
            \int^\infty_a\rho(t,r)\Psi(t,r)r^m\dr-\int^\infty_a\rho^a_0\Psi(0,r)r^m\dr=\int^t_0\int^\infty_a(\rho u\Psi_r+\rho\Psi_t)r^m\dr\ds.
        \end{align*}
        \textbf{2. The weak form of the momentum equation.} For all $\Psi\in\mathscr{H}^a_0$,  there exists $L>0$ such that
        $\supp (\Psi(t,r))\subset[a,L]$  and $\Psi(t,a)=0$ for all $t\in[0,T]$. If $b\ge\max\left\{2^nL^n\max\{C_0,C(a,T)\},\mathscr{B}\right\}$, then $\widetilde{r}^b(t,b),\widetilde{r}^a_0(b)>2L$. Thus 
        \begin{align*}
            &\int^\infty_a\rho^b(t,r) u^b(t,r)\Psi(t,r) r^m\dr-\int^\infty_a\rho^a_0u^a_0\Psi(0,r)r^m\dr\\
            &=\int^{\widetilde{r}^b(t,b)}_a\dfrac{1}{\bar v^b}\bar u^b\Psi r^m\dr-\int^{\widetilde{r}^a_0(b)}_a\rho^{a,b}_0u^{a,b}_0\Psi(0,r)r^m\dr\\
            &=\int^b_0\widetilde{u}^b(t,x)\Psi(t,\widetilde{r}(t,x))\dx-\int^b_0\widetilde{u}^{a,b}_0(x)\Psi(0,\widetilde{r}^a_0(x))\dx\\
            &=\int^t_0\int^b_0\partial_t\widetilde{u}^b\Psi(t,\widetilde{r}(t,x))\dx\ds+\int^t_0\int^b_0\widetilde{u}^b\Psi_t(s,\widetilde{r}^b(s,x))\dx\ds+\int^t_0\int^b_0\Psi_r(s,\widetilde{r}^b(s,x))(\widetilde{u}^b)^2\dx\ds\\
            &:=\sum^3_{i=1}J_i.
        \end{align*}
        Combining the momentum equation of \eqref{approximate problem}, we have
        \begin{align*}
            J_1&=\int^t_0\int^b_0\partial_t\widetilde{u}^b\Psi(t,\widetilde{r}(t,x))\dx\ds\\
            &=\int^t_0\int^b_0\Psi(t,\widetilde{r}(t,x))(\widetilde{r}^b)^m\partial_x\left(\beta\dfrac{\partial_x\left((\widetilde{r}^b)^m\widetilde{u}^b\right)}{\widetilde{v}^b}-\kappa\dfrac{(\widetilde{r}^b)^{2m}\partial_x\widetilde{v}^b}{(\widetilde{v}^b)^5}-P^b-\dfrac{5\kappa}{2}\dfrac{(\widetilde{r}^b)^{2m}|\partial_x\widetilde{v}^b|^2}{(\widetilde{v}^b)^6}\right)\dx\ds\\
            &-\int^t_0\int^b_0\Psi(t,\widetilde{r}(t,x))m\kappa\dfrac{(\widetilde{r}^b)^{2m-1}|\widetilde{v}^b|^2}{(\widetilde{v}^b)^5}\dx\ds\\
            &:=\sum^2_{i=1}J^i_1.
        \end{align*}
        For $J^1_1$, it holds that
        \begin{align*}
            J^1_1&=\int^t_0\int^b_0\Psi(t,\widetilde{r}(t,x))(\widetilde{r}^b)^m\partial_x\left(\beta\dfrac{\partial_x\left((\widetilde{r}^b)^m\widetilde{u}^b\right)}{\widetilde{v}^b}-\kappa\dfrac{(\widetilde{r}^b)^{2m}\partial_x\widetilde{v}^b}{(\widetilde{v}^b)^5}-P^b-\dfrac{5\kappa}{2}\dfrac{(\widetilde{r}^b)^{2m}|\partial_x\widetilde{v}^b|^2}{(\widetilde{v}^b)^6}\right)\dx\ds\\
            &=-\int^t_0\int^b_0\left(\beta\dfrac{\partial_x\left((\widetilde{r}^b)^m\widetilde{u}^b\right)}{\widetilde{v}^b}-\kappa\dfrac{(\widetilde{r}^b)^{2m}\partial_x\widetilde{v}^b}{(\widetilde{v}^b)^5}-P^b-\dfrac{5\kappa}{2}\dfrac{(\widetilde{r}^b)^{2m}|\partial_x\widetilde{v}^b|^2}{(\widetilde{v}^b)^6}\right)\left(\dfrac{\partial_x\Psi(\widetilde{r}^b)^m}{\widetilde{v}^b}+\dfrac{m\Psi}{\widetilde{r}}\right)\dx\ds\\
            &=-\int^t_0\int^{\widetilde{r}^b(t,b)}_a\beta\left(\bar u^b_r+\dfrac{m\bar u^b}{r}\right)\left(\Psi_r+\dfrac{m\Psi}{r}\right)r^m\dr\ds+\int^t_0\int^{\widetilde{r}^b(t,b)}_aA(\bar \rho^b)^\gamma\left(\Psi_r+\dfrac{m\Psi}{r}\right)r^m\dr\ds\\
            &-\kappa\int^t_0\int^{\widetilde{r}^b(t,b)}_a\bar \rho^b\left(\bar \rho^b_{rr}+\dfrac{m\bar \rho^b_r}{r}\right)\left(\Psi_r+\dfrac{m\Psi}{r}\right)r^m\dr\ds+\dfrac{\kappa}{2}\int^t_0\int^{\widetilde{r}^b(t,b)}_a|\bar\rho^b_r|^2\left(\Psi_r+\dfrac{m\Psi}{r}\right)r^m\dr\ds\\
            &=-\int^t_0\int^{\infty}_a\beta\left(u^b_r+\dfrac{m u^b}{r}\right)\left(\Psi_r+\dfrac{m\Psi}{r}\right)r^m\dr\ds+\int^t_0\int^{\infty}_aA(\rho^b)^\gamma\left(\Psi_r+\dfrac{m\Psi}{r}\right)r^m\dr\ds\\
            &-\kappa\int^t_0\int^{\infty}_a \rho^b\left( \rho^b_{rr}+\dfrac{m \rho^b_r}{r}\right)\left(\Psi_r+\dfrac{m\Psi}{r}\right)r^m\dr\ds-\kappa\int^t_0\int^\infty_a\rho^b_r\rho^b_{rr}\Psi r^m\dr\ds.
        \end{align*}
       For $J^2_1$, it follows that
        \begin{align*}
            J^2_1&=-\int^t_0\int^b_0\Psi(t,\widetilde{r}(t,x))m\kappa\dfrac{(\widetilde{r}^b)^{2m-1}|\widetilde{v}^b|^2}{(\widetilde{v}^b)^5}\dx\ds\\
            &=-m\kappa\int^t_0\int^\infty_a|\rho^b_r|^2\Psi r^{m-1}\dr\ds.
        \end{align*}
        Thus,
        \begin{align*}
            J_1&=-\int^t_0\int^{\infty}_a\beta\left(u^b_r+\dfrac{m u^b}{r}\right)\left(\Psi_r+\dfrac{m\Psi}{r}\right)r^m\dr\ds-\kappa\int^t_0\int^{\infty}_a \rho^b\left( \rho^b_{rr}+\dfrac{m \rho^b_r}{r}\right)\left(\Psi_r+\dfrac{m\Psi}{r}\right)r^m\dr\ds\\
            &+\int^t_0\int^{\infty}_aA(\rho^b)^\gamma\left(\Psi_r+\dfrac{m\Psi}{r}\right)r^m\dr\ds-\kappa\int^t_0\int^\infty_a\rho^b_r\rho^b_{rr}\Psi r^m\dr\ds-m\kappa\int^t_0\int^\infty_a|\rho^b_r|^2\Psi r^{m-1}\dr\ds.
        \end{align*}
        With the same argument, it follows that
        \begin{align*}
            J_2+J_3&=\int^t_0\int^b_0\widetilde{u}^b\Psi_t(s,\widetilde{r}^b(s,x))\dx\ds+\int^t_0\int^b_0\Psi_r(s,\widetilde{r}^b(s,x))(\widetilde{u}^b)^2\dx\ds\\
            &=\int^t_0\int^\infty_a\rho^bu^b\Psi_t(s,r)r^m\dr\ds+\int^t_0\int^\infty_a\rho^b(u^b)^2\Psi_r(s,r)r^m\dr\ds.
        \end{align*}
        Combining the above results, we have
        \begin{align*}
             &\int^\infty_a\rho^b(t,r) u^b(t,r)\Psi(t,r) r^m\dr-\int^\infty_a\rho^a_0u^a_0\Psi(0,r)r^m\dr\\
             &=-\int^t_0\int^{\infty}_a\beta\left(u^b_r+\dfrac{m u^b}{r}\right)\left(\Psi_r+\dfrac{m\Psi}{r}\right)r^m\dr\ds-\kappa\int^t_0\int^{\infty}_a \rho^b\left( \rho^b_{rr}+\dfrac{m \rho^b_r}{r}\right)\left(\Psi_r+\dfrac{m\Psi}{r}\right)r^m\dr\ds\\
            &+\int^t_0\int^{\infty}_aA(\rho^b)^\gamma\left(\Psi_r+\dfrac{m\Psi}{r}\right)r^m\dr\ds-\kappa\int^t_0\int^\infty_a\rho^b_r\left(\rho^b_{rr}+\dfrac{m\rho^b_r}{r}\right)\Psi r^m\dr\ds\\
            &+\int^t_0\int^\infty_a\rho^bu^b\Psi_t(s,r)r^m\dr\ds+\int^t_0\int^\infty_a\rho^b(u^b)^2\Psi_r(s,r)r^m\dr\ds.
        \end{align*}
        Based on Lemma \ref{weak convergence} and Lemma \ref{strong convergence of solution}, it follows that
        \begin{align}\label{all convergence results}
            &\rho^b\rightarrow\rho,\quad u^b\rightarrow u,\quad\rho^b_r\rightarrow\rho,\quad(u^b_r,\rho^b_{rr})\rightharpoonup(u_r,\rho_{rr}) \text{ in }L^2(0,T;L^2([a,\infty);r^m\dr)),\nonumber\\
            &C^{-a}(a,T)\le\rho^b+\rho\le C(a,T),\quad\sigma^{\frac14}(t)\left(|\rho^b_r|+|\rho|+|u^b|+|u|\right)\le C(a,T)\text{ for }(t,r)\in[0,T]\times[a,\infty).
        \end{align}
         The Definition $\ref{5.1}_2$ can be proved by \eqref{all convergence results}.
         \par The proof of Lemma \ref{weak solution of exterior problem} is finished.
    \end{proof}
\end{lemma}

\section{Global Weak Solutions Including the Origin}\label{5}

In this section, we establish the global existence of the Cauchy problem by passing the inner radii $a\rightarrow0$. The defect measure, mollification and truncations, and the Hardy inequality are used in this section to obtain the weak form of the momentum equation. Now we denote $(\rho^a,u^a)$ as the solutions derived in \textnormal{Theorem} \ref{exterior existence}, and the corresponding initial data $(\rho^a_0,u^a_0)$ are established as before, therefore the particle path $\widetilde{r}^a(t,x)$ satisfies that
\begin{align*}
    x=\int^{\widetilde{r}^a(t,x)}_a\rho^a(t,r)r^m\dr
\end{align*}
Furthermore, we extend the pair $(u^a,\rho_r)=(0,0) \text{ for } r\in[0,a)$. Hence $\rho^a$ is a constant in $[0,a)$ with the bound $\rho^a(t,r)\le C_0 a^{-\frac{2m}{3}}$.
\begin{lemma}\label{mass  estimate}
    For  $\forall t\in[0,T]$, the following estimates hold 
    \begin{align}
        \sup_{a\in(0,1)}\int_{E}\rho^a(t,r)r^m\dr\le C(T)\left(\int_Er^m\dr\right)^{\frac{\gamma-1}{\gamma}}+C(T)\int_{E}r^m\dr,
    \end{align}
    where $E$ is any measurable set in $[a,\infty)$.  In addition, for each $N>0$ and $a\in(0,1)$, it holds that
    \begin{align*}
        \int^a_0\rho^a r^m\dr\le C(T)a^{n-\frac{2m}{3}}.
    \end{align*}.
\end{lemma}
\begin{proof}
    For $\forall E\subset[a,\infty)$, it follows that
    \begin{align*}
        &\int_E\rho^ar^m\dr\le \left(\int_E(\rho^a)^\gamma r^m\dr\right)^{\frac1\gamma}\left(\int_Er^m\dr\right)^{\frac{\gamma-1}{\gamma}}+C(T)\int_{E}r^m\dr\\
        &\le C(T)\left(\int_Er^m\dr\right)^{\frac{\gamma-1}{\gamma}}+C(T)\int_{E}r^m\dr.
    \end{align*}
    Since $\rho^a$ is a constant in $r\in[0,a]$ with the upper bound: $\rho^a\le C_0a^{-\frac{2m}{3}}$, it holds that
    \begin{align*}
        \int^a_0\rho^a r^m\dr\le C(T)a^{n-\frac{2m}{3}}.
    \end{align*}
    \par The proof of Lemma \ref{mass  estimate} is finished.
\end{proof}
\end{lemma}
\begin{lemma}
    There exist  $u_r,\dfrac{u}{r}\in L^2(0,T;L^2([0,\infty);r^m\dr))$ and $\varrho_r\in L^\infty(0,T;L^2([a,\infty);r^m\dr))$ such that
    \begin{align*}
        \left(u^a_r,\dfrac{u^a}{r}\right)\rightharpoonup \left(u_r,\dfrac{u}{r}\right)\text{ in }L^2(0,T;L^2([0,\infty);r^m\dr)),\quad\rho_r^a\overset{*}{\rightharpoonup}\varrho_r\in L^\infty(0,T;L^2([a,\infty);r^m\dr)).
    \end{align*}
    \begin{proof}
        It follows from the extension $(\rho^a_r,u^a)=(0,0)$ in $r\in[0,a]$ that
        \begin{align*}
            \sup_{t\in[0,T]}\int^\infty_0|\rho^a_r|^2r^m\dr+\int^T_0\int^\infty_0\left|u^a_r+\dfrac{mu^a}{r}\right|^2r^m\dr\ds\le C(T).
        \end{align*}
         There exists a $\varrho_r$ such that $\rho^a_r\overset{*}{\rightharpoonup}\rho$ in $L^\infty(0,T;L^2([a,\infty);r^m\dr))$. Next since 
         \begin{align*}
             \int^T_0\int^\infty_0\left|u^a_r+\dfrac{mu^a}{r}\right|^2r^m\dr\ds&=\int^T_0\int^\infty_0\left(\left|u^a_r\right|^2+2u^a_r\dfrac{mu^a}{r}+\left|\dfrac{mu^a}{r}\right|^2\right)r^m\dr\ds\\
             &=\int^T_0\int^\infty_0\left(\left|u^a_r\right|^2+\left|\dfrac{u^a}{r}\right|^2\right)r^m\dr\ds.
         \end{align*}
         Thus, there exists a $u$ such that $u_r\in L^2(0,T;L^2([0,\infty);r^m\dr)),\dfrac{u}{r}L^2(0,T;L^2([0,\infty);r^m\dr))$.
    \end{proof}
\end{lemma}
\subsection{Particle Path and the Possibly Vacuum Ball}\label{Particle path and the Possibly Vacuum ball}
In this section we establish the limit of approximate particle path $\widetilde{r}^a$.
\begin{lemma}\label{last particle path}
There exists a continuous function $\widetilde{r}(t,x):[0,T]\times(0,\infty)\rightarrow[0,\infty)$ such that 
    \begin{align*}
        \lim_{j\rightarrow\infty}\sup_{(t,x)\in \mathscr{K}}|\widetilde{r}^{a_j}(t,x)-\widetilde{r}(t,x)|=0,
    \end{align*}
    where $\mathscr{K}$ is a arbitrarily compact subset in $[0,T]\times(0,\infty)$ and $\{a_j\}$ is a subsequence satisfies $a_j\rightarrow0$ as $j\rightarrow\infty$. Furthermore, for each $\varepsilon\in(0,1)$, there exists a positive constant $C(\varepsilon,T)$  such that the following H$\ddot{o}$lder continuity holds:
    \begin{align*}
        |\widetilde{r}(t,x_1)-\widetilde{r}(t,x_2)|\le C(\varepsilon,T)|x_1-x_2|\quad \text{ for  }\forall(t,x_1,x_2)\in[0,T]\times[\varepsilon,\infty)^2,\\
        |\widetilde{r}(t_1,x)-\widetilde{r}(t_2,x)|\le C(\varepsilon,T)|t^{\frac{3}{4}}_1-t^{\frac{3}{4}}_2|\quad \text{ for  }\forall(t_1,t_2,x)\in[0,T]^2\times\in[\varepsilon,\infty),
    \end{align*}
 In addition, for each fixed $t\in[0,T]$, the $\widetilde{r}(t,x)$ is strictly increasing in $x$. For  $\forall\varepsilon\in(0,1)$, there exists a positive constant $C(\varepsilon,T)$  such that the following upper and lower bounds hold for all $(t,x)\in[0,T]\times(0,\infty)$
    \begin{align*}
        nx\psi^{-1}_-\left(\dfrac{C_0}{x}\right)\le\widetilde{r}^n(t,x)\le C_0(1+x) .
    \end{align*}
    \begin{proof}
         \textbf{1.} We set $G_i=[0,T]\times\left[\dfrac{1}{i},i\right]$ for each $i\in\mathbb{N}$. It follows from  Theorem \ref{exterior existence} that there exists a positive constant $C(i)$ independent of $a$ such that
         \begin{align*}
             \widetilde{r}^a(t,x)=\widetilde{r}^a_0(x)+\int^t_0u^a(s,\widetilde{r}(s,x))\ds\le C_0+ C(i)\int^T_0\sigma^{-\frac{1}{4}}(s)\ds\le C(i)(1+T).
         \end{align*}
         It also holds that  
         \begin{align*}
             \begin{cases}
                 |\widetilde{r}^a(t,x_1)-\widetilde{r}^a(t,x_2)|\le C(i)|x_1-x_2|\quad \text{ for all }(t,x_1,x_2)\in[0,T]\times\left[\frac{1}{i},i\right]^2,\\
                  |\widetilde{r}^a(t_1,x)-\widetilde{r}^a(t_2,x)|\le C(i)|t_1^{\frac{3}{4}}-t_2^{\frac{3}{4}}|\quad \text{ for all }(t_2,t_2,x)\in[0,T]^2\times\left[\frac{1}{i},i\right].
             \end{cases}
         \end{align*}
         Combining the above and the fact: $[0,T]\times(0,\infty)=\cup_{i\in \mathbb{N}}G_i$, there exists a continuous function $\widetilde{r}(t,x):[0,T]\times(0,\infty)\rightarrow[0,\infty)$ and a subsequence $\{a_j\}\rightarrow0$ such that for each compact subset $\mathscr{K}\Subset[0,T]\times(0,\infty)$ we have
         \begin{align*}
             \lim_{j\rightarrow\infty}\sup_{(t,x)\in \mathscr{K}}|\widetilde{r}^{a_j}(t,x)-\widetilde{r}(t,x)|=0,
         \end{align*}
         where we used Lemma \ref{2.1}. Based on this strong convergence, the H$\ddot{o}$lder continuity of $\widetilde{r}$ holds. \\
         \textbf{2.} We prove the strictly increasing property of $\widetilde{r}$. First, for each $0<x_1<x_2<\infty$ we have
         \begin{align*}
             x_2-x_1&=\int^{\widetilde{r}^a(t,x_2)}_{\widetilde{r}^a(t,x_1)}\rho^a(t,r)r^m\dr\le C(x_1,T)\int^{\widetilde{r}^a(t,x_2)}_{\widetilde{r}^a(t,x_1)}r^m\dr\\
             &=\dfrac{C(x_1,T)}{n}\left(\left[\widetilde{r}^a(t,x_2)\right]^n-\left[\widetilde{r}^a(t,x_1)\right]^n\right),
         \end{align*}
         where we have used Lemma \ref{strong convergence of solution} and Lemma \ref{x}. Let $j\rightarrow\infty$ it holds that
         \begin{align*}
             \left(\widetilde{r}(t,x_2)\right)^n-\left(\widetilde{r}(t,x_1)\right)^n>0
         \end{align*}
         Thus $\widetilde{r}(t,x)$ is strictly increasing in $x$ for each fixed $t\in[0,T]$.\\
         \par The proof of Lemma \ref{last particle path} is finished.
    \end{proof}
\end{lemma}
\subsection{Interface of Vacuum}\label{sec Interface of Vacuum}
In this subsection, the existence of the possible vacuum near the origin is proved via Lemma \ref{mass  estimate} and Lemma \ref{last particle path}. 
\begin{lemma}\label{interface of vacuum}
    There exists a H$\ddot{o}$lder continuous function $t\rightarrow\underline{r}(t)\in C^{0,1/n}([0,T])$ such that
    \begin{align*}
        &\lim_{x\rightarrow0^+}\widetilde{r}(t,x)=\underline{r}(t),\quad\lim_{t\rightarrow0^+}\underline{r}(t)=0,\quad\underline{r}^{\frac{n}{2}}(t)\in H^1(0,T),
    \end{align*}
    and $\mathcal{F}:=\left\{(t,r)\in(0,T)\times(0,\infty):\underline{r}(t)<r\right\}$ is an open set. If there exists some $t_*>0$ such that $\underline{r}(t_*)>0$, then 
    \begin{align*}
        \lim_{j\rightarrow\infty}\int^{\underline{r}(t_*)}_{a_j}\rho^{a_j}(t,r)r^m\dr=0.
    \end{align*}
    In addition,
    \begin{align*}
        \liminf_{j\rightarrow\infty}\rho^{a_j}(t,r)=0 \text{ for  }a.e. ~r\in(0,\underline{r}(t_*)).
    \end{align*}
    \begin{proof}
        \textbf{1.} Since $\widetilde{r}(t,x)$ is a non-negative continuous function with the strict monotonicity in $x$, there exists a function $\underline{r}(t):=\lim_{x\rightarrow0}\widetilde{r}(t,x)$.
        
        Since $\widetilde{r}^a(0,x)=\widetilde{r}^a_0(x)$ for each $x\in[0,\infty)$ and $a\in(0,1)$, with the help of Lemma \ref{x} and \eqref{initial x}  we have the following relation for each $x>0$
        \begin{align*}
            x=\int^{\widetilde{r}(0,x)}_a\rho_0(t,r)r^m\dr\ge\dfrac{\left(\widetilde{r}(0,x)\right)^n}{C_0n}.
        \end{align*}
        Thus $\widetilde{r}(0,x)\le \left(C_0nx\right)^{\frac1n}$. For any given $\nu>0$, we choose $x_\nu$ sufficiently small such that
        \begin{align*}
            \widetilde{r}(x_\nu,0)\le \left(C_0nx_\nu\right)^{\frac1n}\le\dfrac{\nu}{2}.
        \end{align*}
        By the continuity, there exists $t_\nu>0$ such that if $0<t<t_\nu$ then $\widetilde{r}(x_\nu,t)<\nu$. By virtue of the definition of $\underline{r}(t)$, for each $t\in(0,t_\nu)$, we have
        \begin{align}
            \underline{r}(t)\le\widetilde{r}(x_\nu,t)<\nu.
        \end{align}
        Thus $\lim_{t\rightarrow0}\underline{r}(t)=0$. \\ 
        \textbf{2. } 
      Since $u^a(t,a)=0$, thus for each $r>a$ we have $r^mu^a(t,r)=\displaystyle\int^r_a(u^a+\dfrac{mu^a}{s})s^m\ds$. Hence, it holds that
      \begin{align*}
          |u^a(t,r)|\le\dfrac{1}{r^m}\left(\int^r_a\left|u^a+\dfrac{mu^a}{r}\right|^2s^m\ds\right)^\frac{1}{2}\left(\dfrac{\widetilde{r}^n-a^n}{b}\right)^\frac12\le C_0r^{1-\frac n2}\left(\int^\infty_0\left|u^a+\dfrac{mu^a}{r}\right|^2s^m\ds\right)^\frac{1}{2}
      \end{align*}
      Thus, by Lemma \ref{x}, it follows for each $t\in[0,T]$ that $\partial_t\widetilde{r}^a(t,x)=\widetilde{u}^a(\widetilde{r}^a(t,x),t)$. Furthermore, we have
        \begin{align*}
            \left|\partial_t\left(\widetilde{r}^a(t,x)^{\frac n2}\right)\right|=\left|\dfrac{n}{2}\widetilde{r}^a(t,x)^{\frac{n-2}{2}}\widetilde{u}^a(\widetilde{r}^a(t,x),t)\right|\le C_0\left(\int^\infty_0\left|u^a+\dfrac{mu^a}{r}\right|^2s^m\ds\right)^\frac{1}{2}.
        \end{align*}
        Then we integrating in time to obtain that
        \begin{align}\label{holder ra}
            &|\widetilde{r}^a(t,x)^{\frac n2}-\widetilde{r}^a(s,x)^{\frac n2}|\le C(T)|t-s|^\frac12 \text{ for each } s,t\in[0,T],\nonumber\\
            &\int^T_0\left|\partial_t(\widetilde{r}^a(t,x)^{\frac n2})\right|^2\dt\le C(T) \text{ for each } x>0.
        \end{align}
        By virtue of Lemma \ref{last particle path}, it follows by passing $j\rightarrow\infty$ to \eqref{holder ra} that
        \begin{align*}
            |\widetilde{r}(t,x)^{\frac n2}-\widetilde{r}(s,x)^{\frac n2}|\le C(T)|t-s|^\frac12 \text{ for each } s,t\in[0,T].
        \end{align*}
        Let $x\rightarrow0$, the H$\ddot{o}$lder continuity of $\underline{r}$ is proved.
        Next, combining $\Vert\partial_t(\widetilde{r}^a(t,x)^{\frac n2})\Vert_{L^2(0,T)}\le C(T)$ and Lemma \ref{last particle path}, it holds that
        \begin{align*}
            \int^T_0\left|\partial_t(\widetilde{r}(t,x)^{\frac n2})\right|^2\dt\le C(T) \text{ for each } x>0,
        \end{align*}
        Let $x\rightarrow0$, we obtain $\underline{r}^{\frac n 2}(t)\in H^1(0,T)$.\\
        \textbf{3. } If there exists $t_*>0$, then for each $a\in(0,\underline{r}(t_*))$
        \begin{align*}
            \left|\int^{\underline{r}(t_*)}_a\rho^a(t,r)r^m\dr\right|&\le \left|\int^{\widetilde{r}^a(t,x)}_a\rho^a(t,r)r^m\dr\right|+ \left|\int^{\widetilde{r}(t_*,x)}_{\widetilde{r}^a(t_*,x)}\rho^a(t,r)r^m\dr\right|+\left|\int^{\underline{r}(t_*)}_{\widetilde{r}(t_*,x)}\rho^a(t,r)r^m\dr\right|\\
            &\le x+C(T)\left(\left|
            (\widetilde{r}(t_*,x))^n-(\widetilde{r}^a(t_*,x))^n\right|+\left|
            (\widetilde{r}(t_*,x))^n-(\underline{r}(t_*))^n\right|\right).
        \end{align*}
        Thus for each $\eta>0$, there exists $\delta>0$ such that for each $x\in(0,\delta)$
        \begin{align*}
            \left|
            (\widetilde{r}(x,t_*))^n-(\underline{r}(t_*))^n\right|\le\eta,
        \end{align*}
        where $x>0$. Then let $\delta'=\min\{\delta,\eta\}$, we have
        \begin{align*}
            \limsup_{j\rightarrow\infty}\int^{\underline{r}(t_*)}_{a_j}\rho^{a_j}(t,r)r^m\dr\le 2\eta+C(T)\limsup_{j\rightarrow\infty}\left|
            (\widetilde{r}(t_*,x))^n-(\widetilde{r}^a(t_*,x))^n\right|=2\eta.
        \end{align*}
        Now passing limit to $\eta\rightarrow0$, we obtain that
        \begin{align*}
            \limsup_{j\rightarrow\infty}\int^{\underline{r}(t_*)}_{a_j}\rho^{a_j}(t,r)r^m\dr=0.
        \end{align*}
        Combining $\rho^a\ge0$ for each $a\in(0,1)$, then
        \begin{align*}
            \lim_{j\rightarrow\infty}\int^{\underline{r}(t_*)}_{a_j}\rho^{a_j}(t,r)r^m\dr=0.
        \end{align*}
        \par The proof of Lemma \ref{interface of vacuum} is finished via the Fatou's lemma.
    \end{proof}
\end{lemma}
\subsection{Strong Convergence of Solutions}\label{Strong convergence of solutions}
In this subsection we establish the strong convergence of the approximate solutions $(\rho^a,u^a)$. From now on, we denote
\begin{align*}
    \mathcal{F}:=\left\{(t,r)\in[0,T]\times(0,\infty):\underline{r}(t)<r\right\}
\end{align*}
\begin{lemma}\label{last convergence}
There exist a subsequence $a_j\rightarrow0$ as $j\rightarrow\infty$ and H$\ddot{o}$lder continuous functions $(u,\rho_r):\mathcal{F}\backslash\{t=0\}\rightarrow \mathbb{R}\times\mathbb{R}$, $\rho:\mathcal{F}\rightarrow [0,\infty)$ such that 
\begin{align*}
    \lim_{j\rightarrow\infty}\sup_{(t,r)\in \Omega_1}|(u^{a_j}-u,\rho_r^{a_j}-\rho_r)|=0,\quad\lim_{j\rightarrow\infty}\sup_{(t,r)\in \Omega_2}|\rho^{a_j}-\rho|=0,
\end{align*}
where $\Omega_1$ and $\Omega_2$ are compact subset in $\mathcal{F}\backslash\{t=0\}$ and $\mathcal{F}$, respectively. In addition, for each $\varepsilon\in(0,1)$, there exists a positive constant $C(\varepsilon,T)$ such that
\begin{align*}
    &\sigma^{\frac14}(t)\left(|u|+|\rho_r|\right)\le C(\varepsilon,T)\quad \text{ for each }r\in[\widetilde{r}(t,\varepsilon),\infty) \text{ and }t\in[0,T],\\
    &C^{-1}(\varepsilon,T)\le\rho\le C(\varepsilon,T)\quad \text{ for each }r\in[\widetilde{r}(t,\varepsilon),\infty) \text{ and }t\in[0,T].
\end{align*}
Moreover, for $\forall\varepsilon\in(0,1)$, there exists a positive constant $C(\varepsilon,T)$ such that
\begin{itemize}
    \item[1. ] For $\forall r\ge\sup_{t\in[t_1,t_2]}\widetilde{r}(t,\varepsilon)$ and $0\le t_1<t_2\le T$, it holds that
    \begin{align*}
        |\rho(t_1,r)-\rho(t_2,r_2)|+\sigma^{\frac12}(t_1)|u(t_1,r)-u(t_2,r)|+\sigma^{\frac12}(t_1)|\rho_r(t_2,r)-\rho_r(t_1,r)|\le C(\varepsilon,T)|t_2-t_1|^{\frac{1}{4}},
    \end{align*}
    \item[2. ] For  $\forall r_1,r_2\ge\widetilde{r}(t,\varepsilon)$ and $t\in[0,T]$, it holds that
    \begin{align*}
        |\rho(t,r_1)-\rho(t,r_2)|+\sigma^{\frac{1}{2}}(t)|u(t,r_1)-u(t,r_2)|+\sigma(t)|\rho_r(t,r_1)-\rho_r(t,r_2)|\le C(\varepsilon,T)|r_1-r_2|^{\frac{1}{2}}.
    \end{align*}
\end{itemize}

\begin{proof}
    We only prove for $\rho_r$, since the case for $(\rho,u)$ is similar to $\rho_r$.  For any fixed $\varepsilon\in(0,1)$, we define the exterior fluid region by
    \begin{align*}
        \mathcal{F}_\varepsilon:=\left\{(t,r)\in[0,T]\times(0,\infty):t\in[\varepsilon,T],\widetilde{r}(t,\varepsilon)\le r\le \dfrac{1}{\varepsilon}\right\}.
    \end{align*}
    Based on the upper and lower bound  established in Lemma \ref{last particle path}, there exists $\varepsilon'>0$ such that $Int(\mathcal{F}_\varepsilon)\neq\emptyset$ for each $0<\varepsilon<\varepsilon'$. Furthermore, the strict monotonicity holds: $\mathcal{F}_{\varepsilon_1}\subset \mathcal{F}_{\varepsilon_2}$ if $\varepsilon_1>\varepsilon_2$.
    It follows from a straightforward proof that
    \begin{align*}
        \mathcal{F}\backslash\{t=0\}=\underset{0<\varepsilon<1}{\cup}\mathcal{F}_\varepsilon,
    \end{align*}
    and we omit here for simplicity. Since $\widetilde{r}(t,x)$ is strictly increasing in $x$ and $(t,x)\rightarrow\widetilde{r}(t,x)$ is continuous, then we let $\wp_\varepsilon:=\inf_{t\in[0,T]}\left(\widetilde{r}(t,\varepsilon)-\widetilde{r}\left(t,\dfrac{\varepsilon}{2}\right)\right)>0$. Thus there exists $N_\varepsilon\in\mathbb{N}$ such that for $\forall j\ge N_\varepsilon$
    \begin{align*}
        \sup_{t\in[0,T]}|\widetilde{r}^{a_j}\left(t,\frac{\varepsilon}{2}\right)-\widetilde{r}\left(t,\frac{\varepsilon}{2}\right)|\le\frac{\wp_\varepsilon}{2}
    \end{align*}
    Therefore, for $\forall j\ge N_\varepsilon$ and $\forall t\in[0,T]$, it holds that  $\widetilde{r}^{a_j}\left(t,\dfrac{\varepsilon}{2}\right)<\widetilde{r}(t,\varepsilon)$.
    With the help of Lemma \ref{strong convergence of solution}, 
    \begin{align}\label{rhor}
      \begin{cases}
            &\sup_{j\ge N_\varepsilon}|\rho_r^{a_j}(t,r_1)-\rho_r^{a_j}(t,r_2)|\le C(\frac{\varepsilon}{2},T)\sigma^{-\frac{1}{2}}(t)|r_1-r_2|^{\frac{1}{2}}\text{ for each }(t,r_1),(t,r_2)\in \mathcal{F}_\varepsilon,\\
        & \sup_{j\ge N_\varepsilon}|\rho_r^{a_j}(t_1,r)-\rho_r^{a_j}(t_2,r)|\le C(\frac{\varepsilon}{2},T)\sigma^{-\frac{1}{2}}(t_1)|t_1-t_2|^{\frac{1}{4}}\text{ for each }(t_1,r),(t_2,r)\in \mathcal{F}_\varepsilon,\\
        &\sup_{j\ge N_\varepsilon}|\rho_r^{a_j}(t,r)|\le C(\varepsilon,T)\sigma^{-\frac14}(t)\le \varepsilon^{-\frac14}C(\varepsilon,T)\text{ for each }(t,r)\in \mathcal{F}_\varepsilon,
      \end{cases}
    \end{align}
    where we have used the fact that if $(t,r)\in \mathcal{F}_\varepsilon$, then for all $j\ge N_\varepsilon$, $\widetilde{r}^{a_j}\left(t,\dfrac{\varepsilon}{2}\right)<r$. 
     By applying Lemma \ref{2.1} to \ref{rhor}, we obtain that there exists a continuous function $\rho_r:\mathcal{F}\backslash\{t=0\}\rightarrow\mathbb{R}$ such that, for each compact set $\Omega_1\in \mathcal{F}\backslash\{t=0\}$, it holds that
     \begin{align*}
         \lim_{j\rightarrow\infty}\sup_{(t,r)\in\Omega_1}|\rho^a_r-\rho_r|=0.
     \end{align*}
     \par The proof is finished.
\end{proof}
\end{lemma}
There exist some uniform estimates of limit of solution which are independent of inner radii $a\in(0,1)$,
\begin{lemma}\label{last estimate epsilon}
        For each fixed $\varepsilon\in(0,1)$, there exists a positive constant $C(\varepsilon,T)$ such that
        \begin{align*}
        &\sup_{t\in[0,T]}\int^\infty_{\widetilde{r}(t,\varepsilon)}\left(|\rho-1|^2+\sigma |u|^4+\sigma|u_r|^2+\sigma|\rho_{rr}|^2+\sigma|\rho_r|^4\right)r^m\dr\nonumber\\
        &+\int^T_0\int^\infty_{\widetilde{r}(t,\varepsilon)}\left(|u_r|^2+|\rho_{rr}|^2+|\rho_r|^4+\sigma|(u_t,u_{rr},\rho_{rrr},\rho_{tr})|^2+|\rho_t|^2\right)r^m\dr\ds\nonumber\\
        &\le C(\varepsilon,T).
    \end{align*}
    \begin{proof}
        These estimates can be proved by the same argument of \eqref{estimate epsilon}. 
    \end{proof}
\end{lemma}
In addition, we have the following relation,
\begin{lemma}\label{x mass}
    For  $\forall x>0$ and  $\forall t\in[0,T]$
    \begin{align*}
        x=\int^{\widetilde{r}(t,x)}_{\underline{r}(t)}\rho(t,r)r^m\dr.
    \end{align*}
    \begin{proof}
        From Lemma \ref{x}, for each $0<y<x$, we have
        \begin{align*}
            x-y=\int^{\widetilde{r}^{a_j}(t,x)}_{\widetilde{r}^{a_j}(t,y)}\rho^{a_j}r^m\dr
        \end{align*}
        By virtue of Lemma \ref{last convergence} and Lemma \ref{last particle path}, for the right hand side, we have
        \begin{align*}
            \int^{\widetilde{r}^{a_j}(t,x)}_{\widetilde{r}^{a_j}(t,y)}\rho^{a_j}r^m\dr\rightarrow \int^{\widetilde{r}(t,x)}_{\widetilde{r}(t,y)}\rho r^m\dr.
        \end{align*}
        Then we let $y\rightarrow0$,
        \begin{align*}
            x=\int^{\widetilde{r}(t,x)}_{\underline{r}(t)}\rho(t,r)r^m\dr,
        \end{align*}
        where we used the monotone convergence theorem.
    \end{proof}
\end{lemma}

\subsection{Entropy Estimates and  Speed of the Interface}\label{Entropy estimates}
In this subsection, we establish some estimate independent of inner radii $a$. Since $\rho\in L^\infty_{loc}(\mathcal{F})$, now we make the following extension for $\rho$ based on Lemma \ref{x mass} and Lemma \ref{interface of vacuum} 
\begin{align*}
    \rho(t,r)=0  \text{ in }\mathbb{R}^n\times[0,T]\backslash \mathcal{F}.
\end{align*}
In addition, we let $\mathcal{V}:=\{(t,r):0<r<\underline{r}(t)\}$ denote the possibly vacuum ball. To prove the possibly vacuum ball is empty, our method is proof by contradiction. From now on, we assume $\mathcal{V}$ is not empty.
\par First we obtain the following entropy estimates of limit solutions $(\rho,u)$.
\begin{lemma}\label{last entropy estimate}
There exists a positive constant $C(T)$ such that
    \begin{align*}
        \esssup_{t\in[0,T]}\int^\infty_0\left(\rho|u|^2+G(\rho)+|\rho_r|^2\right)r^m\dr+\int^T_0\int^\infty_0\left|u_r+\dfrac{mu}{r}\right|^2r^m\dr\ds\le C(T).
    \end{align*}

    \begin{proof}
    We divide the proof into 3 steps.\\
    \textbf{1. }
        Since $(t,x)\rightarrow\widetilde{r}(t,x)$ is continuous and $\widetilde{r}(t,x)$ is strictly increasing in $x$, there exists \begin{align*}
            \wp_\varepsilon:=\inf_{t\in[0,T]}\left(\widetilde{r}(t,\varepsilon)-\widetilde{r}(t,\frac{\varepsilon}{2})\right)>0.
        \end{align*}
        By virtue of Lemma \ref{last particle path}, there exists $N_\varepsilon$ such that for all $j\ge N_\varepsilon$
        \begin{align*}
            \sup_{t\in[0,T]}\left|\widetilde{r}\left(t,\frac{\varepsilon}{2}\right)-\widetilde{r}^{a_j}\left(t,\frac{\varepsilon}{2}\right)\right|\le\dfrac{\wp_\varepsilon}{2}.
        \end{align*}
        Then we have $\widetilde{r}(t,\varepsilon)>\widetilde{r}^{a_j}(t,\frac{\varepsilon}{2})$. For   each bounded set $E\subset[0,\infty)$, we denote   $\chi_E$ as the indicator function, and  set $E=[0,M]$ for arbitrary fixed $M\in \mathbb{R}_+$. Now for each $0<t_1<t_2\le T$, we have
        \begin{align*}
            \int^{t_2}_{t_1}&\int^\infty_{\widetilde{r}^{a_j}(t,\varepsilon)}\rho^{a_j}\left|u^{a_j}\right|^2\chi_Er^m\dr\ds-\int^{t_2}_{t_1}\int^\infty_{\widetilde{r}(t,\varepsilon)}\rho|u|^2\chi_Er^m\dr\ds\\
            &=\int^{t_2}_{t_1}\int^\infty_{\widetilde{r}(t,\varepsilon)}\bigg[\rho^{a_j}\left|u^{a_j}\right|^2-\rho |u|^2\bigg]\chi_Er^m\dr\ds+\int^{t_2}_{t_1}\int^{\widetilde{r}(t,\varepsilon)}_{\widetilde{r}^{a_j}(t,\varepsilon)}\rho^{a_j}\left|u^{a_j}\right|^2 \chi_Er^m\dr\ds\\
            &=I_1+I_2.
        \end{align*}
        By virtue of  $\widetilde{r}(t,\varepsilon)>\widetilde{r}^{a_j}(t,\frac{\varepsilon}{2})$ for each $j\ge N_\varepsilon$ and the point-wise bound: $\rho+\rho^{a_j}+\sigma^{\frac{1}{4}}(t)|u|+\sigma^{\frac{1}{4}}(t)|u^{a_j}|\le C(\varepsilon,T)$ for all $r\in[\widetilde{r}(t,\varepsilon),\widetilde{r}^{a_j}(t,\varepsilon)]\text{ or }[\widetilde{r}^{a_j}(t,\varepsilon),\widetilde{r}(t,\varepsilon)]$, for $I_2$, it holds that
        \begin{align*}
            I_2&\le C(\varepsilon,T)\int^{t_2}_{t_1}\sigma^{-\frac{1}{2}}(s)\left|\big(\widetilde{r}(\varepsilon,s)\big)^n-\big(\widetilde{r}^{a_j}(\varepsilon,s)\big)^n\right|\ds\\
            &\le C(\varepsilon,T)\sup_{s\in[t_1,t_2]}\left|\big(\widetilde{r}(\varepsilon,s)\big)^n-\big(\widetilde{r}^{a_j}(\varepsilon,s)\big)^n\right|\rightarrow0 \text{ as }j\rightarrow\infty.
        \end{align*}
        Since
        \begin{align*}
            I_1&=\int^{t_2}_{t_1}\int^\infty_{\widetilde{r}(t,\varepsilon)}\bigg[\rho^{a_j}\left|u^{a_j}\right|^2-\rho |u|^2\bigg]\chi_Er^m\dr\ds\\
            &=\int^{t_2}_{t_1}\int^\infty_{\widetilde{r}(t,\varepsilon)}\bigg[(\rho^{a_j}-\rho)\left|u^{a_j}\right|^2-\rho (|u|^2-|u^{a_j}|^2)\bigg]\chi_Er^m\dr\ds
            \\
        &=I_{1,1}+I_{1,2}
        \end{align*}
        With the help of Lemma \ref{last convergence} and $\widetilde{r}(t,\varepsilon)>\widetilde{r}^{a_j}(t,\frac{\varepsilon}{2})$, we use the dominated convergence theorem to check $I_{1,1}$ and $I_{1,2}$ tend to zero as $j\rightarrow\infty$. Now passing limit to $M\rightarrow\infty$, we have
        \begin{align}\label{momen}
            \int^{t_2}_{t_1}\int^M_{\widetilde{r}(t,\varepsilon)}\rho|u|^2r^m\dr\ds=&\lim_{j\rightarrow\infty}\int^{t_2}_{t_1}\int^M_{\widetilde{r}^{a_j}(t,\varepsilon)}\rho^{a_j}|u^{a_j}|^2r^m\dr\ds\le C(T)|t_2-t_1|,
        \end{align}
        where we used Lemma \ref{limit entropy}. Now we apply the Lebesgue differentiation theorem and monotone convergence theorem to obtain that
        \begin{align*}
            \esssup_{t\in[0,T]}\int^\infty_{\underline{r}(t)}\rho|u|^2r^m\dr\le C(T).
        \end{align*}
        Since $\rho=0$ in $r\in[0,\underline{r}(t)]$ we have
        \begin{align*}
            \esssup_{t\in[0,T]}\int^\infty_0\rho|u|^2r^m\dr\le C(T).
        \end{align*}
       \textbf{2. } Since $C^{-1}(\varepsilon,T)\le\rho^{a_j}(t,r)\le C(\varepsilon,T)$ for each $r\ge\widetilde{r}(t,\varepsilon)>\widetilde{r}^{a_j}(t,\frac{\varepsilon}{2})$, thus $|G(\rho^{a_j})|\le C(\varepsilon,T)$. Now for each $0<t_1<t_2\le T$ and each bounded set $E\subset[0,\infty)$
        \begin{align*}
            &\int^{t_2}_{t_1}\int^\infty_{\widetilde{r}^{a_j}(t,\varepsilon)}G(\rho^{a_j})\chi_Er^m\dr\ds-\int^{t_2}_{t_1}\int^\infty_{\widetilde{r}(t,\varepsilon)}G(\rho)\chi_Er^m\dr\ds\\
            &=\int^{t_2}_{t_1}\int^\infty_{\widetilde{r}(t,\varepsilon)}\bigg[G(\rho^{a_j})-G(\rho)\bigg]\chi_Er^m\dr\ds+\int^{t_2}_{t_1}\int^{\widetilde{r}(t,\varepsilon)}_{\widetilde{r}^{a_j}(t,\varepsilon)}G(\rho^{a_j})\chi_Er^m\dr\ds\rightarrow0\text{ as }j\rightarrow\infty,
        \end{align*}
        where we used the same argument mentioned above. By virtue of Lemma \ref{limit entropy}, we have
        \begin{align}\label{G}
            \int^{t_2}_{t_1}\int^\infty_{\widetilde{r}(t,\varepsilon)}G(\rho)\chi_Er^m\dr\ds=\lim_{j\rightarrow\infty}\int^{t_2}_{t_1}\int^\infty_{\widetilde{r}^{a_j}(t,\varepsilon)}G(\rho^{a_j})\chi_Er^m\dr\ds\le C(T)|t_2-t_1|.
        \end{align}
        Then we let $E=[0,M]$ for each fixed $M\in \mathbb{R}$ and apply the Lebesgue differentiation theorem to \eqref{G} to obtain that
        \begin{align*}
            \esssup_{t\in(0,T)}\int^M_{\widetilde{r}(t,\varepsilon)}G(\rho)r^m\dr\le C(T).
        \end{align*}
        Then by the monotone convergence theorem, we pass the limit to $M$ and $\varepsilon$  to obtain that 
        \begin{align*}
            \esssup_{t\in(0,T)}\int^\infty_{\underline{r}(t)}G(\rho)r^m\dr\le C(T).
        \end{align*}
        Since $G=1$ in $\{\rho=0\}$ thus 
        \begin{align*}
            \esssup_{t\in(0,T)}\int^\infty_0G(\rho)r^m\dr=\esssup_{t\in(0,T)}\int^\infty_{\underline{r}(t)}G(\rho)r^m\dr+\int^{\underline{r}(t)}_0r^m\dr\le C(T).
        \end{align*}
\textbf{3. } Since $\rho^a_r\rightarrow\rho_r$ in $\mathcal{F}$, and $\rho^a_r\overset{*}{\rightharpoonup}\varrho_r$ in $L^\infty(0,T;[0,\infty);r^m\dr)$, it holds that $\rho_r=\varrho_r~a.e. $  in the fluid region. Next, for each $N\in\mathbb{R}$,
\begin{align*}
    \int^N_0|\rho^{a_j}-\rho|r^m\dr&=\int^N_{\widetilde{r}(t,x)}|\rho^{a_j}-\rho|r^m\dr+\int^{\widetilde{r}^{a_j}(t,x)}_{a_j}\rho^{a_j}r^m\dr+\int^{\widetilde{r}(t,x)}_{\widetilde{r}^{a_j}(t,x)}\rho^{a_j}r^m\dr+\int^{\widetilde{r}(t,x)}_{\underline{r}(t)}\rho r^m\dr\\
    &\le2x+\int^N_{\widetilde{r}(t,x)}|\rho^{a_j}-\rho|r^m\dr+C(T)\left[(\widetilde{r}^{a_j}(t,x))^n-(\widetilde{r}(t,x))^n\right],
\end{align*}
where $x>0$. We pass the limit to $j$ and $x$ to obtain $\rho^{a_j}\rightarrow\rho$ in $L^\infty(0,T;L^1((0,N);r^m\dr))$. Thus, for each $\Phi\in C^\infty_c((0,N))$, it holds that
\begin{align*}
   \int^N_0\varrho_r\Phi\dr \leftarrow\int^N_0\rho^{a_j}_r\Phi\dr=-\int^N_0\rho^{a_j}\Phi_r\dr\rightarrow-\int^N_0\rho\Phi_r\dr
\end{align*}
Thus $\rho_r=\varrho_r~a.e. $  in $[0,\underline{r}(t)]$. Combining the above, now, we have proved $\rho=\varrho$ for  $a.e. ~r\in[0,\infty)$. Thus,
\begin{align*}
    \sup_{t\in[0,T]}\int^\infty_0|\rho_r|^2r^m\dr\le C(T).
\end{align*}
Next, the limit dissipation estimate is derived from the weakly lower semi-continuity of Sobolev norm straightly.
\par The proof of Lemma \ref{last entropy estimate} is finished.
    \end{proof}
\end{lemma}
In addition, we prove the speed of the interface of vacuum,
\begin{lemma}\label{interface speed}
     The speed of the vacuum interface  satisfies that 
     \begin{align*}
         \dfrac{\mathrm{d}\underline{r}(t)} {\mathrm{d}t}=u(t,\underline{r}(t)) \text{ for } a.e.  t\in \{t\in(0,T):\underline{r}(t)>0\}.
     \end{align*}
    \begin{proof}
        For  $\forall x>0$ and $0<t_1<t_2< T$, it holds that $\{(\widetilde{r}(t,x),t):t_1\le t\le t_2\}$ is a compact subset in $\mathcal{F}$. Furthermore, by virtue of Lemma \ref{last particle path} and Lemma \ref{last convergence}, we have the following convergence for each $t\in[t_1,t_2]$,
        \begin{align*}
            u^{a_j}(t,\widetilde{r}^{a_j}(t,x))\rightarrow u(t,\widetilde{r}(t,x)) \text{ as } j\rightarrow\infty.
        \end{align*}
        Then, it holds from Lemma \ref{x} that
        \begin{align}\label{rminusr}
            \widetilde{r}(t,x)-\widetilde{r}(s,x)=\int^t_su(\tau,\widetilde{r}(\tau,x))\mathrm{d}\tau, \text{ for }0<s<t< T.
        \end{align}
        Since for each $\mathcal{I}\Subset\{t\in(0,T):\underline{r}(t)>0\}$, there exists two positive constants $M_1,M_2$ such that  $M_1\le\underline{r}(t)< M_2$ for each $t\in \mathcal{I}$. Base on Lemma \ref{last entropy estimate}, it holds that 
        \begin{align*}
            u\in L^2\left(0,T;H^1\left(\dfrac{M_1}{2},2M_2;\dr\right)\right)\subset L^2\left(0,T;C\left(\dfrac{M_1}{2},2M_2\right)\right).
        \end{align*}
        Thus, $\Vert u(\cdot,\underline{r}(\cdot))\Vert_{L^2(I)}\le C(M_1,M_2,T)\Vert u\Vert_{L^2\left(0,T;H^1\right)}<\infty$.
        With the help of Cauchy--Schwarz inequality, it holds that
        \begin{align*}
            \int_I|u(t,\widetilde{r}(t,x))-u(t,\underline{r}(t))|^2\dt&=\int_I\left|\int^{\widetilde{r}(t,x)}_{\underline{r}(t)}u_r(t,s)\ds\right|^2\dt\\
            &\le C(M_1,M_2,T)\sup_{t\in I}|\widetilde{r}(t,x)-\underline{r}(t)|\int_I\int^{2M_2}_{\frac{M_1}{2}}|u_r|^2\dr\dt\rightarrow0\text{ as } x\rightarrow0.
        \end{align*}
Now passing limit to $x\rightarrow0$ in \eqref{rminusr}, we obtain
\begin{align*}
    \underline{r}(t)-\underline{r}(s)=\int^t_su(\tau,\underline{r}(\tau))\mathrm{d}\tau.
\end{align*}
Then by the Lebesgue differentiation theorem, for  $a.e. ~t\in \mathcal{I}$, it holds that $\dfrac{\mathrm{d}r(t)}{\mathrm{d}t}=u(t,\underline{r}(t))$. Since $\{t\in(0,T):\underline{r}(t)>0\}$ is an open set,  there exists a union of at most countably many open intervals $(t^j_1,t^j_2)$ such that
$\{t\in(0,T):\underline{r}(t)>0\}=\cup_j(t^j_1,t^j_2)$. We let $\mathcal{I}^j_n=\left[t^j_1+\dfrac{t^j_2-t^j_1}{n+2},t^j_2-\dfrac{t^j_2-t^j_1}{n+2}\right]\Subset(t^j_1,t^j_2)$, then for  $a.e. ~t\in\mathcal{I}^j_n$, it holds that $\dfrac{\mathrm{d}r(t)}{\mathrm{d}t}=u(t,\underline{r}(t))$. Since $(t^j_1,t^j_2)=\cup_n\mathcal{I}^j_n$, we have $\dfrac{\mathrm{d}r(t)}{\mathrm{d}t}=u(t,\underline{r}(t))$ for $a.e. ~t\in(t^j_1,t^j_2)$. With the same argument, it holds that
\begin{align*}
    \dfrac{\mathrm{d}r(t)}{\mathrm{d}t}=u(t,\underline{r}(t))\text{ for }a.e. t\in\{t\in(0,T):\underline{r}(t)>0\}.
\end{align*}
\par The proof of Lemma \ref{interface speed} is finished.
    \end{proof}
\end{lemma}
\subsection{Weak Forms of the Cauchy Problem}\label{Weak forms of the Cauchy problem}
In this section, we prove that the solutions $(\rho,u)$ satisfies the weak form of the Cauchy problem which are specified in Definition \ref{defni}.
\begin{lemma}\label{weak form of continuity equation}
For  $\forall \Psi\in C^1\left([0,T];C^1_c([0,\infty)\right)$, it follows that
\begin{align*}
    \int^\infty_0\rho(t,r)\Psi(t,r)r^m\dr-\int^\infty_0\rho_0(r)\Psi(0,r)r^m\dr=\int^t_0\int^\infty_0(\rho\Psi_t+\rho u\Psi_r)r^m\dr\ds.
\end{align*}
    \begin{proof}
    Since $(t,y)\rightarrow\widetilde{r}(t,y)$ is continuous and $\widetilde{r}(t,y)$ is strictly increasing in $y$, there exists 
    \begin{align*}
            \wp_x:=\inf_{t\in[0,T]}\left(\widetilde{r}(t,x)-\widetilde{r}\left(t,\frac{x}{2}\right)\right)>0
        \end{align*}
        By virtue of Lemma \ref{last particle path}, there exists $N_x$ such that for all $j\ge N_x$ we have
        \begin{align*}
            \sup_{t\in[0,T]}\left|\widetilde{r}(t,\frac{x}{2})-\widetilde{r}^{a_j}\left(t,\frac{x}{2}\right)\right|\le\dfrac{\wp_x}{2}.
        \end{align*}
        Then we $\widetilde{r}(t,x)>\widetilde{r}^{a_j}(t,\frac{x}{2})$.
        From Definition \ref{5.1}, it follows that
        \begin{equation*}
        \int^\infty_{a_j}\rho^{a_j}(t,r)\Psi (t,r)r^m\dr= \int^\infty_{a_j}\rho_0^{a_j}(r)\Psi (0,r)r^m\dr+\int^t_0\int^\infty_{a_j}\big(\rho^{a_j}\Psi_t +\rho^{a_j}u^{a_j}\Psi_r\big)r^m\dr\ds.
    \end{equation*}
    Since for each $x>0$ we have
    \begin{align*}
        &\left|\int^\infty_{a_j}\rho^{a_j}(t,r)\Psi (t,r)r^m\dr-\int^\infty_{\underline{r}(t)}\rho(t,r)\Psi (t,r)r^m\dr\right|\\
        &\le\left|\int^\infty_{\widetilde{r}^{a_j}(t,x)}\rho^{a_j}(t,r)\Psi(t,r) r^m\dr-\int^\infty_{\widetilde{r}(t,x)}\rho(t,r)\Psi(t,r)r^m\dr\right|\\
        &+\Vert\Psi\Vert_{L^\infty}\bigg[\int^{\widetilde{r}(t,x)}_{\underline{r}(t)}\rho(t,r)r^m\dr+\int^{\widetilde{r}(t,x)}_{a_j}\rho^{a_j}(t,r)r^m\dr+\int^{\widetilde{r}^{a_j}(t,x)}_{\widetilde{r}(t,x)}\rho^{a_j}(t,r)r^m\dr\bigg]\\
        &:=\sum^4_{i=1}I_1.
    \end{align*}
    For $I_1$,
    \begin{align*}
        I_1=&\left|\int^\infty_{\widetilde{r}^{a_j}(t,x)}\rho^{a_j}(t,r)\Psi(t,r) r^m\dr-\int^\infty_{\widetilde{r}(t,x)}\rho(t,r)\Psi(t,r)r^m\dr\right|\\
        \le&\left|\int^\infty_{\widetilde{r}(t,x)}\left(\rho^{a_j}-\rho\right)\Psi r^m\dr\right|+\left|\int^{\widetilde{r}(t,x)}_{\widetilde{r}^{a_j}(t,x)}\rho^{a_j}\Psi r^m\dr\right|\rightarrow0 \text{ as }j\rightarrow\infty.
    \end{align*}
    It follows from Lemma \ref{x} and Lemma \ref{x mass} that $I_2+I_3\le Cx$. Thus for each $\delta>0$, there exists $x'>0$ such that for each $x\in(0,x')$, $I_2+I_3\le \delta$. By virtue of Lemma \ref{interface of vacuum}, $I_4$ tends to zero as $j\rightarrow\infty$. Thus
    \begin{align*}
        \lim_{j\rightarrow\infty}\int^\infty_{a_j}\rho^{a_j}(t,r)\Psi (t,r)r^m\dr=\int^\infty_{\underline{r}(t)}\rho(t,r)\Psi (t,r)r^m\dr.
    \end{align*}
    By the same argument, we have
    \begin{align*}
        \int^t_0\int^\infty_{\underline{r}(t)}\rho\Psi_tr^m\dr\ds=\lim_{j\rightarrow\infty}\int^t_0\int^\infty_{a_j}\rho^{a_j}\Psi_tr^m\dr\ds.
    \end{align*}
    Since $\left(\rho^a_0,u^a_0\right)=(1,0)$ in $[0,a]$,
    \begin{align*}
        &\left|\int^\infty_a\rho^a_0\Psi(0,r)r^m\dr-\int^\infty_0\rho_0\Psi(0,r)r^m\dr\right|\\
        &\le\Vert\Psi(\cdot,0)\Vert_{L^\infty
        ([0,\infty);r^m\dr)}\Vert\rho^a_0-\rho_0\Vert_{L^2([0,\infty);r^m\dr)}+\int^a_0|\Psi(0,r)|r^m\dr\rightarrow0\text{ as }a\rightarrow0.
    \end{align*}
    Thus, $\displaystyle\int^\infty_0\rho_0\Psi(0,r)r^m\dr=\displaystyle\lim_{j\rightarrow\infty}\int^\infty_a\rho^a_0\Psi(0,r)r^m\dr$.
    
    Next, since 
    \begin{align*}
        &\left|\int^t_0\int^\infty_{a_j}\rho^{a_j}u^{a_j}\Psi_r(s,r)r^m\dr\ds-\int^t_0\int^\infty_{\underline{r}(s)}\rho u\Psi_r(s,r)r^m\dr\ds\right|\\
        &\le\left|\int^t_0\int^\infty_{\widetilde{r}^{a_j}(s,x)}\rho^{a_j}u^{a_j}\Psi_r(s,r)r^m\dr\ds-\int^t_0\int^\infty_{\widetilde{r}(s,x)}\rho u \Psi_r(t,r)r^m\dr\ds\right|\\
        &+\Vert\Psi_r\Vert_{L^\infty}\bigg(\int^t_0\int^{\widetilde{r}^{a_j}(t,x)}_{a_j}|\rho^{a_j}u^{a_j }|r^m\dr\ds+\int^t_0\int^{\widetilde{r}(s,x)}_{\underline{r}(t)}|\rho u|r^m\dr\ds\bigg)\\
        &=\sum^7_{i=5}I_i.
    \end{align*}
    For $I_5$, it follows that
    \begin{align*}
        I_5&=\left|\int^t_0\int^\infty_{\widetilde{r}^{a_j}(s,x)}\rho^{a_j}u^{a_j}\Psi_r(s,r)r^m\dr\ds-\int^t_0\int^\infty_{\widetilde{r}(s,x)}\rho u \Psi(t,r)r^m\dr\ds\right|\\
        &\le\left|\int^t_0\int^\infty_{\widetilde{r}(s,x)}\big(\rho^{a_j}u^{a_j}-\rho u\big)\Psi_rr^m\dr\ds\right|+\left|\int^t_0\int^{\widetilde{r}(s,x)}_{\widetilde{r}^{a_j}(s,x)}\rho^{a_j}u^{a_j}\Psi_rr^m\dr\ds\right|\\
        &=I_{5,1}+I_{5,2}.
    \end{align*}
    By virtue of Lemma \ref{last convergence}, we have the point wise bounds: $\rho^{a_j}(t,r)+\rho(t,r)+\sigma^{-\frac{1}{4}}(t)\left(|u^{a_j}(t,r)|+|u(t,r)|\right)\le C(t,x)$ for $r\in[\widetilde{r}(t,x),\infty)$ and the strong convergence 
    \begin{align*}
        \rho^{a_j}\rightarrow\rho,u^{a_j} \rightarrow u \text{ in }[\widetilde{r}(t,x),\infty)\cap\supp\Psi.
    \end{align*}
    Thus we have
    \begin{align*}
        I_{5,1}+I_{5,2}&=\left|\int^t_0\int^\infty_{\widetilde{r}(s,x)}\big(\rho^{a_j}u^{a_j}-\rho u\big)\Psi_rr^m\dr\ds\right|+\left|\int^t_0\int^{\widetilde{r}(s,x)}_{\widetilde{r}^{a_j}(s,x)}\rho^{a_j}u^{a_j}\Psi_rr^m\dr\ds\right|\\
        &\le\left|\int^t_0\int^\infty_{\widetilde{r}(s,x)}\big((\rho^{a_j}-\rho)u^{a_j}-\rho (u-u^{a_j})\big)\Psi_rr^m\dr\ds\right|\\
        &+C(t,x)\sup_{s\in[0,T]}\left|\left(\widetilde{r}(s,x)\right|^n-\left(\widetilde{r}^{a_j}(s,x)\right)^n\right|\cdot\left|\int^t_0\sigma^{-\frac{1}{4}}\ds\right|\\
        &\le\left|\int^t_0\int^\infty_{\widetilde{r}(s,x)}\big((\rho^{a_j}-\rho)u^{a_j}\Psi_rr^m\dr\ds\right|+\left|\int^t_0\int^\infty_{\widetilde{r}(s,x)}\rho (u-u^{a_j})\Psi_rr^m\dr\ds\right|\\
        &+C(\varepsilon,T)\sup_{s\in[0,T]}\left|\left(\widetilde{r}(s,x)\right|^n-\left(\widetilde{r}^{a_j}(s,x)\right)^n\right|\rightarrow0 \text{ as }j\rightarrow\infty.
    \end{align*}
    Hence, $I_5\rightarrow0$ as $j\rightarrow\infty$. Furthermore,
    \begin{align*}
        I_6+I_7&\le\Vert\Psi_r\Vert_{L^\infty}\int^t_0\left(\int^{\widetilde{r}^{a_j}(s,x)}_{a_j}\rho^{a_j}|u^{a_j}|^2r^m\dr\right)^{\frac{1}{2}}\left(\int^{\widetilde{r}^{a_j}(s,x)}_{a_j}\rho^{a_j}r^m\dr\right)^{\frac{1}{2}}\ds\\
        &+\Vert\Psi_r\Vert_{L^\infty}\int^t_0\left(\int^{\widetilde{r}(s,x)}_{\underline{r}(s)}\rho|u|^2r^m\dr\right)^{\frac{1}{2}}\left(\int^{\widetilde{r}(s,x)}_{\underline{r}(s)}\rho r^m\dr\right)^{\frac{1}{2}}\ds\\
        &\le\Vert\Psi_r\Vert_{L^\infty}\sqrt{x}\cdot C(T)t.
    \end{align*}
    Hence for each $\nu>0$, there exists $x''$ sufficient small such that if $x\in(0,x'')$ we have
    \begin{align*}
        I_6+I_7\le\nu.
    \end{align*}
    Thus it follows that
    \begin{align*}
        \int^t_0\int^\infty_{\underline{r}(s)}\rho u\Psi_r(s,r)r^m\dr\ds=\lim_{j\rightarrow\infty}\int^t_0\int^\infty_{a_j}\rho^{a_j}u^{a_j}\Psi_r(s,r)r^m\dr\ds.
    \end{align*}
    Now combining the above results, we have
    \begin{align*}
        \int^\infty_{\underline{r}(t)}\rho(t,r)\Psi(t,r)r^m\dr-\int^\infty_0\rho_0(r)\Psi(0,r)r^m\dr=\int^t_0\int^\infty_{\underline{r}(t)}(\rho\Psi_t+\rho u\Psi_r)r^m\dr\ds.
    \end{align*}
    With the help of our extension $\rho(t,r)=0$ for $r\in[0,\underline{r}(t)]$, we have
    \begin{align*}
        \int^\infty_{0}\rho(t,r)\Psi(t,r)r^m\dr-\int^\infty_0\rho_0(r)\Psi(0,r)r^m\dr=\int^t_0\int^\infty_{0}(\rho\Psi_t+\rho u\Psi_r)r^m\dr\ds.
    \end{align*}
    \par The proof of Lemma \ref{weak form of continuity equation} is finished.
 \end{proof}
\end{lemma}
Now, we use the defect measure, the Hardy's inequality, mollification, and truncations to prove the weak form of the momentum equation.
\begin{lemma}\label{momentum equation}
    For $\forall\Phi\in C^1\left([0,T];C^2_c([0,\infty))\right)$ with $\Phi(t,0)=0$ for each $t\in[0,T]$, 
    \begin{align*}
        &\int^\infty_0\left(\rho u\right)(t,r)\Phi(t,r)r^m\dr-\int^\infty_0\left(\rho_0 u_0\right)(r)\Phi(0,r)r^m\dr\\
        &=\int^t_0\int^\infty_0\left\{\rho u\Phi_t+\rho|u|^2\Phi_r+\left(A\rho^\gamma-\beta(u_r+\dfrac{mu}{r})\right)\left(\Phi_r+\dfrac{m\Phi}{r}\right)\right\}r^m\dr\ds\\
        &+\kappa\int^t_0\int^\infty_0\rho\rho_r\left(\Phi_r+\dfrac{m\Phi}{r}\right)_rr^m\dr\ds+\dfrac{\kappa}{2}\int^t_0\int^\infty_0|\rho_r|^2(3\Phi_r+\dfrac{m\Phi}{r})r^m\dr\ds.
    \end{align*}
    \begin{proof}
    We divide the proof into 11 steps.\\
    \textbf{1. }
        Since $\Psi\in C^1([0,T];C^2_c((0,\infty)))$ there exists $a_0,N>0$ such that $\supp\Psi\subset[a_0,N]$, thus for each $a<a_0$, $\Psi(t,a)=0$. Then for each $a<a_0$, it follows from Definition \ref{5.1} that
        \begin{align*}
       & \int^\infty_a(\rho^au^a)(t,r)\Psi(t,r) r^m\dr\\
       &=\int^\infty_a(\rho^a_0u^a_0)(r)\Psi (0,r)r^m\dr+\int^t_0\int^\infty_a\rho^au^a(\Psi_t+u^a\Psi_r)r^m\dr\ds\\
        &+\int^t_0\int^\infty_a\left(A(\rho^a)^\gamma-\beta\left(u^a_r+\dfrac{mu^a}{r}\right)\right)\left(\Psi_r+\dfrac{m\Psi}{r}\right)r^m\dr\ds\\
        &+\kappa\int^t_0\int^\infty_0\rho^a\rho^a_r\left(\Psi_r+\dfrac{m\Psi}{r}\right)_rr^m\dr\ds+\dfrac{\kappa}{2}\int^t_0\int^\infty_0|\rho^a_r|^2(3\Psi_r+\dfrac{m\Psi}{r})r^m\dr\ds.
    \end{align*}
     By the same argument, the convergence of  $\rho^au^a\Psi_t$ can be proved. Next, for each $r\in[a_0,N]$, it holds that
    \begin{align*}
        |u^a(t,r)|^2=-2\int^\infty_ru^a_r(t,s)u^a(t,s)\ds\le C(T)a_0^{1-m}\int^\infty_r\left(u_r^2+\dfrac{u^2}{s^2}\right)s^m\ds,
    \end{align*}
and the above bound also holds for $|u|^2$. Thus
\begin{align*}
    \left|\int^t_0\int^\infty_0\left(\rho^a|u^a|^2-\rho|u|^2\right)\Psi_rr^m\dr\ds\right|&\le\left|\int^t_0\int^\infty_{\widetilde{r}(t,x)}\left(\rho^a|u^a|^2-\rho|u|^2\right)\Psi_rr^m\dr\ds\right|\\
    &+\left|\int^t_0\int_{a_0}^{\widetilde{r}(t,x)}\left(\rho^a|u^a|^2-\rho|u|^2\right)\Psi_rr^m\dr\ds\right|\\
    &:=\mathcal{I}_1+\mathcal{I}_2.
\end{align*}
By the same argument mentioned before, $\mathcal{I}_1\rightarrow0$ as $a\rightarrow0$. Next, for each $r\in[a_0,N]$, it holds that
\begin{align*}
    \rho^a(t,r)-\rho^a(t,1)&\le\int^{\max\{1,r\}}_{\min\{1,r\}}|\rho^a_r(t,s)|\ds\\
    &\le\left(\int^{\max\{1,r\}}_{\min\{1,r\}}|\rho^a_r|^2s^m\ds\right)^{\frac12}\left(\int^{\max\{1,r\}}_{\min\{1,r\}}\dfrac{1}{s^m}\ds\right)^{\frac12}\\
    &\le C(a_0,T).
\end{align*}
Thus $\rho^a\le C(a_0,T)$, with the same argument, $\rho\le C(a_0,T)$. Now for $\mathcal{I}_2$, it holds that
\begin{align*}
    \mathcal{I}_2&=\left|\int^t_0\int_{a_0}^{\widetilde{r}(t,x)}\left(\rho^a|u^a|^2-\rho|u|^2\right)\Psi_rr^m\dr\ds\right|\\
    &\le \int^t_0\Vert\Psi_r\Vert_{L^\infty}\Vert u^a\Vert_{L^\infty([a_0,\widetilde{r}(t,x)])}\left(\int^\infty_0\rho^a( u^a)^2r^m\dr\right)^{\frac12}\left(\int^{\widetilde{r}(t,x)}_{a_0}\rho^ar^m\dr\right)^{\frac12}\ds\\
    &+\int^t_0\Vert\Psi_r\Vert_{L^\infty}\Vert u\Vert_{L^\infty([a_0,\widetilde{r}(t,x)])}\left(\int^\infty_0\rho u^2r^m\dr\right)^{\frac12}\left(\int^{\widetilde{r}(t,x)}_{a_0}\rho r^m\dr\right)^{\frac12}\ds\\
    &\le C(a_0,T)\left[(\widetilde{r}^a(t,x))^n-(\widetilde{r}(t,x))^n+x\right]^{\frac12} \le C(a_0,T)x^\frac12.
\end{align*}
Thus,$\lim_a\mathcal{I}_2\le C(a_0,T)x^{\frac12}$ as $a\rightarrow0$. Now we let $x\rightarrow0$ to obtain the convergence
\begin{align*}
    \lim_{a\rightarrow0}\int^t_0\int^\infty_0\rho^a|u^a|^2\Psi_rr^m\dr\ds=\int^t_0\int^\infty_0\rho|u|^2\Psi_rr^m\dr\ds.
\end{align*}
The convergence of viscosity term is straightforward and we omit the proof here.\\
\textbf{2. The convergence of nonlinear pressure.}
By virtue of Lemma \ref{last entropy estimate} step 3 and the uniform bound $(\rho^a,\rho)\le C(a_0,T)$ for each $r\in[a_0,N]$, it follows that
\begin{align*}
    \int^N_{a_0}|(\rho^a)^\gamma-\rho^\gamma|r^m\dr\le C(a_0,T)\int^N_{a_0}|\rho^a-\rho|r^m\dr\rightarrow0\text{ as }a\rightarrow0.
\end{align*}
\textbf{3. The convergence of the Korteweg term.}
First, we consider the low-order Korteweg term $\displaystyle\kappa\int^t_0\int^\infty_0\rho^a\rho^a_r\left(\Psi_r+\dfrac{m\Psi}{r}\right)_rr^m\dr\ds$. It is obvious that $\Vert\rho ^a\Vert_{L^\infty_tH^1([a_0,N];\dr)}\le C(a_0,T)$. Since $H^1\hookrightarrow L^\infty$, $\Vert\rho^a\Vert_{L^\infty_tL^\infty_r([a_0,N];\dr)}\le C(a_0,T)$. By virtue of Definition \ref{5.1}, the mass equation holds in the sense of distribution,
\begin{align}\label{mass weak form}
    \rho^a_t+(\rho^au^a)_r+\dfrac{m\rho^au^a}{r}=0\text{ in }\mathcal{D}'((0,T)\times[a,\infty)).
\end{align}
Since $\displaystyle\int^N_{a_0}|\rho^au^a|^2\dr\le C(a_0,T)$, it holds that
\begin{align*}
    \sup_{a<a_0}\Vert\rho^a_t\Vert_{L^\infty(0,T;H^{-1}([a_0,N]))}\le C(a_0,T).
\end{align*}
Thus, by Aubin--Lions--Simon lemma, $\rho^a\rightarrow \rho$ in $C(0,T;L^2([a_0,N]))$. Thus combining the $\rho_r^a{\rightharpoonup}\rho_r$ in $L^2(0,T;L^2((a_0,N);\dr))$, we obtain that
\begin{align*}
\lim_{a\rightarrow0}\int^t_0\int^\infty_0\rho^a\rho^a_r\left(\Psi_r+\dfrac{m\Psi}{r}\right)_rr^m\dr\ds=\int^t_0\int^\infty_0\rho\rho_r\left(\Psi_r+\dfrac{m\Psi}{r}\right)_rr^m\dr\ds.
\end{align*}
Next, we consider the last term $\displaystyle\dfrac{\kappa}{2}\int^t_0\int^\infty_0|\rho^a_r|^2(3\Psi_r+\dfrac{m\Psi}{r})r^m\dr\ds$.\\

\textbf{4. The defect measure of the gradient of density.}
Now we let $Q\Subset(0,T)\times(0,\infty)$, there exists some $a_0,N$, $0<t_1<t_2<T$ such that
$Q\subset(t_1,t_2)\times\left(\dfrac{a_0}{2},2N\right)$. It is straightforward to check that
\begin{align*}
    \int^{t_2}_{t_1}\int^{2N}_{\frac{a_0}{2}}|\rho^{a_j}_r|^2\dr\ds\le C(a_0,T).
\end{align*}
Thus $|\rho^{a_j}_r|^2\dr\ds$ is a non-negative measure in $(t_1,t_2)\times\left(\dfrac{a_0}{2},2N\right)$. Now we define the functional for each $j\in \mathbb{N }$
\begin{align*}
    \mathcal{L}_j(\phi):=\int^{t_2}_{t_1}\int^{2N}_{\frac{a_0}{2}}|\rho^{a_j}_r|^2\Phi\dr\ds\text{ for each }\phi\in C^\infty_c(Q)
\end{align*}
Thus it is easy to check $\mathcal{L}_j$ is a positive bounded linear functional on $C_c(Q)$ equipped with the supremum norm. Since $C_c(Q)$ is dense in $C_0$, thus for each fixed $\phi\in C_0$, for $\forall\varepsilon>0$ there exists a $\phi_\varepsilon\in C_c$ such that $|\phi-\phi_\varepsilon|<\varepsilon$. Since 
\begin{align*}
    \Vert\mathcal{L}_j(\phi_{\varepsilon_1})-\mathcal{L}_j(\phi_{\varepsilon_2})\Vert\le C(T)\Vert\phi_{\varepsilon_1}-\phi_{\varepsilon_2}\Vert_{L^\infty}\rightarrow0\text{ as }\varepsilon_1,\varepsilon_2\rightarrow0. 
\end{align*}
Thus $\{\mathcal{L}_j(\phi_{\varepsilon})\}$ is a Cauchy sequence. Now we let $\mathcal{L}_j(\phi):=\lim_{\varepsilon\rightarrow0}\mathcal{L}_j(\phi_\varepsilon)$, thus, $\mathcal{L}_j$ is a positive bounded linear functional on $C_0(Q)$. By the Banach--Alaoglu theorem, there exist a subsequence $\{j_n\}$ and a $\mathcal{L}\in (C_0(Q))^*$ such that
\begin{align*}
    \mathcal{L}(\phi)=\lim_{n\rightarrow\infty}\mathcal{L}_{j_n}(\phi)\text{ for each }\phi\in C_0(Q),\quad0\le\mathcal{L}(\phi) \text{ for each }0\le\phi\in C_0(Q).
\end{align*}
By the Riesz representation theorem, there exists a non-negative finite Radon measure $\nu$ such that
\begin{align*}
    \mathcal{L}(\phi)=\int_Q\phi\mathrm{d}\nu,\quad0\le\nu(Q)\le C(a_0,T).
\end{align*}
Now, for each $0\le\phi\in C_0(Q)$, it holds that $\sqrt{\phi}\rho^{a_j}_r\rightharpoonup\sqrt{\phi}\rho_r$ in $L^2(Q)$. With the help of the weakly lower semi-continuity of Sobolev norm, the following inequality holds
\begin{align*}
    \int^{t_2}_{t_1}\int^{2N}_{\frac{a_0}{2}}\phi|\rho_r|^2\dr\ds\le\liminf_{j\rightarrow\infty}\int^{t_2}_{t_1}\int^{2N}_{\frac{a_0}{2}}\phi|\rho^{a_j}_r|^2\dr\ds=\int_Q\phi\mathrm{d}\nu.
\end{align*}
Now, for each Borel set $E\subset Q$, we define $\mathfrak{d}(E)=\nu(E)-\displaystyle\iint_E|\rho_r|^2\dr\ds$. Thus  $\mathfrak{d}$ is a non-negative measure on $Q$.  Since 
\begin{align}\label{limit density deirivative}
    \iint_Q|\rho^{a}_r-\rho_r|^2\Phi\dr\ds&=\iint_Q(|\rho^{a}_r|^2-2\rho^a_r\rho_r+|\rho_r|^2)\Phi\dr\ds\nonumber\\
    &\rightarrow\int_Q\phi\mathrm{d}\nu-\iint_Q|\rho_r|^2\phi\dr\ds=\int_Q\phi
    \mathrm{d}\mathfrak{d}\text{ as }j\rightarrow\infty,
\end{align}
where we have used $\rho^{a_j}_r\overset{*}{\rightharpoonup}\rho_r$ in $L^\infty(0,T;L^2((\frac{a_0}{2},2N);\dr))$. This shows that
\begin{align*}
    |\rho^{a_j}_r-\rho_r|^2\dr\ds\overset{*}{\rightharpoonup}\mathfrak{d}\text{ in }\mathcal{M}_{loc}((0,T)\times(0,\infty)).
\end{align*}
Based on the fact $\rho_r^{a_j}\rightarrow\rho_r$ locally uniformly in $\{\mathcal{F}\backslash\{t=0\}\}_{loc}$, $\mathfrak{d}=0$ in $\mathcal{F}$. In the following steps, we prove that $\mathfrak{d} $ vanish on $(0,T)\times(0,\infty)\backslash\{\mathcal{F}\backslash\{t=0\}\}$, hence the strong convergence holds. In the next step, we prove  $\mathfrak{d}=0$ on the vacuum interface $\{(t,r):r=\underline{r}(t)\}$. \\
\textbf{5. $\mathfrak{d}=0$ on the vacuum interface.} For  $\chi\in C_c(0,T)$ and $\xi_l\in C^\infty_c(0,\infty)$, where $\xi_l\rightarrow1 $ for each $r>0$  as $l\rightarrow\infty$ and $0\le\xi_l\le1$. Since
\begin{align*}
0\le\int_Q\chi(t)\xi_l(r)r^m\mathrm{d}\mathfrak{d}\int_Q\chi(t)\xi_l(r)r^m\mathrm{d}\nu=\lim_{j\rightarrow\infty}\int^{t_2}_{t_1}\chi(t)\int^{2N}_{\frac{a_0}{2}}\xi_l|\rho^{a_j}_r|^2r^m\dr\ds\le C(T)\int^T_0\chi(t)\dt.
\end{align*}
Passing $l\rightarrow\infty$, it holds that
\begin{align*}
    \int_Q\chi(t)r^m\mathrm{d}\nu\le C(T)\int^T_0\chi(t)\dt.
\end{align*}
For each Borel set $E\subset(0,T)$, we define $\dot{\pi}(E)=\displaystyle\int_{ E\times(0,\infty) }r^m\mathrm{d}\mathfrak{d}$. Thus $\dot{\pi}(E)\le C(T)|E|$, where $|\cdot|$ denote the Lebesgue measure. Thus by the Radon-Nikodym theorem, there exists a non-negative measurable function $h(t)$ such that $\dfrac{\mathrm{d}\dot{\pi}}{\dt}=h(t)$. Furthermore, it holds that
\begin{align*}
    \int_Eh(t)\dt\le C(T)|E| \text{ for each Borel set }E\subset(0,T).
\end{align*}
Thus, by the Lebesgue differentiation theorem, it holds that $h(t)\le C(T)$ for  $a.e. ~t\in(0,T)$. With the help of Lemma \ref{Disintergration theorem}, there exists a unique family of probability measure kernel $\{\pi_t\}_{t\in(0,T)}$ called the time-disintegration of $\mathfrak{d}$ such that $r^m\mathrm{d}\mathfrak{d}=\mathrm{d}\pi_t\mathrm{d}\dot{\pi}$. Now we define $\mathrm{d}\mathfrak{d}_t=h(t) r^{-m}\mathrm{d}\pi_t$, it holds that
\begin{align}\label{time slices}
    \int^\infty_0r^m\mathrm{d}\mathfrak{d}_t=h(t)\int^\infty_0 \mathrm{d}\pi_t=h(t)\le C(T)\text{ for }a.e. t\in(0,T)
\end{align}
Now, for each $\phi\in C^\infty_c((0,T);C^\infty_c(0,\infty))$, we define 
\begin{align*}
    <\mathfrak{d},\Phi>=\int\Phi\mathrm{d}\mathfrak{d},\quad <\mathfrak{d}_r,\Phi>=-\int\Phi_r\mathrm{d}\mathfrak{d}.
\end{align*}
Based on step 1--3 and \eqref{limit density deirivative}, it holds that for each $\phi\in C^\infty_c((0,T);C^\infty_c(0,\infty))$ that
\begin{align}\label{weak forms without time terms}
    0&=\int^T_0\int^\infty_0\left(\rho u\phi_t+\rho u^2\phi_r+A\rho^\gamma\left(\phi_r+\dfrac{m\phi}{r}\right)-\beta\left(u_r+\dfrac{mu}{r}\right)\left(\phi_r+\dfrac{m\phi}{r}\right)\right)r^m\dr\ds\nonumber\\
    &+\kappa\int^T_0\int^\infty_0\rho\rho_r\left(\phi_r+\dfrac{m\phi}{r}\right)_rr^m\dr\ds+\dfrac{\kappa}{2}\int^T_0\int^\infty_0|\rho_r|^2\left(3\phi_r+\dfrac{m\phi}{r}\right)r^m\dr\ds\nonumber\\
    &+\dfrac{\kappa}{2}\int^T_0\int^\infty_0\left(3\phi_r+\dfrac{m\phi}{r}\right)r^m\mathrm{d}\mathfrak{d}.
\end{align}
Thus, it holds for the possibly vacuum ball that
\begin{align}\label{the defect measure and limit velocity}
    \beta\left(u_r+\dfrac{mu}{r}\right)_r=\dfrac{\kappa}{2}\left(3\mathfrak{d}_r+\dfrac{2m\mathfrak{d}}{r}\right)\text{ in }\mathcal{D}'(\mathcal{V}).
\end{align}
Now let $I=[t_1,t_2]\Subset\{t:\underline{r}(t)>0\}$, there exist two positive constants $M_1,M_2$ such that
\begin{align*}
    0<M_1\le \underline{r}(t)\le M_2\text{ for each }t\in[t_1,t_2].
\end{align*}
Let $\chi_1\in C^\infty_c((-1,1))$ and $\chi_1=1$ on $\left[-\dfrac{1}{2},\dfrac{1}{2}\right]$. We define 
\begin{align}
    \varphi_h(t,r)=\eta(t)(r-\underline{r}(t))\chi_1\left(\dfrac{r-\underline{r}(t)}{h}\right),
\end{align}
where $0<h<\min\left\{\dfrac{1}{2},\dfrac{M_1}{4}\right\}$ and $\eta\in C^\infty_c((t_1,t_2)),0\le\eta\le1$. Thus it follows that
\begin{align}\label{properties of varphi}
    &\supp \varphi_h\subset\{(t,r):t\in\supp\eta,|r-\underline{r}(t)|\le h\},\quad\varphi_h(t,\underline{r}(t))=0,\quad(\varphi_h)_r(t,\underline{r}(t))=\eta(t)\nonumber\\
    &\varphi_h(t,r)=\eta(t)[r-\underline{r}(t)]\text{ if }|r-\underline{r}(t)|\le \dfrac{h}{2},\quad|\varphi_h|\le Ch,\quad|(\varphi_h)_r|\le C\Vert\eta\Vert_{L^\infty},\nonumber\\
    &(\varphi_h)_t\in L^2\left(I;L^\infty\left(\dfrac{M_1}{2},M_2\right)\right),\quad(\varphi_h)_t(t,\underline{r}(t))=-\eta(t)\underline{r}'(t)\text{ for each }t\in(t_1,t_2). 
\end{align}
Since $\underline{r}(t)$ is not smooth, we need to mollify the vacuum interface. First we extend $\underline{r}(t)$ by
\begin{align}\label{extension of interface}
\underline{r}^{ext}(t)=
    \begin{cases}
        \underline{r}(t_1)\quad t<t_1,\\
        \underline{r}(t)\quad t_1\le t\le t_2,\\
        \underline{r}(t_2)\quad t>t_2.
    \end{cases}
\end{align}
It follows from \eqref{extension of interface} that $ M_1\le\underline{r}^{ext}(t)\le M_2$. Now we mollify the extended vacuum interface by $\underline{r}^\epsilon=\vartheta_\epsilon*\underline{r}^{ext}(t)$, where $\vartheta_\epsilon$ is the standard mollifier. Thus the $\underline{r}^\epsilon$ satisfies 
\begin{align*}
    \underline{r}^\epsilon\rightarrow\underline{r} \text{ in }H^1((t_1,t_2))\text{ as } \epsilon\rightarrow0,\quad M_1\le\underline{r}^\epsilon\le M_2.
\end{align*}
Based on the construction of $\underline{r}^\epsilon$, we define the mollified cut-off function by
\begin{align*}
    \varphi_{h,\epsilon}=\eta(t)(r-\underline{r}^\epsilon(t))\chi_1\left(\dfrac{r-\underline{r}^\epsilon(t)}{h}\right).
\end{align*}
With the help of the properties of mollifier, it holds that
\begin{align*}
    &(\varphi_{h,\epsilon},(\varphi_{h,\epsilon})_r,(\varphi_{h,\epsilon})_{rr})\rightarrow(\varphi_{h},(\varphi_{h})_r,(\varphi_{h})_{rr})\quad a.e.  \text{ as }\epsilon\rightarrow0,\\
    &(\varphi_{h,\epsilon})_t\rightarrow(\varphi_h)_t\text{ in }L^2(I;L^\infty).
\end{align*}
Thus we take the mollified cut-off function $\varphi_{h,\epsilon}$ into \eqref{weak forms without time terms}  and let $\epsilon\rightarrow0$ to obtain  that
\begin{align*}
    0&=\int^T_0\int^\infty_0\left(\rho u(\varphi_h)_t+\rho u^2(\varphi_h)_r+A\rho^\gamma\left((\varphi_h)_r+\dfrac{m\varphi_h}{r}\right)-\beta\left(u_r+\dfrac{mu}{r}\right)\left((\varphi_h)_r+\dfrac{m\varphi_h}{r}\right)\right)r^m\dr\ds\nonumber\\
    &+\kappa\int^T_0\int^\infty_0\rho\rho_r\left((\varphi_h)_r+\dfrac{m\varphi_h}{r}\right)_rr^m\dr\ds+\dfrac{\kappa}{2}\int^T_0\int^\infty_0|\rho_r|^2\left(3(\varphi_h)_r+\dfrac{m\varphi_h}{r}\right)r^m\dr\ds\nonumber\\
    &+\dfrac{\kappa}{2}\int_{(t_1,t_2)\times(0,\infty)}\left(3(\varphi_h)_r+\dfrac{m\varphi_h}{r}\right)r^m\mathrm{d}\mathfrak{d}.
\end{align*}
Now we let $h\rightarrow0$. With the help of  \eqref{properties of varphi} and $\rho(t,\underline{r}(t))=0$, it follows from the dominated convergence theorem that
\begin{align*}
    \int^T_0\int^\infty_0\left(\rho u(\varphi_h)_t+\rho u^2(\varphi_h)_r+A\rho^\gamma\left((\varphi_h)_r+\dfrac{m\varphi_h}{r}\right)-\beta\left(u_r+\dfrac{mu}{r}\right)\left((\varphi_h)_r+\dfrac{m\varphi_h}{r}\right)\right)r^m\dr\ds\overset{h\rightarrow0}{\longrightarrow}0.
\end{align*}
Next for each $r\in(\underline{r}(t),\underline{r}(t)+h)$, it holds by virtue of $\rho(t,\underline{r}(t))=0$ that
\begin{align*}
    |\rho(t,r)|^2\le|r-\underline{r}(t)|\int^{\underline{r}(t)+h}_{\underline{r}(t)}|\rho_r(t,s)|^2\ds
\end{align*}
Thus 
\begin{align*}
    \dfrac{1}{h^2}\int^{\underline{r}(t)+h}_{\underline{r}(t)}|\rho(t,r)|^2\dr\le C(M_1,M_2,T)\int^{\underline{r}(t)+h}_{\underline{r}(t)}|\rho_r(t,s)|^2s^m\ds.
\end{align*}
Since $|(\varphi_{h})_{rr}|\le \dfrac{C}{h}$, it holds that
\begin{align*}
    \left|\kappa\int^T_0\int^\infty_0\rho\rho_r\left((\varphi_h)_r+\dfrac{m\varphi_h}{r}\right)_rr^m\dr\ds\right|&\le \dfrac{C(M_1,M_2,T)}{h}\int^{t_2}_{t_1}\int^{\underline{r}(t)+h}_{\underline{r}(t)}|\rho\rho_r|r^m\dr\ds\\
    &\le C(M_1,M_2,T)\int^{t_2}_{t_1}\int^{\underline{r}(t)+h}_{\underline{r}(t)}|\rho_r|^2r^m\dr\ds\rightarrow0\text{ as }h\rightarrow0.
\end{align*}
By the same argument, $\displaystyle\int^T_0\int^\infty_0|\rho_r|^2\left(3(\varphi_h)_r+\dfrac{m\varphi_h}{r}\right)r^m\dr\ds\rightarrow0 $ as $h\rightarrow0$. Using \eqref{properties of varphi} to the defect measure term, we obtain that
\begin{align*}
    \lim_{h\rightarrow0}\dfrac{\kappa}{2}\int_{(t_1,t_2)\times(0,\infty)}\left(3(\varphi_h)_r+\dfrac{m\varphi_h}{r}\right)r^m\mathrm{d}\mathfrak{d}=\dfrac{3\kappa}{2}\int_{\{(t,r):t\in(t_1,t_2),r=\underline{r}(t)\}}\eta(t)r^m\mathrm{d}\mathfrak{d}
\end{align*}
Combining above, it holds that
\begin{align*}
    \int_{\{(t,r):t\in(t_1,t_2),r=\underline{r}(t)\}}\eta(t)r^m\mathrm{d}\mathfrak{d}=0.
\end{align*}
Thus $\mathfrak{d}(\{(t,r):r=\underline{r}(t),t\in(t_1,t_2)\})=0$. By virtue of the argument mentioned in Lemma \ref{interface speed}, it is obvious that $\mathfrak{d}(\{(t,r):r=\underline{r}(t)\})=0$.\\
\textbf{6. $\mathfrak{d}$ is locally square-integrable function.}
Since 
\begin{align*}
    \left(r^{\frac{2m}{3}}\mathfrak{d}\right)_r=\left[\dfrac{2\beta}{3\kappa}\left(u_r+\dfrac{mu}{r}\right)r^{\frac{2m}{3}}\right]_r-\dfrac{4m\beta}{9\kappa}\left(u_r+\dfrac{mu}{r}\right)r^{\frac{2m}{3}-1}\text{ in }\mathcal{D}'(\mathcal{V}).
\end{align*}
For $I=[t_1,t_2]\Subset\{t:\underline{r}(t)>0\}$, we choose $0<r_0<\min_{t\in[t_1,t_2]}\underline{r}(t)$ and define 
\begin{align}\label{local vi impulse}
    \mathcal{W}(t,r):=\int^r_{r_0}s^{\frac{2m}{3}-1}\left(u_r(t,s)+\dfrac{mu(t,s)}{s}\right)\ds\text{ for }r\in(0,r_0).
\end{align}
and we set 
\begin{align}\label{distri}
    \mathcal{S}:=r^{\frac{2m}{3}}\left[\mathfrak{d}-\dfrac{2\beta}{3\kappa}\left(u_r+\dfrac{mu}{r}\right)\right]+\dfrac{4m\beta}{9\kappa}\mathcal{W}\text{ in }\mathcal{D}'((t_1,t_2)\times(0,r_0)) 
\end{align}
Then it holds that $\mathcal{S}_r=0$. Thus by Lemma \ref{distribution derivatives}, there exists $\mathscr{S}\in\mathcal{D}'((t_1,t_2))$ such that $\mathcal{S}=1\Box \mathscr{S}$. Now we prove that $\mathscr{S}$ is in fact a square-integrable function. For $\forall \chi\in C^\infty_c\left(\left(\dfrac{r_0}{2},\dfrac{3r_0}{4}\right)\right)\ge0$ with $\int\chi=1$, we define
\begin{align*}
    \nu_\chi(E):=\int_{ E\times(0,r_0)}\chi(r)r^{\frac{2m}{3}}\mathrm{d}\mathfrak{d}\text{ for each E }\subset(t_1,t_2)
\end{align*}
Then $0\le\nu_\chi(E)\le\Vert\chi r^{-\frac{m}{3}}\Vert_{L^\infty}\int_{(0,r_0)\times E}r^m\mathrm{d}\mathfrak{d}\le C(T)\Vert\chi r^{-\frac{m}{3}}\Vert_{L^\infty}|E|$, where $|\cdot|$ denotes the Lebesgue measure. By the Radon-Nikodym theorem, there exists a  function $A(t)\in L^\infty(I)$ such that
\begin{align*}
    \mathrm{d}\nu_{\chi}=A(t)\mathrm{d}t,\quad 0\le A(t)\le C(T)\Vert\chi r^{-\frac{m}{3}}\Vert_{L^\infty}.
\end{align*}
By virtue of \eqref{local vi impulse} and \eqref{distri}, it holds for each $\eta\in  C^\infty_c(t_1,t_2)$ that
\begin{align*}
    \langle\mathscr{S},\eta\rangle=\dfrac{4m\beta}{9\kappa}\int_I\eta\int^{r_0}_0\chi(r)\mathcal{W}\dr\dt+\int_I\eta(t)A(t)\dt-\dfrac{2\beta}{3\kappa}\int_I\eta\int^{r_0}_0\chi r^{\frac{2m}{3}}\left(u_r+\dfrac{mu}{r}\right)\dr\dt.
\end{align*}
By the definition of \eqref{local vi impulse},  it holds that 
\begin{align*}
    \dfrac{4m\beta}{9\kappa}\int_I\eta\int^{r_0}_0\chi(r)\mathcal{W}\dr\dt+\int_I\eta(t)A(t)\dt-\dfrac{2\beta}{3\kappa}\int_I\eta\int^{r_0}_0\chi r^{\frac{2m}{3}}\left(u_r+\dfrac{mu}{r}\right)\dr\dt\le C(r_0,T)\Vert\eta\Vert_{L^2(I)}.
\end{align*}
Thus $\mathscr{S}\in L^2(I)$. Since $\mathcal{S}=1\Box\mathscr{S}$, thus for each $\psi\in C^\infty_c((t_1,t_2)\times(0,r_0))$ such that
\begin{align*}
    \langle\mathcal{S},\psi\rangle=\int^{t_2}_{t_1}\mathscr{S}(t)\int^{r_0}_0\psi\dr\dt=\int^{t_2}_{t_1}\int^{r_0}_0\mathscr{S}(t)\psi(t,r)\dr\dt.
\end{align*}
Based on the definition of \eqref{distri}, it follows for each $\psi\in C^\infty_c((t_1,t_2)\times(0,r_0))$that
\begin{align}\label{distri defect}
    \langle \mathfrak{d},\psi\rangle=\left\langle \dfrac{2\beta}{3\kappa}\left(u_r+\dfrac{mu}{r}\right),\psi\right\rangle+\int^{t_2}_{t_1}\int^{r_0}_0\mathscr{S}(t)\psi(t,r)r^{-\frac{2m}{3}}\dr\dt-\left\langle\dfrac{4m\beta}{9\kappa}r^{-\frac{2m}{3}}\mathcal{W},\psi\right\rangle
\end{align}
Then it holds that,
\begin{align*}
    \mathfrak{d}=\dfrac{2\beta}{3\kappa}\left(u_r+\dfrac{mu}{r}\right)+r^{-\frac{2m}{3}}\mathscr{S}(t)-\dfrac{4m\beta}{9\kappa}r^{-\frac{2m}{3}}\mathcal{W}\text{ in }\mathcal{D}'((0,r_0)\times(t_1,t_2)).
\end{align*}
Since 
\begin{align*}
    \int^{t_2}_{t_1}\int^{r_0}_0\left|r^{-\frac{2m}{3}}\mathscr{S}(t)\right|^2r^m\dr\ds=\int^{t_2}_{t_1}|\mathscr{S}|^2\ds\int^{r_0}_0r^{-\frac{m}{3}}\dr<\infty,
\end{align*}
where we used $m<3$. Now we estimate for the last term. Since $u_r+\dfrac{mu}{r}\in L^2(t_1,t_2;L^2((r,r_0);r^m\dr))$, thus for  $a.e. ~t\in(t_1.t_2)$, $\displaystyle\int|u_r+\dfrac{mu}{r}|^2r^m\dr<\infty$, thus
\begin{align*}
    &\int^{\max\{r,r_0\}}_{\min\{r,r_0\}}s^{\frac{2m}{3}-1}\left|u_r(t,s)+\dfrac{mu(t,s)}{s}\right|\ds\\
    &\le\left(\int^{\max\{r,r_0\}}_{\min\{r,r_0\}}\left|u_s+\dfrac{mu}{s}\right|^2s^m\ds\right)^\frac{1}{2}\left(\int^{\max\{r,r_0\}}_{\min\{r,r_0\}}s^{\frac{m}{3}-2}\ds\right)^{\frac{1}{2}}<\infty
\end{align*}
since for each $\epsilon>0$, it follows that
\begin{align*}
    \int^{r_0}_\epsilon r^{-\frac{m}{3}}|\mathcal{W}|^2\dr&=\dfrac{r^{1-\frac{m}{3}}}{1-\frac{m}{3}}|\mathcal{W}|^2\bigg|^{r_0}_\epsilon-\dfrac{2}{1-\frac{m}{3}}\int^{r_0}_\epsilon r^{1-\frac{m}{3}}\mathcal{W}\mathcal{W}_r\dr\\
    &\le\dfrac{2}{1-\frac{m}{3}}\left(\int^{r_0}_\epsilon r^{-\frac{m}{3}}|\mathcal{W}|^2\dr\right)^{\frac12}\left(\int^{r_0}_\epsilon r^{m}\left|u_r+\dfrac{mu}{r}\right|^2\dr\right)^{\frac12}.
\end{align*}
Thus, $\displaystyle\int^{r_0}_\epsilon r^{-\frac{m}{3}}|\mathcal{W}|^2\dr\le C\int^{r_0}_0\left|u_r+\dfrac{mu}{r}\right|^2r^m\dr$. By the monotone convergence theorem, it holds that$ \displaystyle\int^{r_0}_0 r^{-\frac{m}{3}}|\mathcal{W}|^2\dr\le C\int^{r_0}_0\left|u_r+\dfrac{mu}{r}\right|^2r^m\dr$. Thus,
\begin{align*}
   \int^{t_2}_{t_1}\int^{r_0}_0r^{-\frac m3}|\mathcal{W}|^2\dr\ds\le C(T).
\end{align*}
As a result, it holds that $\mathfrak{d}\in L^2((t_1,t_2);L^2(0,r_0;r^m\dr))$. Next, we define $\mathcal{V}_I:=\{(t,r):t_1<t<t_2,0<r<\underline{r}(t)\}$ and extend \eqref{distri} to $\mathcal{V}_I$ where the corresponding \eqref{local vi impulse} is also extended to $\mathcal{V}_I$. Thus it holds that $\mathcal{S}_r=0$ in $\mathcal{D}'(\mathcal{V}_I)$.

For $\forall(r_*,t_*)\in \mathcal{V}_I$, there exists a constant $c_2$ such that $\max\{r_*,r_0\}<c_2<\underline{r}(t_*)$. Since $\underline{r}$ is continuous, there exists an open interval $J$ such that $t_*\in J\Subset(t_1,t_2)$ and $c_2<\underline{r}(t)$ for each $t\in J$. Furthermore, it also holds that $ J\times\left(\dfrac{r_0}{2},c_2\right)\subset\mathcal{V}_I$. In addition, we have $\mathcal{S}_r=0$ on $J\times\left(\dfrac{r_0}{2},c_2\right)$. By the same argument, there exists $\mathscr{S}_J$ such that $\mathcal{S}=1\Box \mathscr{S}_J$. Let $\chi\in C^\infty_c\left(\dfrac{r_0}{2},r_0\right),\eta(t)\in C^\infty_c(J)$ and $\displaystyle\int^{r_0}_{\frac{r_0}{2}}\chi(r)\dr=1$, it holds that 
\begin{align*}
    \langle  \mathscr{S}_J,\eta\rangle=\langle \mathcal{S},\chi\eta\rangle=\langle \mathscr{S},\eta\rangle.
\end{align*}
Thus $\mathscr{S}_J=\mathscr{S}|_J$. Since these local representations agree with $1\Box\mathscr S$
and the corresponding rectangles, together with
$(t_1,t_2)\times(0,r_0)$, cover $\mathcal V_I$, we conclude that
\[
    S=1\Box\mathscr S
    \qquad\text{in }\mathcal D'(\mathcal V_I).
\] 
With the same argument in \eqref{distri defect}, it holds for each $\psi\in C^\infty_c\left((t_1,t_2)\times\left(\frac{r_0}{2},c_2\right)\right)$ that
\begin{align}
    \langle \mathfrak{d},\psi\rangle=\iint_{\mathcal{V}_I} \dfrac{2\beta}{3\kappa}\left(u_r+\dfrac{mu}{r}\right)\psi\dr\dt+\iint_{\mathcal{V}_I}\mathscr{S}(t)\psi(t,r)r^{-\frac{2m}{3}}\dr\dt-\iint_{\mathcal{V}_I}\dfrac{4m\beta}{9\kappa}r^{-\frac{2m}{3}}\mathcal{W}\psi\dr\dt
\end{align}
Thus, with the help of Cauchy--Schwarz inequality, it is straightforward to check $\mathfrak{d}\in L^2(\mathcal{V}_I;r^m\dr\dt)$.\\
\textbf{7. Local relative energy and its limit form.}
Now we define the local relative energy $e^a=\dfrac{1}{2}\rho^a\left|u^a\right|^2+\dfrac{A}{\gamma-1}(\rho^a)^\gamma+\dfrac{\kappa}{2}\left|\rho^a_r\right|^2$ for $r\in[a,\infty)$. Thus, it holds in the sense of distribution that
\begin{align}\label{local relative energy}
    e^a_t+\dfrac{1}{r^m}\left(r^mu^a(e^a+A(\rho^a)^\gamma)\right)_r&=\dfrac{\beta}{r^m}(r^mu^au^a_r)_r-\beta\left(|u^a_r|^2+\dfrac{m|u^a|^2}{r^2}\right)\nonumber\\
    &+\dfrac{\kappa}{r^m}\left[r^m\left((\rho^au^a\rho^a_r)_r-\dfrac{3}{2}|\rho^a_r|^2u^a-2\rho^a\rho^a_ru^a_r\right)\right]_r
\end{align}
For $\forall\Omega\Subset(0,\infty)$, it holds that
\begin{align*}
    \Vert e^a\Vert_{L^\infty((t_1,t_2);L^1(\Omega;\dr))}+\Vert \rho^a\Vert_{L^\infty((t_1,t_2);H^1(\Omega;\dr))}+\Vert u^a\Vert_{L^2((t_1,t_2);H^1(\Omega;\dr))}\le C(\Omega,a_1,t_2)
\end{align*}
Thus for each $\zeta\in C^\infty_c(\Omega)$, it holds that $\Vert\zeta e^a_t\Vert_{L^1((t_1,t_2);H^{-3}(\Omega;\dr))}\le C(\Omega,a_1,t_2)<\infty$. Since 
\begin{align*}
    L^1\hookrightarrow \hookrightarrow H^{-s}\hookrightarrow H^{-3} \text{ for each } \dfrac{1}{2}<s<1,
\end{align*}
where the detail proof we refer the reader to Brezis \cite{brezis} \S 6. With the help of Aubin--Lions--Simon lemma, it holds that
\begin{align*}
    \zeta e^a\overset{a\rightarrow0}{\longrightarrow}\zeta\bar e\text{ in }L^2((t_1,t_2);H^{-s}(\Omega;\dr)).
\end{align*}
We let $I\times J_1\Subset\mathcal{V}$. Since $\rho^a\rightarrow\rho$ in $C(0,T;L^2_{loc}(0,\infty))$ and $\Vert\rho^a\Vert_{L^\infty_tH^1(J_1;\dr)}\le C(J_1,T)$, we have $\Vert\rho^a\Vert_{L^\infty(I\times J_1)}\rightarrow0$ as $a\rightarrow0$. Thus,
\begin{align*}
    \dfrac{\kappa}{2}|\rho^a_r|^2\zeta\rightarrow\zeta\bar e\text{ in  }L^2((t_1,t_2);H^{-s}(J_1;\dr)).
\end{align*}
In addition, we have $\rho^a|u^a|^2\rightarrow0$ in $L^1(I\times J_1)$, thus Lemma \ref{last entropy estimate}, $\rho^a|u^a|^2\rightarrow0$ in $L^2(I\times J_1)$.  Since $|\rho^a_r|^2\rightarrow\mathfrak{d}$ in $\mathcal{D}'(\mathcal{V})$, by the uniqueness of limit, it follows that $\mathfrak{d}=\dfrac{2}{\kappa}\bar e$. Thus,
\begin{align*}
    |\rho^a_r|^2\rightarrow\mathfrak{d}\text{ in }L^2(I;H^{-s}(J_1)).
\end{align*}
Hence $u^a|\rho^a_r|^2\rightarrow\mathfrak{d}u$ in $\mathcal{D}'(\mathcal{V})$, where we used $\mathfrak{d}\in L^2$. 
\par Since for each $I\times J_1\Subset\mathcal{V}$, it holds that
\begin{align*}
    \sup_j\iint_{I\times J_1}\left(|u^{a_j}_r|^2+|u^{a_j}|^2\right)\dr\ds<\infty.
\end{align*}
With the same argument mentioned in step 4, there exist non-negative finite Radon measure $\overline{|u|^2}$ and $\overline{|u_r|^2}$ such that
\begin{align*}
    |u^{a_j}|^2\dr\ds\overset{*}{\rightharpoonup}\overline{|u|^2},\quad |u^{a_j}_r|^2\dr\ds\overset{*}{\rightharpoonup}\overline{|u_r|^2}\text{ in }\mathcal{M}_{loc}(\mathcal{\mathcal{V}}).
\end{align*}
Now we define $\mathcal{Z}_1=\overline{|u|^2}-|u|^2\dr\ds$ and  $\mathcal{Z}_2=\overline{|u_r|^2}-|u_r|^2\dr\ds$. It is  clear  that $\mathcal{Z}_i,i=1,2$ are non-negative. Now we pass the limit to \eqref{local relative energy},
\begin{align*}
    \mathfrak{d}_t+\dfrac{4}{r^m}(r^mu\mathfrak{d})_r=\dfrac{\beta}{\kappa r^m }\left(r^m(\overline{|u|^2})_r\right)_r-\dfrac{2\beta}{\kappa}\left(\overline{|u_r|^2}+\dfrac{m\overline{|u|^2}}{r^2}\right)\text{ in }\mathcal{D}'(\mathcal{V}).
\end{align*}
Using \eqref{the defect measure and limit velocity}, $\mathcal{Z}_1=\overline{|u|^2}-u^2$, and $\mathcal{Z}_2=\overline{|u_r|^2}-|u_r|^2$, it holds that
\begin{align}\label{limit local relative energy}
    \mathfrak{d}_t+\dfrac{1}{r^m}(r^mu\mathfrak{d})_r=\dfrac{\beta}{\kappa r^m}\left(r^m(\mathcal{Z}_1)_r\right)_r-\dfrac{2\beta}{\kappa}\left(\mathcal{Z}_2+\dfrac{m\mathcal{Z}_1}{r^2}\right)-\left(3u_r+\dfrac{mu}{r}\right)\mathfrak{d}.
\end{align}
\textbf{8. Modified velocity.}
Let $w_j=u^{a_j}-u$, it follows from step 3 and step7 that
\begin{align*}
    |w_j|^2\dr\ds\overset{*}{\rightharpoonup}\mathcal{Z}_1\text{ in }\mathcal{M}_{loc}(\mathcal{V}),
\end{align*}
Let $I=[t_1,t_2]\Subset\{t:\underline{r}(t)>0\}$, for each $h>0$, it holds that
\begin{align*}
    \int_I|u^{a_j}(t,\underline{r}(t))-u(t,\underline{r}(t))|^2\dt&\le C(t_1,t_2)h\int_I\int^{\underline{r}(t)+h}_{\underline{r}(t)}\left(|u_r|^2+|u^{a_j}_r|^2\right)\dr\ds\\
    &+C\int_I|u^{a_j}(t,\underline{r}(t)+h)-u(t,\underline{r}(t)+h)|^2\dt.
\end{align*}
We first let $j\rightarrow$ then pass $h\rightarrow0$ to obtain that
\begin{align*}
    \lim_{j\rightarrow\infty}\int_I|u^{a_j}(t,\underline{r}(t))-u(t,\underline{r}(t))|^2\dt=0.
\end{align*}
Thus $u^{a_j}(t,\underline{r}(t))\rightarrow u(t,\underline{r}(t))$ in $L^2(I)$. Now we consider the modified velocity
\begin{align}\label{modified velocity}
    m^{a_j}(t,r):=u^{a_j}(t,r)-u(t,r)-\dfrac{r}{\underline{r}(t)}[u^{a_j}(t,\underline{r}(t))-u(t,\underline{r}(t))]=w_j(t,r)-\frac{r}{\underline{r}(t)}w_j(t,\underline{r}(t)).
\end{align}
It is obvious that $m^{a_j}(t,\underline{r}(t))=0$, thus
\begin{align*}
    &\int^{\underline{r}(t)}_0\left|\dfrac{r}{\underline{r}(t)}\left(u^{a_j}(t,\underline{r}(t))-u(t,\underline{r}(t))\right)\right|^2r^m\dr=\dfrac{|u^{a_j}(t,\underline{r}(t))-u(t,\underline{r}(t))|^2}{|\underline{r}(t)|^2}\int^{\underline{r}(t)}_0r^{m+2}\dr\\
    &=\dfrac{|\underline{r}(t)|^{m+1}}{m+3}|u^{a_j}(t,\underline{r}(t))-u(t,\underline{r}(t))|^2
\end{align*}
Thus, 
\begin{align*}
    \left\Vert\dfrac{r}{\underline{r}(t)}\left(u^{a_j}(t,\underline{r}(t))-u(t,\underline{r}(t))\right)\right\Vert_{L^2(\mathcal{V}_I;r^m\dr\dt)}\rightarrow0 \text{ as }  j\rightarrow\infty.
\end{align*}
Furthermore, it follows that
\begin{align*}
    \left\Vert\dfrac{1}{\underline{r}(t)}w_j(t,\underline{r}(t))\right\Vert_{L^2(\mathcal{V}_I;r^m\dr\dt)}\rightarrow0\text{ as }j\rightarrow\infty.
\end{align*}
Since $|m^{a_j}|^2-|w_j|^2=\left|\dfrac{r}{\underline{r}(t)}w_j(t,\underline{r}(t))\right|^2-2w_j(t,r)\dfrac{r}{\underline{r}(t)}w_j(t,\underline{r}(t))$, we have
\begin{align*}
    \int_{\mathcal{V}_I}\left||m^{a_j}|^2-|w_j|^2\right|r^m\dr\ds\rightarrow0\text{ as }j\rightarrow\infty.
\end{align*}
Thus $|m^{a_j}|^2r^m\dr\ds\overset{*}{\rightharpoonup}r^m\mathcal{Z}_1$ in $\mathcal{M}_{loc}(\mathcal{V})$. This shows that the modification preserves the defect measure associated with the squared velocity. In the next two steps, we use the defect measure identity \eqref{limit local relative energy}, obtained by passing to the limit in the local relative energy identity \eqref{local relative energy} as $a\to0$, to prove that the defect measure $\mathfrak{d}$ vanishes in the vacuum region.\\
\textbf{9. Inequality of the defect momentum equation.}
By virtue of Lemma \ref{Hardy}, it follows that: $\displaystyle\sup_j\iint_{\mathcal{V}_I}\left(|m^{a_j}_r|^2+\dfrac{m|m^{a_j}|^2}{r^2}\right)r^m\dr\dt\le C$, $m^{a_j}(t,\underline{r}(t))=0$. Therefore,
\begin{align*}
    \iint_{\mathcal{V}_I}\left(\dfrac{1}{(\underline{r}(t)-r)^2}+\dfrac{1}{r^2}\right)|m^{a_j}|^2r^m\dr\dt\le C.
\end{align*}
By choosing $\xi_n\in C_c(\mathcal{V}_I)$, $0\le\xi_n\le1$, and $\xi_n\rightarrow1$ as $n\rightarrow\infty$, we have
\begin{align*}
    \int_{\mathcal{V}_I}\xi_n\left(\dfrac{1}{(\underline{r}(t)-r)^2}+\dfrac{1}{r^2}\right)r^m\mathrm{d}\mathcal{Z}_1=\lim_{j\rightarrow\infty}\iint_{\mathcal{V}_I}\xi_n\left(\dfrac{1}{(\underline{r}(t)-r)^2}+\dfrac{1}{r^2}\right)|m^{a_j}|^2r^m\dr\dt\le C.
\end{align*}
With the help of the monotone convergence theorem, it holds that $\displaystyle\int_{\mathcal{V}_I}\left(\dfrac{1}{(\underline{r}(t)-r)^2}+\dfrac{1}{r^2}\right)r^m\mathrm{d}\mathcal{Z}_1\le C$. Now we define $\mathcal{A}_h:=\{(t,r)\in\mathcal{V}_I:0<\underline{r}(t)-r<2h\}$ and $\mathfrak{A}_\epsilon=\{(t,r)\in\mathcal{V}_I:0<r<2\epsilon\}$. Then, by virtue of the dominated convergence theorem, we have
\begin{align*}
    \int_{\mathcal{A}_h}\dfrac{r^m}{(\underline{r}(t)-r)^2}\mathrm{d}\mathcal{Z}_1\rightarrow0\text{ as }h\rightarrow0,\quad\int_{\mathfrak{A}_\epsilon}r^{m-2}\mathrm{d}\mathcal{Z}_1\rightarrow0\text{ as }\epsilon\rightarrow0.
\end{align*}
Now we define the cut-off function $\chi_{h,\epsilon}(t,r)$ by
\begin{align*}
    \chi_{h,\epsilon}(t,r)=\vartheta\left(\dfrac{r}{\epsilon}\right)\vartheta\left(\dfrac{\underline{r}(t)-r}{h}\right),
\end{align*}
where $\vartheta$ is a smooth function with the following properties: $0\le\vartheta\le1$,   $\vartheta=0 $ for $r\le1$, and $\vartheta=1$ for $r\ge2$. Let $\min_{t\in[t_1,t_2]}\underline{r}(t):= c_I$, and we choose $\epsilon$ and $h$ sufficiently small $2\epsilon+2h<\underline{r}(t)$ such that $\supp\vartheta'\left(\dfrac{r}{\epsilon}\right)\cap\supp\vartheta'\left(\dfrac{\underline{r}(t)-r}{h}\right)=\emptyset$. By a straightforward calculation,
\begin{align*}
    \int_{\mathfrak{A}_\epsilon}\left|\vartheta\left(\dfrac{r}{\epsilon}\right)_{rr}+\dfrac{m\vartheta\left(\dfrac{r}{\epsilon}\right)_r}{r}\right|r^m\mathrm{d}\mathcal{Z}_1\le C\epsilon^{-2}\int_{\mathfrak{A}_\epsilon}r^m\mathrm{d}\mathcal{Z}_1+C\epsilon^{-1}\int_{\mathfrak{A}_\epsilon}r^m\mathrm{d}\mathcal{Z}_1\le C\int_{\mathfrak{A}_\epsilon}r^{m-2}\mathrm{d}\mathcal{Z}_1\rightarrow0\text{ as }\epsilon\rightarrow0.
\end{align*}
With the same argument,
\begin{align*}
    \int_{\mathcal{A}_h}\left|\vartheta\left(\dfrac{\underline{r}(t)-r}{h}\right)_{rr}+\dfrac{m\vartheta\left(\dfrac{\underline{r}(t)-r}{h}\right)_r}{r}\right|r^m\mathrm{d}\mathcal{Z}_1\le Ch^{-2}\int_{\mathcal{A}_h}r^m\mathrm{d}\mathcal{Z}_1+Ch^{-1}\int_{\mathcal{A}_h}r^m\mathrm{d}\mathcal{Z}_1\rightarrow0\text{ as }h\rightarrow0.
\end{align*}
Combining above, it holds for each $\eta(t)\in C^\infty_c(I)\ge0$ that
\begin{align}\label{second order of cutoff}
    \int_{\mathcal{V}_I}\eta(t)\left(\left(\chi_{\epsilon,h}\right)_{rr}+\dfrac{m(\chi_{\epsilon,h})_r}{r}\right)r^m\mathrm{d}\mathcal{Z}_1\rightarrow0.
\end{align}
We temporarily assume $\underline{r}(t)$ is smooth, the agrument  can be prove rigorously by mollification, we omit here for simplicity. Multiplying $\eqref{limit local relative energy}$ with $r^m\chi_{\epsilon,h}$ and integrating in $\mathcal{V}_I$,
\begin{align}\label{integration defect}
    &-\iint_{\mathcal{V}_I}\mathfrak{d}\eta'\chi_{\epsilon,h}r^m\dr\dt+\iint_{\mathcal{V}_I}\eta\chi_{\epsilon,h}\left(3u_r+\dfrac{mu}{r}\right)\mathfrak{d}r^m\dr\dt\nonumber\\
    &=\iint_{\mathcal{V}_I}\eta\mathfrak{d}\left((\chi_{\epsilon,h})_t+u(\chi_{\varepsilon,h})_r\right)r^m\dr\dt+\int_{\mathcal{V}_I}\eta(t)\dfrac{\beta}{\kappa}\left(\left(\chi_{\epsilon,h}\right)_{rr}+\dfrac{m(\chi_{\epsilon,h})_r}{r}\right)r^m\mathrm{d}\mathcal{Z}_1\nonumber\\
    &-\dfrac{2\beta}{\kappa}\int_{\mathcal{V}_I}r^m\eta\chi_{\epsilon,h}\mathrm{d}\mathcal{Z}_2-\dfrac{2m\beta}{\kappa}\int_{\mathcal{V}_I}\eta\chi_{\epsilon,h}r^{m-2}\mathrm{d}\mathcal{Z}_1
\end{align}
For the first term in the right hand side, 
\begin{align*}
    &\iint_{\mathcal{V}_I}\eta\mathfrak{d}\left((\chi_{\epsilon,h})_t+u(\chi_{\varepsilon,h})_r\right)r^m\dr\dt\\&=\iint_{\mathcal{V}_I}\eta\mathfrak{d}\left[\vartheta\left(\dfrac{r}{\epsilon}\right)\dfrac{-u+\underline{r}'(t)}{h}\vartheta'\left(\dfrac{\underline{r}(t)-r}{h}\right)+\dfrac{u}{\epsilon}\vartheta'\left(\dfrac{r}{\epsilon}\right)\vartheta\left(\dfrac{\underline{r}(t)-r}{h}\right)\right]r^m\dr\dt\\
    &:=\mathscr{R}_1+\mathscr{R}_2
\end{align*}
\par For $\mathscr{R}_1$, we define $d(t,r)=u(t,r)-\dfrac{\underline{r}'(t)}{\underline{r}(t)}r$, thus it follows from Lemma \ref{interface speed} that $d(t,\underline{r}(t))=0$. In addition, 
\begin{align*}
    \iint_{\mathcal{V}_I}|d_r|^2r^m\dr\dt<\infty.
\end{align*}
Thus by Lemma \ref{Hardy}, $\displaystyle\iint_{\mathcal{V}_I}\dfrac{|d|^2}{(\underline{r}(t)-r)^2}r^m\dr\dt<\infty$. Furthermore, it follows that
\begin{align*}
\mathscr{R}_1&=\iint_{\mathcal{V}_I}\eta\mathfrak{d}\vartheta\left(\dfrac{r}{\epsilon}\right)\dfrac{-u+\underline{r}'(t)}{h}\vartheta'\left(\dfrac{\underline{r}(t)-r}{h}\right)r^m\dr\dt\\
    &\le C\iint_{\mathcal{A}_h}|\mathfrak{d}|\left(\dfrac{|d|}{\underline{r}-r}+\dfrac{|\underline{r}'|}{\underline{r}}\right)r^m\dr\dt\rightarrow0\text{ as }h\rightarrow0.
\end{align*}
Since $\supp{\vartheta'\left(\dfrac{r}{\epsilon}\right)}\subset[\epsilon,2\epsilon]$, for $\mathscr{R}_2$, we have
\begin{align*}
    \mathscr{R}_2&=\iint_{\mathcal{V}_I}\eta\mathfrak{d}\dfrac{u}{\epsilon}\vartheta'\left(\dfrac{r}{\epsilon}\right)\vartheta\left(\dfrac{\underline{r}(t)-r}{h}\right)r^m\dr\dt\\
    &\le C\left\Vert\dfrac{u}{r}\right\Vert_{L^2(\mathcal{V}_I;r^m\dr\dt)}\left(\iint_{\mathfrak{A}_\epsilon}|\mathfrak{d}|^2r^m\dr\dt\right)^{\frac12}\rightarrow0\text{ as }\epsilon\rightarrow0.
\end{align*}
With the help of the dominated convergence theorem, it follows that
\begin{align*}
    &\iint_{\mathcal{V}_I}\mathfrak{d}\eta'\chi_{\epsilon,h}r^m\dr\dt\rightarrow\iint_{\mathcal{V}_I}\mathfrak{d}\eta'r^m\dr\dt\text{ as }\epsilon,h\rightarrow0,\\
    &\iint_{\mathcal{V}_I}\eta\mathfrak{d}\chi_{\epsilon,h}\left(3u_r+\dfrac{mu}{r}\right)r^m\dr\dt\rightarrow\iint_{\mathcal{V}_I}\eta\mathfrak{d}\left(3u_r+\dfrac{mu}{r}\right)r^m\dr\dt\text{ as }\epsilon,h\rightarrow0,
\end{align*}
Using the monotone convergence theorem, it holds that
\begin{align}\label{cut off defect measure}
    \int_{\mathcal{V}_I}r^m\eta\chi_{\epsilon,h}\mathrm{d}\mathcal{Z}_2+\int_{\mathcal{V}_I}\eta\chi_{\epsilon,h}r^{m-2}\mathrm{d}\mathcal{Z}_1\rightarrow\int_{\mathcal{V}_I}r^m\eta\mathrm{d}\mathcal{Z}_2+\int_{\mathcal{V}_I}\eta r^{m-2}\mathrm{d}\mathcal{Z}_1\ge0
\end{align}
Taking limit to $\epsilon,h$, it follows from  \eqref{integration defect} that
\begin{align}\label{defect momentum inequality}
    \dfrac{\mathrm{d}}{\dt}\int^{\underline{r}(t)}_0\mathfrak{d}r^m\dr+\int^{\underline{r}(t)}_0\left(3u_r+\dfrac{mu}{r}\right)\mathfrak{d}r^m\dr\le0 \text{ in }\mathcal{D}'((t_1,t_2)).
\end{align}
\textbf{10. $\mathfrak{d}=0$.}
 Let $d_{\eta,\epsilon,h}=\eta(t)\vartheta\left(\dfrac{r}{\epsilon}\right)\vartheta\left(\dfrac{\underline{r}(t)-r}{h}\right)d(t,r)$, where $\eta(t)\in C^\infty_c((t_1,t_2))\ge0$. By the same argument in step 8, it holds that
\begin{align*}
    \iint_{\mathcal{V_I}}|(d_{\eta,\epsilon,h})_r-d_r\eta|^2r^m\dr\dt&\le\iint_{\mathcal{V_I}}\left|\vartheta\left(\dfrac{r}{\epsilon}\right)\vartheta\left(\dfrac{\underline{r}(t)-r}{h}\right)-1\right|^2|\eta d_r|^2r^m\dr\dt\\
    &+\iint_{\mathcal{V_I}}\left|\vartheta'\left(\dfrac{r}{\epsilon}\right)\dfrac{1}{\epsilon}\vartheta\left(\dfrac{\underline{r}(t)-r}{h}\right)\eta d\right|^2 r^m\dr\dt\\
    &+\iint_{\mathcal{V_I}}\left|\vartheta\left(\dfrac{r}{\epsilon}\right)\dfrac{1}{h}\vartheta'\left(\dfrac{\underline{r}(t)-r}{h}\right)\eta d\right|^2 r^m\dr\dt\\
    &\rightarrow0 \text{ as }\epsilon,h\rightarrow0.
\end{align*}
Thus we obtain the conclusion: $(d_{\eta,\epsilon,h})_r\rightarrow \eta d_r$ in $L^2(\mathcal{V}_I)$. By the same argument, it holds that $\dfrac{d_{\eta,\epsilon,h}}{r}\rightarrow \dfrac{d\eta}{r}$ in $L^2(\mathcal{V}_I)$. Then we extend $d_{\eta,\epsilon,h}$ to zero in $(0,\infty)\backslash(0,\underline{r}(t))$ and denote this extended function by $\bar{d}_{\eta,\epsilon,h}$. 
\par Now we mollify the extended function as $d_{\eta,\epsilon,h,\sigma}=\chi_\sigma*\bar{d}_{\eta,\epsilon,h}$, where $\chi_\sigma$ is the standard 2D mollifier. It is obvious that
$\supp d_{\eta,\epsilon,h}\subset\{(t,r):t\in\supp\eta,\epsilon\le r\le\underline{r}(t)-h\}:=\bar{\mathcal{K}}$. Let $0<\sigma<\dfrac{1}{2}\text{dist}\{\bar{\mathcal{K}},\mathbb{R}^2\backslash\mathcal{V}_I\}$, it is straightforward to check $d_{\eta,\epsilon,h,\sigma}\in C^\infty_c(\mathcal{V}_I)$, and $\left((d_{\eta,\epsilon,h,\sigma})_r,d_{\eta,\epsilon,h,\sigma}\right)\rightarrow \left((d_{\eta,\epsilon,h})_r,d_{\eta,\epsilon,h}\right)$ in $L^2(\mathcal{V}_I;r^m\dr\dt))$. Based on $\mathfrak{d}\in L^2$,  we multiply   \eqref{the defect measure and limit velocity} with $d_{\eta,\epsilon,h,\sigma}$ and  pass the limit to $\sigma,\epsilon,h\rightarrow0$, in that order, then it follows that
\begin{align*}
    \beta\int^{t_2}_{t_1}\eta(t)\int^{\underline{r}(t)}_0\left(u_r+\dfrac{mu}{r}\right)\left(d_r+\dfrac{md}{r}\right)r^m\dr\dt=\dfrac{\kappa}{2}\int^{t_2}_{t_1}\eta(t)\int^{\underline{r}(t)}_0\mathfrak{d}(3d_r+\dfrac{md}{r})r^m\dr\dt.
\end{align*}
Since 
\begin{align*}
    d_r+\dfrac{md}{r}=u_r+\dfrac{mu}{r}-(m+1)\dfrac{\underline{r}'}{\underline{r}},\quad 3d_r+\dfrac{md}{r}=3u_r+\dfrac{mu}{r}-(m+3)\dfrac{\underline{r}'}{\underline{r}}.
\end{align*}
It holds that
\begin{align}\label{urd}
    \int^{\underline{r}(t)}_0\left(3u_r+\dfrac{mu}{r}\right)\mathfrak{d}r^m\dr=\dfrac{2\beta}{\kappa}\int^{\underline{r}(t)}_0\left(u_r+\dfrac{mu}{r}\right)\left(d_r+\dfrac{md}{r}\right)r^m\dr+(m+3)\dfrac{\underline{r}'(t)}{\underline{r}(t)}\int^{\underline{r}(t)}_0\mathfrak{d}r^m\dr.
\end{align}
By virtue of Lemma \ref{interface speed}, $r^mu\rightarrow0$ as $r\rightarrow0$, and $\displaystyle\int^{\underline{r}(t)}_0\left(u_r+\dfrac{mu}{r}-(m+1)\dfrac{\underline{r}'}{\underline{r}}\right)r^m\dr=0$, we have
\begin{align}\label{urdr}
    \displaystyle\int^{\underline{r}(t)}_0\left(u_r+\dfrac{mu}{r}\right)\left(d_r+\dfrac{md}{r}\right)r^m\dr&=\int^{\underline{r}(t)}_0\left|u_r+\dfrac{mu}{r}-(m+1)\dfrac{\underline{r}'}{\underline{r}}\right|^2r^m\dr\nonumber\\
    &+(m+1)\dfrac{\underline{r}'}{\underline{r}}\int^{\underline{r}(t)}_0\left(u_r+\dfrac{mu}{r}-(m+1)\dfrac{\underline{r}'}{\underline{r}}\right)r^m\dr\nonumber\\
    &=\int^{\underline{r}(t)}_0\left|u_r+\dfrac{mu}{r}-(m+1)\dfrac{\underline{r}'}{\underline{r}}\right|^2r^m\dr.
\end{align}
Combining \eqref{urd} and \eqref{urdr}, it holds that
\begin{align}\label{ulti urd}
    \int^{t_2}_{t_1}\eta\int^{\underline{r}(t)}_0\left(3u_r+\dfrac{mu}{r}\right)\mathfrak{d}r^m\dr\dt&=\dfrac{2\beta}{\kappa}\int^{t_2}_{t_1}\eta\int^{\underline{r}(t)}_0\left|u_r+\dfrac{mu}{r}-(m+1)\dfrac{\underline{r}'}{\underline{r}}\right|^2r^m\dr\dt\nonumber\\
    &+(m+3)\int^{t_2}_{t_1}\eta\dfrac{\underline{r}'(t)}{\underline{r}(t)}\int^{\underline{r}(t)}_0\mathfrak{d}r^m\dr\dt.
\end{align}
Taking \eqref{ulti urd} into \eqref{defect momentum inequality}, it follows that
\begin{align}\label{mathfrak d inequality}
    \dfrac{\mathrm{d}}{\dt}\int^{\underline{r}(t)}_0\mathfrak{d}r^m\dr+\dfrac{2\beta}{\kappa}\int^{\underline{r}(t)}_0\left|u_r+\dfrac{mu}{r}-(m+1)\dfrac{\underline{r}'}{\underline{r}}\right|^2r^m\dr+(m+3)\dfrac{\underline{r}'(t)}{\underline{r}(t)}\int^{\underline{r}(t)}_0\mathfrak{d}r^m\dr\le0 \text{ in }\mathcal{D}'((t_1,t_2)).
\end{align}
Since for each $0\le\psi\in C^\infty_c(t_1,t_2)$, 
\begin{align*}
    \left\langle\dfrac{\mathrm{d}}{\dt}\int^{\underline{r}(t)}_0\mathfrak{d}r^m\dr+\dfrac{2\beta}{\kappa}\int^{\underline{r}(t)}_0\left|u_r+\dfrac{mu}{r}-(m+1)\dfrac{\underline{r}'}{\underline{r}}\right|^2r^m\dr+(m+3)\dfrac{\underline{r}'(t)}{\underline{r}(t)}\int^{\underline{r}(t)}_0\mathfrak{d}r^m\dr,\psi\right\rangle\le0 
\end{align*}
If we let 
\begin{align*}
    \mathscr{T}:=-\displaystyle\dfrac{\mathrm{d}}{\dt}\int^{\underline{r}(t)}_0\mathfrak{d}r^m\dr-\dfrac{2\beta}{\kappa}\int^{\underline{r}(t)}_0\left|u_r+\dfrac{mu}{r}-(m+1)\dfrac{\underline{r}'}{\underline{r}}\right|^2r^m\dr-(m+3)\dfrac{\underline{r}'(t)}{\underline{r}(t)}\int^{\underline{r}(t)}_0\mathfrak{d}r^m\dr,
\end{align*}
there exists a non-negative Radon measure $\nu_2$ such that $\mathscr{T}=\nu_2$. It follows from \eqref{time slices} that
\begin{align*}
    \int_{(0,T)\times(0,\infty)}f(t,r)\mathrm{d}\mathfrak{d}=\int^T_0\left(\int^\infty_0f(t,r)\mathrm{d}\mathfrak{d}_t\right)\dt,
\end{align*}
where $f$ is measurable function. Combining $\mathfrak{d}\in L^2$. For $\forall E\subset(t_1,t_2)$, we choose 
\begin{align*}
    f_E=1_E(t)1_{\{0<r<\underline{r}(t)\}}r^m,
\end{align*}
so it holds that
\begin{align*}
    \int_E\int^{\underline{r}(t)}_0\mathfrak{d}r^m\dr\dt=\int_{(0,T)\times(0,\infty)}1_E(t)1_{\{0<r<\underline{r}(t)\}}r^m\mathrm{d}\mathfrak{d}=\int_E\left(\int_{(0,\underline{r}(t))}r^m\mathrm{d}\mathfrak{d}_t\right)\dt.
\end{align*}
Thus $\displaystyle\int^{\underline{r}(t)}_0\mathfrak{d}r^m\dr=\int_{(0,\underline{r}(t))}r^m\mathrm{d}\mathfrak{d}_t\le C(T)$ for $a.e. ~t\in(t_1,t_2)$. In addition, it holds that the total variation measure of $\displaystyle-\dfrac{2\beta}{\kappa}\int^{\underline{r}(t)}_0\left|u_r+\dfrac{mu}{r}-(m+1)\dfrac{\underline{r}'}{\underline{r}}\right|^2r^m\dr-(m+3)\dfrac{\underline{r}'(t)}{\underline{r}(t)}\int^{\underline{r}(t)}_0\mathfrak{d}r^m\dr$ is finite. Since $\nu_2$ is a non-negative Radon measure, it holds for each $J\Subset(t_1,t_2)$ that $\nu_2(J)<\infty$. Combining the above it follows that $\left|\displaystyle\dfrac{\mathrm{d}}{\dt}\int^{\underline{r}(t)}_0\mathfrak{d}r^m\dr\right|(J)<\infty$, therefore,  $\displaystyle\int^{\underline{r}(t)}_0\mathfrak{d}r^m\dr\in BV_{loc}((t_1,t_2))$. Next, since $\underline{r}(t)$ is absolutely continuous, we can multiply \eqref{mathfrak d inequality} with $\underline{r}^{m+3}(t)$ to obtain that
\begin{align}\label{last mathfrak d}
    \dfrac{\mathrm{d}}{\dt}\left(\underline{r}^{m+3}\int^{\underline{r}(t)}_0\mathfrak{d}r^m\dr\right)+\dfrac{2\beta}{\kappa}\underline{r}^{m+3}(t)\int^{\underline{r}(t)}_0\left|u_r+\dfrac{mu}{r}-(m+1)\dfrac{\underline{r}'}{\underline{r}}\right|^2r^m\dr\le0 \text{ in }\mathcal{D}'((t_1,t_2)).
\end{align}
Now we are able to prove $\mathfrak{d}=0$ in $\mathcal{V}_I$. First we mollify $\displaystyle\underline{r}^{m+3}\int^{\underline{r}(t)}_0\mathfrak{d}r^m\dr$,
\begin{align*}
    \left(\underline{r}^{m+3}\int^{\underline{r}(t)}_0\mathfrak{d}r^m\dr\right)_\epsilon(t)=\int^{t_2}_{t_1}\left(\underline{r}^{m+3}\int^{\underline{r}(t)}_0\mathfrak{d}r^m\dr\right)(\tau)\vartheta_\epsilon(t-\tau)\mathrm{d}\tau,
\end{align*}
where $\vartheta_\epsilon$ is the non-negative  standard mollifier. Thus,
\begin{align*}
    \left(\underline{r}^{m+3}\int^{\underline{r}(t)}_0\mathfrak{d}r^m\dr\right)_\epsilon'(t)&=\int^{t_2}_{t_1}\left(\underline{r}^{m+3}\int^{\underline{r}(t)}_0\mathfrak{d}r^m\dr\right)(\tau)\vartheta'_\epsilon(t-\tau)\mathrm{d}\tau\\
    &=\left\langle\left(\underline{r}^{m+3}\int^{\underline{r}(t)}_0\mathfrak{d}r^m\dr\right)_t,\vartheta_\epsilon(t-)\right\rangle\le0\text{ for each }t\in(t_1+\epsilon,t_2-\epsilon).
\end{align*}
Now we let $\epsilon$ sufficiently small, it holds for each $t_1<s<t<t_2$ that 
\begin{align*}
    \left(\underline{r}^{m+3}\int^{\underline{r}(t)}_0\mathfrak{d}r^m\dr\right)_\epsilon(s)\ge\left(\underline{r}^{m+3}\int^{\underline{r}(t)}_0\mathfrak{d}r^m\dr\right)_\epsilon(t),
\end{align*} 
and with the help of the Lebesgue differentiation theorem, we have
\begin{align}
    \left(\underline{r}^{m+3}\int^{\underline{r}(t)}_0\mathfrak{d}r^m\dr\right)(s)\ge\left(\underline{r}^{m+3}\int^{\underline{r}(t)}_0\mathfrak{d}r^m\dr\right)(t)
\end{align}
Thus for  $a.e. ~t\in(t_1,t_2)$, there exists a subsequence $s_j:s_j\rightarrow t_1$ such that $0\le\left(\underline{r}^{m+3}\int^{\underline{r}(t)}_0\mathfrak{d}r^m\dr\right)(t)\le\left(\underline{r}^{m+3}\int^{\underline{r}(t)}_0\mathfrak{d}r^m\dr\right)(s_j)\le C(T)\underline{r}^{m+3}(s_j)$. Since $\{t:\underline{r}(t)>0\}$ is an open set, there exists at most countably open interval $(t^n_1,t^n_2)$ such that $\{t:\underline{r}(t)>0\}=\cup_n(t^n_1,t^n_2)$ and $\underline{r}(t^n_1)=0$. Now we define the approximate time sequence 
\begin{align*}
    I^j_n=\left(t^n_1+\dfrac{t^n_2-t^n_1}{j+2},t^n_2-\dfrac{t^n_2-t^n_1}{j+2}\right)\subset\left[t^n_1+\dfrac{t^n_2-t^n_1}{j+2},t^n_2-\dfrac{t^n_2-t^n_1}{j+2}\right]\Subset(t^n_1,t^n_2),
\end{align*}
Since $\mathcal{V}_{I^j_n}\subset\mathcal{V}_{I^{j+1}_{n}}$ for each $n\in\mathbb{N}$ and $\cup_{j=1}\mathcal{V}_{I^j_n}=\mathcal{V}_{(t^n_1,t^n_2)}$, thus for $a.e. ~t\in(t^n_1,t^n_2)$ there exists a subsequence $s_j$ such that $s_j\rightarrow t^n_1$ with $0\le\left(\underline{r}^{m+3}\int^{\underline{r}(t)}_0\mathfrak{d}r^m\dr\right)(t)\le C(T)\underline{r}^{m+3}(s_j)\rightarrow0$. Thus , it holds that
\begin{align*}
    \left(\int^{\underline{r}(t)}_0\mathfrak{d}r^m\dr\right)(t)=0\text{ for } a.e.  t\in(t^n_1,t_2^n).
\end{align*}
This fact illustrates that $\mathfrak{d}=0~a.e. $  in $\mathcal{V}_{(t_1^n,t_2^n)}$, and thus, $\mathfrak{d}=0~a.e. $  in $\mathcal{V}$. Then for each $\Psi\in C^1([0,T];C^2_c(0,\infty))$, it holds that
\begin{align}\label{intermid weak form of momentum}
        &\int^\infty_0\left(\rho u\right)(t,r)\Psi(t,r)r^m\dr-\int^\infty_0\left(\rho_0 u_0\right)(r)\Psi(0,r)r^m\dr\nonumber\\
        &=\int^t_0\int^\infty_0\left\{\rho u\Psi_t+\rho|u|^2\Psi_r+\left(A\rho^\gamma-\beta(u_r+\dfrac{mu}{r})\right)\left(\Psi_r+\dfrac{m\Psi}{r}\right)\right\}r^m\dr\ds\nonumber\\
        &+\kappa\int^t_0\int^\infty_0\rho\rho_r\left(\Psi_r+\dfrac{m\Psi}{r}\right)_rr^m\dr\ds+\dfrac{\kappa}{2}\int^t_0\int^\infty_0|\rho_r|^2(3\Psi_r+\dfrac{m\Psi}{r})r^m\dr\ds.
    \end{align}
\textbf{11. Weak forms for $\Phi$.} Now for each $\Phi\in C^1([0,T];C^2_c[0,\infty))$ with $\Phi(t,0)=0$. Let $\chi^h$ be a smooth, increasing cut-off function satisfying $\chi^h=0$ for $0\le r\le h$ and $\chi^h(r)=1$ for $r\ge2h$, also, it holds that $|\chi^h_r|\le\dfrac{C}{h},|\chi^h_{rr}|\le\dfrac{C}{h^2}$.
Set $\chi^h\Phi$ be the test function in \eqref{intermid weak form of momentum}, with the help of the dominated convergence theorem, we can conclude the desired results except for the term $\kappa\displaystyle\int^t_0\int^\infty_0\rho\rho_r\left((\chi^h\phi)_r+\dfrac{m(\chi^h\phi)}{r}\right)_r r^m\dr\ds$. 
Now , we start the proof with three-dimensional case. For some fixed  $R_0>0$ and each $t\in[0,T]$,   by virtue of  Lemma \ref{last entropy estimate}, there exists the corresponding $s_t\in[R_0,2 R_0]$ such that $\rho(s_t,t)\le C(T)$. Therefore, it holds that 
\begin{align*}
    |\rho(R_0,t)-\rho(s_t,t)|\le \left(\int^{2R_0}_{R_0}\rho^2_r r^2\dr\right)^{\frac12}\left(\int^{2R_0}_{R_0}\dfrac{1}{r^2}\dr\right)^{\frac12}\le C(T).
\end{align*}
Thus $\rho(R_0,t)\le C(T)$. For each $0<r<R_0$, it holds that
\begin{align*}
    |\rho(t,r)-\rho(R_0,t)|&\le\int^{R_0}_r|\rho_r(t,s)|\dr\le \left(\int^{R_0}_r|\rho_r|^2s^2\ds\right)^{\frac12}\left(\int^{R_0}_r\dfrac{1}{s^2}\ds\right)^{\frac12}\\
    &\le C(T)\left(\dfrac{1}{r}-\dfrac{1}{R_0}\right)^{\frac12}\le C(T)\left(1+r^{-\frac12}\right).
\end{align*}
Thus $\rho(t,r)\le C(T)\left(1+r^{-\frac12}\right)$ for $0<r<R_0$. Since
\begin{align*}
    &\kappa\int^t_0\int^\infty_0\rho\rho_r\left((\chi^h\Phi)_r+\dfrac{m(\chi^h\Phi)}{r}\right)_rr^2\dr\ds\\
    &=\kappa\int^t_0\int^\infty_0\rho\rho_r\chi^h\left(\Phi_r+\dfrac{m\Phi}{r}\right)_rr^2\dr\ds+\kappa\int^t_0\int^\infty_0\left(\chi^h_{rr}\Phi+2\chi^h_r\Phi_r+\dfrac{m\Phi\chi^h_r}{r}\right)\rho\rho_rr^2\dr\ds.
\end{align*}
It follows from $\Phi(t,0)=0$ that
\begin{align*}
    &\left|\kappa\int^t_0\int^\infty_0\left(\chi^h_{rr}\Phi+2\chi^h_r\Phi_r+\dfrac{m\Phi\chi^h_r}{r}\right)\rho\rho_rr^2\dr\ds\right|\\
    &\le C(T)\left(\int^t_0\int^{2h}_h|\rho_r|^2r^2\dr\ds\right)^{\frac12}\left(\dfrac{1}{h^2}\int^t_0\int^{2h}_h\rho^2r^2\dr\ds\right)^{\frac12}\\
    &\le C(T)\left(\int^{2h}_h|\rho_r|^2r^2\dr\right)^{\frac12}\rightarrow0\text{ as }h\rightarrow0.
\end{align*}
Therefore the desired conclusion has been obtained. Next we prove the weak form of two-dimensions. With a similar argument with respect to three-dimensional case, it follows that
\begin{align*}
    \rho^2(t,r)\le C(T)\left(1+\log\dfrac{R_0}{r}\right)\text{ for }0<r<R_0,
\end{align*}
where we choose $R_0$ sufficiently large. Let $\chi^h$ be the smooth cut-off function mentioned above, with the same argument, it holds that
\begin{align}\label{h term of korteweg}
    &\kappa\int^t_0\int^\infty_0\left(\chi^h_{rr}\Phi+2\chi^h_r\Phi_r+\dfrac{m\Phi\chi^h_r}{r}\right)\rho\rho_rr\dr\ds\nonumber\\
    &\le C(T)\left(\int^t_0\int^{2h}_h|\rho_r|^2r\dr\ds\right)^{\frac12}\left(\dfrac{1}{h^2}\int^t_0\int^{2h}_h\rho^2r\dr\ds\right)^{\frac12}\nonumber\\
    &\le C(T)\left(1+\log\dfrac{R_0}{h}\right)^{\frac12}\left(\int^{2h}_h|\rho_r|^2r\dr\right)^{\frac12}.
\end{align}
Now we choose $h_j=R_02^{-j}$ and let $f_j=\displaystyle\int^T_0\int^{2h_j}_{h_j}|\rho_r|^2 r\dr\dt$,  it holds that $\Sigma^\infty_{j=1}f_j<\infty$. Thus $\liminf_{j\rightarrow\infty}j f_j=0$, so there exists a subsequence $\{j_k f_{j_k}\}^\infty_{k=1}$ such that $j_k f_{j_k}<\dfrac{1}{k^2}$. Since $\log \dfrac{R_0}{h_{j_k}}=j\log 2$, it holds that
\begin{align*}
    &\left|\kappa\int^t_0\int^\infty_0\left(\chi^{h_{j_k}}_{rr}\Phi+2\chi^{h_{j_k}}_r\Phi_r+\dfrac{m\Phi\chi^{h_{j_k}}_r}{r}\right)\rho\rho_rr\dr\ds\right|\\
&\le C(T)\left(\left(1+j_k\right)\int^T_0\int^{2h_{j_k}}_{h_{j_k}}|\rho_r|^2r\dr\dt\right)^{\frac12}\le \dfrac{C(T)}{k}\rightarrow0\text{ as }k\rightarrow\infty.
\end{align*}
Thus we let $h=h_k$ in \eqref{h term of korteweg} and let $k\rightarrow\infty$ to obtain the desired results.
\par The proof of Lemma \ref{momentum equation} is finished.
    \end{proof}
\end{lemma}
\section{Non-Formation of Vacuum States  in Two Spatial Dimensions}\label{6}
 In this section, we rule out the possibly vacuum ball via the momentum-impulse logarithmic contradiction. Specifically, the unweighted viscosity impulse tends to infinity as $r\rightarrow0$ with speed $\log\dfrac{1}{r}$ in the vacuum ball, however, the unweighted momentum tends to infinity as $r\rightarrow0$ only with speed $\sqrt[4]{\left(\log\dfrac{1}{r}\right)^3}$, this leads to the contradiction if other terms can be controlled by a finite bound. Furthermore, the decay rates of $\rho$ as $r\rightarrow0$ are characterized. 
\begin{lemma}\label{Ruling out the possible vacuum}
    \begin{align}
        \underline{r}(t)=0\text{ for }\forall t\in[0,T].
    \end{align}
    \begin{proof}
        By virtue of Lemma \ref{momentum equation}, it holds that
\begin{align*}
    \left(u_r+\dfrac{mu}{r}\right)_r=0\text{ in }\mathcal{D}'(\mathcal{V}).
\end{align*}
With the help of Lemma \ref{last entropy estimate} there exists a $c(t)$ such that $u(t,r)=c(t)r$ for $0<r<\underline{r}(t)$. It is obvious that $c(t)\in L^2_{loc}(\{\underline{r}(t)>0\})$. Since the following inequality holds for  $a.e. ~t\in[0,T]$
\begin{align*}
    |u(t,r)|^2-|u(t,\epsilon)|=2\int^r_\epsilon u_ru\ds\le\int^r_\epsilon\left(|u_r|^2s+\dfrac{|u|^2}{s}\right)\ds<\infty
\end{align*}
Thus $\lim_{\epsilon\rightarrow0}|u(\epsilon,t)|^2<\infty$ and the entropy estimates Lemma \ref{last entropy estimate} force that $u(t,0)=0$ for  $a.e. ~t\in[0,T]$. Let $I=[t_1,t_2]\Subset\{\underline{r}(t)>0\}$, it follows from Lemma \ref{interface speed} that
\begin{align*}
    \underline{r}'(t)=u(t,\underline{r}(t))=c(t)\underline{r}(t)\text{ for each }t\in[t_1,t_2].
\end{align*}
Integrating $c(t)$ from $s\rightarrow t$, we obtain that
\begin{align}\label{infinite vis impulse}
    \int^t_sc(\tau)\mathrm{d}\tau=\log\dfrac{\underline{r}(t)}{\underline{r}(s)}.
\end{align}
Let $r_*\ge\sup_{t\in[0,T]}\underline{r}(t)$ be a fixed positive constant, by virtue of $\rho^\gamma\le C(1+G(\rho))$, it holds that
 $\displaystyle\int^{2r_*}_{r_*}\rho^\gamma r\dr\le C(T)$. With the help of the H$\ddot{o}$lder inequality and the mean value theorem, there exists $s_t\in[r_*,2r_*]$ for each $t\in[0,T]$ such that $\rho(s_t,t)\le C(T)$. Thus for each $0<r<r_*$
 \begin{align*}
     \rho(t,r)-\rho(t,s_t)=\int^r_{s_t}\rho_r(t,s)\ds\le C(T)\left(1+\sqrt{\log\dfrac{2r_*}{r}}\right),
 \end{align*}
 and this estimate holds for $r\in[r_*,2r_*]$ obviously. Let the cut-off function $\chi$ satisfy $0\le\chi\le1,\chi=1\text{ for }[0,r_*],\chi=0\text{ for }[2r_*,\infty),\chi'\le0$. Now we define the unweighted momentum in $\{t:\underline{r}(t)>0\}$ 
 \begin{align}\label{unweighted momentum}
     \mathscr{M}_\chi(t)=\int^\infty_0\chi\rho u\dr
 \end{align}
 It follows from the H$\ddot{o}$lder inequality that
 \begin{align}\label{Linfty cutoff momentum}
     |\mathscr{M}_\chi|\le C(T)\left(\int^{2r_*}_{\underline{r}(t)}\dfrac{\rho}{r}\dr\right)^{\frac12}\le C(T)\left(1+\left(\log\dfrac{2r_*}{\underline{r}(t)}\right)^{\frac34}\right).
 \end{align}
 By virtue of Lemma \ref{momentum equation}, for each $\psi\in C^\infty_c((0,T)\times(0,\infty))$, it implies that
 \begin{align}\label{unweighted momentum equation}
     0&=\int^T_0\int^\infty_0\left(\rho u\psi_t+\left(\rho u^2+A\rho^\gamma-\beta\left(u_r+\dfrac{mu}{r}\right)\right)\psi_r-\dfrac{\rho u^2+\kappa|\rho_r|^2}{r}\psi\right)\dr\dt\nonumber\\
     &+\kappa\int^T_0\int^\infty_0\rho\rho_r(\psi_{rr}-\dfrac{\psi_r}{r})\dr\dt+\dfrac{3\kappa}{2}\int^T_0\int^\infty_0|\rho_r|^2\psi_r\dr\dt.
 \end{align}
 Now we consider the following cut-off function
 \begin{align*}
     \vartheta(t)\in C^\infty_c(t_1,t_2),\eta_\delta\in C^\infty(0,\infty),\quad\eta_{\delta}=0 \text{ for }r\le\delta,\eta_{\delta}=1 \text{ for }r\ge2\delta,
 \end{align*}
where $0<2\delta<\min_{s\in[t_1,t_2]}\underline{r}(t)$. Taking $\psi(t,r)=\chi(r)\vartheta(t)\eta_{\delta}(r)$ into \eqref{unweighted momentum equation}, it holds that
\begin{align*}
    0&=\int^T_0\int^\infty_0\left(\rho u(\vartheta'\chi\eta_\delta)+\left(\rho u^2+A\rho^\gamma-\beta\left(u_r+\dfrac{mu}{r}\right)\right)(\chi'\vartheta\eta_\delta+\chi\vartheta\eta_\delta')-\dfrac{\rho u^2+\kappa|\rho_r|^2}{r}(\chi\vartheta\eta_\delta)\right)\dr\dt\nonumber\\
     &+\kappa\int^T_0\int^\infty_0\rho\rho_r\left((\chi''\vartheta\eta_\delta+2\chi'\vartheta\eta'_{\delta}+\chi\vartheta\eta_\delta'')+\dfrac{\chi'\vartheta\eta_\delta+\chi\vartheta\eta_\delta'}{r}\right)\dr\dt\\
     &+\dfrac{3\kappa}{2}\int^T_0\int^\infty_0|\rho_r|^2(\chi'\vartheta\eta_\delta+\chi\vartheta\eta_\delta')\dr\dt.
\end{align*}
It is straightforward to check $\supp\eta_\delta'\subset\mathcal{V}$, thus 
\begin{align*}
    0&=\int^{t_2}_{t_1}\mathscr{M}_\chi(s)\vartheta'(s)\ds+\int^{t_2}_{t_1}\vartheta(s)\left(-2\beta c(s)-\int^\infty_0\dfrac{\chi}{r}(\rho u^2+\kappa|\rho_r|^2)\dr\right)\ds\\
    &+\dfrac{3\kappa}{2}\int^{t_2}_{t_1}\vartheta(s)\int^{2r_*}_{r_*}|\rho_r|^2\chi'\dr\ds+\int^{t_2}_{t_1}\vartheta(s)\int^{2r_*}_{r_*}\chi\left(\rho u^2+A\rho^{\gamma}-\beta\left(u_r+\dfrac{mu}{r}\right)\right)\dr\ds\\
    &+\int^{t_2}_{t_1}\vartheta(s)\int^{2r_*}_{r_*}\rho\rho_r\left(\chi''+\dfrac{\chi'}{r}\right)\dr\ds.
\end{align*}
Thus, it holds that
\begin{align}\label{last equality of chi momen}
    (\mathscr{M}_\chi)_t&=-2\beta c(t)-\int^\infty_0\dfrac{\chi}{r}(\rho u^2+\kappa|\rho_r|^2)\dr+\kappa\int^{2r_*}_{r_*}\rho\rho_r\left(\chi''-\dfrac{\chi'}{r}\right)\dr+\dfrac{3\kappa}{2}\int^{2r_*}_{r_*}|\rho_r|^2\chi'\dr\nonumber\\
    &+\int^{2r_*}_{r_*}\chi'\left(\rho u^2+A\rho^{\gamma}-\beta\left(u_r+\dfrac{mu}{r}\right)\right)\dr\text{ in }\mathcal{D}'(t_1,t_2).
\end{align}
Since $c(t)\in L^2(t_1,t_2)$, it holds that $(\mathscr{M}_\chi)_t\in L^1(t_1,t_2)$. Combining \eqref{Linfty cutoff momentum}, $\mathscr{M}_\chi\in W^{1,1}(t_1,t_2)$. Thus, $\mathscr{M}_\chi$ is absolutely continuous and $\displaystyle\mathscr{M}_\chi(t_2)-\mathscr{M}_\chi(t_1)=\int^{t_2}_{t_1}\mathscr{M'_\chi}\dt$. It is straightforward to check 
\begin{align*}
    &\int^{t_2}_{t_1}\int^{2r_*}_{r_*}\chi\left(\rho u^2+A\rho^{\gamma}-\beta\left(u_r+\dfrac{mu}{r}\right)\right)\dr\dt\\
    &+\int^{2r_*}_{r_*}\rho\rho_r\left(\chi''+\dfrac{\chi'}{r}\right)\dr+\dfrac{3\kappa}{2}\int^{2r_*}_{r_*}|\rho_r|^2\chi'\dr\dt\le C(r_*,T)\le C(T).
\end{align*}
Thus we integrate \eqref{last equality of chi momen} to obtain
\begin{align*}
    \displaystyle\mathscr{M}_\chi(t_2)-\mathscr{M}_\chi(t_1)+2\beta\log\dfrac{\underline{r}(t_2)}{\underline{r}(t_1)}\le C(T).
\end{align*}
Since $\{t:\underline{r}(t)>0\}=\cup^\infty_{i=1}(t^i_1,t^i_2)$. For each fixed $i$, we choose $t_0\in(t^i_1,t^i_2)$, thus there exists a subsequence $\{s_j\}:s_j\rightarrow t^i_1$ as $j\rightarrow\infty$ such that
\begin{align}\label{contradiction}
    2\beta\log\dfrac{\underline{r}(t_0)}{\underline{r}(s_j)}\le C(T)+|\mathscr{M}_\chi(t_0)|+|\mathscr{M}_\chi(s_j)|\le C(T)\left(1+\left(\log\dfrac{\underline{r}(t_0)}{\underline{r}(s_j)}\right)^{\frac34}\right).
\end{align}
Now we let $j\rightarrow\infty$, \eqref{contradiction} force $2\beta\le0$ and this is contradict to $\beta>0$. Thus, $\underline{r}(t)\equiv0~a.e. $  in $(t^i_1,t^i_2)$, and therefore, $\underline{r}(t)\equiv0$  for  $a.e. ~t\in[0,T]$. By virtue of the continuity of $\underline{r}(t)$, the desired conclusion is obtained.
\par The proof of Lemma \ref{Ruling out the possible vacuum} is finished.
    \end{proof}
\end{lemma}

Now, we are ready to characterize the pointwise bound of $\rho$ for each positive radius $r>0$. The upper bound of $\rho$ has already been obtained in the proof of Lemma \ref{Ruling out the possible vacuum}, thus the proof of the following lemma focuses on the lower bound.
\begin{lemma}\label{characterize point wise bound of rho}
    For $\forall r>0$, there exist two  positive constants $C(T)$ and $C(p,T)$ such that 
    \begin{align*}
        &\min\{C^{-1}(p,T)r^{\frac{2p}{1-p}},C^{-1}(T)\}\le\rho(t,r) \text{ for each }r>0,t\in[0,T],p\in\left(\frac{2}{3},1\right),\\
        &\rho(t,r)\le C(T)\text{ for each }r\ge1,\quad\rho(t,r)\le C(T)\left(1+\log\dfrac{1}{r}\right)\text{ for each } 0<r<1.
    \end{align*}
    \begin{proof}
        \textbf{1.} Set $\widetilde{H}(t,r)=\displaystyle\int^{\rho(t,r)}_1\sqrt{G(t,\zeta)}\mathrm{\zeta}$. It is obvious that $s\mapsto\displaystyle\int^{s}_1\sqrt{G(t,s)}\ds$ is monotone increasing in $s$. Since 
        \begin{align*}
            \left|\int^{\rho(t,r)}_1\sqrt{G(t,s)}\ds\right|&\le\left|\int^\infty_r\sqrt{G(\rho(t,z))}\rho_r(t,z)\dz\right|\\
            &\le\dfrac{C(T)}{r}.
        \end{align*}
        Then we let $h_*=\min\left\{\displaystyle\int^{1}_{\frac12}\sqrt{G(t,s)}\ds,\int^2_1\sqrt{G(t,s)}\ds\right\}$, and we choose a sufficiently large $R_0$ such that  $\dfrac{C(T)}{r}<h_*$ for each $r\ge R_0$.
        Thus for each $r\ge R_0$, we have $\dfrac{1}{2}\le\rho(t,r)\le2$. With the help of the Newton--Leibniz formula, it holds that
        \begin{align}\label{log upp bound of rho}
            \rho(t,r)\le2+\left|\int^r_{R_0}\rho_r(t,s)\ds\right|\le2 +C(T)\left|\log \dfrac{R_0}{r}\right|^{\frac12}.
        \end{align}
        Since $|u(t,r)|^2\le\displaystyle\int^\infty_0\left|u_r+\dfrac{u}{r}\right|^2r\dr$, it holds that
        \begin{align*}
            |r\rho(t,r)u(t,r)|\le r\left(2+ C(T)\left|\log \dfrac{R_0}{r}\right|^{\frac12}\right)\left(\int^\infty_0\left|u_r+\dfrac{u}{r}\right|^2r\dr\right)^{\frac12}\rightarrow0\text{ as }r\rightarrow0.
        \end{align*}
        Thus by the dominated convergence theorem, $r\rho(\cdot,r)u(\cdot,r)\rightarrow0$ in $L^1(0,T)$. By virtue of Lemma \ref{weak form of continuity equation}, it follows that
        \begin{align*}
            (\rho r)_t+(\rho ur)_r=0\text{ in }\mathcal{D}'((0,T)\times(0,\infty)).
        \end{align*}
         Let $\chi_\varepsilon$ satisfy 
        \begin{align*}
            0\le\chi_\varepsilon\le1,\quad\chi_\varepsilon=0\text{ for }r\le\varepsilon,\quad\chi_\varepsilon=1\text{ for }r\ge2\varepsilon,\quad|\chi'_\varepsilon|\le\dfrac{C}{\varepsilon}.
        \end{align*}
Now we define the mass function $x(t,r)=\displaystyle\int^r_0\rho(t,s)s\ds$, it holds that
\begin{align*}
    0&=\int^T_0\int^\infty_0\rho(t,s)s\chi_\varepsilon(s)\eta_t(t,s)\ds\dt-\int^T_0\int^\infty_0u(t,s)\rho(t,s)s\chi_\varepsilon(s)\psi(t,s)\ds\dt\\
    &+\int^T_0\int^\infty_0u(t,s)\rho(t,s)s\chi'_\varepsilon(s)\eta(t,s)\ds\dt,
\end{align*}
where $\eta(t,r)=\displaystyle\int^\infty_r\psi(t,s)\ds$ and $\psi\in C^\infty_c((0,T)\times(0,\infty))$. For the third term in right hand side, we have
\begin{align*}
    &\left|\int^T_0\int^\infty_0\rho(t,s)s\chi'_\varepsilon(s)\eta(t,s)u(t,s)\ds\dt\right|\\
    &\le\dfrac{C\Vert\eta\Vert\sqrt{T}}{\varepsilon}\left(\int^T_0\int^\infty_0\left|u_r+\dfrac{u}{r}\right|^2r\dr\ds\right)^{\frac{1}{2}}\int^{2\varepsilon}_{\varepsilon}\left(2+C(T)\left|\log\dfrac{R_0}{r}\right|^{\frac12}\right)r\dr\\
    &\le C(T)\varepsilon\left(1+\left|\log\dfrac{R_0}{\varepsilon}\right|^{\frac12}\right)\rightarrow0\text{ as }\varepsilon\rightarrow0.
\end{align*}
Since $\eta$ has compact support, it follows that 
\begin{align*}
    \int^T_0\int^\infty_0\rho(t,s)s\chi_\varepsilon(s)\eta_t(t,s)\ds\dt\overset{\varepsilon\rightarrow0}{\longrightarrow}\displaystyle\int^T_0\int^\infty_0\rho(t,s)s\eta_t(t,s)\ds\dt.
\end{align*}
With the same argument,
\begin{align*}
    \int^T_0\int^\infty_0u(t,s)\rho(t,s)s\chi_\varepsilon(s)\psi(t,s)\ds\dt\rightarrow\displaystyle\int^T_0\int^\infty_0u(t,s)\rho(t,s)s\psi(t,s)\ds\dt.
\end{align*}
Thus,
\begin{align*}
    \int^T_0\int^\infty_0\rho(t,s)s\eta_t(t,s)\ds\dt=\int^T_0\int^\infty_0u(t,s)\rho(t,s)s\psi(t,s)\ds\dt.
\end{align*}
By the Fubini's theorem, $\displaystyle\int^T_0\int^\infty_0\rho(t,s)s\eta_t(t,s)\ds\dt=\int^T_0\int^\infty_0x(t,r)\psi_t(t,r)\dr\dt$, thus, $x_t(t,r)=-r\rho u$ for $a.e. ~(t,r)\in(0,T)\times(0,\infty)$. It is straightforward to obtain $x_r(t,r)=\rho r$. Since $x$ is a strictly increasing function, by virtue of the inverse function theorem, there exists a function $\widetilde{r}(t,x)$ such that $\widetilde{r}(t,x)$ is strictly increasing in $x$. Moreover, it follows from \eqref{log upp bound of rho} that $x(t,r)\le C(T)R_0 r \text{ for }r\le R_0$. we have $x(t,2R_0)\le C(T)R^2_0$, thus, we can choose a sufficiently large $X_0> C(T)R^2_0$ such that $x(t,2R_0)<X_0$ and the corresponding particle path $\widetilde{r}(t,X_0)>2R_0$. In addition, $x(t,r)\le\dfrac{X_0}{R_0}r$ for $0\le r\le R_0$. Since $\dfrac{1}{2}\le\rho(t,r)\le 2$ for $r\ge R_0$, we have
\begin{align*}
    2R_0\le \widetilde{r}(t,X_0)\le (R_0^2+4X_0)^{\frac12}.
\end{align*}
Now we consider the particle path $\widetilde{r}(t,x)$ in $x\in[0,X_0]$. If $\widetilde{r}(t,x)<R_0$, we have $x=x(\widetilde{r}(t,x),t)\le\dfrac{X_0}{R_0}\widetilde{r}(t,x)$, thus $\widetilde{r}(t,x)\ge\dfrac{R_0x}{X_0}$. If $R_0\le\widetilde{r}(t,x)$, $\dfrac{R_0x}{X_0}\le R_0\le \widetilde{r}(t,x)$. Combing above, it holds that
\begin{align}\label{ inequality of x}
    \dfrac{R_0x}{X_0}\le\widetilde{r}(t,x)\text{ for all }x\in[0,X_0].
\end{align}
\textbf{2. }Now we define the corresponding unweighted momentum $\mathscr{M}^{\chi _{R_0}}(t,r):=\displaystyle\int^\infty_r\chi(s)\rho(t,s)u(t,s)\ds$, where $0\le\chi_{R_0}\le1$, $\chi_{R_0}=1$ for $r\in[0,R_0]$ and $\chi_{R_0}=0$ for $r\in[2R_0,\infty)$. It follows that
\begin{align}\label{unweighted momentum of chi R0}
    \mathscr{M}_t^{\chi _{R_0}}+u\mathscr{M}_r^{\chi _{R_0}}&=\chi_{R_0} \left(A\rho^\gamma-\beta\left(u_r+\dfrac{u}{r}\right)-\kappa\left(\rho\left(\rho_{rr}+\dfrac{\rho_r}{r}\right)-\dfrac{|\rho_r|^2}{2}\right)\right)\nonumber\\
    &+\int^\infty_r\chi'_{R_0}(s)\left(\rho u^2+A\rho^\gamma-\beta\left(u_r+\dfrac{u}{r}\right)-\kappa\left(\rho\left(\rho_{rr}+\dfrac{\rho_r}{r}\right)-\dfrac{|\rho_r|^2}{2}\right)\right)\ds\nonumber\\
    &-\int^\infty_r\chi_{R_0}\dfrac{\rho u^2+\kappa|\rho_r|^2}{s}\ds.
\end{align}
By virtue of Lemma \ref{last entropy estimate}, it follows that $|\mathscr{M}^{\chi _{R_0}}(t,r)|\le C(T)\left(1+\left|\log\dfrac{R_0}{r}\right|\right)^{\frac34}$. Thus it follows from \eqref{ inequality of x} that
\begin{align*}
|\mathscr{M}^{\chi _{R_0}}(t,\widetilde{r}(t,x))|\le C(T)\left(1+\log\dfrac{X_0}{x}\right)^{\frac34}\text{ for }0<x\le X_0.
\end{align*}
Since $\supp\chi'_{R_0}\subset\left[R_0,2R_0\right]$,
\begin{align*}
    \left|-\kappa\int^{2R_0}_{R_0}\chi'_{R_0}\left(\rho\left(\rho_{rr}+\dfrac{\rho_r}{r}\right)-\dfrac{|\rho_r|^2}{2}\right)\dr\right|=\left|\dfrac{3\kappa}{2}\int^{2R_0}_{R_0}\chi'_{R_0}|\rho_r|^2\dr+\kappa\int^{2R_0}_{R_0}\left(\chi''_{R_0}-\dfrac{\chi'_{R_0}}{s}\right)\rho\rho_r\ds\right|\le C(T).
\end{align*}
Thus,
\begin{align*}
    &\int^T_0\left|\int^{2R_0}_{R_0}\chi'_{R_0}(s)\left(\rho u^2+A\rho^\gamma-\beta\left(u_r+\dfrac{u}{r}\right)-\kappa\left(\rho\left(\rho_{rr}+\dfrac{\rho_r}{r}\right)-\dfrac{|\rho_r|^2}{2}\right)\right)\ds\right|\dt\\
    &\le C(T)\left(1+\left(\int^T_0\int^{2R_0}_{R_0}\left|u_r+\dfrac{u}{r}\right|^2r\dr\ds\right)^{\frac12}\right)\le C(T).
\end{align*}
If we let $I=[a,b]\Subset(0,\infty)$, there exists a positive constant $C_I$ such that $0<(C_I)^{-1}\le\rho\le C_I$. With the help of Lemma \ref{last estimate epsilon}, we have $\rho_r\in L^2(0,T;L^\infty(I;\dr)),u_r\in L^2(0,T;L^2(I;\dr))$. Thus it follows from Lemma \ref{weak form of continuity equation} that $\rho_t\in L^2(0,T;L^2(I;\dr))$, in addition, $\mathscr{M}^{\chi_{R_0}}_t\in L^2(0,T;L^2(I;\dr))$. Furthermore, since $\mathscr{M}^{\chi_{R_0}}_{rr}=-\chi'\rho u -\chi\rho_ru-\chi\rho u_r\in L^2(0,T;L^2(I;\dr)) $.

\textbf{3. } We set $0<\rho_*<\min\left\{\dfrac{1}{4},\left(\dfrac{\beta^2}{8\kappa}\right)^{\frac13}\right\}$ and  define 
$\widetilde{S}'(\rho)=\dfrac{\beta}{2\rho}\left(1+\sqrt{1-\dfrac{4\kappa \rho^3}{\beta^2}}\right)$ for $\rho<\rho_*$. Then we extend $\widetilde{S}'(\rho)$ by
\begin{align*}
    S'(\rho)=\begin{cases}
        \dfrac{\beta+\eta_{S'}\sqrt{\beta^2-4\kappa\rho^3}}{2\rho}\quad&\text{ for } 0<\rho<\rho_1,\\
        \dfrac{\beta}{2\rho}\quad&\text{ for }\rho\ge\rho_1,
    \end{cases}
\end{align*}
where $\rho_1:=\left(\dfrac{3}{2}\right)^{\frac13}\rho_*$. It is obvious that $\dfrac{\beta}{2}\le \rho S'(\rho)\le\beta$. 
Then we set 
\begin{align*}
    S_1(\rho)=\int^{\rho_*}_\rho S'(\sigma)\mathrm{d}\sigma, 
\end{align*}
where the cut-off function $\eta_{S'}$ satisfies
\begin{align*}
    0\le \eta_{S'}\le 1,\quad\eta_{S'}=1 \text{ for }\rho\le\rho_*,\quad\eta_{S'}=0\text{ for }\rho\ge\left(\dfrac{3}{2}\right)^{\frac13}\rho_*.
\end{align*}
Moreover, it holds that $\dfrac{1}{2}-\dfrac{\rho S''(\rho)}{S'(\rho)}=\dfrac{3\beta}{2\sqrt{\beta^2-4\kappa\rho^3}}$ for $0<\rho\le\rho_*$. Since $\beta-S'(\rho)\rho=\dfrac{\kappa\rho^2}{S'(\rho)}$ for $\rho<\rho_*$, thus for each $\rho<\rho_*$ it holds that
\begin{align*}
    \beta\log\dfrac{\rho_*}{\rho}-S_1(\rho)=\int^{\rho_*}_\rho\dfrac{\beta-S'(s)s}{s}\mathrm{d}s\le\dfrac{2\kappa}{\beta}\int^{\rho_*}_0\rho^2\mathrm{d}\rho=\dfrac{2\kappa\rho^3_*}{3\beta}.
\end{align*}
Thus $\beta\log\dfrac{\rho_*}{\rho}-\dfrac{2\kappa\rho^3_*}{3\mu}\le S_1(\rho)\le\beta\log\dfrac{\rho_*}{\rho}$ for all $\rho\le\rho_*$. By the construction of $S_1$, it follows that
\begin{align*}
    \left((S_1)_t(\rho),(S_1(\rho))_{rr}\right)\in L^2(0,T;L^2(I;\dr)).
\end{align*}
\textbf{4.} Now we define $z(t,r)=S_1(\rho(t,r))+\mathscr{M}^{\chi_{R_0}}(t,r)$,  it is straightforward to obtain $z\in L^2(0,T;H^2(I;\dr))\cap L^\infty(0,T;H^1(I;\dr))$ and $z_t\in L^2(0,T;L^2(I;\dr))$. Since $\rho(t,r)\ge \dfrac{1}{2}$ for $r\ge R_0$,  $\{\rho<\rho_*\}\subset\{r<R_0\}$. Furthermore, it can be proved that the following equation holds for $\{\rho<\rho_*\}$
\begin{align}\label{equation of z}
    &z_t+\left(u+\dfrac{\kappa\rho_r}{S'}\right)z_r-\kappa\dfrac{\rho}{S'}\left(z_{rr}+\dfrac{z_r}{r}\right)\\
    &=A\rho^\gamma-\int^\infty_r\chi_{R_0}\dfrac{\rho u^2+\kappa|\rho_r|^2}{s}\ds-\dfrac{3\kappa\beta|\rho_r|^2}{2\sqrt{\beta^2-4\kappa\rho^3}}\nonumber\\
    &+\int^{2R_0}_{R_0}\chi'_{R_0}(s)\left(\rho u^2+A\rho^\gamma-\beta\left(u_r+\dfrac{u}{r}\right)-\kappa\left(\rho\left(\rho_{rr}+\dfrac{\rho_r}{r}\right)-\dfrac{|\rho_r|^2}{2}\right)\right)\ds.
\end{align}
By a standard argument we have $\partial_x\widetilde{r}(t,x)=\dfrac{1}{\rho(\widetilde{r}(t,x))\widetilde{r}(t,x)}$ and $\partial_t\widetilde{r}(t,x)=u(t,\widetilde{r}(t,x))$, thus we have $\partial_t\widetilde{r}(t,x)\in{L^2((0,T);L^\infty_{loc}((0,\infty);\dx))}$. Now set $\widetilde{z}(t,x)=z(t,\widetilde{r}(t,x))$, we have
\begin{align}\label{regularity of z}
    \widetilde{z}\in L^\infty(0,T;H^1_{loc}(0,\infty);\dx)\cap L^2(0,T;H^2_{loc}(0,\infty);\dx)\cap H^1(0,T;L^2_{loc}(0,\infty);\dx).
\end{align}
\textbf{5. } In this step we consider \eqref{equation of z} in Lagrangian coordinates, and we omit the superscript for simplicity in the step 5, step 6, and step 7. In addition, we temporarily fix the positive constant $C(T)$. Derivatives in Lagrangian coordinates are denoted by $\partial$. First, we reformulate the \eqref{equation of z} in Lagrangian coordinates 
\begin{align}\label{equation of z in lagrange}
    \partial_tz-\dfrac{\kappa\rho^3r^2}{S'}\partial_{xx}z-\dfrac{2\kappa\rho^2}{S'}\partial_xz&=A\rho^\gamma-\int^\infty_r\chi_{R_0}\dfrac{\rho u^2+\kappa|\rho_r|^2}{s}\ds-\dfrac{3\kappa\beta|\partial_x\rho\rho r|^2}{2\sqrt{\beta^2-4\kappa\rho^3}}\nonumber\\
    &+\int^{2R_0}_{R_0}\chi'_{R_0}(s)\left(\rho u^2+A\rho^\gamma-\beta\left(u_r+\dfrac{u}{r}\right)-\kappa\left(\rho\left(\rho_{rr}+\dfrac{\rho_r}{r}\right)-\dfrac{|\rho_r|^2}{2}\right)\right)\ds.
\end{align}
Now we define $\psi_\delta(x)=-\beta \log\left(\left(\dfrac{x}{X_0}\right)^{\theta}-\delta\right)$ in $x\in(x_\delta,X_0]$, where $\dfrac{2}{3}<\theta<1$, $0<\delta<1$, and $x_\delta=X_0\delta^{\frac1\theta}$. It also holds that $\psi_\delta(x)\rightarrow\infty$ as $x\rightarrow x_\delta$ and $\theta\beta\log\dfrac{X_0}{x}\le \psi_\delta(x)$. By straightforward calculations, it holds that
\begin{align*}
    &\psi'_\delta(x)=-\dfrac{\beta\theta}{x\left(1-\delta\left(\dfrac{X_0}{x}\right)^\theta\right)}<0,\quad \psi''_\delta(x)=\dfrac{\beta\theta\left(1-(1-\theta)\delta\left(\dfrac{X_0}{x}\right)^\theta\right)}{x^2\left(1-\delta\left(\dfrac{X_0}{x}\right)^\theta\right)^2}>0,\\
    &\dfrac{\rho S'(\rho)\psi''_\delta}{(\psi'_\delta)^2}\le\dfrac{1}{\theta}<\dfrac{3\beta}{2\sqrt{\beta^2-4\kappa\rho^3}}\text{ for all }\rho<\rho_*.
\end{align*}
Set  $H_\theta=\sup_{\log\frac{X_0}{x}\ge0}\left(C(T)\left(1+\log\dfrac{X_0}{x}\right)^{\frac34}-\theta\beta \log\dfrac{X_0}{x}\right)$, $S_0=\beta\log_+\dfrac{\rho_*}{C^{-1}_0}$, and $M_0=S_0+H_\theta+\beta$, where $C^{-1}_0$ is the lower bound of initial density $\rho_0$. There exists a positive constant $\theta\beta+\dfrac{27 C^4(T)}{256(\theta\beta)^3}$ such that $H_\theta\le\theta\beta+\dfrac{27 C^4(T)}{256(\theta\beta)^3}$. Let
\begin{align*}
    b(t)&=M_0+\displaystyle\int^t_0\left[A\rho^\ga_*+\left|\int^{2R_0}_{R_0}\chi'_{R_0}(s)\left(\rho u^2+A\rho^\gamma-\beta\left(u_r+\dfrac{u}{s}\right)-\kappa\left(\rho\left(\rho_{rr}+\dfrac{\rho_r}{s}\right)-\dfrac{|\rho_r|^2}{2}\right)\right)\ds\right|\right]\mathrm{d}\tau\\
    &+\dfrac{12\kappa\rho^4_*}{\beta^2(3\theta-2)}\int^t_0\int^\infty_0\left(u^2_r+\dfrac{u^2}{r^2}\right)r\dr\mathrm{d}\tau.
\end{align*}
It is straight to check $b(t)\in W^{1,1}(0,T)$ and $b$ is monotone increasing in $t$. Furthermore, it holds that
\begin{align*}
    &A\rho^\gamma+\int^{2R_0}_{R_0}\chi'_{R_0}(s)\left(\rho u^2+A\rho^\gamma-\beta\left(u_r+\dfrac{u}{r}\right)-\kappa\left(\rho\left(\rho_{rr}+\dfrac{\rho_r}{r}\right)-\dfrac{|\rho_r|^2}{2}\right)\right)\ds\\
    &-\int^\infty_r\chi_{R_0}\dfrac{\rho u^2+\kappa|\rho_r|^2}{s}\ds-b'+\dfrac{12\kappa\rho^4_*}{\beta^2(3\theta-2)}|u|^2\le0\text{ for all }\rho<\rho_*.
\end{align*}\\
\textbf{6. }  Now we define the auxiliary function $W_\delta(t,x)=z-\psi_\delta-b(t)$. Since $z_r=-S'\rho_r-\rho u$ and $z_r=z_x\rho r$, it holds that
\begin{align*}
    \rho_r=\dfrac{-\rho}{S'}\left(r\psi'_\delta(x)+u+r(W_\delta)_x\right).
\end{align*}
Next, by virtue of \eqref{equation of z in lagrange}, it holds that
\begin{align*}
    &\partial_tW_\delta-\dfrac{\kappa\rho^3 r^2}{S'}\partial_{xx}W_\delta-2\kappa\dfrac{\rho^2}{S'}\partial_xW_\delta\\
    &=A\rho^\gamma+\int^{2R_0}_{R_0}\chi'_{R_0}(s)\left(\rho u^2+A\rho^\gamma-\beta\left(u_r+\dfrac{u}{r}\right)-\kappa\left(\rho\left(\rho_{rr}+\dfrac{\rho_r}{r}\right)-\dfrac{|\rho_r|^2}{2}\right)\right)\ds\\
    &-\int^\infty_r\chi_{R_0}\dfrac{\rho u^2+\kappa|\rho_r|^2}{s}\ds-b'-\dfrac{3\beta\kappa|\rho_x\rho r|^2}{2\sqrt{\beta^2-4\kappa\rho^3}}+\dfrac{\kappa\rho^3 r^2}{S'}\psi''_\delta+2\kappa\dfrac{\rho^2}{S'}\psi'_\delta.
\end{align*}
Since 
\begin{align*}
    &-\dfrac{3\beta\kappa|\rho_x\rho r|^2}{2\sqrt{\beta^2-4\kappa\rho^3}}+\dfrac{\kappa\rho^3 r^2}{S'}\psi''_\delta+2\kappa\dfrac{\rho^2}{S'}\psi'_\delta\\
    &=\dfrac{\kappa\rho^2}{(S')^2}\left(\dfrac{\rho S'\psi''_\delta}{(\psi'_\delta)^2}r^2|\psi'_\delta|^2-(r\psi'_\delta+u)^2\cdot\dfrac{3\beta}{2\sqrt{\beta^2-4\kappa\rho^3}}\right)-\dfrac{3\beta}{2\sqrt{\beta^2-4\kappa\rho^3}}\cdot\dfrac{2\kappa\rho^2r}{(S')^2}\partial_xW_\delta(r\psi'_\delta+u)\\
    &-\dfrac{\kappa\rho^2r^2|\partial_xW_\delta|^2}{(S')^2}\dfrac{3\beta}{2\sqrt{\beta^2-4\kappa\rho^3}}+\dfrac{2\kappa\rho^2}{S'}\psi'_\delta.
\end{align*}
It is obvious that the last two terms are negative, thus, 
\begin{align}\label{inequality of W}
    &\partial_tW_\delta-\dfrac{\kappa\rho^3r^2}{S'}\partial_{xx}W_\delta-2\kappa\dfrac{\rho^2}{S'}\partial_xW_\delta\nonumber\\
    &\le A\rho^\gamma+\int^{2R_0}_{R_0}\chi'_{R_0}(s)\left(\rho u^2+A\rho^\gamma-\beta\left(u_r+\dfrac{u}{r}\right)-\kappa\left(\rho\left(\rho_{rr}+\dfrac{\rho_r}{r}\right)-\dfrac{|\rho_r|^2}{2}\right)\right)\ds-b'\nonumber\\
    &+\dfrac{\kappa\rho^2}{(S')^2}\left(\dfrac{\rho S'\psi''_\delta}{(\psi'_\delta)^2}r^2|\psi'_\delta|^2-(r\psi'_\delta+u)^2\cdot\dfrac{3\beta}{2\sqrt{\beta^2-4\kappa\rho^3}}\right)\\
    &-\dfrac{3\beta}{2\sqrt{\beta^2-4\kappa\rho^3}}\cdot\dfrac{2\kappa\rho^2r}{(S')^2}\partial_xW_\delta(r\psi'_\delta+u).
\end{align}
Since $\theta\beta\log\dfrac{X_0}{x}\le\psi_\delta(x)$, it holds that $ W_\delta\le S_1(\rho)-S_0-\beta$. Thus $W_\delta(0,x)\le-\beta<0$ and  $W_\delta\le-\beta<0$ for $\{\rho\ge\rho_*\}$. In addition, if $W_\delta>0$, it must satisfy $S_1>0$. By the construction of $S_1$, it holds that $\rho<\rho_*$. Specifically, it holds that
\begin{align}\label{W subset}
   \{W_\delta>0\}\subset\left\{\rho<\rho_*\exp\left(-1-\dfrac{S_0}{\beta}\right)\right\}\subset\{\rho<\dfrac{\rho_*}{e}\}.
\end{align}
In addition, since $r(X_0,t)>2R_0$, $\mathscr{M}^{\chi_{R_0}}(X_0,t)=0$, $0<\psi_\delta(X_0,t)$, and $\dfrac{1}{2}\le\rho(X_0,t)$, we have $W_\delta(t,X_0)<0$. For $\forall \delta>0$, there exists a positive constant $C_\delta(T)$ such that
\begin{align*}
    z(t,x)\le C_\delta(T)\text{ for each }\dfrac{x_\delta}{2}\le x\le X_0,0\le t\le T,
\end{align*}
where we used $\rho(t,x)$ has   positive uniform upper and lower bound in $x\in\left[\dfrac{x_\delta}{2},X_0\right]$. 
\par We set $h_\delta=\min\left\{\dfrac{1-\delta}{2},\dfrac{e^{-\frac{C_{\delta}(T)}{\beta}-1}}{2}\right\}$ and $l_\delta=X_0(\delta+h_\delta)^{\frac1\theta}$, it follows that
\begin{align*}
    W_\delta(t,x)\le0 \text{ for }x_\delta<x\le l_\delta,0\le t\le T,
\end{align*}
where we used the lower bound $C_\delta(T)+\beta\le\psi_\delta(l_\delta)$. Since $\rho S'(\rho)\ge\dfrac{\beta}{2}$ for $\rho<\rho_*$, it holds that
\begin{align*}
    \dfrac{\kappa\rho^2}{(S')^2}\left(\dfrac{\rho S'\psi''_\delta}{(\psi'_\delta)^2}r^2|\psi'_\delta|^2-(r\psi'_\delta+u)^2\cdot\dfrac{3\beta}{2\sqrt{\beta^2-4\kappa\rho^3}}\right)\le \dfrac{12\kappa\rho^4_*}{\beta^2(3\theta-2)}|u|^2\text{ for all }\rho\le\rho_*,
\end{align*}
thus, by virtue of \eqref{inequality of W}, it holds that
\begin{align*}
    \partial_tW_\delta-\dfrac{\kappa\rho^3r^2}{S'}\partial_{xx}W_\delta-\left(\dfrac{2\kappa\rho^2}{S'}-\dfrac{2\kappa\rho^2r}{(S')^2}(r\psi'_\delta+u)\dfrac{3\beta}{2\sqrt{\beta^2-4\kappa\rho^3}}\right)\partial_xW_\delta\le0.
\end{align*}
Based on the fact: $\dfrac{2\kappa\rho^2}{S'}=\partial_x\left(\kappa\dfrac{\rho^3r^2}{S'}\right)-\dfrac{\kappa r\rho}{S'}\left(3-\dfrac{\rho S''}{S'}\right)\rho_x\rho r$, the following inequality holds in $\{\rho<\rho_*\}$ that
\begin{align}\label{equation of div W}
    \partial_tW_\delta-\partial_x\left(\dfrac{\kappa\rho^3 r^2}{S'}\partial_xW_\delta\right)+\left(\dfrac{\kappa r\rho}{S'}\left(3-\dfrac{\rho S''}{S'}\right)\partial_x\rho\rho r+\dfrac{2\kappa\rho^2r}{(S')^2}(r\psi'_\delta+u)\left(\dfrac{1}{2}-\dfrac{\rho S''(\rho)}{S'(\rho)}\right)\right)\partial_xW_\delta\le0
\end{align}
Thus \eqref{equation of div W} holds in $\{W_\delta>0\}$ by virtue of \eqref{W subset}.\\
\textbf{7. } For $\forall \delta>0$, there exist two positive constants $C^{-1}_{\delta,l_\delta}(T)$ and $C_{\delta,l_\delta}(T)$ such that
\begin{align*}
    C^{-1}_{\delta,l_\delta}(T)\le\rho(t,x)\le C_{\delta,l_\delta}(T),\quad C^{-1}_{\delta,l_\delta}(T)\le r(t,x)\le C_{\delta,l_\delta}(T),\text{ for }x\in(l_\delta,X_0)\times(0,T).
\end{align*}
Thus $0<\dfrac{\kappa (C^{-1}_{\delta,l_\delta}(T))^{6}}{\beta}\le\dfrac{\kappa\rho^3r^2}{S'}\le\dfrac{2\kappa(C_{\delta,l_\delta}(T))^{6}}{\beta}$. Since $\psi'_\delta$ is bounded on $[l_\delta,X_0]$, it can be proved that 
\begin{align*}
    \left(\dfrac{\kappa\rho^3r^2}{S'},\dfrac{\kappa r\rho}{S'}\left(3-\dfrac{\rho S''}{S'}\right)\partial_x\rho\rho r+\dfrac{2\kappa\rho^2r}{(S')^2}(r\psi'_\delta+u)\left(\dfrac{1}{2}-\dfrac{\rho S''(\rho)}{S'(\rho)}\right)\right)\in L^2(0,T;L^\infty((l_\delta,X_0);\dx)).
\end{align*}
By virtue of \eqref{regularity of z}, we have
\begin{align*}
    W_\delta\in L^\infty(0,T;H^1(l_\delta,X_0))\cap L^2(0,T;H^2(l_\delta,X_0))\cap W^{1,1}(0,T;L^2(l_\delta,X_0)).
\end{align*}
Now we let 
\begin{align*}
    \mathscr{W}_\delta=\partial_tW_\delta-\partial_x\left(\dfrac{\kappa\rho^3 r^2}{S'}\partial_xW_\delta\right)+\left(\dfrac{\kappa r\rho}{S'}\left(3-\dfrac{\rho S''}{S'}\right)\partial_x\rho\rho r+\dfrac{2\kappa\rho^2r}{(S')^2}(r\psi'_\delta+u)\left(\dfrac{1}{2}-\dfrac{\rho S''(\rho)}{S'(\rho)}\right)\right)\partial_xW_\delta.
\end{align*}
By virtue of  $\partial_x\left(\dfrac{\kappa\rho^3 r^2}{S'}\partial_xW_\delta\right)\in L^2(0,T;L^2(l_\delta,X_0))$, it follows that $\mathscr{W}_\delta\in L^1(0,T;L^2(l_\delta,X_0))$. Furthermore, we have $\mathscr{W}_\delta\le0~a.e. $  in $\{W_\delta>0\}$, thus, $\displaystyle\int^{X_0}_{l_\delta}\mathscr{W}_\delta W^+_\delta\dx\le0$. Since $\partial_xW_\delta=0~a.e. $  in $\{W_\delta=0\}$ and $W_\delta(t,X_0)<0,W_\delta(t,l_\delta)\le 0$ for each $t\in[0,T]$, we have $W_\delta^+\in L^2(0,T;H^1_0(l_\delta,X_0))\cap L^\infty(0,T;L^2(l_\delta,X_0))$. Since $W^{1,1}\hookrightarrow C$, by a standard argument, we have
\begin{align*}
    \Vert W^+_\delta(t)\Vert^2_{L^2(l_\delta,X_0)}-\Vert W^+_\delta(s)\Vert^2_{L^2(l_\delta,X_0)}=2\int^t_s\int^{X_0}_{l_\delta}W^+_\delta\partial_tW_\delta\dx\mathrm{d}\tau\quad \text{ for each }0\le s\le t\le T.
\end{align*}
Thus we have
\begin{align*}
    &\dfrac{1}{2}\dfrac{\mathrm{d}}{\dt}\Vert W^+_\delta(t)\Vert^2_{L^2(l_\delta,X_0)}+\dfrac{\kappa (C^{-1}_{\delta,l_\delta}(T))^{6}}{\beta}\Vert \partial_xW^+_\delta(t)\Vert^2_{L^2(l_\delta,X_0)}\\
    &\le C\Vert W^+_\delta(t)\Vert^2_{L^2(l_\delta,X_0)}\left\Vert\dfrac{\kappa r\rho}{S'}\left(3-\dfrac{\rho S''}{S'}\right)\partial_x\rho\rho r+\dfrac{2\kappa\rho^2r}{(S')^2}(r\psi'_\delta+u)\left(\dfrac{1}{2}-\dfrac{\rho S''(\rho)}{S'(\rho)}\right)\right)\Vert^2_{L^\infty(l_\delta,X_0)}\\
    &+\dfrac{\kappa (C^{-1}_{\delta,l_\delta}(T))^{6}}{2\beta}\Vert \partial_xW^+_\delta(t)\Vert_{L^2(l_\delta,X_0)}^2.
\end{align*}
Since $W_\delta(0,x)\le0$, we have $W^+_\delta(0,x)=0$. Thus with the help of the Gronwall's inequality, it holds that $\Vert W^+_\delta(t,\cdot)\Vert_{L^2(l_\delta,X_0)}=0$, therefore, $W_\delta\le0~a.e. $  in $(l_\delta,X_0)\times(0,T)$. By virtue of $z\in H^1(0,T;L^2(I))\cap L^2(0,T;H^2(I))$, we have $z\in C([0,T];H^1(I))$, thus, $W_\delta\in C([0,T]; H^1(I))$. Since $W_\delta(t,x)\le0$ for $x_\delta<x\le l_\delta,t\in[0,T]$, we have $S_1(\rho(t,r(t,x)))+\mathscr{M}^{\chi_{R_0}}(t,r(t,x))\le \psi_\delta(x)+b(t)$ for each $(x_\delta,X_0]\times[0,T]$. For $\forall x>0$, there exists a subsequence $\{\delta_n\}\rightarrow0$ such that $x>x_{\delta_n}$
 and $\psi_{\delta_n}\rightarrow\beta\theta\log\dfrac{X_0}{x}$. Let $n\rightarrow\infty$ we have
 \begin{align*}
     S_1(\rho)+\mathscr{M}^{\chi_{R_0}}\le\beta\theta\log\dfrac{X_0}{x}+ C(T) \text{ for }0<x\le X_0, 0\le t\le T.
 \end{align*}
 Thus for the region $\{\rho<\rho_*\}\cap(x,X_0]$, it holds that
 \begin{align*}
     \beta\log\dfrac{\rho_*}{\rho}\le\theta\beta\log\dfrac{X_0}{x}+C(T)\left(1+\left(\log\dfrac{X_0}{x}\right)^{\frac34}\right)\le C(T)+\beta p\log\dfrac{X_0}{x},
 \end{align*}
 where $\theta<p<1$. Thus we have 
 \begin{align}\label{rho inequality}
     \exp\{-C(p,T)\}X^{-p}_0x^p\rho_*\le\rho(t,r(t,x)).
 \end{align}
 Then for each $\rho>\rho_*$, it also holds that $\exp\{-C(p,T)\}X^{-p}_0x^p\rho_*\le\rho(t,r(t,x))$. \\
 \textbf{8. } Based on our Lemma \ref{Ruling out the possible vacuum}, it follows that $\widetilde{r}(t,x)\rightarrow0$ as $x\rightarrow0$. Thus we have
 \begin{align}
     \widetilde{r}^2(t,x)=2\int^x_0\dfrac{1}{\rho(t,\widetilde{r}(t,y))}\dy\le C(p,T)\int^x_0y^{-p}\dy\le C(p,T)x^{1-p}.
 \end{align}
 Consequently, $x\ge C(T)\left(\widetilde{r}(t,x)\right)^{\frac{2}{1-p}}$, substituting this inequality into \eqref{rho inequality}, we have 
 \begin{align*}
     \rho(t,r)\ge C^{-1}(p,T)r^{\frac{2p}{1-p}}\text{ for }r\le 2R_0,
 \end{align*}
 where we used $2R_0<\widetilde{r}(t,X_0)$. Since for  $r\ge R_0$ we have $\dfrac{1}{2}\le\rho(t,r)$, we obtain the conclusion that $\min\{C^{-1}(p,T)r^{\frac{2p}{1-p}},C^{-1}(T)\}\le\rho(t,r)$.
    \end{proof}
\end{lemma}
\begin{remark}
    This \text{Theorem} not only characterizes the lower bound of $\rho$ but also excludes some decay rate as $r\rightarrow0$. It implies that if $\rho$ vanishes algebraically at the origin, then its algebraic decay exponent could not exceed $\dfrac{2p}{1-p}$. However, the lower bound remains compatible with slower algebraic decay and with logarithmic decay, such as $\rho(r)\sim|\log r|^{-\beta}$ with $\beta>0$.
\end{remark}

Now we transform the weak form obtained above into $(t,x)\in \mathbb{R}^n\times[0,T]$. For  $\forall\Phi\in C^1(0,T;C^1_c(\mathbb{R}^n))$, we define 
\begin{align*}
    \rho(t,\boldsymbol{x}):=\rho(t,|\boldsymbol{x}|),\quad U(t,\boldsymbol{x}):=u(t,|\boldsymbol{x}|)\dfrac{\boldsymbol{x}}{|\boldsymbol{x}|},\quad \Psi(t,r)=\int_{S^{n-1}}\Phi(r\boldsymbol{x},t)\mathrm{d}S_x.
\end{align*}
where $S^{n-1}$ denotes the $n-1$ dimensional sphere and $\mathrm{d}S_x$ is the surface measure on $S^{n-1}$. For the continuity equation, we have
\begin{align*}
    \int^\infty_0\rho(t,r)\Psi(t,r)r^m\dr&=\int^\infty_0\int_{S^{n-1}}\rho(t,r)\Phi(r\boldsymbol{x},t)r^m\mathrm{d}S_x\dr
    =\int_{\mathbb{R}^n}\rho(t,|\boldsymbol{x}|)\Phi(t,\boldsymbol{x})\dx.
\end{align*}
By the same argument,
\begin{align*}
    \int^t_0\int^\infty_0\rho(s,r)\Psi_t(s,r)r^m\dr\ds=\int^t_0\int_{\mathbb{R}^n}\rho(s,|\boldsymbol{x}|)\Phi_t(s,\boldsymbol{x})\dx\ds.
\end{align*}
\begin{align*}
    \int^t_0\int^\infty_0\left(\rho u\right)(s,r)\Psi_r(s,r)r^m\dr\ds&=\int^t_0\int^\infty_0\int_{S^{n-1}}\left(\rho u\right)(s,r) \nabla\Phi(r\boldsymbol{x},s)\cdot \boldsymbol{x}r^m\mathrm{d}S_x\dr\ds\\
    &=\int^t_0\int_{\mathbb{R}^n}\rho(s,|\boldsymbol{x}|)u(s,|\boldsymbol{x}|)\dfrac{x}{|\boldsymbol{x}|}\cdot\nabla\Phi(s,\boldsymbol{x})\dx\ds\\
    &=\int^t_0\int_{\mathbb{R}^n}\rho(s,|\boldsymbol{x}|)U(s,\boldsymbol{x})\cdot\nabla\Phi(s,\boldsymbol{x})\dx\ds
\end{align*}
Now we have obtained the weak form of the continuity equation in Eulerian coordinates $(t,\boldsymbol{x})\in[0,T]\times\in \mathbb{R}^n$. Next we prove the momentum equation. For the momentum equation, we only prove the Korteweg term, since the proofs for the remaining terms are similar to the proof of the continuity equation. Let $\Phi\in C^1(0,T;C^2_c(\mathbb{R}^n);\mathbb{R}^n)$,  we let the function $\Psi$ be that
\begin{align*}
    \Psi(t,r):=\int_{S^{n-1}}\Phi(r\boldsymbol{y},t)\cdot\boldsymbol{y}\mathrm{d}S_y
\end{align*}
It is straightforward to check $\Psi(t,0)=0$. Then we take $\Psi$ as a test function mentioned in Lemma \ref{momentum equation}.  First, we have
\begin{align*}
    \dfrac{\kappa}{2}\int^\infty_0|\rho_r|^2\left(3\Psi_r+\dfrac{\Psi}{r}\right)r\dr=\dfrac{\kappa}{2}\int^\infty_0|\rho_r|^2\left(\Psi_r+\dfrac{\Psi}{r}\right)r\dr+\kappa\int^\infty_0|\rho_r|^2\Psi_rr\dr.
\end{align*}
It follows from $(r^m\Psi)_r=\displaystyle r^m\int_{S^{n-1}}\partial_{x_i}\Phi^i(r\boldsymbol{y},t)\mathrm{d}S_y$ that
\begin{align*}
   \dfrac{\kappa}{2}\int^\infty_0|\rho_r|^2\left(\Psi_r+\dfrac{\Psi}{r}\right)r\dr=\dfrac{\kappa}{2}\int_{\mathbb{R}^{n}}|\nabla\rho(t,x)|^2\partial_{x_i}\Phi^i(t,x)\dx.
\end{align*}
Since $(\Phi(r\boldsymbol{y},t)\boldsymbol{y})_r=\boldsymbol{y}\otimes \boldsymbol{y}:\nabla\Phi(r\boldsymbol{y},t)$, it follows that $\Psi_r=\displaystyle\int_{S^{n-1}}\boldsymbol{y}\otimes \boldsymbol{y}:\nabla\Phi(r\boldsymbol{y},t)\mathrm{d}S_y$, thus for the rest term,  it holds that
\begin{align*}
    \kappa\int^\infty_0|\rho_r|^2\Psi_rr\dr=\kappa\int^\infty_0\int_{S^{n-1}}\nabla_x\rho\otimes\nabla_x\rho:\nabla_x\Phi(r\boldsymbol{y},t)r\mathrm{d}S_y\dr=\kappa\int_{\mathbb{R}^n}\nabla\rho\otimes\nabla\rho:\nabla\Phi(t,\boldsymbol{x})\dx.
\end{align*}
With the same argument, it holds that
\begin{align*}
    \int_{\mathbb{R}^n}\rho\nabla\rho\cdot\nabla\div\Phi\dx=\int^\infty_0\rho\rho_r\left(\Psi_r+\dfrac{m\Psi}{r}\right)_rr\dr.
\end{align*}
Now we complete proofs for the momentum  equation in Eulerian coordinates $(t,\boldsymbol{x})\in [0,T]\times\mathbb{R}^n$. For 3-D case, the proof is similar. Now the proof is entirely finished.


\appendix
\section{Some Basic Lemmas}\label{appendix}

For the convenience of readers, we list some basic lemmas that have been used frequently in this paper. 
The first lemma is a variant of the Arzela--Ascoli theorem. 
\begin{lemma}\cite{CHZ}\label{2.1}
    Let $\mathbb{D}\subseteq \mathbb{R}^2$ be a domain with Lipschitz boundary such that $\operatorname{Int}(\mathbb{D})\neq\varnothing$, where $\operatorname{Int}(\cdot)$ denotes the set
of interior points. Assume that $\{g_i:\mathbb{D}\to \mathbb{R}\}_{i\in N}$ is a collection
of functions, and $\{\mathbb{D}_\alpha,C_\alpha\}_{\alpha\in\mathbb{N}}$ is a collection of
compact subsets $\mathbb{D}_\alpha \Subset \mathbb{D}$ and constants $C_\alpha\in(0,\infty)$ such that, 
for each $\alpha\in \mathbb{N}$,
\begin{equation}
\begin{gathered}\label{2.1.1}
\mathbb{D}_\alpha \subseteq \operatorname{Int}(\mathbb{D}_{\alpha+1}),\qquad
\operatorname{Int}(\mathbb{D}_\alpha)\neq \varnothing,\qquad
\mathbb{D}=\bigcup_{\alpha\in N}\mathbb{D}_\alpha, \\
\sup_{i\in\mathbb{N}}\sup_{x\in \mathbb{D}_\alpha} |g_i(x)|\le C_\alpha, \\
\sup_{i\in N} |g_i(x)-g_i(y)|
\le C_\alpha |x-y|
\qquad \text{for all } x,y\in \mathbb{D}_\alpha .
\end{gathered}
\end{equation}
Then there exist both a continuous function $g:\mathbb{D}\to\mathbb{R}$ and a subsequence
$\{i_j\}_{j\in \mathbb{N}}$ with $i_j\to\infty$ as $j\to\infty$ such that
\[
\lim_{j\to\infty}\sup_{x\in K}|g_{i_j}(x)-g(x)|=0
\quad \text{ for all compact subsets } K\Subset \mathbb{D},
\]
and
\[
|g(x)-g(y)|\le C_\alpha |x-y|
\qquad \text{for all } x,y\in \mathbb{D}_\alpha,\ \text{for each } \alpha\in \mathbb{N}.
\]
\end{lemma}
Then we prove a disintegration theorem, and this lemma is crucial to derive the limiting radial momentum equation. 
\begin{lemma}\label{Disintergration theorem}
    Let $\mathbb{X}:=(0,\infty),\mathbb{Y}:=(0,T)$, $\nu$ be a non-negative finite Borel measure, $\Lambda(B)=\nu(\mathbb{X}\times B)$. Then  there exists a probability measure kernel $\{\Pi_t\}_{t\in\mathbb Y}$ satisfying the following properties.

    \begin{itemize}
        \item [1.] for each fixed $A\in \mathcal{B}(\mathbb{X})$, $t\rightarrow\Pi_t(A)$ is a measurable mapping.
        \item [2.] for each $A\in \mathcal{B}(\mathbb{X}),~B\in \mathcal{B}(\mathbb{Y})$, it holds that
        \begin{align*}
            \nu(A\times B)=\int_B\Pi_t(A)\mathrm{d}\Lambda(t),
        \end{align*}
        thus for any non-negative Borel function $\Phi:\mathbb{X}\times\mathbb{Y}\rightarrow[0,\infty]$,
        \begin{align*}
            \int_{\mathbb{X\times\mathbb{Y}}}\Phi(x,t)\mathrm{d}\nu=\int_\mathbb{Y}\left(\int_\mathbb{X}\Phi(x,t)\mathrm{d}\Pi_t\right)\mathrm{d}\Lambda(t).
        \end{align*}
        Moreover, this kernel is unique for $\Lambda$-a.e.  $t\in\mathbb Y$.
    \end{itemize}
    \begin{proof}
        We divide the proof into five steps.\\
        \textbf{1. } First, for each $q\in\mathbb{Q}_+$ and each $B\in\mathcal{B}(\mathbb{Y})$ we define
        \begin{align*}
            \nu_q(B):=\nu((0,q]\times B),
        \end{align*}
        since $0\le\nu_q\le\Lambda$, by Radon--Nikodym theorem, there exists a non-negative Borel measurable function $F_q$ such that
        \begin{align*}
            \nu((0,q]\times B)=\int_B F_q(t)\mathrm{d}\Lambda(t),\quad0\le F_q\le1 \text{ for }\Lambda ~a.e. 
        \end{align*}
        Combing the monotonicity of $\nu_q,q\in\mathbb{Q}_+$ and the uniqueness of Radon--Nikodym derivative, we have
        \begin{align}\label{properties of F}
            &F_q\le F_p ~a.e. \text{  for }\Lambda~a.e. \text{ and each }q<p,\quad p,q\in\mathbb{Q}_+,\nonumber\\
            &F_{q+\frac{1}{k}}\downarrow F_q,\quad F_k\uparrow1,\quad F_{\frac{1}{k}}\downarrow0,\text{ as k}\rightarrow\infty \text{ for }\Lambda~a.e. ,
        \end{align}
where $k\in\mathbb{N}$. Thus there exists a zero measure set $N$ such that $\Lambda(N)=0$ and \ref{properties of F} hold on $\mathbb{Y}\backslash N$, if we define $F_q$ on $N$ as follows: 
\begin{align*}
    F_q(t)=\begin{cases}
        1\quad q\ge1,\\
        0\quad q<1.
    \end{cases}
\end{align*}
Then \ref{properties of F} holds for each $t\in\mathbb{Y}$.\\
\textbf{2.} For $\forall t\in\mathbb{Y},s\in(0,1)$ we define 
\begin{equation*}
    \mathcal{R}_t(s):=\inf\{q\in\mathbb{Q}_+:F_q(t)\ge s\}.
\end{equation*}
It is straightforward to check $\{q\in\mathbb{Q}_+:F_q(t)\ge s\}$ is not empty and $0<\mathcal{R}_t(s)<\infty$ by \ref{properties of F}, in addition, the proof of the equivalence relation: $\mathcal{R}_t(s)\le q\Longleftrightarrow s\le F_q(t)$ for each $q\in\mathbb{Q}_+$ is also easy to check. Since $\{(t,s):\mathcal{R}_t(s)\le q\}=\{(t,s):s\le F_q(t)\}$ is a Borel measurable set in $\mathbb{Y}\times(0,1)$, thus $(t,s)\rightarrow\mathcal{R}_t(s)$ is a jointly measurable mapping.\\
\textbf{3.} For $\forall A\in\mathcal{B}(\mathbb{X})$, we define
\begin{align*}
\Pi_t(A):=\int^1_0I_A(\mathcal{R}_t(s))\mathrm{d}s=\left|\{s\in(0,1):\mathcal{R}_t(s)\in A\}\right|,
\end{align*}
where $|\cdot|$ denotes the Lebesgue measure. By virtue of the definition, $\Pi_t$ is a probability measure. Since $(t,s)\rightarrow\mathcal{R}_t(s)$ is non-negative jointly measurable, the mapping $t\rightarrow\Pi_t(A)$ is measurable for each $A\in\mathcal{B}(\mathbb{X})$ with the help of Tonelli theorem. Furthermore, for each $q\in \mathbb{Q}_+$,
\begin{align*}
    \Pi_t((0,q])=\left|\{s\in(0,1):\mathcal{R}_t(s)\le q\}\right|=F_q(t).
\end{align*}
\textbf{4.} For each Borel set $E\subset\mathbb{X}\times\mathbb{Y}$, we consider the non-negative finite measure 
\begin{align*}
    \widetilde{\nu}(E):=\int_\mathbb{Y}\int^1_0I_E(\mathcal{R}_t(s),t)\ds\mathrm{d}\Lambda
\end{align*}
It is straightforward to verify that $\widetilde{\nu}(\mathbb{X}\times\mathbb{Y})=\Lambda(\mathbb{Y})$. Moreover, for each $q\in\mathbb{Q}_+$ and $B\in\mathcal{B}(\mathbb{Y})$, we have
\begin{align*}
    &\widetilde{\nu}((0,q]\times B)=\int_B\Pi_t((0,q])\mathrm{d}\Lambda=\int_B F_q(t)\mathrm{d}\Lambda=\nu((0,q]\times B),\\
    &\widetilde{\nu}(\mathbb{X}\times B)=\Lambda(B)=\nu(\mathbb{X}\times B).
\end{align*}
Thus $\widetilde{\nu}=\nu$. And it holds for each Borel non-negative function $\Phi$ that 
\begin{align*}
    \int_{\mathbb{X}\times\mathbb{Y}}\Phi(x,t)\mathrm{d}\nu=\int_\mathbb{Y}\int^1_0\Phi(\mathcal{R}_t(s),t)\ds\mathrm{d}\Lambda=\int_\mathbb{Y}\left(\int_\mathbb{X}\Phi(x,t)\mathrm{d}\Pi_t\right)\mathrm{d}\Lambda.
\end{align*}
\textbf{5.} Assuming there exists another probability measure kernel $\{\pi_t\}_{t\in\mathbb{Y}}$ satisfying the same properties as $\{\Pi_t\}_{t\in\mathbb{Y}}$, then for each $q\in\mathbb{Q}_+$ and each $B\in\mathcal{B}(\mathbb{Y})$, we have
\begin{align*}
    \int_B\Pi_t((0,q])\mathrm{d}\Lambda=\nu((0,q]\times B)=\int_B\pi_t((0,q])\mathrm{d}\Lambda.
\end{align*}
Thus $\Pi_t=\pi_t ~\Lambda$-a.e.  $t$.

    \end{proof}
\end{lemma}
The Hardy's inequality is used to deal with the modified velocity near the vacuum interface.
\begin{lemma}\label{Hardy}
    For  $\forall a\in[0,R)$ and $v(r)\in H^1(a,R)$ with $v(R)=0$, the following inequality holds
    \begin{align*}
        \int^R_a\dfrac{|v(r)|^2}{(R-r)^2}r^m\dr\le4\int^R_a|v_r|^2r^m\dr.
    \end{align*}
    \begin{proof}
        Since
        \begin{align*}
            \int^R_a\dfrac{|v(r)|^2}{(R-r)^2}r^m\dr=\dfrac{|v|^2r^m}{R-r}\bigg|^R_a-\int^R_a\dfrac{m|v|^2r^{m-1}+2vv_rr^m}{R-r}\dr.
        \end{align*}
        Thus by virtue of H$\ddot{o}$lder inequality,
        \begin{align*}
             \int^R_a\dfrac{|v(r)|^2}{(R-r)^2}r^m\dr\le2\int^R_a\left|\dfrac{vv_rr^m}{R-r}\right|\dr\le 2\left(\int^R_a\dfrac{|v(r)|^2}{(R-r)^2}r^m\dr\right)^{\frac12}\left(\int^R_a|v_r|^2r^m\dr\right)^{\frac{1}{2}}.
        \end{align*}
    \end{proof}
\end{lemma}
The following Lemma is crucial to prove the defect measure of $|\rho^a_r|^2$ is  a square-integrable function.
\begin{lemma}\label{distribution derivatives}
Let $\mathbb{X}=(c,d),\mathbb{Y}=(0,t)$ with $0<c<d,0<t$.  If $f\in\mathcal{D}'((c,d)\times(0,t))$
    satisfies $f_r=0$ in $\mathcal{D'}((c,d)\times(0,t))$, then there exists a unique $g\in\mathcal{D}'((0,t))$ such that $f=1\Box g$, where $1\Box g$ is defined by
    \begin{align*}
        \langle1\Box g,\psi\rangle:=\langle g,\int^d_c\psi(t,r)\dr\rangle\text{ for each }\psi\in C^\infty_c((c,d)\times(0,t)).
    \end{align*}
    \begin{proof}
        Let $\chi\in C^\infty_c((c,d))$ satisfy $\int^d_c\chi(r)\dr=1$. For $\forall\eta\in C^\infty_c((0,t))$, we define $<g,\eta>:=\langle f,\chi\eta\rangle$. For  $\forall\psi\in C^\infty_c((c,d)\times(0,t))$,  we let $\bar\psi=\displaystyle\int^d_c\psi(t,r)\dr$, $\Psi(t,r)=\displaystyle\int^r_c(\psi(t,s)-\chi(s)\bar\psi(t))\ds$, and thus $\Psi\in C^\infty_c((c,d)\times(0,t))$. Since $\psi=\chi\bar\psi+\Psi_r$, it holds that
        \begin{align*}
            \langle f,\psi\rangle=\langle f,\chi\bar\psi\rangle+\langle f,\Psi_r\rangle=\langle f,\chi\bar\psi\rangle-\langle f_r,\Psi\rangle=\langle g,\bar\psi\rangle.
        \end{align*}
        The uniqueness is immediately.
    \end{proof}
\end{lemma}
Next, we show that the spherical symmetry of the initial data
is preserved under mollification.

\begin{lemma}\label{mollify initial data still spher}
    Let $f(\boldsymbol{x})=f(|\boldsymbol{x}|)$ be a local integrable radial  function in $\mathbb{R}^n$ and $j_a$ be a standard mollifier, then $j_a*f$ is still a radial function.
    Moreover, if $f$ is a vector radial  function, then there must exist a smooth scalar function $g_a$ such that
    \begin{align*}
        f^a(\boldsymbol{x})=g_a(|\boldsymbol{x}|^2)\boldsymbol{x}.
    \end{align*}
    \begin{proof}
    \textbf{1. }
        Let $f^a(\boldsymbol{x})=(j_a*f)(\boldsymbol{x})$. For each orthogonal matrix $Q\in O(n)$, where $O(n)$ denotes the orthogonal group of order $n$. We have
        \begin{align*}
            f^a(Q\boldsymbol{x})&=\int_{\mathbb{R}^n}j_a(Q\boldsymbol{x}-\boldsymbol{y})f(\boldsymbol{y})\mathrm{d}\boldsymbol{y}=\int_{\mathbb{R}^n}j_a(Q(\boldsymbol{x}-\boldsymbol{z}))f(Q\boldsymbol{z})\mathrm{d}\boldsymbol{z}\\
            &=\int_{\mathbb{R}^n}j_a(\boldsymbol{x}-\boldsymbol{y})f(\boldsymbol{z})\mathrm{d}\boldsymbol{z}=f^a(\boldsymbol{x}).
        \end{align*}
\textbf{2. } For each $\boldsymbol{x}\neq0$, let $e=\dfrac{\boldsymbol{x}}{|\boldsymbol{x}|}$, it is straightforward to check 
\begin{align*}
    Q=2e\otimes e-\mathbb{I}_n\in O(n),\quad Qe=e.
\end{align*}
For each $w$: $w\cdot e=0$, it holds that $Qw=-w$. Since $Qx=x$,
\begin{align*}
    f^a(\boldsymbol{x})=f^a(Q\boldsymbol{x})=Qf^a(\boldsymbol{x}).   
\end{align*}
We decompose $f^a(\boldsymbol{x})$ into its components parallel and orthogonal to $e$,
\begin{align*}
f^a(\boldsymbol{x})=\alpha e+w,\qquad w\cdot e=0,
\end{align*}
where $\alpha= f^a(\boldsymbol{x})\cdot e$. Thus, $Q f^a(\boldsymbol{x})=\a e-w$. Since $\a e+w=\a e-w$, we have $w=0$. Thus, $f^a$ is parallel to $\boldsymbol{x}$, so there exists a scalar radial function $g_a(|\boldsymbol{x}|)$ such that
\begin{align*}
    f^a(\boldsymbol{x})=\widetilde{f}(|\boldsymbol{x}|)\dfrac{\boldsymbol{x}}{|\boldsymbol{x}|} \text{ for each }x\neq0.
\end{align*}
Let $e_1=(1,0,\cdots,0)$, we define $g(s):=f^a(se_1)\cdot e_1$ for $s\in\mathbb{R}$. It is obvious that $g\in C^\infty(\mathbb{R})$, $g(-s)=-g(s)$, and $g(0)=0$. Let $h(s):=\int^1_0 g'(\theta s)\mathrm{d}\theta$, we have $g(s)=sh(s)$ and $h(-s)=h(s)$. Since $h$ is a smooth even function, there exists  $g_a\in C^\infty([0,\infty))$ such that  $h(s)=g_a(s^2)$. Thus, $g(s)=sg_a(s^2)$. Moreover, since $\widetilde{f}(r)=g(r)$
\begin{align*}
    f^a(\boldsymbol{x})=\widetilde{f}(|\boldsymbol{x}|)\dfrac{\boldsymbol{x}}{|\boldsymbol{x}|}=|\boldsymbol{x}|g_a(|\boldsymbol{x}|^2)\dfrac{\boldsymbol{x}}{|\boldsymbol{x}|}=\boldsymbol{x}g_a(|\boldsymbol{x}|^2).
\end{align*}
       
    \end{proof}
\end{lemma}

\bigskip



\noindent{\bf{Acknowledgments:}}
This research is partially supported by National Key R$\&$D Program of China (No. 2022YFA1007300). The research of Yachun Li was also supported in part by the National Natural Science Foundation of China (No. 12371221), the National Key R\&D Program of China (N0. 2024YFA1013302), the Shanghai Center for Mathematical Sciences, and the Fundamental Research Funds of the Central Universities.
The research of Shengguo Zhu was also supported in part by 
the National Natural Science Foundation of China under the Grant  12471212 
and the Royal Society (UK)-Newton International 
Fellowships NF170015.

\bigskip 
 
\noindent{\bf Data Availability Statements:}
Data sharing not applicable to this article as no datasets were generated or analyzed during the current study.

\bigskip

\noindent{\bf Conflict of interests:}
The authors declare that they have no competing interests.

\end{document}